\documentclass{article}
\usepackage[utf8]{inputenc}

\usepackage{amsmath, amssymb, stmaryrd}
\usepackage{amsthm}
\theoremstyle{plain}
\newtheorem{assumption}{Assumption}
\newtheorem{proposition}{Proposition}

\usepackage{ mathrsfs }
\usepackage{graphics}
\usepackage{epsfig}
\usepackage{graphicx}
\usepackage{tabularx}
\usepackage{hyperref}
\usepackage{booktabs}
\usepackage{systeme}
\usepackage{hyperref}
\usepackage{caption}
\usepackage{url}
\usepackage{subcaption}
\usepackage[shortlabels]{enumitem}
\usepackage{geometry}
 \usepackage{listings}
\usepackage{color} 
\definecolor{mygreen}{RGB}{28,172,0} 
\definecolor{mylilas}{RGB}{170,55,241}
\hypersetup{colorlinks=true, urlcolor=blue}

\newcommand\Ls{\mathscr{L}}

\newcommand\lam{{\bs \lambda}}

\newcommand\p{{\bs p}}

\newcommand\R{\mathbb{R}}

\renewcommand\u{{\bs u}}
\renewcommand\v{{\bs v}}

\newcommand\x{{\bs x}}

\newcommand\grad{\nabla}

\newcommand{\bs}{\boldsymbol}

\renewcommand{\exp}{\mathop{\rm exp}}

\renewcommand{\min}{\mathop{\rm min}}

\renewcommand{\max}{\mathop{\rm max}}
\renewcommand{\sup}{\mathop{\rm sup}}

\newcommand{\norm}[1]{\ensuremath{\left\| #1 \right\|}} 

\newtheorem{theorem}{Theorem}[section]

\title{Solving Singular Control Problems in Mathematical Biology, Using PASA}
\author{Summer Atkins, Mahya Aghaee, Maia Martcheva, William Hager}
\date{October 2020}
\begin{document}
\maketitle
\begin{abstract}
In this paper, we will demonstrate how to use a nonlinear polyhedral constrained optimization solver called the Polyhedral Active Set Algorithm (PASA) for solving a general singular control problem. We present methods of discretizing a general optimal control problem that involves the use of the gradient of the Lagrangian for computing the gradient of the cost functional so that PASA can be applied.  When a numerical solution contains artifacts that resemble ``chattering'', a phenomenon where the control oscillates wildly along the singular region, we recommend a method of regularizing the singular control problem by adding a term to the cost functional that measures a scalar multiple of  the total variation of the control, where the scalar is viewed as a tuning parameter. We then demonstrate PASA's performance on three singular control problems that give rise to different applications of mathematical biology. We also provide some exposition on the heuristics that we use in determining an appropriate size for the tuning parameter. 
\end{abstract}
\section{Introduction}

Optimal control theory is a tool that is used in mathematical biology for observing how a dynamical system  behaves when employing one or many  variables that can be controlled outside of that system. 
Mathematical biologists apply optimal control theory to disease models of immunologic  and epidemic types \cite{Joshi, lenhartHIV, LenhartTB, Vaccinemobile}, to management decisions in harvesting \cite{Neubert2003, Ding2009}, and to resource allocation models \cite{JoshiPlant, King1982}. 
In practice, mathematical biologists tend to construct optimal control problems with quadratic dependence on the control. 
Problems of this structure are well-behaved in the sense that there are established methods  of proving existence and uniqueness of an optimal control \cite{Filippov,Pontryagin, fleming}.
In addition, for numerically solving problems of this form, many employ the forward-backward sweep method, a numerical method presented in Lenhart and Workman's book \cite{Lenhart2007} that involves combinations of the forward application and the backward application of a fourth-order Runge-Kutta method. 
We direct the reader to the following references  \cite{gradientalgorithms,quasinewton,BettsSurvey,Raosurvey} which are excellent surveys of other numerical methods, such as gradient methods, quasi-Newton methods, shooting methods, and collocation methods,  that are used within the optimization community for solving optimal control problems.

Control problems in biology tend to depend quadratically with respect to the control due to the construction of the objective or cost functional, which is the functional that is being optimized with respect to the control variables.  
The construction of the objective functional is an essential component to optimal control theory because it measures our criteria for determining what control strategy is deemed ``best". 
In mathematical biology, the costs are frequently nonlinear and depend on the states and the controls, and a cost term for a particular control may be the sum of a bilinear term in that control and one for the states and a quadratic term in the control. Frequently, the quadratic term has a lower coefficient.
With regard to the principle of parsimony, it is difficult to justify the use of a quadratic term for representing the cost of administering a control. 
A linear term  would be a more realistic representation of the cost of applying a control.  
For optimal control problems with linear dependence on the control, it is possible to obtain a solution that is piecewise constant where the constant values correspond to the bounds of the control.
An optimal control of this structure, which is often called a ``bang-bang" control,   can be readily interpreted and implemented. 
These characteristics compel many (see \cite{Mahya,Ledzewicz2011,Rongpaper, Neubert2003,Ding2009, King1982, JoshiHarvest}), to use control problems with linear dependence on the control for biological models. 
There are however evident setbacks to using  optimal control problems in which the control appears linearly. 
For one thing these problems are much more difficult to solve analytically due to the potential existence of a singular subarc. 
As demonstrated in \cite{Neubert2003, King1982,Ledzewicz2011}, procedures for obtaining an explicit formula for the singular case involve taking one or many time derivatives of the switching function as a means to gain a system of equations. 
Often the first order and second order necessary conditions for optimality, which are respectively named Pontryagin's minimum principle \cite{Pontryagin}  and the Generalized  Legendre-Clebsch Condition/Kelley's Condition \cite{Powers, Robbins,Lewis, Zelikin1994}, are checked. 
Naturally, procedures for explicitly solving for singular control problems increase in difficulty when  multiple state variables and multiple control variables are involved. 

Additionally, numerical procedures for solving singular control problems are inevitably problematic. 
If the parameters of an optimal control problem with linear dependence are set to where all optimal control variables are bang-bang, then the forward backward sweep method \cite{Lenhart2007} can successfully run.
However, if the presence of a singular control is a possibility, then  forward backward sweep is not advisable. 
In \cite{AUV}, Foroozandeh and De Pinho test four numerical methods, including the Imperial College London Optimal Control Software (ICLOCS) \cite{ICLOCS} and the Gauss Pseudospectral Optimization Software (GPOPS)\cite{GPOPS},  by solving a singular optimal control problem for Autonomous Underwater Vehicles (AUV). 
They find that three of these methods have difficulty in detecting the structure of the optimal control  and in accurately computing the swtiching points without a priori information. 
Switching points are the corresponding points in time when an optimal control switches from singular to non-singular and vice versa.
When solving for the AUV problem, both ICLOCS and GPOPS obtain a control that ehxibitted oscillations within the singular region, causing both methods to be unable to provide direct information about the switching points. 
Foroozandeh and De Pinho conclude that only the mixed binary non-linear programming method (MBNLP) \cite{Mixedbin} is successful in accurately approximating the optimal control to the AUV problem and its switching points.

One of  the predominate issues associated with singular control problems is the concept of ``chattering", which is also known as the  ``Fuller Phenomenon"  (see \cite{Zelikin1994}).
As mentioned in Zelikin and Borisov's book \cite{Zelikin1994}, an optimal control is said to be \emph{chattering} if the control  oscillates infinitely many times between the bounds of the control over a finite region. 
It is thought that such an event occurs in singular control problems when the optimal control consists of singular and non-singular (bang) subarcs and those regions cannot be directly joined.
In \cite{Powers}, MacDanell and Powers present some necessary conditions for joining singular and nonsingular subarcs. 
And Zelikin and Borisov \cite{Zelikin1994} mention other theorems that can be used to verify when chattering is present. 
It is possible that the discretization of an opitmal control problem causes a numerical method to generate numerical artifacts that resemble chattering even though the optimal control does not chatter. 
However, it is difficult to determine when a numerical solution is exhibiting many oscillations due to chattering or due to numerical artifacts. 

Regardless of the situation, it is evident that a chattering optimal control would be an unrealistic procedure to implement. 
A way to bypass this issue is to solve for a penalized version of the optimal control problem. 
In \cite{regularize}, Yang et al. present methods of regularizing optimal control problems by adding a penalty term to the cost functional of the original problem. 
The penalization terms suggested  in \cite{regularize} are: a weighted parameter times the $L^2$ norm of the control, a weighted parameter times the $L^2$ norm of the derivative of the control, and a weighted parameter times the $L^2$ norm of the second derivative of the control.  
In \cite{Ding2009}, Ding and Lenhart employ a penalty term to a harvesting optimal control problem to avoid a potential chattering result found in the control $h$. 
The penalty term  that was applied for this problem consisted of  a weighted parameter times $\norm{\grad h}^2_{L^2}$. 
The penalty term that is used  in \cite{Ding2009} adds convexity properties to the  control problem  which  can be beneficial for verifying existence and uniqueness of an optimal control. 
Additionally,  $\norm{\grad h}^2_{L^2}$ can be discretized to where it can serve as a crude estimate for the total variation of the control. 
However,  Lenhart and Ding \cite{Ding2009} need to use variational inequalities to solve for their problem, and incorporating such a penalty term restricts their set of admissible controls into a functional space that requires its controls to be differentiable. 
The penalty term that we believe to show the most promise has been recently suggested in Capognigro et al.'s work \cite{Caponigro2018}. 
 In \cite{Caponigro2018}, Caponigro et. al. recommends a method of regularizing chattering  in optimal control problems by adding a penalty term that represents a  penalty  parameter times the total variation of the control. 
 Such penalization will influence the numerical solution in a way that will reduce the number of oscillations. 
In \cite{Caponigro2018}, this penalty is applied to the Fuller's problem, which is the classical example that introduced the concept of chattering, and they  obtain a ``quasi" optimal solution to the Fuller problem that does not chatter. 
Caponigro et. al also prove that the optimal value of the penalized problem converges to the associated optimal value of the original problem as the penalty weight parameter $p$ tends to zero.

In this paper, we will demonstrate how to use a nonlinear polyhedral constrained optimization solver called 
PASA\footnote{To access PASA software that can be used on MATLAB for Linux and Unix operating systems, download the SuiteOPT Version 1.0.0 software given on \url{http://users.clas.ufl.edu/hager/papers/Software/ } . For future reference, any updates to the software will be uploaded to this link. }, developed by Hager and Zhang \cite{Hager2016}, to solve a general singular control problem that is being regularized by use of a total variation term \cite{Caponigro2018}.
We recommend PASA because it is user friendly  to those who are not as acquainted with optimization techniques for optimal control problems, and it is freely accessible to 
 use on MATLAB for Linux and Unix operating systems. 
According to Hager and Zhang \cite{Hager2016}, PASA consists of two phases, with the first phase being the gradient projection algorithm and the second phase being  any algorithm that is used for solving a linearly constrained optimization problems. 
The gradient projection algorithm is an optimization solver commonly used for bounded constrained optimization problems. 
When applying gradient descent to a bounded constrained optimization problem, it is possible to obtain an iterate that lies outside of the feasible set due to the negative direction of the gradient. 
Projected gradient method takes this issue into consideration by adding additional steps involving projecting points outside the feasible set onto the feasible set. 
For more information on projected gradient methods we direct the reader to the following works \cite{Bertsekas1976, Calamai1987, Goldstein1964, Levitin1966, McCormick1972, Hager2016}.
Using PASA for solving optimal control problems will  involve converting optimal control problems into discretized optimization problem. 
 In this paper, the discretization for the general singular control problem will involve using explicit Euler's method for the state equations, and left-rectangular integral approximation for the original cost functional. 
 Additionally, we need to ensure that the discretized and penalized cost functional is differentiable, which will require performing a decomposition of the terms associated with the total variation penalty. 
 We will use the gradient of Lagrangian function of the discretized optimal control problem  for computing the gradient to the discretized cost functional. 
 Conveniently, the process of ensuring that the gradient of the Lagrangian is equal to the gradient of the discretized objective functional will yield a discretization procedure for the adjoint equations. 

 Further we demonstrate PASA's performance on three  singular control problems that are being regularized via bounded variation. 
 The order of these examples will increase in difficulty based upon the number of state variables and the number of control variables. 
 Additionally, each example gives rise to different applications of mathematical biology. 
 Explicit formulas for the singular case and for the switching points are obtained for the first two examples, allowing us to compare PASA's numerical results with the exact solution.
 For each problem, we will illustrate the discretization process and then present numerical results that were obtained when solving for both the unpenalized and penalized problem. 
 We also provide some exposition on the heuristics that we use  in determining an appropriate size for the tuning parameter $\rho$.

 The first example is a fishery harvesting problem that was  presented in Clark's \cite{Clark} and later restated in  Lenhart and Workman's \cite{Lenhart2007}.
 The fishery problem consisted of one state variable and one control variable, where the state variable represents the fish population and the control variable represents harvesting effort. 
 In \cite{Lenhart2007}, an explicit formula for the singular case can be obtained through using Pontryagin's maximum principle \cite{Pontryagin}; however, forward-backward sweep is unable to solve for the problem whenever parameters are set to ensure existence of a singular subarc.  
 The second example is from King and Roughgarden's work \cite{King1982}, and the optimal control problem is a resource allocation model for studying an annual plant's allocation procedure in distributing photosynthate.
 The control problem consists of two state variables  and one control variable. 
 One  variable measures the weight of the components of a plant that correspond to vegetative growth while the other variable measures the weight of the components of a plant that correspond to reproductive growth. 
 The control variable used for this problem represents the fraction of photosynthate being reserved for vegetative growth. 
 King and Roughgarden use Pontryagin's maximum principle \cite{Pontryagin} and find conditions based upon the parameters of the problem for determining when the optimal control would be bang-bang or concatenations of singular and bang control. 
 They verify that their explicit formula for the singular subarc satisfy the generalalized Legendre-Clebsch Condition and the strengthened Legendre-Clebsch Condition \cite{Powers, Robbins,Lewis, Zelikin1994}.
 Additionally, King and Roughgarden use one of  MacDanell and Power's Junction theorems \cite{Powers} to show the optimal control to the problem satisfies the necessary conditions for joining singular and non-singular subarcs. 
 When using PASA to solve for this plant problem, we only need to use the regularization term for a degenerate case of the problem. 
 
 The final example, which is from  Ledzewicz, Aghaee, Sch\"{a}ttler's \cite{Mahya}, is an optimal control problem where the three state variables involved correspond to an SIR model with demography.  
 An SIR model is a compartmental model that is used for modeling the spread of an infectious disease in a population, where the population is divided into the following three classes: 1. $S$ is the class of individuals who are susceptible to the disease; 2. $I$ is the class of infected individuals who are assumed to be infectious;  and 3. $R$ is the class of individuals who have recovered from the disease and are considered immune to the disease.
 For books covering mathematical models for epidemiology,  we recommend Brauer and Castillo-Chavez's \cite{Chavez} and Martcheva's \cite{MartchevaBook}. 
 The optimal control problem used in  Ledzewicz et al. \cite{Mahya} consists of two control variables where one control represents vaccination while the other represents treatment. 
Ledzewicz et al. numerically solve for this problem with parameters set to where a singular subarc is present in the optimal vaccination strategy while the optimal treatment strategy obtained appears to be bang-bang. 
 When using PASA, we regularize only the vaccination control via bounded variation since the treatment control contains no oscillations. 
 In the Appendix section, we provide the MATLAB code that was used for solving the last example to illustrate how to write up these optimal control problems.

\section{Discretization of Penalized Control Problem via Bounded Variation}

The following problem is the optimal control problem of interest:
\begin{equation}\label{eqn: generalOCP}
    \begin{array}{rl}
      \min\limits_{u\in\mathscr{A}} &J(\u) = \int_0^T g(\x(t),\u(t))dt\\
      \text{s.t.} & \dot x_i(t) = f_i(\x(t),\u(t)) \; \text{ for  all } i=1,\dots, n,\\
      & x_{i}(0)=x_{i,0}, \text{ for all } i=1,\dots, n, \\
      &{\xi_j\leq u_j(t)\leq \omega_j }\; \text{ for all } t\in[0,T] \text{ and for all }  j=1,\dots, m,
    \end{array}
\end{equation}
where functions $f_1$,$\dots$, $f_n$, and $g$ are assumed to be continuously differentiable in all arguments. 
In addition we assume that  if any of the $m$ control variables $u_i$  appear in functions $f_1,\dots,f_n$ and $g$, then $u_i$ will appear linearly.  
The class of admissible controls,  will be as follows:
\[\mathscr{A}=\{\u\in (L^1(0,T))^m\;|\;  {\bs{\xi}\leq \u(t)\leq \bs{\omega}}\text { for all } t\in[0,T] \}.\] 
We assume that the conditions for Filippov-Cesari Existence Theorem \cite{Filippov} hold for problem (\ref{eqn: generalOCP}).
State vector $\x\in\R^n$ consists of $n$ state variables that satisfy the state equations and initial values given in problem (\ref{eqn: generalOCP}). 

The common procedure for solving problem (\ref{eqn: generalOCP}) is employing Pontryagin's Minimum Principle \cite{Pontryagin} to generate the first order necessary conditions for optimality. 
We first define the Hamiltonian function to problem (\ref{eqn: generalOCP}), $H(\x,\u,\lam)$, as 
\begin{align*}
    H(\x,\u,\lam)&=g(\x, u)+\lam^{T}\bs{f}(\x,\u)\\
                &= g(\x,\u)+\sum\limits_{i=0}^{n}\lambda_{i}(t)f_i(\x,\u),
\end{align*}
where $\lam\in\R^n$ is the adjoint vector. 
We have from Pontryagin's minimum principle \cite{Pontryagin} that if $\u^*$ is the optimal control to problem (\ref{eqn: generalOCP}) with corresponding trajectory $\x^*$, then there exists a non-zero adjoint vector function $\lam^*$ that is a solution to the following adjoint system
\begin{align*}
    \dot\lambda_{\ell}(t) &=-\frac{\partial H(\x,\u,\lam)}{\partial x_{\ell}},\;\text{for all } \ell=1,\dots,n,\\
    \lambda_{\ell}(T)&=0, \text{ for all } \ell=1,\dots, n,
    \end{align*}
   and satisfies 
   \begin{align*}
       H(\u^*, \x^*,\lam^*)&=\min\limits_{\u\in\mathscr{A}}H(\u,\x^*, \lam^*).
   \end{align*}
   Using our definition of the Hamiltonian function the adjoint equations will be 
\begin{align}\label{eqn: adjoint}
   \dot\lambda_{\ell}(t)&= -\frac{\partial g(\x,\u)}{\partial x_{\ell}} -\sum\limits_{i=1}^{n}\lambda_i(t) \frac{\partial f_i(\x,\u)}{\partial x_{\ell}} \; \text{ for all } \ell=1,\dots, n,
\end{align}
 with transversality conditions being
 \begin{equation}\label{eqn: transversality}
     \lambda_{\ell}(T)=0, \;\text{ for all } \ell =1,\dots, n. 
 \end{equation}
Based on assumptions of functions $g$ and $f_1,\dots, f_n$ in problem (\ref{eqn: generalOCP}), the Hamiltonian depends linearly with respect to each  control. 

{For demonstration purposes, we assume that every component of control vector $\u$ in problem (\ref{eqn: generalOCP}) needs to be regularized via bounded variation \cite{Caponigro2018}. 
This means that when numerically solving for problem (\ref{eqn: generalOCP}) without this regularization term, we obtain oscillatory numerical artifacts or other unwarranted numerical artifacts in each control.}
To regularize problem (\ref{eqn: generalOCP}) via bounded variation we introduce a tuning vector $\bs{\rho}$ where $0\leq \rho_j< 1$  for all $j=1,\dots,m$, and we present Conway's \cite{Conway} definition of the total variation function $V$ of a real or complex valued function  $u$  that is defined on the interval $[a,b]$.
Let $\mathbb{K}$ be the field of complex numbers or the field of real numbers. Given $u:[a,b]\to \mathbb{K}$ the total variation of $u$ on $[a,b]$ is defined to be
\begin{equation}\label{eqn: totalvariation}
    V(u) = \sup\limits_{\mathcal{P}}\sum\limits_{k=0}^{N_P-1}|u(t_{k+1})-u(t_k)|,
\end{equation}
where $\mathcal{P}=\{P=\{t_0,t_1,\dots,t_{N_P}\}:\; a\leq t_0\leq t_1\leq \dots\leq t_{N_p}\leq b\}$.
Notice that if we assume that function $u$ is real-valued and piecewise constant on $[a,b]$, then 
the total variation of $u$ on $[a,b]$ is the sum of the absolute value of the jumps in $u$. 
The penalization of problem (\ref{eqn: generalOCP}) via bounded variation will be as follows:
\begin{equation}\label{eqn: generalOCPpen}
    \begin{array}{rl}
        \min\limits_{\u\in\mathscr{A}} &J_{\bs{\rho}}(\u)=\int_0^Tg(\x(t)\u(t))dt+\sum\limits_{j=1}^{m}\rho_j V(u_j)  \\
         \text{s.t.} & \dot x_i(t) = f_i(\x(t),\u(t)), \; \text{for all } i=1,\dots, n,\\
         &x_i(0)=x_{i,0},\;\text{ for all } i=1,\dots,n, \\
         & \xi_j \leq u_j(t) \leq \omega_j \; \text{ for all } t\in [0,T] \; \text{and }j=1,\dots, m,\\
    \end{array}
\end{equation}
where  $V(u_j)$ is the total variation of $u_j$ on interval $[0,T]$  and $0\leq \rho_j<1$ is the bounded variation penalty parameter associated with control variable $u_j$ for all $j=1,\dots,m$.
{For this problem,  we assume that all control variables need to be penalized.  
However, if based upon observations of the numerical solutions for  the unregularized problem (\ref{eqn: generalOCP}), we notice that one of the control variables, $u_{\ell}$, exhibits no numerical artifacts or unusual oscillations, then we recommend solving for problem (\ref{eqn: generalOCPpen}) with the corresponding tuning parameter $\rho_{\ell}$ set to being zero.  }
We can construct the Hamiltonian function that corresponds to problem (\ref{eqn: generalOCPpen}), and the Hamiltonian gives the same adjoint equations (\ref{eqn: adjoint}) and transversality conditions (\ref{eqn: transversality}). 

For using PASA to numerically solve for the penalized problem, we need to discretize  problem (\ref{eqn: generalOCPpen}) and discretize the adjoint equations (\ref{eqn: adjoint}). 
We present the method of discretizing the penalized problem only, but we emphasize that we can use PASA to solve for the discretized unpenalized problem by numerically solving for the discretized penalized problem when  the penalty vector, 
$\bs{\rho}$, is set to being the zero vector. 
We begin by partitioning time interval $[0,T]$, by using $N+1$ equally spaced nodes. 
For all $i=1,\dots, n$ and for all $k=0,\dots, N$ we denote $x_{i,k}=x_{i}(t_k)$, and we emphasize that component $x_{i,0}$ will be the initial value given in problem (\ref{eqn: generalOCPpen}).
So for each $i=1,\dots,n$,  we have that $\x_i\in\R^{N+1}$  with $x_{i,0}$ being the initial value given. 
We will assume that for all $j=1,\dots,m$, control variable $u_j$ is constant over each mesh interval. 
For all $j=1,\dots, m$ and for all $k=0,\dots, N-2$ we denote $u_{j,k} = u_j(t)$ for all $t_k\leq t<t_{k+1}$, and for all $j=1,\dots, m$ we denote $u_{j,N-1}=u_j(t)$ for all $t_{N-1}\leq t\leq t_N=T$.
So for each $j=1,\dots,m$ we have that $\u_j\in\R^{N}$.
The assumption of each control variable $u_j$ being constant on each mesh interval allows us to express the total variation of $u_j$ on $[0,T]$ in terms of the particular partition of $[0,T]$ that is used for the discretization: 
\[
   V(u_j) = \sum\limits_{k=0}^{N-2} |u_{j,k+1}-u_{j,k}|.
\]
For simplicity of discussion, we use left-rectangular integral approximation to discretize the integral used in objective functional $J_{\bs{\rho}}$, and we use forward Euler's method to discretize the state equations given in problem (\ref{eqn: generalOCPpen}).
In general, we recommend using an explicit scheme for discretizing the state equations if the dynamics of the system have only initial conditions involved. 
We then have the following 
\begin{equation}\label{eqn: generaldiscretepen}
    \begin{array}{rl}
    \min &J_{\bs{\rho}}(\u_1,\dots, \u_m)= \sum\limits_{k=0}^{N-1}hg(\x_{\bs{\cdot},k},\u_{\bs{\cdot},k}) + \sum\limits_{j=1}^{m}\rho_j\sum\limits_{k=0}^{N-2}|u_{j,k+1}-u_{j,k}|\\
    & x_{i,k+1} = x_{i,k}+ hf_i(\x_{\bs{\cdot},k},\u_{\bs{\cdot},k})  \; \text{ for all } i=1,\dots, n \; \text {and for all } k=0,\dots, N-1,\\
    & \xi_j\leq u_{j,k} \leq \omega_j, \; \text{ for all } j=1,\dots, m\; \text{ and for all } k=0,\dots, N-1,
    \end{array}
\end{equation}
where $h=\frac{T}{N}$ is the mesh size, $\x_{\bs{\cdot},k}=[x_{1,k}, x_{2,k}, \dots, x_{n,k}]$, and $\u_{\bs{\cdot},k}= [u_{1,k},u_{2,k},\dots, u_{m,k}]$ for all $k=0,\dots,N-1$.

Since PASA consists of a phase that uses a projected gradient method, we need the objective functional in problem (\ref{eqn: generaldiscretepen}) to be differentiable, which is not the case due to the absolute value terms that correspond to the discretization of the total variation function. 
We need to perform a decomposition of each absolute value term so that $J_{\rho}$ can be differentiable. 
For each $j=1,\dots, m$, we introduce two $N-1$ dimensional vectors $\bs{\zeta}_j$ and $\bs{\iota}_j$ whose entries are non-negative, and every entry of $\bs{\zeta}_j$ and  $\bs{\iota}_j$ will be defined as 
\[
{|u_{j,k+1}-u_{j,k}|=\zeta_{j,k}+\iota_{j,k}\; \text{ for all } k=0,\dots, N-2.}\]
Another way of viewing $\bs{\zeta}_j$ and  $\bs{\iota}_j$, is that each component $\zeta_{j,k}$ and  $\iota_{j,k}$ will be defined based upon the following conditions:
\begin{align*}
    \textbf{Condition 1:}& \text{ If } u_{j,k+1}-u_{j,k}>0\text{, then }{ \zeta_{j,k}=u_{j,k+1}-u_{j,k} \text{ and } \iota_{j,k}= 0,}\\
    \textbf{Condition 2:}& \text{ If } u_{j,k+1}-u_{j,k}\leq0\text{, then }{ \zeta_{j,k}=0 \text{ and } \iota_{j,k}= -(u_{j,k+1}-u_{j,k})}.
\end{align*}
Employing this decomposition to problem (\ref{eqn: generaldiscretepen}) yields
\begin{equation}\label{eqn: generalpendisc}
\begin{array}{rl}
    \min &J_{\bs{\rho}}(\u_1,\bs{\zeta}_1,\bs{\iota}_1, \dots,\u_m,\bs{\zeta}_m,\bs{\iota}_m)=\sum\limits_{k=0}^{N-1}hg(\x_{\cdot,k},\u_{\cdot,k}) + \sum\limits_{j=1}^{m}\left(\rho_j\sum\limits_{k=0}^{N-2}(\zeta_{j,k}+\iota_{j,k}) \right) \\
     & x_{i,k+1} = x_{i,k}+ hf_i(\x_{\cdot,k},\u_{\cdot,k}),  \; \text{for all } i=1,\dots, n \text { and for all } k=0,\dots, N-1,\\
     &\xi_j\leq u_{j,k} \leq \omega_j, \; \text{ for all } j=1,\dots, m\; \text{ and for all } k=0,\dots, N-1,\\
     & u_{j,k+1}-u_{j,k}=\zeta_{j,k}-\iota_{j,k},\;\text{ for all } j=1,\dots, m \text{ and for all } k=0,\dots, N-2,\\
     & 0\leq\bs{\zeta}_j, \; 0\leq \bs{\iota}_j \text{ for all } j=1,\dots, m.
\end{array}
\end{equation}
For the above problem, we are now minimizing $J_{\bs{\rho}}$ with respect to vectors $\u_1,\bs{\zeta}_1,\bs{\iota}_1, \dots,\u_m,\bs{\zeta}_m,$ and $\bs{\iota}_m$. 
The constraints associated with each $\bs{\zeta}_j$ and $\bs{\iota}_j$, are constraints that PASA can interpret. 
For all $j=1,\dots,m$, the equality constraints associated with each $\bs{\zeta}_j$ and $\bs{\iota}_j$ in problem (\ref{eqn: generalpendisc}) are linear and can be expressed as: 
\begin{equation}
  \left[
      \begin{array}{c | c | c}
      \bs{A_j} &-\bs{I}_{N-1}&\bs{I}_{N-1}
      \end{array}
   \right]
   \begin{bmatrix}
   \u_j\\
   \hline
   \bs{\zeta}_j\\
   \hline
   \bs{\iota}_j
   \end{bmatrix}
   =
   \bs{0},
  \end{equation}
  where $\bs{I}_{N-1}$ is the identity matrix of dimension $N-1$, $\bs{0}$ is an $N-1$ dimensional all zeros vector, and $\bs{A_j}$ is an $N-1$ by $N$ matrix defined as 
  \begin{equation}\label{eqn: sparsematrix}
\bs{A_j} = 
\begin{bmatrix}
-1 	& 1         & 0 		& \cdots & 0 \\
 0 	& \ddots & \ddots    & \ddots & \vdots\\
\vdots &\ddots &  \ddots   & \ddots & 0\\
0	& \cdots& 0 		& -1	&   1
\end{bmatrix}.
\end{equation}
Moreover, the equality constraints of the decomposition vectors given in problem (\ref{eqn: generalpendisc}) can be written as 
\begin{equation}
  \left[
      \begin{array}{c | c | c | c | c |c | c|c|c|c|c}
      \bs{A_1} &-\bs{I}_{N-1}&\bs{I}_{N-1} & \cdots &\bs{A_j}& -\bs{I}_{N-1}&\bs{I}_{N-1}&\cdots&\bs{A_{m}} & -\bs{I}_{N-1}&\bs{I}_{N-1}
      \end{array}
   \right]
   \begin{bmatrix}
   \u_1\\
   \hline
   \bs{\zeta}_1\\
   \hline
   \bs{\iota}_1\\
   \hline
   \vdots\\
   \hline
   \u_j\\
   \hline
   \bs{\zeta}_j\\
   \hline
   \bs{\iota}_j\\
   \hline
   \vdots\\
   \hline
   \u_m\\
   \hline
   \bs{\zeta}_m\\
   \hline
   \bs{\iota}_m
   \end{bmatrix}
   =
   \bs{0},
  \end{equation}
  where $\bs{A_j}$ is defined as (\ref{eqn: sparsematrix}) for all $j=1,\dots, m$ and $\bs{0}\in\R^{N-1}$.

From problem (\ref{eqn: generalpendisc}), we wish to find the gradient of  $J_{\p}$ with respect to $\u_1,\bs{\zeta}_1,\bs{\iota}_1, \dots,\u_m,\bs{\zeta}_m,$ and $\bs{\iota}_m$. 
We will compute the gradient of the Lagrangian of Problem (\ref{eqn: generalpendisc}) with respect to $\u_1,\bs{\zeta}_1,\bs{\iota}_1, \dots,\u_m,\bs{\zeta}_m,$ and $\bs{\iota}_m$ to find $\nabla J_{\bs{\rho}}$. 
This is necessary because based upon problem (\ref{eqn: generalpendisc}), state variables $\x_1,\dots, \x_n$ can be viewed as functions of $\u_1,\dots, \u_m$. 
So when computing the gradient of $J_{\p}$ with respect to $\u_1,\bs{\zeta}_1,\bs{\iota}_1, \dots,\u_m,\bs{\zeta}_m,$ and $\bs{\iota}_m$, we should consider that vectors $\x_1,\dots, \x_m$ depend on the controls. 

The discretized problem can be put in the following general form for which the technique of using the Lagrangian to compute the gradient of a cost functional can be used: 
\begin{equation}\label{eqn: optproblem}
{ \begin{array}{rl}
     \min & G(\x,\u)  \\
      & \bs{F}(\x,\u)=\bf{0},
 \end{array}}
\end{equation}
where $\u\in \R^m$, $\x\in\R^n$, $G:\R^n\times\R^m\to\R$, and $\textbf{F}:\R^n\times\R^m\to\R^n$ are differentiable; moreover, it is assumed in problem (\ref{eqn: optproblem}) that we can uniquely solve for $\x$ in term of $\u$. 
Based on these assumptions, we can rewrite problem (\ref{eqn: optproblem}) as 
\begin{equation}\label{eqn: 4lagrangethm}
 {  
 \begin{array}{rl}
      \min   & \mathscr{J}(\u)=G(\x(\u),\u) \\
         & \bs{F}(\x(\u),\u)=\bf{0}, 
    \end{array}
}
\end{equation}
where $\x=\x(\u)$ denotes the unique solution of $\bs{F}(\x,\u)=\bf{0}$ for a given $\u\in\R^m$. 
The following result can be deduced from the implicit function theorem and the chain rule
(see \cite[Remark 3.2]{Hager2000}, \cite{Hager1987}):
  \begin{theorem}\label{thm: lagrange}
   If the Jacobian $\nabla_{\x}\bs{F}(\x,\u)$ is invertible for each $\u\in\R^m$ and $\x=\x(\u)$, then for each $\u\in\R^m$ the gradient of $\mathscr{J}$ in problem (\ref{eqn: 4lagrangethm}) is 
   \begin{equation}\label{eqn: thmgrad}
       \nabla_{\u}\mathscr{J}(\u)=\nabla_{\u} \Ls(\x,\u,\lam)_{|_{\x=\x(\u)}},
   \end{equation}
   where $\Ls:\R^n\times\R^m\times\R^n\to \R$ is the Lagrangian of problem (\ref{eqn: 4lagrangethm}), 
   \begin{equation*}
       \Ls(\x,\u,\lam)=G(\x,\u)+\lam^T\bs{F}(\x,\u),
   \end{equation*}
   and $\lam$ is chosen such that 
   \begin{equation}\label{eqn: thmadjoint}
       \grad_{\x}\Ls(\x,\u,\lam)=\grad_{\x}G(\x,\u)+\lam^T\grad_{\x}\bs{F}(\x,\u)=\bf{0}.
   \end{equation}
  \end{theorem}
  
   Before computing the Lagrangian to problem (\ref{eqn: generalpendisc}) we will rewrite the state equations accordingly: 
   \[
   -x_{i,k+1}+x_{i,k} + hf_i(\x_{\cdot,k},\u_{\cdot,k})=0,\; \text{ for all } i=1,\dots,n \text{ and for all }  k=0,\dots,N-1.  
   \]
   The Lagrangian to problem (\ref{eqn: generalpendisc}) will then be 
   \begin{align}\label{eqn: lagrangegeneral}
       \Ls(\u_1,\bs{\zeta}_1,\bs{\iota}_1,\dots, \u_m, \bs{\zeta}_m,\bs{\iota}_m)=& \sum\limits_{k=0}^{N-1}hg(\x_{\bs{\cdot},k},\u_{\bs{\cdot},k}) + \sum\limits_{j=1}^{m}\rho_j\sum\limits_{k=0}^{N-2}(\zeta_{j,k}+\iota_{j,k})\nonumber \\
       &\;\;\;+\sum\limits_{i=1}^{n}\sum\limits_{k=0}^{N-1}\lambda_{i,k}(-x_{i,k+1}+x_{i,k}+hf_i(\x_{\bs{\cdot},k},\u_{\bs{\cdot},k})),
   \end{align}
   where $\lam_1,\dots, \lam_n\in\R^{N-1}$ are the Lagrange multiplier vectors.   
   For all $j=1,\dots, m$ and for all $k=0,\dots, N-1$, we compute the partial derivative of $\Ls$ with respect to $u_{j,k}$ and obtain
   \begin{align*}
       \frac{\partial }{\partial u_{j,k}}\Ls &= h\frac{\partial}{\partial u_{j,k}}g(\x_{\bs{\cdot},k},\u_{\bs{\cdot},k})+\sum\limits_{i=1}^n \lambda_{i,k}\left(h\frac{\partial}{\partial u_{j,k}}f_i(\x_{\bs{\cdot},k},\u_{\bs{\cdot},k})\right) .
   \end{align*}
   For all $j=1,\dots,m$, we have that  
   \begin{equation}\label{eqn: gradLuj}
       \bs{\grad_{\u_j} \Ls} = 
       \begin{bmatrix}
        h\frac{\partial}{\partial u_{j,0}}g(\x_{\bs{\cdot},0},\u_{\bs{\cdot},0})+\sum\limits_{i=1}^n \lambda_{i,0}\left(h\frac{\partial}{\partial u_{j,0}}f_i(\x_{\bs{\cdot},0},\u_{\bs{\cdot},0})\right)\\
        \vdots\\
         h\frac{\partial}{\partial u_{j,k}}g(\x_{\bs{\cdot},k},\u_{\bs{\cdot},k})+\sum\limits_{i=1}^n \lambda_{i,k}\left(h\frac{\partial}{\partial u_{j,k}}f_i(\x_{\bs{\cdot},k},\u_{\bs{\cdot},k})\right)\\
         \vdots \\
          h\frac{\partial}{\partial u_{j,N-1}}g(\x_{\bs{\cdot},N-1},\u_{\bs{\cdot},N-1})+\sum\limits_{i=1}^n \lambda_{i,N-1}\left(h\frac{\partial}{\partial u_{j,N-1}}f_i(\x_{\bs{\cdot},N-1},\u_{\bs{\cdot},N-1})\right),
       \end{bmatrix}
   \end{equation}
   where $\bs{\grad_{\u_j}\Ls}\in\R^{N}$.
   For all $j=1,\dots, m$ and for all $k=0,\dots, N-2$, we compute the partial derivative of $\Ls$ with respect to $\zeta_{j,k}$ and the partial derivative of $\Ls$ with respect to $\iota_{j,k}$: 
   \begin{align*}
       \frac{\partial }{\partial \zeta_{j,k}} \Ls&= \rho_j, \\
       \frac{\partial }{\partial \iota_{j,k}} \Ls &= \rho_j.
   \end{align*}
   We can then say for all $j=1,\dots, m$ that 
   \begin{equation}\label{eqn: gradLzj}
     \bs{ \grad_{\bs{\zeta}_j}\Ls}= 
      \begin{bmatrix}
      \rho_j\\
      \vdots\\
      \rho_j
      \end{bmatrix},
   \end{equation}
   and 
    \begin{equation}\label{eqn: gradLij}
      \bs{\grad_{\bs{\iota}_j}\Ls}= 
      \begin{bmatrix}
      \rho_j\\
      \vdots\\
      \rho_j
      \end{bmatrix}
   \end{equation}
   where $\bs{\grad_{\bs{\zeta}_j}\Ls}\in\R^{N-2}$ and $\bs{\grad_{\bs{\iota}_j}\Ls}\in\R^{N-2}$. 
   By Theorem \ref{thm: lagrange}, provided that Lagrange multiplier vectors $\lam_1,\dots, \lam_n$ satisfy theorem condition (\ref{eqn: thmadjoint}), we have that  
   \begin{equation}
       \bs{\grad J} =\bs{\grad \Ls} =
    \begin{bmatrix}
    \bs{\grad_{\u_1}\Ls}\\
   \bs{ \grad_{\bs{\zeta}_1}\Ls}\\
    \bs{\grad_{\bs{\iota}_1}\Ls}\\
    \hline
    \vdots\\
    \hline
     \bs{\grad_{\u_k}\Ls}\\
    \bs{\grad_{\bs{\zeta}_k}\Ls}\\
   \bs{ \grad_{\bs{\iota}_k}\Ls}\\
    \hline
    \vdots\\
    \hline
    \bs{\grad_{\u_{m}}\Ls}\\
    \bs{\grad_{\bs{\zeta}_{m}}\Ls}\\
    \bs{\grad_{\bs{\iota}_m}\Ls}\\
    \end{bmatrix},
   \end{equation}
    where entries of $\bs{\grad\Ls}$ are defined in equations (\ref{eqn: gradLuj})-(\ref{eqn: gradLzj}).
    Conveniently, our method of finding vectors $\lam_1,\dots,\lam_n$ that satisfy condition  (\ref{eqn: thmadjoint}) will simultaneously produce a discretization of the adjoint equations (\ref{eqn: adjoint}).
    For all $i=1,\dots, n$ and for all $k=1,\dots, N$ we compute the partial derivative of $\Ls$ with respect to $x_{i,k}$:
    \begin{align*}
        \frac{\partial}{\partial x_{i,k}}\Ls &=h\frac{\partial }{\partial x_{i,k}} g(\x_{\bs{\cdot},k},\u_{\bs{\cdot},k}) -\lambda_{i,k-1} +\lambda_{i,k}+\sum\limits_{\ell=1}^{n}\lambda_{\ell,k}\left(h\frac{\partial }{\partial x_{i,k}}f_{\ell}(\x_{\bs{\cdot},k}, \u_{\bs{\cdot},k})\right)\text{ for } k=1,\dots,N-1,\\
        \frac{\partial}{\partial x_{i,N}}\Ls&= -\lambda_{i,N-1}.
    \end{align*}
    We did not take the partial derivative of $\Ls$ with respect to $x_{i,0}$ since $x_{i,0}$ is a known value for all $i=1,\dots,n$.
   
    To satisfy theorem condition (\ref{eqn: thmadjoint}),  we set $\frac{\partial}{\partial x_{i, k}}\Ls$ equal to zero for all $i=1,\dots,n$ and for all $k=1,\dots,N,$ and solve for $\lambda_{i,k-1}$. 
    We obtain the following for all $i=1,\dots,n$
    \begin{align}
        \lambda_{i,k-1}&= h\frac{\partial }{\partial x_{i,k}} g(\x_{\bs{\cdot},k},\u_{\bs{\cdot},k}) +\lambda_{i,k}+\sum\limits_{\ell=1}^{n}\lambda_{\ell,k}\left(h\frac{\partial }{\partial x_{i,k}}f_{\ell}(\x_{\bs{\cdot},k}, \u_{\bs{\cdot},k})\right)\text{ for } k=1,\dots,N-1,\text{ and }\\
        \lambda_{i,N-1}&= 0\label{eqn: discadjointN}. 
    \end{align}
    The above equations not only allow us to use the gradient of the Lagrangian of problem (\ref{eqn: generaldiscretepen}) to find the gradient of $J_{\bs{\rho}}$, but also give us a discretization for adjoint equations (\ref{eqn: adjoint}). 
    Additionally for all $i=1,\dots,n$, equation (\ref{eqn: discadjointN}) will be analogous to the transversality conditions given in (\ref{eqn: transversality}).

\section{Example 1: Fishery Problem}
{In this section, we focus on a basic resource model on harvesting, which was first presented in Clark's \cite{Clark}. We present the fishery problem as stated in  Lenhart and Workman  \cite[Example 17.4]{Lenhart2007}, where a logistic growth function was utilized within the fishery model. }
\begin{equation}\label{eqn: maxharvest}
    \begin{array}{rl}
     \max\limits_{u\in\mathscr{A}} & J(u) =\int\limits_0^T (pqx(t)-c)u(t) dt\\
   \textrm{s.t.} & x'(t) =x(t)(1-x(t))-qu(t)x(t),\\
    &x(0)=x_0>0,\\
    & 0\leq u(t)\leq M
    \end{array}
\end{equation}
In problem (\ref{eqn: maxharvest}), $u(t)$ is the control variable that measures the effort put into harvesting fish at time $t$, and $x(t)$ measures the total population of fish at time $t$.
We are assuming that there is a maximum harvesting rate, $M$. 
Parameter $q$ represents the ``catchability'' of a fish and parameter $c$ represents the cost of harvesting one unit of fish. 
Parameter $p$ is the selling price for one unit of harvested fish. 
The objective function $J$ is constructed to represent the total profit, revenue less cost, of harvesting fish over time interval $[0,T]$.
The set of admissible controls, $\mathscr{A}$, for problem (\ref{eqn: maxharvest}) will be defined as follows:
\[\mathscr{A}= \{u\in L^1(0,T):  \text{ for all } t\in[0,T], \; u(t)\in A, \; A=[0,M]\}.\] 
Existence of an optimal control for problem (\ref{eqn: maxharvest}) follows from Filippov-Cesari Existence Theorem \cite{Filippov}.

According to Lenhart and Workman \cite{Lenhart2007}, the standard forward-backward sweep method will not converge if parameters in problem (\ref{eqn: maxharvest}) were set to where the optimal control $u^*$ contained a singular region. 
 {In this section, we present an explicit formula for the singular case, which was obtained via Pontryagin's minimum principle \cite{Pontryagin}, and  we show that the singular case satisfies the generalized Legendre-Clebsch Condition\cite{Powers, Robbins,Lewis, Zelikin1994}.}
 {Additionally, we present a set of assumptions on the parameters to the problem in order to obtain an optimal harvesting strategy that begins singular and switches to the maximum harvesting rate. 
 For this scenario, an explicit formula for the switching point is obtained.
With parameters set to meet those particular assumptions, we will use the explicit solution to problem (\ref{eqn: maxharvest}) to test PASA's accuracy in  solving for the regularized problem for varying values of the tuning parameter.}
We also discuss how to discretize for both the fishery problem and the regularized version of the problem in which a bounded variation regularization term is applied.
And finally, we present some empirical evidence for convergence between the numerical solution obtained by PASA for the regularized fishery problem and the exact solution to the fishery problem. 
\subsection{Explicitly Solving Fishery Problem}
  In problem \cite{Lenhart2007}, Lenhart and Workman demonstrate a method of using Pontryagin's maximum principle \cite{Pontryagin} and properties of the switching function to solve for the singular case to problem (\ref{eqn: maxharvest}). Lenhart and Workman also discuss conditions for existence of a singular case. 
 Since we will be using a numerical solver that is used for solving minimization problems, we provide an analytical solution to the minimization problem that is equivalent to problem (\ref{eqn: maxharvest}). 
 The equivalent minimization problem is obtained by negating the objective functional $J(u)$ in problem (\ref{eqn: maxharvest}): 
 \begin{equation}\label{eqn: minprofit}
    \begin{array}{rl}
     \min\limits_{u} & J(u) =\int\limits_0^T -(pqx(t)-c)u(t) dt\\
   \textrm{s.t.} & x'(t) =x(t)(1-x(t))-qu(t)x(t),\\
    &x(0)=x_0>0,\\
    & 0\leq u(t)\leq M.
    \end{array}
\end{equation}
 Our procedure for solving for problem (\ref{eqn: minprofit}) will be analogous to \cite{Lenhart2007}'s methods. 
  We will be using Pontryagin's Minimum Principle \cite{Pontryagin} to solve for problem (\ref{eqn: minprofit}). 
 The Hamiltonian for the above problem is 
 \begin{equation}\label{eqn: profithamiltonian}
     H(x,u, \lambda) = (c-pqx)u+\lambda(x-x^2-qux),
 \end{equation}
 where $\lambda$ is the adjoint variable. 
 By taking the partial derivative of the Hamiltonian with respect to state variable $x$ we obtain the adjoint equation associated with adjoint variable $\lambda$ which is 
 \begin{equation}\label{eqn: profitadjoint}
     \lambda'(t) = -\frac{\partial H}{\partial x} = pqu-\lambda +2\lambda x+q\lambda u, 
 \end{equation}
 with the transversality condition being 
 \begin{equation}\label{eqn: profitadjointTC}
     \lambda(T) = 0.
 \end{equation}
 We also use the Hamiltonian given in (\ref{eqn: profithamiltonian}) to compute the switching function corresponding to problem (\ref{eqn: minprofit})
 \begin{equation}\label{eqn: profitswitchfunction}
     \psi(t) =\frac{\partial H}{\partial u} = c-pqx(t)-q\lambda(t) x(t).
 \end{equation}
 Based on Pontryagin's Minimum Principle \cite{Pontryagin}, if there exists an optimal pair $(u^*,x^*)$ for problem (\ref{eqn: minprofit}) then there exists 
 $\lambda^*$, satisfying adjoint equation (\ref{eqn: profitadjoint}) and terminal condition (\ref{eqn: profitadjointTC}), where $H(x^*,u^*, \lambda^*)\leq H(x^*, u, \lambda^*)$ for all admissible controls $u$.
 Additionally, if $u^*$ is the optimal control then $u^*$ must have the following form:   
 \begin{equation}\label{eqn: profitupsi}
     u^*(t) = 
     \begin{cases}
     0 & \textrm{ whenever }\psi(t) >0,\\
     \textrm{singular} & \textrm{ whenever } \psi(t) = 0,  \\
     M & \textrm{ whenever } \psi(t) <0.
     \end{cases}
 \end{equation}
 To solve for the singular case we suppose $\psi(t) \equiv 0$ on some subinterval $I\subset [0,T]$. 
 Assuming that parameter $c>0$, then we have from the switching function given in (\ref{eqn: profitswitchfunction}) that both $x$ and $pq+q\lambda$ are nonzero on interval $I$. 
 By setting the $\psi$ equal to zero and solving for $\lambda$ we have that 
 \begin{equation}\label{eqn: profitadjoint1}
     \lambda^*(t) = \frac{c-pqx(t)}{qx(t)}
 \end{equation}
 on the interval $I$. 
 We differentiate equation (\ref{eqn: profitadjoint1}) and use the state equation given in problem (\ref{eqn: minprofit}) to obtain the following:
 \begin{equation}\label{eqn: profitadjoint2}
 \begin{split}
     \lambda'(t) &=-\frac{c}{qx^2}(x(1-x)-qux)\\
     &=-\frac{c}{qx}+\frac{c}{q}+\frac{cu}{x}.
  \end{split}
 \end{equation}
 
 We rewrite equation (\ref{eqn: profitadjoint}) by using expression (\ref{eqn: profitadjoint1}), and we get 
 \begin{equation}\label{eqn: profitadjoint3}
 \lambda'(t) =-\frac{c}{qx}+p+\frac{2c}{q}-2px+\frac{cu}{x}.
 \end{equation}
 Equating expressions (\ref{eqn: profitadjoint2}) and (\ref{eqn: profitadjoint3}) gives us a solution for $x$ on the singular interval which is 
 \begin{equation}\label{eqn: profitxsing}
     x^*(t) = \frac{c+pq}{2pq}.
 \end{equation}

 Since $x^*$ is constant on $I$ we have that $x'(t) =0$
 on $I$. 
 We use the singular solution for $x^*$ and set the state equation found in problem (\ref{eqn: minprofit}) equal to zero to solve for $u^*$. 
 On the interval $I$ we have that 
 \begin{align*}
     0 &= x'(t) =x^*(1-x^*)-qux^*\\
     u^*(t) &= \frac{1-x^*}{q}\\
     u^*(t) &= \frac{pq-c}{2pq^2}.
 \end{align*}
 Additionally, $x^*$ being constant on the singular region and expression (\ref{eqn: profitadjoint1}) would imply that $\lambda^*$ is also constant on the singular region with constant value being as follows
 \begin{equation*}
     \lambda^*(t)=\frac{p(c-pq)}{c+pq}.
 \end{equation*}
 By substituting in the constant solutions for $u^*$, $x^*$, and $\lambda^*$ into adjoint equation (\ref{eqn: profitadjoint}), we have that the right hand side of the adjoint equation is zero, as desired. 
 To conclude, we found that if a singular region, $I$, exists, then $u^*$, $x^*$, and $\lambda^*$ are all constant on $I$ where 
 \begin{equation}\label{eqn: profitsing}
     u^* = \frac{pq-c}{2pq^2},
     \quad x^*  = \frac{c+pq}{2pq},
     \quad \lambda^* = \frac{p(c-pq)}{c+pq}. 
\end{equation}
{
We wish to show that the singular case solution satisfies the second order necessary condition of optimality, which is referred as the generalized Legendre-Clebsch Condition \cite{Powers, Robbins} or Kelley's condition  \cite{Zelikin1994, Lewis, Kelley}).
 Before showing that the singular cases given in (\ref{eqn: profitsing}) satisfies the Legendre-Clebsch Condition and/or Kelley's condition \cite{Powers, Robbins, Zelikin1994, Lewis, Kelley}, we would like to present some parameter assumptions on problem (\ref{eqn: minprofit}).
}
{\begin{assumption}\label{as: assump1}
	Parameters $p, c, M, q >0$ are set to satisfy $0<pq-c<2pq^2M$.
\end{assumption}
\begin{assumption}\label{as: assump2}
	Initial value $x_0$  is set to equal $\frac{c+pq}{2pq}$, which is the constant value that is associated to the state solution corresponding to singular $u^*$. 
\end{assumption}}
{
Note that Assumption \ref{as: assump1} implies that the singular case solution $u^*=\frac{pq-c}{2pq^2}$ satisfies the boundary constraints that are assumed on the control. 
Assumption \ref{as: assump2}  dynamically forces the problem to yield an optimal control that begins singular. 
 The generalized Legendre-Clebsch Condition involves finding what is called the \emph{order} of a singular arc which is defined as the being the integer $q$ such that $\left(\frac{d^{2q}}{dt^{2q}}\frac{\partial H}{\partial u}\right)$ is the lowest order total derivative of the partial derivative of the Hamiltonian with respect to $u$, in which control $u$ appears explicitly. 
 We use the state equations given in problem (\ref{eqn: minprofit}) and the adjoint equations given in (\ref{eqn: profitadjoint}) to find the first and second time derivative of the switching function (\ref{eqn: profitswitchfunction}):
 \begin{align}
 \frac{d}{dt}\psi &= pq (x^2-x)-q\lambda x^2, \label{eqn: profitsdswitch}\\
 \frac{d^2}{dt^2}\psi &= pq(2x-1)(x-x^2)-q\lambda x^2+[pq^2x+q^2\lambda x^2-3pq^2x^2]u.\label{eqn: profitsdswitch2}
 \end{align}
 From above, we have that the order of the singular arc is $q=1$. 
 We need to show that if $u^*$ is an optimal singular control on some interval of order $q$, then it is necessary that 
 \begin{equation}\label{eqn: legendre}
\left . (-1)^{q}\frac{\partial}{\partial u }\left[\frac{d^2}{dt^2}\left(\frac{\partial H}{\partial u}\right)\right]\right |_{x=x^*,\lambda=\lambda^*}\geq 0. 
 \end{equation}
 We take the partial derivative of (\ref{eqn: profitsdswitch2}) with respect to $u$ and evaluate at the corresponding singular case solutions for $x$ and $\lambda$ given in (\ref{eqn: profitsing}): 
\begin{equation*}
\left .\frac{\partial}{\partial u }\left[\frac{d^2}{dt^2}\left(\frac{\partial H}{\partial u}\right)\right]\right |_{x=\frac{c+pq}{2pq},\lambda=\frac{p(c-pq)}{2pq^2}}=\frac{q(c+pq)}{2}-\frac{3(c+pq)^2}{4p}-\frac{(pq-c)(c+pq)^2}{8p^2q^2}.
\end{equation*}
Note that the term $\frac{(pq-c)(c+pq)^2}{8p^2q^2}$ is positive by Assumption \ref{as: assump1}. 
Using algebra, we combine the first two terms in the above equation to obtain
\begin{equation*}
\left .\frac{\partial}{\partial u }\left[\frac{d^2}{dt^2}\left(\frac{\partial H}{\partial u}\right)\right]\right |_{x=\frac{c+pq}{2pq},\lambda=\frac{p(c-pq)}{2pq^2}}=-\frac{3c^2}{4p}-cq-\frac{pq^2}{4}-\frac{(pq-c)(c+pq)^2}{8p^2q^2}<0.
\end{equation*}
}
By multiplying the above inequality by $(-1)^q$ where $q=1$,  we then have the second order necessary condition of optimality (\ref{eqn: legendre}) being satisfied.

{Notice that only Assumption \ref{as: assump1} is needed in proving that the singular case satisfies the Legendre-Clebsch condition; however, we can use both assumptions to obtain a control that begins singular and switches to the maximum harvesting rate, where an explicit formula for the switching point can be obtained. 
This particular scenario for problem (\ref{eqn: minprofit}), makes it an excellent candidate problem to use for testing PASA's accuracy in solving for the regularized variant of this problem.}
We first verify Assumptions \ref{as: assump1} and \ref{as: assump2} imply that $u^*$ will not be singular on the entire time interval $[0,T]$. 
 By using the transversality condition, $\lambda(T)=0$, we recognize that the singular case solution for $\lambda$ given in (\ref{eqn: profitadjoint1}) is 0 if and only if parameters are set to either satisfy $p=0$ or $c-pq=0$, and Assumption \ref{as: assump1} ensures that both cases are not possible. 
 Let $0<t^*<T$ be the time when $u^*$ switches from being singular to non-singular, and let $I=[0,t^*)$ be the interval corresponding to when the optimal control $u^*$ is singular.
 By looking at the objective functional for problem (\ref{eqn: minprofit}), intuition tells us that $u^*\neq 0$ on $[t^*,T]$, and  we can prove this by proof by contradiction. 
 Assume that $u^*(t)=\begin{cases}
                \frac{pq-c}{2pq^2} & 0\leq t<t^*\\
                0 & t^*\leq t\leq T
            \end{cases}$
 is the optimal control to problem (\ref{eqn: minprofit}). 
 Consider the following admissible control $\hat{v}$ where $\hat{v} $ is singular over the entire interval, i.e. 
 \[\hat{v}(t) = \frac{pq-c}{2pq^2} \; \quad \text{ for all } t\in[0,T].
     \]
 By assumption of $u^*$ being optimal for problem (\ref{eqn: minprofit}), we have $J(\hat{v})\geq J(u^*)$. 
    Since $\hat{v}=u^*$ on the interval $[0,t^*]$ and $u^*\equiv 0$ on interval $[t^*,T]$, we obtain the following: 
    \[J(\hat{v})=J(u^*)-\int\limits_{t^*}^{T}((pqx_{\hat{v}}(t)-c)\hat{v}(t))dt, \]
    where $x_{\hat{v}}(t)$ is the corresponding solution to the state equation given in problem (\ref{eqn: minprofit}). 
    A contradiction is obtained if we can show that $pqx_{\hat{v}}-c>0$ for all $t\in[t^*,T]$ because the following implies $J(\hat{v})\leq J(u^*)$. 
    Since $\hat{v}$ is singular on the interval $[t^*,T]$, we have that $x_{\hat{v}}$ will be the corresponding singular case solution given in (\ref{eqn: profitsing}).
    We then have that 
    \begin{align*}
        pqx_{\hat{v}}-c &= pq\left(\frac{c+pq}{2pq}\right)-c\\
                        &= \frac{1}{2}(pq-c)\\
                        &>0,
    \end{align*}
    where the above inequality holds by Assumption \ref{as: assump1}. 
    Therefore, we have our contradiction. 
   
 It then follows that the optimal harvesting policy for problem (\ref{eqn: minprofit}) with parameters set to satisfying Assumptions \ref{as: assump1} and \ref{as: assump2} will be a control that begins singular and switches once to the maximal harvesting effort, meaning $u^*\equiv M$ on the interval $I=[t^*,T]$. 
 We can find an explicit expression for $t^*$. 
 Assume that $u^*(t) = M$ for all $t\in [t^*,T]$, then the adjoint equation (\ref{eqn: profitadjoint}) becomes 
 \begin{equation}\label{eqn: profitadjointeqM}
     \lambda'(t) = pqM -\lambda +2\lambda x +q\lambda M. 
 \end{equation}
 Now, $\lambda$ is continuous at $t^*$. 
 Hence, $\lambda(t^*)= \lim\limits_{t\to t^{*^-}}\lambda(t)=\frac{p(c-pq)}{c+pc}$, which is the solution for equation (\ref{eqn: profitadjoint}) when $u$ is singular. 
 We can then solve for equation (\ref{eqn: profitadjointeqM}) to where $\lambda$ must satisfy the terminal condition, $\lambda(T)=0$, and the condition that $\lambda(t^*)=\frac{p(c-pq)}{c+pc}$.
 The condition that $\lambda(t^*)=\frac{p(c-pq)}{c+pq}$ allows us to obtain an explicit solution for $t^*$.
 When solving differential equation (\ref{eqn: profitadjointeqM}), we will use a standard method for solving linear differential equations.
 We will rewrite equation (\ref{eqn: profitadjointeqM}) as
 \begin{equation}\label{eqn: lambdapqM}
     \lambda'(t) +\lambda(t)(\alpha-2x(t))=pqM,
 \end{equation}
 where $\alpha=1-qM$. 
{ We let $v(t)$ to be the appropriate integrating factor which is defined as}
 \begin{equation}\label{eqn: profitv}
     v(t) = e^{\int_{t^*}^t(\alpha-2x(\tau))d\tau}.
 \end{equation}
{Multiplying $v(t)$ to both sides of equation (\ref{eqn: lambdapqM}), integrating over the interval $(t^*,t)$, and rearranging terms yields: }
\begin{equation}\label{eqn: profitadjoint5}
    \lambda(t) =\frac{1}{v(t)}\left[v(t^*)\lambda(t^*)+pqM\int_{t^*}^tv(\tau)d\tau\right].
\end{equation}
Now in order to evaluate $v(t)$ and $\int_{t^*}^tv(\tau)d\tau$, we will need to obtain an explicit solution for $x(t)$ over the interval $[t^*,T]$. 

Given that $u^*(t)=M$ for all $t\in [t^*,T]$, the state equation given in problem (\ref{eqn: minprofit}) over the specified interval becomes
 \begin{equation}\label{eqn: stateEqM}
     \frac{dx}{dt} =x(1-x)-qMx,
 \end{equation}
 which is separable. 
To solve the above equation we separate variables $x$ and $t$
\begin{align*}
    \frac{dx}{x(\alpha-x)} &= dt,
\end{align*}
where $\alpha=1-qM$. 
{We perform partial fraction decomposition on the left hand side of the above equation, integrate both sides, and exponentiate to obtain the following}
 \begin{equation}\label{eqn: stateimplicit}
     \frac{x}{|\alpha-x|}= Ke^{\alpha t},
 \end{equation}
 where $K$ is some constant. 
 We have that state variable $x$ is continuous at switching point $t^*$. 
 By continuity and Assumption \ref{as: assump2} we have $x^*(t^*)=\lim\limits_{t\to t^{*^-}}x^*(t)=x_0=\frac{c+pq}{2pq}$, which is the state solution value associated with the singular case. 
 We can use $x^*(t^*)=x_0=\frac{c+pq}{2pq}$ to solve for constant value $K$ found in equation (\ref{eqn: stateimplicit}). 
 Evaluating equation (\ref{eqn: stateimplicit}) at $t=t^*$ yields the following:
\begin{align}\label{eqn: profitK}
K= \frac{x_0e^{-\alpha t^*}}{|\alpha-x_0|}=\frac{(c+pq)e^{-\alpha t^*}}{\gamma},
\end{align}
where $\gamma=|2\alpha pq-c-pq|$. 
 {We obtain an explicit solution of equation (\ref{eqn: stateEqM}), by rewriting equation (\ref{eqn: stateimplicit}) as }

\begin{equation*}
\frac{x}{\alpha -x} =\hat{K}e^{\alpha t} \quad \text{ where }  \hat{K}= \begin{cases}
														K	&\alpha-x\geq 0,\\
														-K 	& \alpha-x<0
													\end{cases}.
\end{equation*}

{Solving for the above equation yields }
 {\begin{equation}\label{eqn: stateExplicit}
     x(t)= \frac{\alpha \hat{K} e^{\alpha t}}{ 1+\hat{K}e^{\alpha t}}\quad \text{ for } t\in[t^*,T].
 \end{equation}
 Note that the value of $\hat{K}$ seems to depend on the sign of $\alpha-x(t)$.
 However, we can use continuity of $x$ and the structure of the solutions for $x$ on $[t^*, T]$ to show that the value of $\hat{K}$ only depends on the sign of $\alpha-x_0$.  
 But first, we would like to prove the following
}
{
\begin{proposition}\label{thm: alphax0}
  Assumptions \ref{as: assump1} and \ref{as: assump2} imply $\alpha-x_0<0$ where $\alpha=1-qM$. 
 \end{proposition}
}
{
\begin{proof}
 By assumption \ref{as: assump1}, parameters $c,p,q,M>0$ are chosen to satisfy $0<pq-c$ and $pq-c<2pq^2M$. 
Dividing both sides of the second inequality by $2pq$ yields $\frac{pq-c}{2pq}<qM$. 
Negating the inequality and adding one to both sides yields
$\alpha<1+\frac{c-pq}{2pq}=\frac{c+pq}{2pq}$, and the right hand side of the inequality is $x_0$, by Assumption \ref{as: assump2}. 
Therefore, we have $\alpha-x_0<0$. 
 \end{proof}
}
 {
Recall that by Assumption \ref{as: assump2} and equations (\ref{eqn: profitsing}) and (\ref{eqn: stateExplicit}), we have the following solution for the state variable:  
 \begin{equation}\label{eqn: profitstatekhat}
 x(t)= 
 \begin{cases}
 		x_0=\frac{c+pq}{2pq} 	& 0\leq t\leq t^*,\\
		\frac{\alpha \hat{K}e^{\alpha t}}{1+\hat{K}e^{\alpha t}}& t^*\leq t\leq T
\end{cases},
 \end{equation}
 where the sign of $\hat{K}$ is determined by the sign $\alpha-x(t)$.
 Since $x$ is continuous on the entire time interval and since $x$ is constant on the singular region, we have that $x(t^*)=x_0$. 
 Hence at $t=t^*$,  $\hat{K}$ will be determined by $\alpha-x(t^*)=\alpha-x_0$.
}
{
 Differentiating $x(t)$ along the non-singular region yields
 \begin{equation*}
 x'(t) =  \frac{\alpha^2\hat{K}e^{\alpha t}}{(1+\hat{K}e^{\alpha t})^2},
 \end{equation*}
 which is either strictly positive or strictly negative based upon the sign of $\hat{K}$.
 Since $\hat{K}$ is negative at $t=t^*$, there is some open interval containing $t^*$ such that the function $x$ will be a non-increasing function. 
}
{
 Note also that, 
 \begin{equation*}
 \lim\limits_{t\to\infty} \frac{\alpha \hat{K}e^{\alpha t}}{1+\hat{K}e^{\alpha t}}= \alpha, 
 \end{equation*}
 so this function has a horizontal asymptote being $x_{hor}=\alpha$ on the $tx$-plane. 
 This implies that $x(t)$ given in (\ref{eqn: profitstatekhat}) remains above the horizontal asymptote on the interval $[t^*,T]$ even though the function is non-increasing on $(t^*,T)$.
 In conclusion, the sign of $\alpha-x_0$ determines the structure of $x(t)$ on the non-singular region. 
 Using proposition \ref{thm: alphax0} and equation (\ref{eqn: profitstatekhat}) we then conclude $x(t)$ will be of the following form: 
 \begin{equation*}
  x(t)= 
 \begin{cases}
 		x_0=\frac{c+pq}{2pq} 	& 0\leq t\leq t^*,\\
		\frac{\alpha Ke^{\alpha t}}{-1+Ke^{\alpha t}}& t^*\leq t\leq T
\end{cases},
\end{equation*}
where $K$ is defined on equation (\ref{eqn: profitK}).
}

 We now use the explicit solution for variable $x(t)$  to evaluate $v(t)$ given in (\ref{eqn: profitv}):
\begin{align}
 v(t) 
      &= \exp{\left(\alpha(t-t^*)-2\int_{t^*}^t \frac{\alpha Ke^{\alpha \tau}}{- 1+Ke^{\alpha \tau}}d\tau \right)}.
\end{align}
 {We can use a $u$-substitution to evaluate the integral term in $v$ and after applying logarithm rules we obtain the following:  }
\begin{align}\label{eqn: v}
v(t) &=
\left(\frac{- 1+Ke^{\alpha t^*}}{- 1+Ke^{\alpha t}}\right)^2 e^{\alpha(t-t^*)}.
\end{align}

 To evaluate $\int_{t^*}^t v(\tau)d\tau$, we will need to use $u$-substitution method. 
  Let $\sigma(\tau)=- 1+Ke^{\alpha \tau} $, then $\int_{t^*}^t v(\tau)d\tau$ with $v$ given in (\ref{eqn: v}) becomes
 \begin{align*}
     \int_{t^*}^t v(\tau)d\tau&= \int_{t^*}^t \left(\frac{- 1+Ke^{\alpha t^*}}{- 1+Ke^{\alpha \tau}}\right)^2e^{\alpha(\tau-t^*)}d\tau\\
     &= \int_{\sigma(t^*)}^{\sigma(t)}\frac{(-1+Ke^{\alpha t^*})^2e^{-\alpha t^*}}{\alpha K \sigma^2}d\sigma \\
     &= \frac{(- 1+Ke^{\alpha t^*})^2}{\alpha K}e^{-\alpha t^* }\left[\frac{1}{- 1+Ke^{\alpha t^*}}-\frac{1}{-1+Ke^{\alpha t}}\right]\\
     &= \frac{e^{-\alpha t^*}(-1+Ke^{\alpha t^*)}}{\alpha K}\left[1-\frac{- 1+Ke^{\alpha t^*}}{- 1+Ke^{\alpha t}}\right]\\
     &= \frac{e^{-\alpha t^*}(- 1+Ke^{\alpha t^*})}{\alpha (- 1+Ke^{\alpha t})}\left[e^{\alpha t}-e^{\alpha t^*}\right]\\
    \int_{t^*}^t v(\tau)d\tau &= \frac{(-1+Ke^{\alpha t^*})}{\alpha (- 1+Ke^{\alpha t})}\left[e^{\alpha(t-t^*)}-1\right].
 \end{align*}

 Now we use equation (\ref{eqn: profitadjoint5}) to evaluate $\lambda(t)$: 
\begin{equation}\label{eqn: profitadjoint6}
     \lambda(t)=\left(\frac{- 1+Ke^{\alpha t}}{- 1+Ke^{\alpha t^*}}\right)^2e^{\alpha(t^*-t)}\left[\frac{p(c-pq)}{c+pq}+\frac{pqM}{\alpha}\left(\frac{- 1+Ke^{\alpha t^*}}{- 1+Ke^{\alpha t}}\right)\left(e^{\alpha(t-t^*)}-1\right)\right]. 
 \end{equation}
 To find $t^*$ set $\lambda(T)=0$ and solve for $t^*$, but note that the term $-1+Ke^{\alpha t}$ that is  used in equation (\ref{eqn: profitadjoint6}) must be non-zero for all $t\in[t^*,T]$, otherwise the state variable solution given in equation (\ref{eqn: stateExplicit}) would not be defined for all $t\in[t^*,T]$.  
\begin{align}
     0&=\left(\frac{- 1+Ke^{\alpha T}}{- 1+Ke^{\alpha t^*}}\right)^2e^{\alpha(t^*-T)}\left[\frac{p(c-pq)}{c+pq}+\frac{pqM}{\alpha}\left(\frac{- 1+Ke^{\alpha t^*}}{- 1+Ke^{\alpha T}}\right)\left(e^{\alpha(T-t^*)}-1\right)\right]\nonumber\\
     0&=\frac{c-pq}{c+pq} +\frac{qM}{\alpha}\left(\frac{- 1+Ke^{\alpha t^*}}{- 1+Ke^{\alpha T}}\right)e^{\alpha(T-t^*)}- \frac{qM}{\alpha}\left(\frac{- 1+Ke^{\alpha t^*}}{- 1+Ke^{\alpha T}}\right)
 \end{align}
 Multiply both sides of the above equation by $\alpha(c+pq)(- 1+Ke^{\alpha T})$ to obtain

 \begin{equation*}
     0=\alpha(c-pq)(- 1+Ke^{\alpha T})+qM(c+pq)(- 1+Ke^{\alpha t^*})e^{\alpha(T-t^*)}-qM(c+pq)(- 1+Ke^{\alpha t^*}).
 \end{equation*}
 Substituting $K=\frac{(c+pq)e^{-\alpha t^*}}{\gamma}$ into the above equation and multiplying everything by $\gamma$ yields
  \begin{equation*}
     0=\alpha(c-pq)\left(- \gamma+(c+pq)e^{\alpha (T-t^*)}\right)+qM(c+pq)(- \gamma+(c+pq))e^{\alpha(T-t^*)}-qM(c+pq)(- \gamma+(c+pq)).
 \end{equation*}
 We rearrange terms from the above equation to isolate expression $e^{\alpha t^*}$ 
 \begin{align*}
     \left[\alpha(c-pq)(c+pq)+qM(c+pq)(- \gamma +(c+pq))\right]e^{\alpha(T-t^*)}&=  \gamma \alpha(c-pq)+qM(c+pq)(- \gamma +(c+pq))\\
     e^{\alpha(T-t^*)} &= \frac{ \gamma \alpha(c-pq)+qM(c+pq)(- \gamma +(c+pq))}{(c+pq)[\alpha(c-pq)+qM(- \gamma +(c+pq))]}.
 \end{align*}
  We take the natural logarithm of both sides of the above equation and rearrange terms to find that
 \begin{equation*}
      t^* =T-\frac{1}{\alpha}\ln{\left[\frac{ \gamma \alpha(c-pq)+qM(c+pq)(- \gamma +(c+pq))}{(c+pq)[\alpha(c-pq)+qM(- \gamma+(c+pq))]}\right]},
  \end{equation*}
  where $\alpha = 1-qM$ and $\gamma = |pq-c-2pq^2M|$.
  We simplify $t^*$ more by substituting in $\alpha=1-qM$ into the above expression:
   \begin{equation}\label{eqn: switchprofit}
      t^* = T-\frac{1}{1-qM}\ln{\left[\frac{-\gamma(pq-c)-2\gamma cqM+qM(c+pq)^2}{(c+pq)[(c-pq)+2pq^2M-\gamma qM]}\right]}.
  \end{equation}
  
  To summarize,  Assumptions \ref{as: assump1} and \ref{as: assump2} imply that the optimal control, $u^*$, to problem (\ref{eqn: minprofit}) must begin singular and switch once to the non-singular case where $u^*\equiv M$ on $[t^*,T]$. 
 Additionally, the solutions for $u^*$,  $x^*$ and $\lambda^*$ are the following: 
  \begin{equation}\label{eqn: fullprofitu}
    u^*(t) =\begin{cases}
                \frac{pq-c}{2pq^2} & 0\leq t<t^*,\\
                M & t^*\leq t\leq T,
            \end{cases}
  \end{equation}
  \begin{equation}\label{eqn: fullprofitx}
      x^*(t) =\begin{cases}
                    \frac{c+pq}{2pq} & 0\leq t\leq t^*,\\
                    \frac{\alpha Ke^{\alpha t}}{- 1+Ke^{\alpha t}} & t^*\leq t\leq T,
                \end{cases}
  \end{equation}
  and 
  \begin{equation}\label{eqn: fullprofitlam}
      \lambda^*(t) = \begin{cases}
      \frac{p(c-pq)}{c+pq} & 0\leq t\leq t^*,\\
      \left(\frac{-1+Ke^{\alpha t}}{- 1+Ke^{\alpha t^*}}\right)^2e^{\alpha(t^*-t)}\left[\frac{p(c-pq)}{c+pq}+\frac{pqM}{\alpha}\left(\frac{- 1+Ke^{\alpha t^*}}{- 1+Ke^{\alpha t}}\right)(e^{\alpha(t-t^*)}-1)\right] & t^*\leq t\leq T,
                    \end{cases}
  \end{equation}
  where $\alpha=1-qM$, $K= \dfrac{x_0e^{-\alpha t^*}}{|\alpha-x_0|}=\dfrac{(c+pq)e^{-\alpha t^*}}{\gamma}$, and $\gamma=|2\alpha pq-c-pq|$.


\subsection{Discretization of Fishery Problem}\label{subsec: profitdiscrete}
 For numerically solving problem (\ref{eqn: minprofit}), we first discretize and then optimize. 
 We will be using the polyhedral active set algorithm (PASA), which was developed by Hager and Zhang \cite{Hager2016}, to find an optimal solution to the discretized problem.
 Additionally, we  need to discretize the adjoint equation associated with problem (\ref{eqn: minprofit}) which is 
 \begin{equation}\label{eqn: profitadjoint4}
     \lambda'(t) =pqu-\lambda+2\lambda x + q\lambda u, 
 \end{equation}
 with the  transversality condition being 
 \begin{equation*}
     \lambda(T)=0. 
 \end{equation*}
 
 
 For discretizing problem (\ref{eqn: minprofit}) we assume that control $u$ is constant over each mesh interval. 
 We partition time interval $[0, T]$, by using $N+1$ equally spaced nodes, $0=t_0<t_1<\cdots<t_N=T$.  
 For all $k=0,1,\dots,N$ we assume that the $x_k=x(t_k)$. For the control, we denote $u_k =u(t)$ for all $t_k\leq t< t_{k+1}$ when $k=0,\dots, N-2$ and $u_{N-1}=u(t)$ for all $t_{N-1}\leq t\leq t_N$.
 So we have $\x\in\R^{N+1}$ while $\u\in \R^{N}$.
 We use a left-rectangular integral approximation for objective function $J$ in (\ref{eqn: minprofit}), and we use forward Euler's method to approximate the state equation in (\ref{eqn: minprofit}). 
 The discretization of problem (\ref{eqn: minprofit}) is then 
 \begin{equation}\label{eqn: profitdiscrete}
     \begin{array}{rl}
         \min  & J(\u)=\sum\limits_{k=0}^{N-1} h(c-pqx_k)u_k \\
          & x_{k+1}=x_k+h(1-x_k-qu_k)x_k  \textrm{ for all } 0\leq k\leq N-1,\\
          & x_{0}>0,\\
          & 0\leq u_k\leq M \textrm{ for all } 0\leq k\leq N-1,
     \end{array}
 \end{equation}
 where $h=\frac{T}{N}$ is the mesh size and the first component of state vector, $x_0$, is set to being the initial condition associated with the state equation given in problem (\ref{eqn: minprofit}).\\
\indent Since PASA uses the gradient projection algorithm for one of its phases, we need to compute the gradient of the cost functional for problem (\ref{eqn: profitdiscrete}). 
 We use Theorem \ref{thm: lagrange} to find $\grad_{\u}J$ which requires finding the Lagrangian to problem (\ref{eqn: profitdiscrete}) and its gradient. 
 Additionally, we need to construct a Lagrange multiplier vector $\lam$ that satisfies equation (\ref{eqn: thmadjoint}). 
 Consequently, the Lagrange multiplier vector that satisfies equation (\ref{eqn: thmadjoint}) produces the numerical scheme that is used for discretizing adjoint equation (\ref{eqn: profitadjoint4}) and produces the transversality condition (\ref{eqn: profitadjointTC}). 
 To compute the Lagrangian to problem (\ref{eqn: profitdiscrete}), we first need to arrange the discretized state equations accordingly
 \begin{equation}\label{eqn: profitstatediscrete}
     -x_{k+1}+x_k+h(x_k-x_k^2-qu_kx_k)=0\; \text{for all } k=0,1,\dots, N-1.
 \end{equation}
 The Lagrangian to problem (\ref{eqn: profitdiscrete}) is 
 \begin{equation*}
     \Ls(\x,\u,\lam) = \sum\limits_{k=0}^{N-1}(h(c-pqx_k)u_k)+\sum\limits_{k=0}^{N-1} \lambda_k(-x_{k+1}+x_k+h(x_k-x_k^2-qu_kx_k)), 
 \end{equation*}
  where $\lam\in \R^{N-1}$ is the Lagrange multiplier vector. 
  Note that we need not worry about the inequality constraints associated with the bounds of the control when computing the Lagrangian to problem (\ref{eqn: profitdiscrete}) because these bounds are not being entered into cost function $J$. 
  By taking the partial derivative of $\Ls$ with respect to $u_k$, we obtain: 
  \begin{equation*}
      \frac{\partial \Ls}{\partial u_k} = h(c-pqx_k-q\lambda_kx_k)\; \textrm{for all } k = 0,\dots, N-1.
  \end{equation*}
  By Theorem \ref{thm: lagrange}, we have that 
  \begin{equation}\label{eqn: profitgradJ}
      \grad_{\u}{J}= \grad_{\u}{\Ls}=
      \begin{bmatrix}
         h(c-pqx_0-q\lambda_0x_0)\\
         \vdots\\
         h(c-pqx_k-q\lambda_kx_k)\\
         \vdots\\
         h(c-pqx_{N-1}-q\lambda_{N-1}x_{N-1})
      \end{bmatrix},
  \end{equation}
   provided that equation (\ref{eqn: thmadjoint}) is satisfied. 
   To satisfy equation (\ref{eqn: thmadjoint}) we take the partial derivative of $\Ls$ with respect to $x_k$ for all $k=1,\dots, N$, and note that we do not take the partial derivative of $\Ls$ with respect to $x_0$ because $x_0$ is a known value. 
   Taking the partial derivative of $\Ls$ with respect to the state vector components yield the following expressions:
   \begin{align}
       \frac{\partial \Ls}{\partial x_k}&= h(-pqu_k+\lambda_k(1-2x_k-qu_k))+\lambda_k-\lambda_{k-1} \; \textrm{ for all } k=1,\dots, N-1, \text{ and}  \label{eqn: profitadjointdiscrete}\\
       \frac{\partial \Ls}{x_N}& = -\lambda_{N-1}\label{eqn: profitadjointdiscreteN}
   \end{align}
   To align with equation (\ref{eqn: thmadjoint}) we set expressions (\ref{eqn: profitadjointdiscrete}) and (\ref{eqn: profitadjointdiscreteN}) equal to zero and solve for $\lambda_{k-1}$ for all $k=1,\cdots, N$: 
   \begin{align}
        \lambda_{k-1}&=\lambda_k+h(-pqu_k+\lambda_k(1-2x_k-qu_k))\; \textrm{for all } k=1,\dots, N-1, \textrm{ and}\label{eqn: profitadjointdiscrete1}\\
        \lambda_{N-1} &= 0.\label{eqn: profitadjointdiscreteN1}
   \end{align}
    Expressions (\ref{eqn: profitadjointdiscrete1}) and (\ref{eqn: profitadjointdiscreteN1}) will serve as the discretization for the costate equation (\ref{eqn: profitadjoint4}) and the transversality condition (\ref{eqn: profitadjointTC}).
   
 We use PASA to solve for the penalized version of problem (\ref{eqn: minprofit}) where the penalty applied to the problem will be a bounded variation penalty, as suggested in Capognigro et al. \cite{Caponigro2018}. 
 The penalized version of problem (\ref{eqn: minprofit}) is as follows: 
 \begin{equation}\label{eqn: profitpen}
     \begin{array}{rl}
         \min & J_{\rho}(u) = \int\limits_0^T (c-pqx)udx + \rho V(u)  \\
          & x'(t) = (1-x)x-qux, \; x(0)=x_0>0,\\
          & 0\leq u(t) \leq M, 
     \end{array}
 \end{equation}
 where $0\leq \rho <1$, is a penalty parameter and $V(u)$ measures the total variation of $u$ which is 
 \begin{equation}\label{eqn: profitV}
     V(u) = \sup\limits_{\mathcal{P}} \sum\limits_{i=0}^{n_P-1} |u(t_{i+1})-u(t_i)|
 \end{equation}
 where $\mathcal{P}=\{P=\{t_0,t_1,\dots, t_{n_P}\}: \textrm{ P is a partition of } [0,T]  \}$.
 Notice that if we numerically solve for problem (\ref{eqn: profitpen}) with $\rho=0$, then the problem is not being penalized via bounded variation. 
 We use the same procedure as before to discretize problem (\ref{eqn: profitpen}).
 Assuming control $u$ to be constant over each mesh interval allows us to express the total variation of $u$ as being the sum of the absolute value of the jumps of $u$.
 So for a sufficiently small mesh size $h$ we have 
 \begin{equation*}
     V(u) = \sum\limits_{k=0}^{N-1}|u(t_{k+1})-u(t_k)|.
 \end{equation*}
 The discretized version of problem (\ref{eqn: profitpen}) is 
 \begin{equation}\label{eqn: profitpendiscrete}
     \begin{array}{rl}
         \min & J_{\rho}(\u)=\sum\limits_{k=0}^{N-1}(h(c-pqx_k)u_k) +\rho \sum\limits_{k=0}^{N-1}|u_{k+1}-u_k| \\
          & x_{k+1}= x_k +h(1-x_k-qu_k)x_k \; \textrm{for all } 0\leq k\leq N-1,\\
         & x_0>0,\\
          &0\leq u_k\leq M \; \textrm{for all }0\leq k\leq N-1.
     \end{array}
 \end{equation}
 
  Because PASA involves a gradient scheme we should be concerned about the absolute value terms that are used in problem (\ref{eqn: profitpendiscrete})'s cost functional.
  We suggest a decomposition of each absolute value term in $J_{\rho}$ to ensure that $J_{\rho}$ is differentiable. 
  We introduce two $N-1$ vectors $\bs{\zeta}$ and $\bs{\iota}$ whose entries are non-negative. 
  Each entry of $\bs{\zeta}$ and $\bs{\iota}$ will be defined as: 
{  
\begin{align*}
      |u_{k+1}-u_k|=\zeta_k+\iota_k, \; \textrm{for all } k =0,\dots, N-2. 
  \end{align*}
}
  An equivalent way of expressing the above equation is to assign values to the components of $\bs{\zeta}$ and $\bs{\iota}$ based upon the following conditions. 
  \begin{align*}
     & \textrm{Condition 1: If } u_{k+1}-u_k > 0, \textrm{ then } \zeta_k= u_{k+1}-u_k \textrm{ and } \iota_k = 0;\\
     & \textrm{Condition 2: If } u_{k+1}-u_k \leq 0, \textrm{ then } \zeta_k = 0 \textrm{ and } \iota_k = -(u_{k+1}-u_k).
  \end{align*}
  This decomposition will convert problem (\ref{eqn: profitpendiscrete}) into the following: 
  \begin{equation}\label{eqn: profitpendiscrete1}
     \begin{array}{rl}
      \min &J_{\rho}(\u,\bs{\zeta},\bs{\iota})=\sum\limits_{k=0}^{N-1}(h(c-pqx_k)u_k)+ \rho \sum\limits_{k=0}^{N-2}(\zeta_k+\iota_k)\\
      & x_{k+1}=x_k+h(1-x_k-qu_k)x_k \textrm{ for all } k = 0,\dots, N-1,\\
      & x_0>0,\\
      & 0\leq u_k\leq M \; \textrm{for all } k= 0,\dots, N-1,\\
      & u_{k+1}-u_k= \zeta_k-\iota_k \; \textrm{for all } k=0,\dots, N-2,\\
      & \zeta_k\geq 0 \textrm{ and } \iota_k\geq 0 \textrm{ for all } 1, \dots, N-1. 
    \end{array}
  \end{equation}
  For problem (\ref{eqn: profitpendiscrete1}), we are minimizing the penalized objective function with respect to vectors $\u,\bs{\zeta},$ and $\bs{\iota}$. 
  The constraints associated with $\bs{\zeta}$ and $\bs{\iota}$, are constraints that PASA can interpret. 
  The equality constraints associated with $\bs{\zeta}$ and $\bs{\iota}$ can be written accordingly: 
  \begin{equation}
  \left[
      \begin{array}{c | c | c}
      \bs{A} &-\bs{I}_{N-1}&\bs{I}_{N-1}
      \end{array}
   \right]
   \begin{bmatrix}
   \u\\
   \hline
   \bs{\zeta}\\
   \hline
   \bs{\iota}
   \end{bmatrix}
   =
   \bs{0},
  \end{equation}
  where $\bs{I}_{N-1}$ is the identity matrix with dimension $N-1$, $\bs{0}$ is the $N-1$ dimensional all zeros vector, and {$\bs{A}$ is the $N-1\times N$ dimensional  sparse matrix given in (\ref{eqn: sparsematrix}).}
  Since problem (\ref{eqn: profitpendiscrete1}) is optimizing $J_{\rho}$ with respect to $\u$, $\bs{\zeta}$, and $\bs{\iota}$, the Lagrangian of problem (\ref{eqn: profitpendiscrete1}) will be
  \begin{equation*}
      \Ls(\x,\u,\bs{\zeta},\bs{\iota},\lam)=\sum\limits_{k=0}^{N-1}(h(c-pqx_k)u_k) + \rho\sum\limits_{k=0}^{N-2}(\zeta_k+\iota_k)+\sum\limits_{k=0}^{N-1} \lambda_k(-x_{k+1}+x_k+h(x_k-x_k^2-qu_kx_k)).
  \end{equation*}
  If we use Theorem (\ref{thm: lagrange}) to find the gradient of $J_\rho$, we would have that 
  \begin{equation}\label{eqn: profitgradJpen}
      \grad_{\u,\bs{\zeta},\bs{\iota}}J_{\rho}=\grad_{\u,\bs{\zeta},\bs{\iota}}\Ls = 
      \begin{bmatrix}
        \grad_{\u}J\\
        \hline
        \rho\\
        \vdots\\
        \rho\\
        \hline
        \rho\\
        \vdots\\
        \rho
      \end{bmatrix},
  \end{equation}
  where $\grad_{\u}J$ is defined in equation (\ref{eqn: profitgradJ}), provided that condition (\ref{eqn: thmadjoint}) is satisfied. 
  As before, we satisfy condition (\ref{eqn: thmadjoint}) by taking the partial derivative of $\Ls$ with respect to each $x_k$ for $k=1,\dots,N$ and set each partial derivative equal to 0 and solve for the components of $\lam$.
  Performing this procedure yields equations (\ref{eqn: profitadjointdiscrete1}) and (\ref{eqn: profitadjointdiscreteN1}), which serves as our discretization of adjoint equation (\ref{eqn: profitadjoint}) and our generalization of the transversality condition  (\ref{eqn: profitadjointTC}). 
  \subsubsection{Summary of Discretization}\label{subsec: profitdiscretesummary}
   We converted the penalized harvesting problem (\ref{eqn: profitpen}) into problem (\ref{eqn: profitpendiscrete1}) by performing the following steps. 
   \begin{enumerate}
       \item \textbf{Discretization of State Equation:} We use explicit Euler's method to discretize the state equation which is given in equation (\ref{eqn: profitstatediscrete}).  
       \item \textbf{Discretization of Objective Functional and Decomposition of Absolute Value Terms:} We use a left-rectangular integral approximation for discretizing the integral shown in (\ref{eqn: profitpen}). 
       We also assume $u$ as being piecewise constant over each mesh interval, which allows us to convert the total variation term given in (\ref{eqn: profitV}) into the finite series of absolute value terms that is used in problem (\ref{eqn: profitpendiscrete}). 
       Because we would like to use a gradient scheme for solving the discretized problem, we need to decompose the absolute value terms by introducing two $N-1$ vectors $\bs{\zeta}$ and $\bs{\iota}$ that satisfy the constraints that were included in problem (\ref{eqn: profitpendiscrete1}). 
       \item \textbf{Finding the Gradient of the Penalized Objective Functional:} By Theorem \ref{thm: lagrange}, we used the  gradient of the Lagrangian to problem (\ref{eqn: profitpendiscrete1}) to find $\grad_{[\u,\bs{\zeta},\bs{\iota}]}J_{\rho}$.
       The formula for  $\grad_{[\u,\bs{\zeta},\bs{\iota}]}J_{\rho}$ is given in equation (\ref{eqn: profitgradJpen}) where $\grad_{\u}J$ is given in equation (\ref{eqn: profitgradJ}). 
       \item \textbf{Discretization of Adjoint Variable:} In order to apply Theorem \ref{thm: lagrange} for computing  the gradient of the penalized objective function, adjoint vector $\bs{\lambda}$ needs to satisfy condition (\ref{eqn: thmadjoint}). 
       We set each partial derivative of the Lagrangian $\Ls$ with respect to $x_k$ equal to 0 and solve for $\lambda_{k-1}$. 
       This results in discretizing adjoint equation (\ref{eqn: profitadjoint}) and the transversality condition $\lambda(T)=0$. 
       The discretized equations are given in equations (\ref{eqn: profitadjointdiscrete1}) and (\ref{eqn: profitadjointdiscreteN1}). 
   \end{enumerate}
  \subsection{Numerical Results of Fishery Problem}\label{subsec: numericalfish}
\begin{figure}[htbp]
      \centering
      \includegraphics[width=0.65\linewidth]{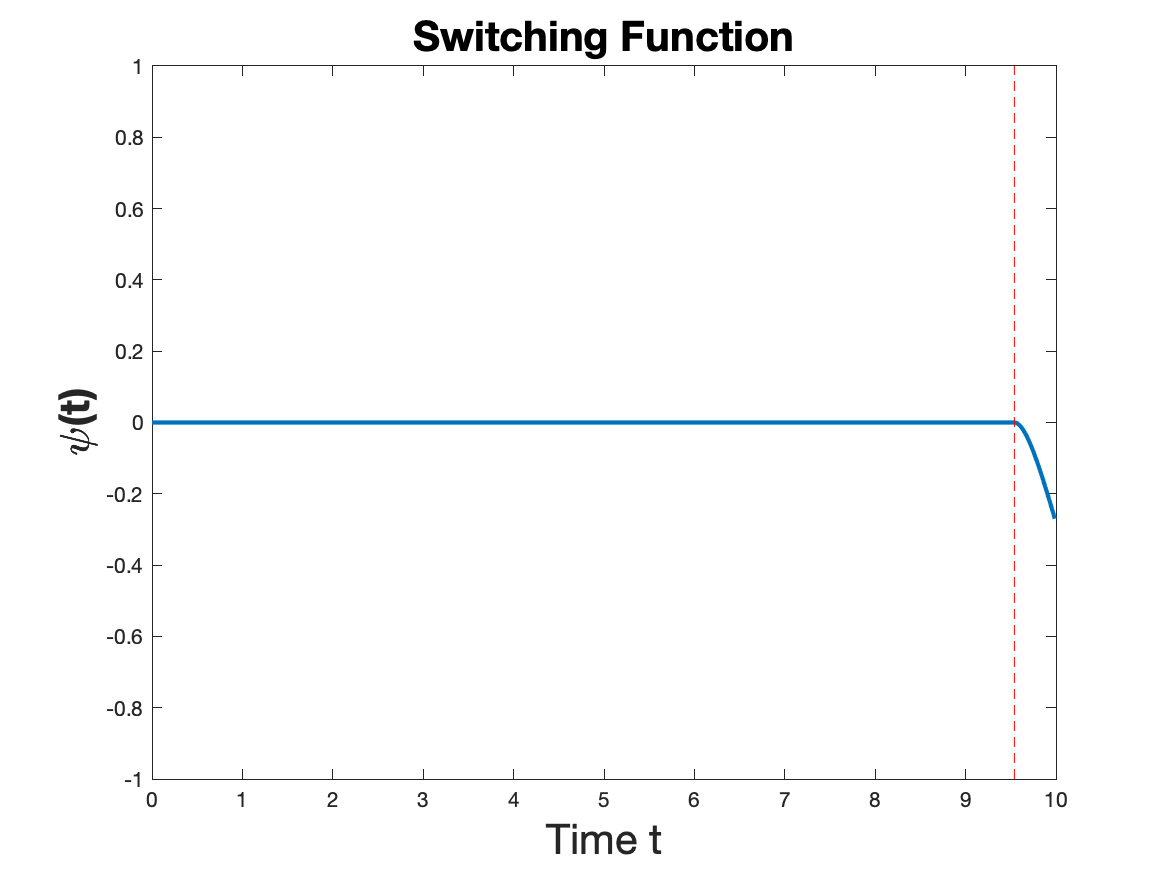}
      \caption{Plot of the Switching Function. The red dotted line is the vertical line $t=t^*\approx 9.5392$.}
      \label{fig: profitswitch}
  \end{figure}

  We wish to numerically solve for problem (\ref{eqn: minprofit}) with the following parameters defined in Table \ref{tab: profitparameters}. 
  Note that these parameter values satisfy Assumptions \ref{as: assump1} and \ref{as: assump2}, so the optimal harvesting strategy control will be of the form given in (\ref{eqn: fullprofitu}).
  Based on Table \ref{tab: profitparameters} and equation (\ref{eqn: switchprofit}), 
  {$t^*\approx 9.5392$} which is between values 0 and $T=10$. 
  Additionally, from these parameter values the singular control solution is $u^*=0.1875$  which will be  between the bounds $0$ and $M=1$, as desired.
  After using Table \ref{tab: profitparameters} to evaluate $t^*$ from (\ref{eqn: switchprofit}), $u^*$ from (\ref{eqn: fullprofitu}), $x^*$ from (\ref{eqn: fullprofitx}), and $\lambda^*$ from (\ref{eqn: fullprofitlam}), we substitute these solutions into the switching function (\ref{eqn: profitswitchfunction}) to see if $u^*$ satisfies (\ref{eqn: profitupsi}). 
  In Figure \ref{fig: profitswitch}, we plot the switching function $\psi$ that we obtained. Regard that the switching function is zero when $u^*$ is singular and becomes negative at the interval $[t^*,T]$, which is when $u^*=M$. 
  So we have that harvesting policy $u^*$ for (\ref{eqn: fullprofitu}) satisfying (\ref{eqn: profitupsi}). 
  This is equivalent to saying that $u^*$ satisfies Pontryagin's Minimum Principle \cite{Pontryagin}, which is the first order necessary condition of optimality. 
   \begin{table}[ht]
  \centering
  \begin{tabular}{|c| l | c|}
     \hline
    \textbf{Parameter} &\textbf{Description}& \textbf{Value}\\
    \hline
     $T$  & Terminal time                         & 10 \\
     $p$  & selling price of one unit of fish     & 2\\
     $q$  & ``catchability" of the fish           & 2\\
     $c$  & cost of harvesting one unit of fish   & 1\\
     $M$  & Maximum harvest effort                & 1\\
     $x_0$& initial population size of fish       &$\frac{c+pq}{2pq}=0.625$ \\
     \hline
  \end{tabular}
  \caption{Numerical Values of Parameters Used in Computations}
  \label{tab: profitparameters}
  \end{table}

\begin{figure}[htbp]
  \begin{subfigure}[t]{0.33\textwidth}
    \includegraphics[width=\linewidth]{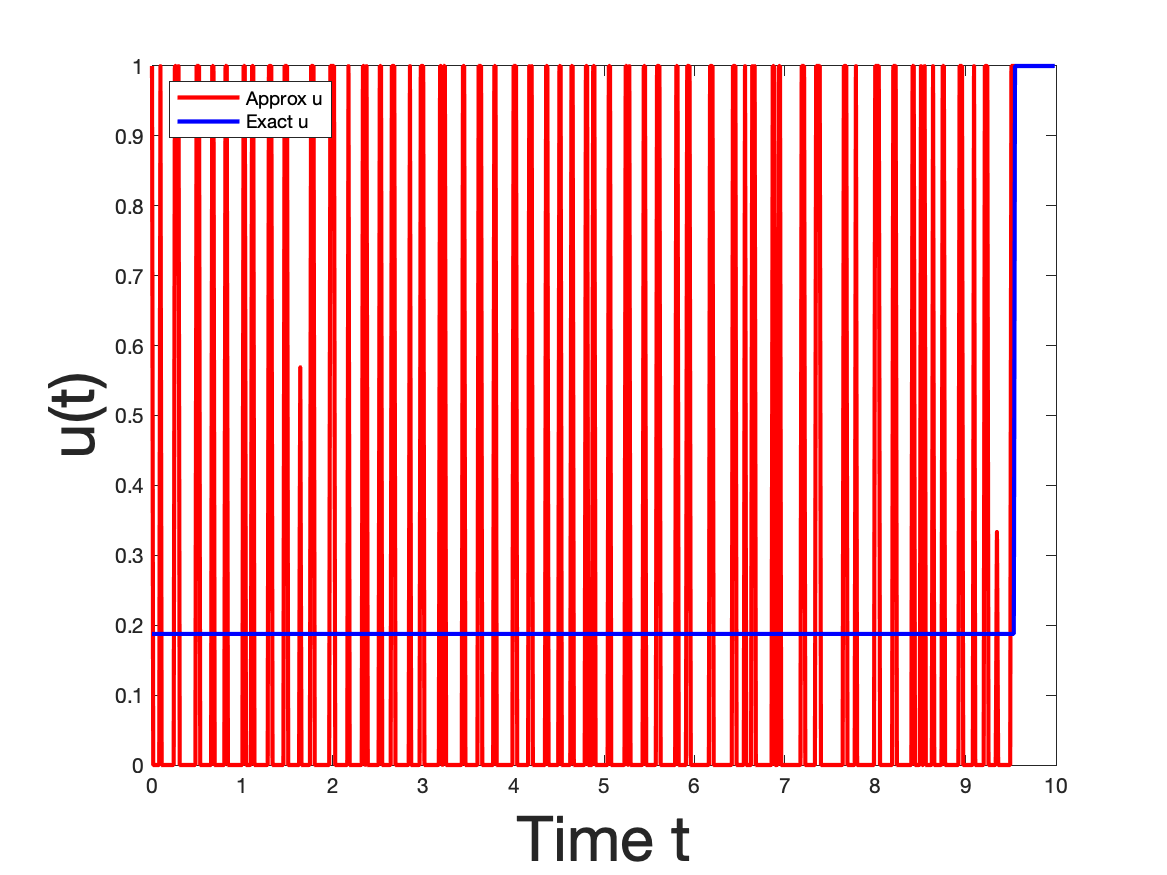}
    \caption{Unpenalized Control $\hat{u}$ (red) vs Optimal Control $u^*$ (blue)}
    \label{fig: ProfitUnpenA}
  \end{subfigure}\hfill
	\begin{subfigure}[t]{0.33\textwidth}
  	  \includegraphics[width=\linewidth]{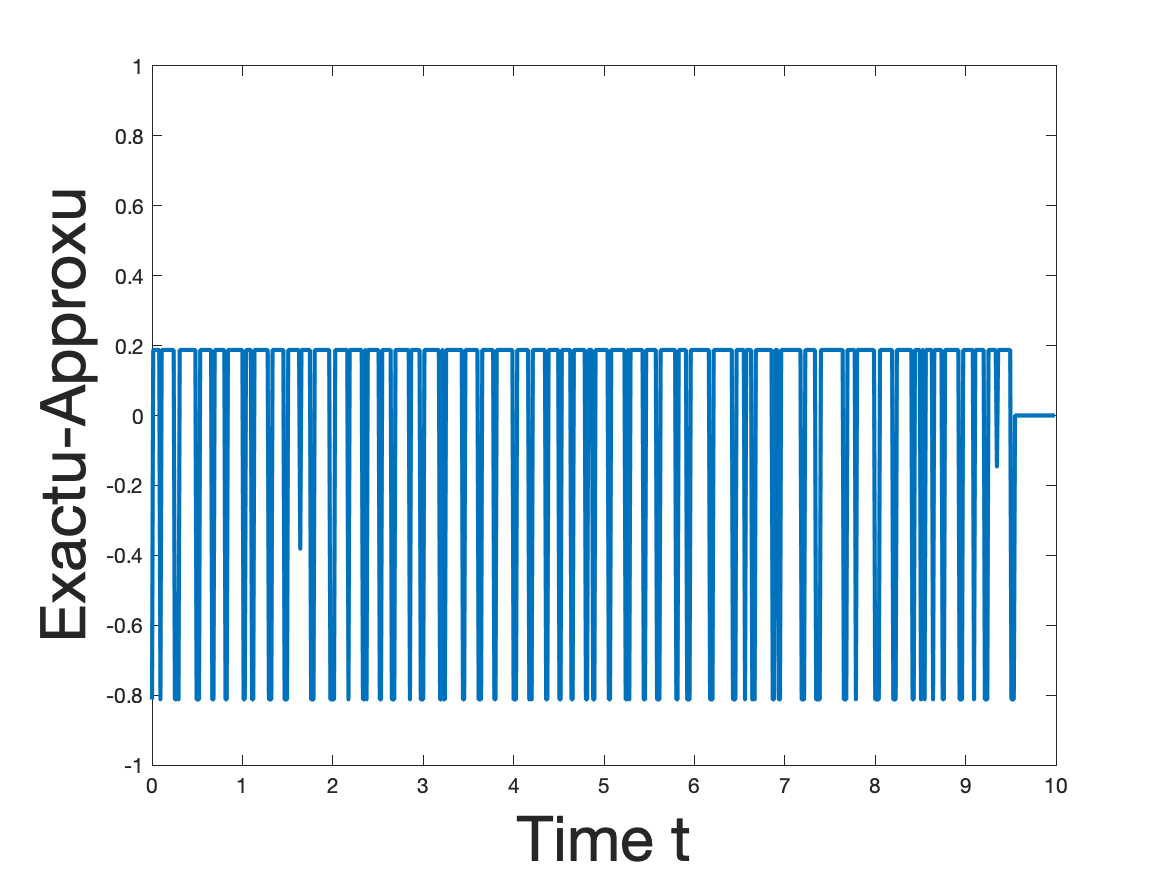}
   	 \caption{$u^*-\hat{u}$ (cyan).}
  	  \label{fig: profftdiffnopen}
	\end{subfigure}
 \begin{subfigure}[t]{0.33\textwidth}
    \includegraphics[width=\linewidth]{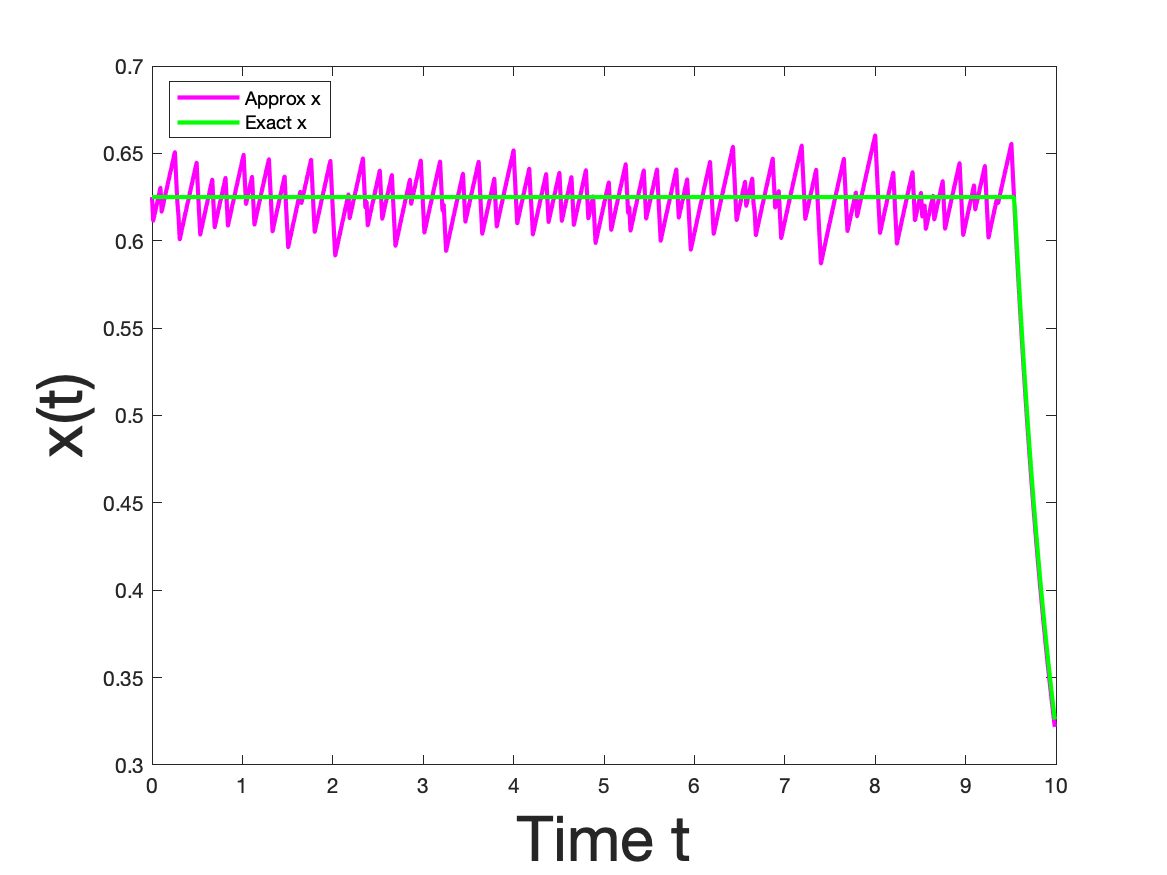}
    \caption{Fish Population Solution $\hat{x}$ Corresponding to Unpenalized $\hat{u}$ (pink) vs Exact State Solution $x^*$ (green) }
     \label{fig: ProfitUnpenB}
  \end{subfigure}\hfill
 \caption{Unpenalized  Results: Time interval $[0,T]$ was partitioned to have $N=750$ mesh intervals and  PASA stopping tolerance was tol=10$^{-10}$. }
  \end{figure}

  We used PASA to solve for problem (\ref{eqn: minprofit}) with parameters settings given in Table \ref{tab: profitparameters}.   
  Our initial guess for control $u$ is  $u(t)=0$ for all $t\in [0,T]$, and the stopping tolerance is set to $10^{-10}$.  
  We partition $[0,T]$ to where there are $N=750$ mesh intervals, and we use the discretization method that is presented in Section \ref{subsec: profitdiscrete}.  
  We first wanted to observe PASA approximation for problem (\ref{eqn: minprofit}) without any penalty being applied, which we denote as being $\hat{u}$.
  In Figure \ref{fig: ProfitUnpenA}, a plot of $\hat{u}$ is shown in red and is compared to the exact solution, $u^*$ which is shown in blue. 
  In Figure \ref{fig: ProfitUnpenA}, $\hat{u}$  resembles chattering on the singular region.
  Additionally, we reran experiments with a tighter stopping tolerance and a finer partition of $[0,T]$, and PASA still obtained an oscillatory solution.  
  {It is likely that the discretization of problem (\ref{eqn: minprofit}) causes PASA's to generate a solution that resembles chattering.} 
  We computed a left-rectangular integral approximation of the total profit of the harvested fish over time interval $[0,T]$ (i.e. the cost functional for the equivalent maximization problem (\ref{eqn: maxharvest})) when employing harvesting policy $\hat{u}$ and and optimal harvesting policy, $u^*$, and  found that $J(\hat{u})\approx 3.074$ and $J(u^*)\approx 3.0575$.  

  The major concern associated with the approximate solution obtained in Figure \ref{fig: ProfitUnpenA} is that it is an unrealistic harvesting strategy to use. 
  We wish to penalize problem (\ref{eqn: minprofit}) by adding a bounded variation term that will reduce the number of oscillations.
  Before getting into the results, we would like to discuss our methods for determining which penalty parameter values gave the best solution.
  We have some advantage for choosing an appropriate penalty parameter value due to knowing the analytic solution. 
  However, if we did not know what $u^*$ looked like we would recommend first looking at the plots of PASA's approximation of the penalized harvesting policy, denoted as $u_{\rho}$.
  Firstly, we would rule out a penalty parameter value if the plot of $u_{\rho}$ contained any unusual jumps and/or oscillations.
  If such a solution occurred, we would suggest that the penalty parameter value is too small and is generating a solution that is similar to the unpenalized solution. Secondly, we would rule out a $\rho$  value if the corresponding penalized solution $u_{\rho}$ does not closely align with Pontryagin's Minimum Principle \cite{Pontryagin}. 
  For example, based upon Pontryagin's Minimum Principle, $u^*$ is piecewise constant where either $u^*(t)\in\{0,\frac{pq-c}{2pq^2},M\}$.
  So if the penalized solution  $u_{\rho}$ ever took on values that were not $0,M,$ or $\frac{pq-c}{2pq^2}$, then we would find value $\rho$ to be suspect.  
  Thirdly, one might find it useful to look at the plots of the switching function that are associated with the penalized solution $u_{\rho}$. 
  This third suggestion might be the most useful if one were penalizing an optimal control problem where the singular case solution cannot be found explicitly.
  The sign of the switching function can help one verify whether or not $u_{\rho}$ aligns with Pontryagin's Minimum Principle.

 \begin{figure}
 \centering
    \begin{subfigure}[t]{0.25\textwidth}
     \includegraphics[width=\linewidth]{images/profitnopen.png}
     \caption{$u_{\rho}$ vs $u^*$ for $\rho=0$}
        \label{fig: profitvary0}
    \end{subfigure}\hfill
    \begin{subfigure}[t]{0.25\textwidth}
     \includegraphics[width=\linewidth]{images/profitdiff.png}
     \caption{$u^*-u_{\rho}$ for $\rho=0$}
     \label{fig: profitvarydiff0}
    \end{subfigure}\hfill
    \begin{subfigure}[t]{0.25\textwidth}
     \includegraphics[width=\linewidth]{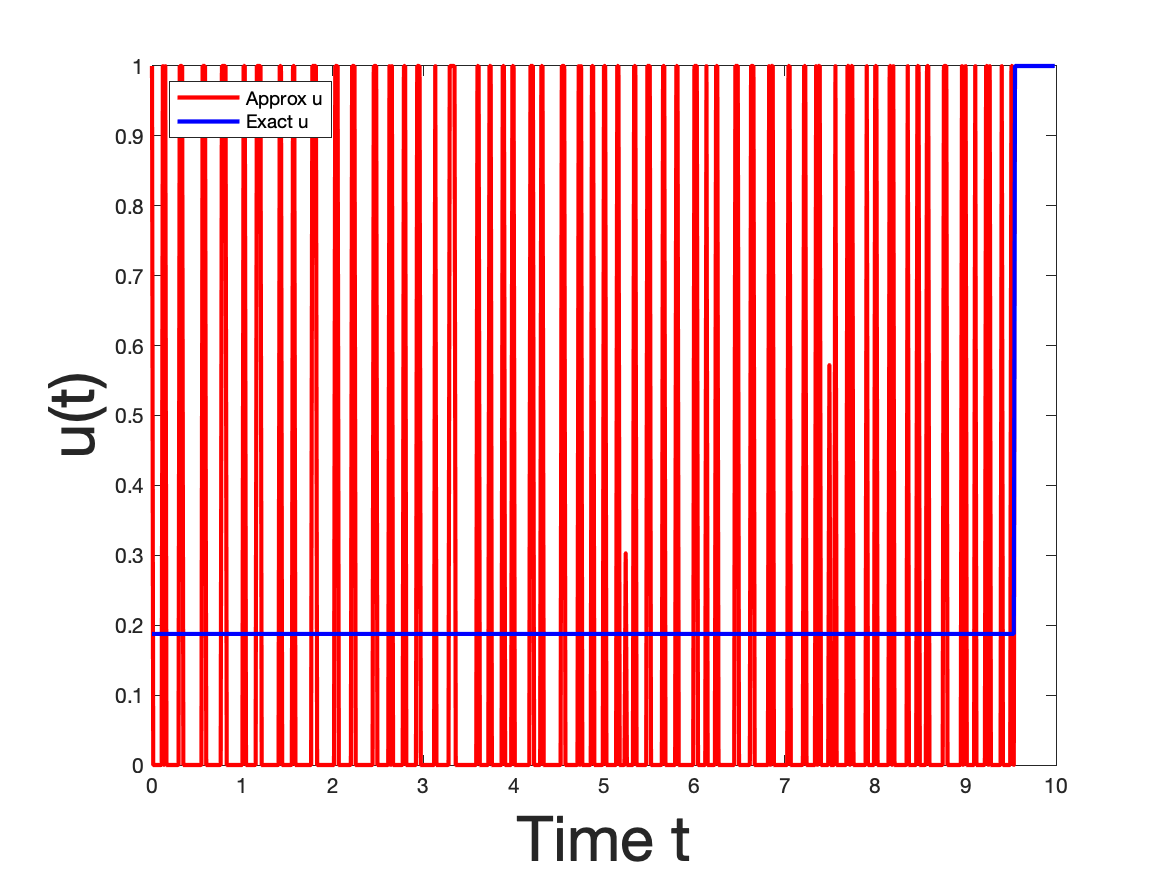}
      \caption{$u_{\rho}$ vs $u^*$ for $\rho=10^{-8}$}
     \label{fig: profitvary8}
    \end{subfigure}\hfill
    \begin{subfigure}[t]{0.25\textwidth}
     \includegraphics[width=\linewidth]{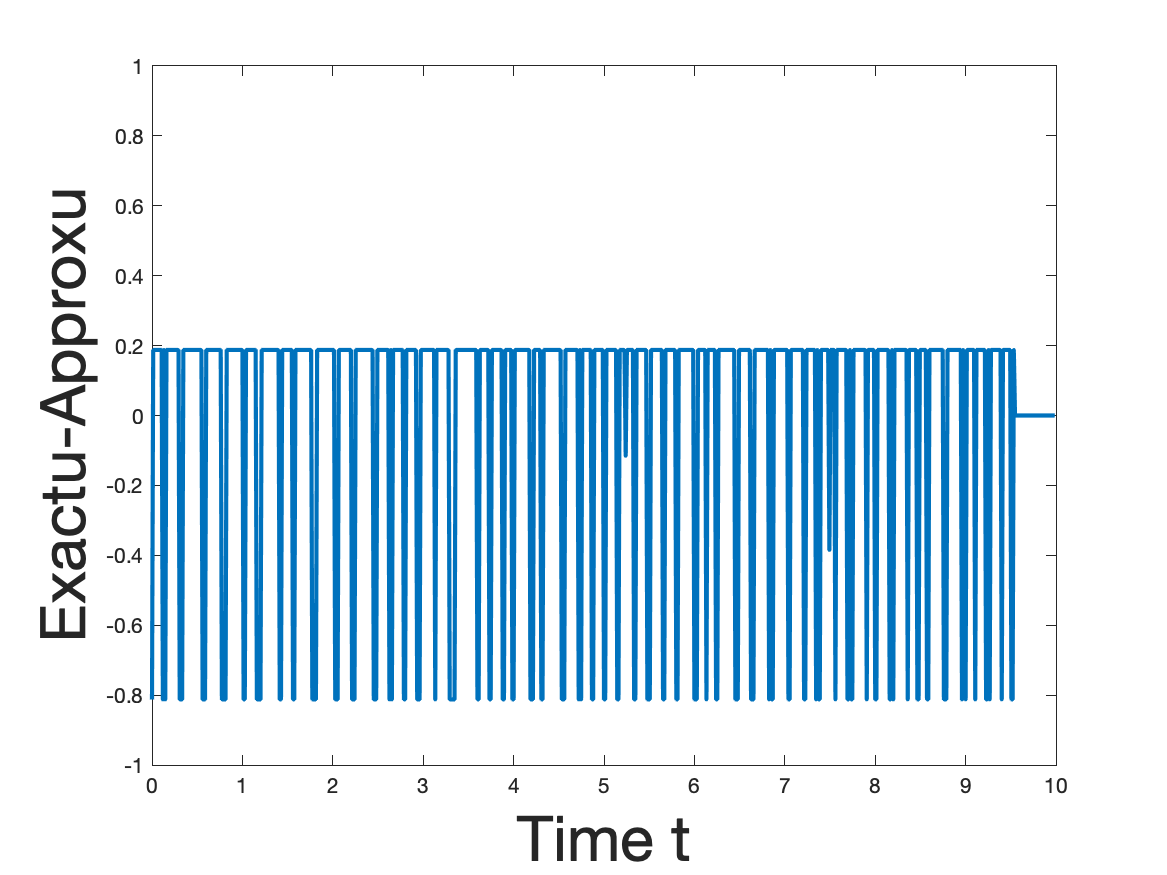}
      \caption{$u^*-u_{\rho}$ for $\rho=10^{-8}$}
     \label{fig: profitvarydiff8}
    \end{subfigure}\hfill
    
    \begin{subfigure}[t]{0.25\textwidth}
     \includegraphics[width=\linewidth]{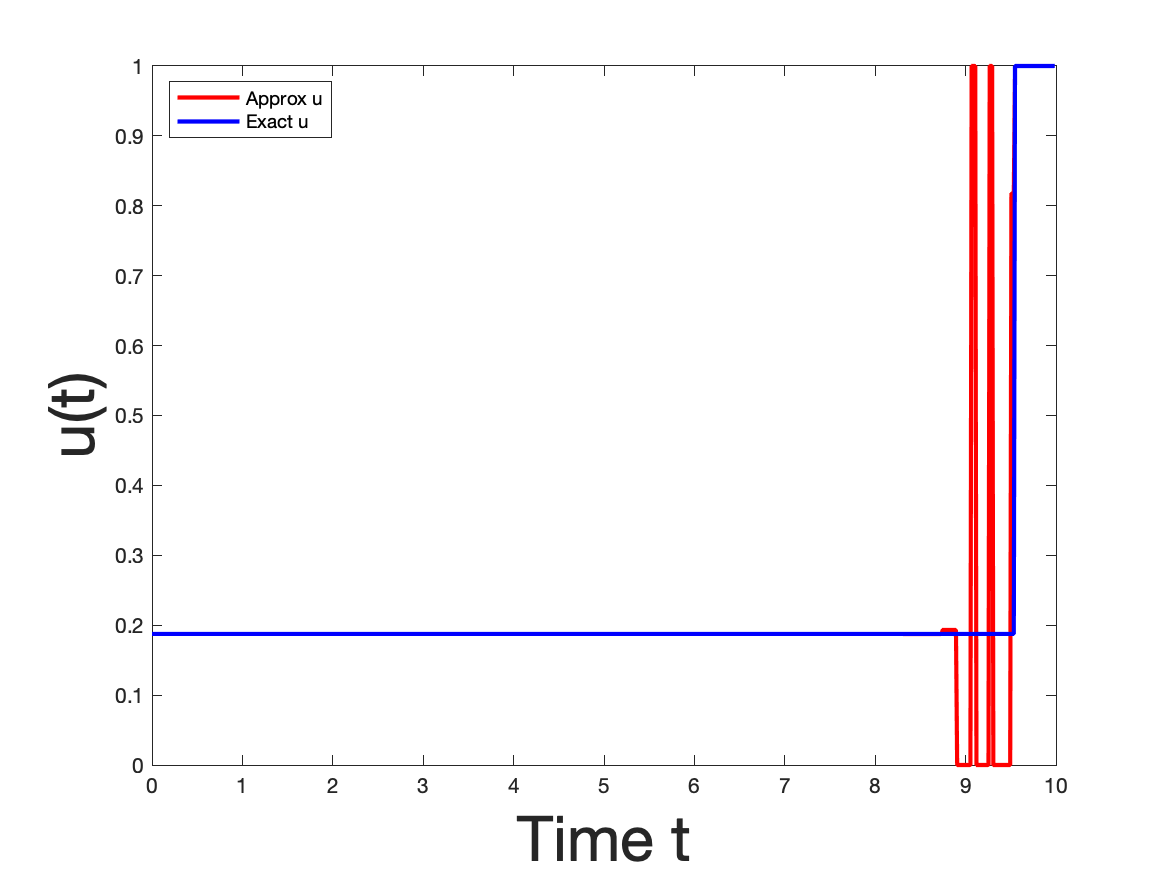}
     \caption{$u_{\rho}$ vs $u^*$ for $\rho=10^{-5}$}
     \label{fig: profitvary5}      
    \end{subfigure}\hfill
    \begin{subfigure}[t]{0.25\textwidth}
     \includegraphics[width=\linewidth]{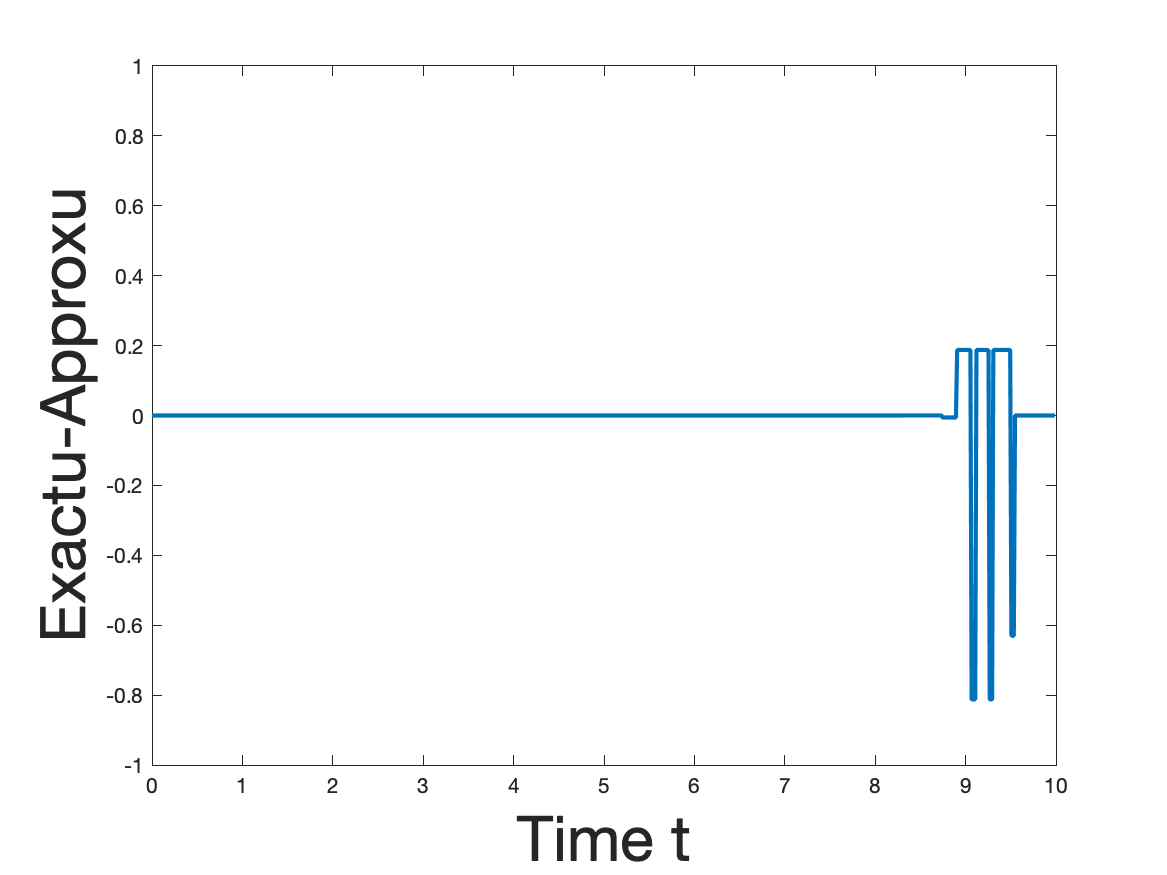}
     \caption{$u^*-u_{\rho}$ for $\rho=10^{-5}$}
     \label{fig: profitvarydiff5}      
    \end{subfigure}\hfill
    \begin{subfigure}[t]{0.25\textwidth}
     \includegraphics[width=\linewidth]{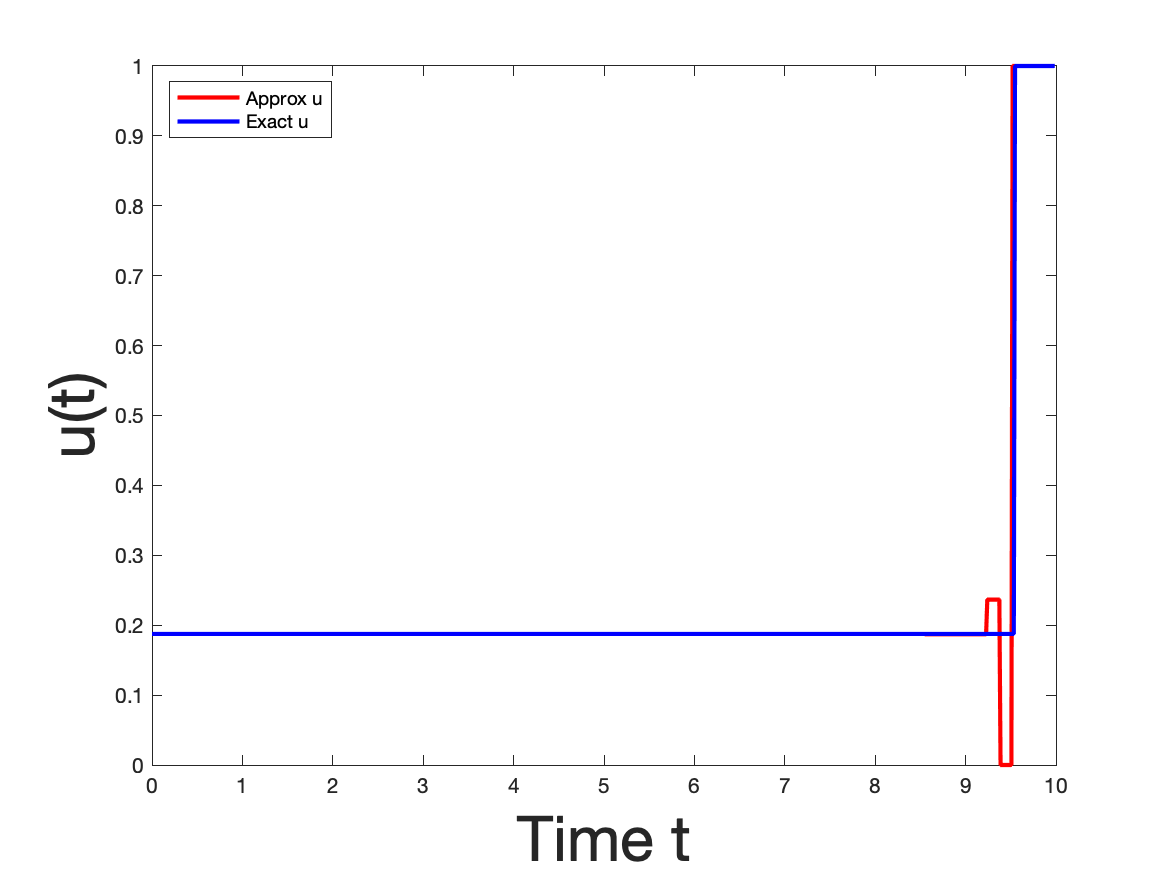}
     \caption{$u_{\rho}$ vs $u^*$ for $\rho=10^{-4}$}
     \label{fig: profitvary4}      
    \end{subfigure}\hfill
    \begin{subfigure}[t]{0.25\textwidth}
     \includegraphics[width=\linewidth]{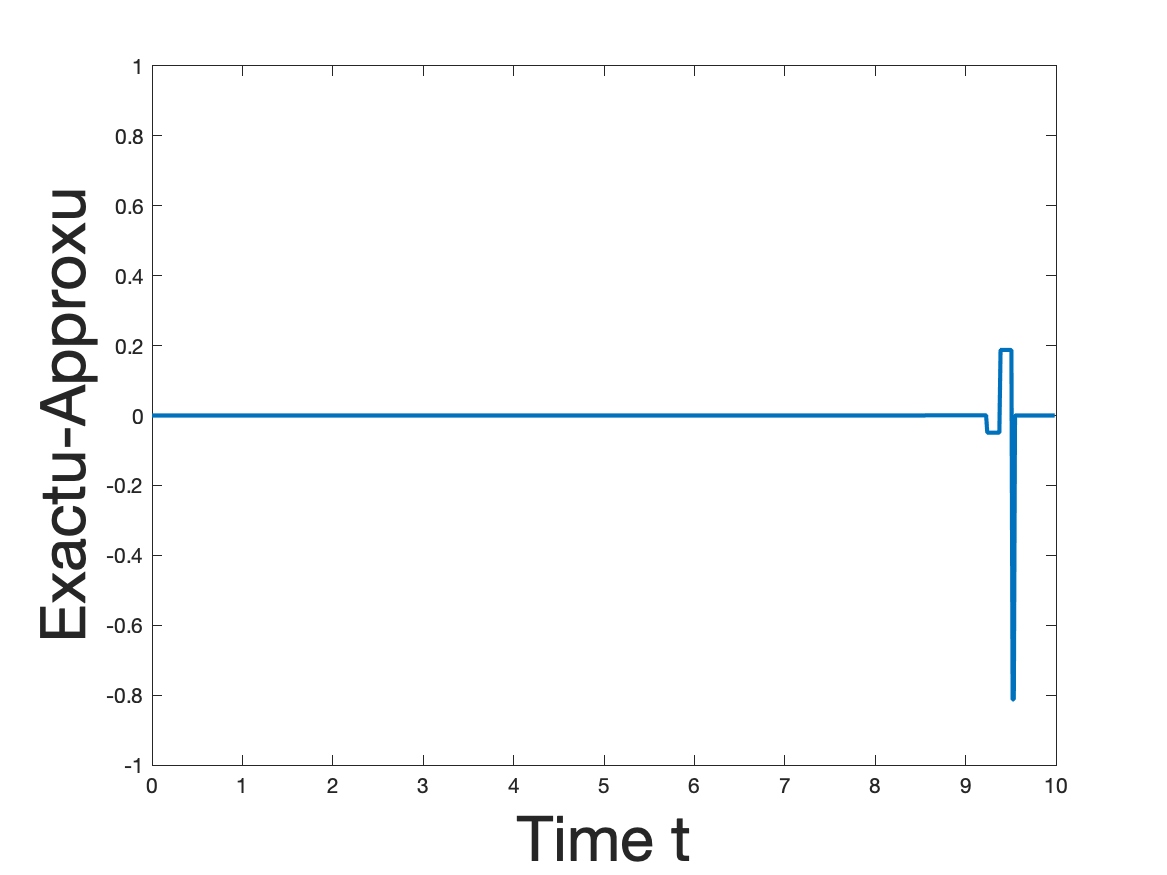}
     \caption{$u^*-u_{\rho}$ for $\rho=10^{-4}$}
     \label{fig: profitvarydiff4}      
    \end{subfigure}\hfill
    
    \begin{subfigure}[t]{0.25\textwidth}
     \includegraphics[width=\linewidth]{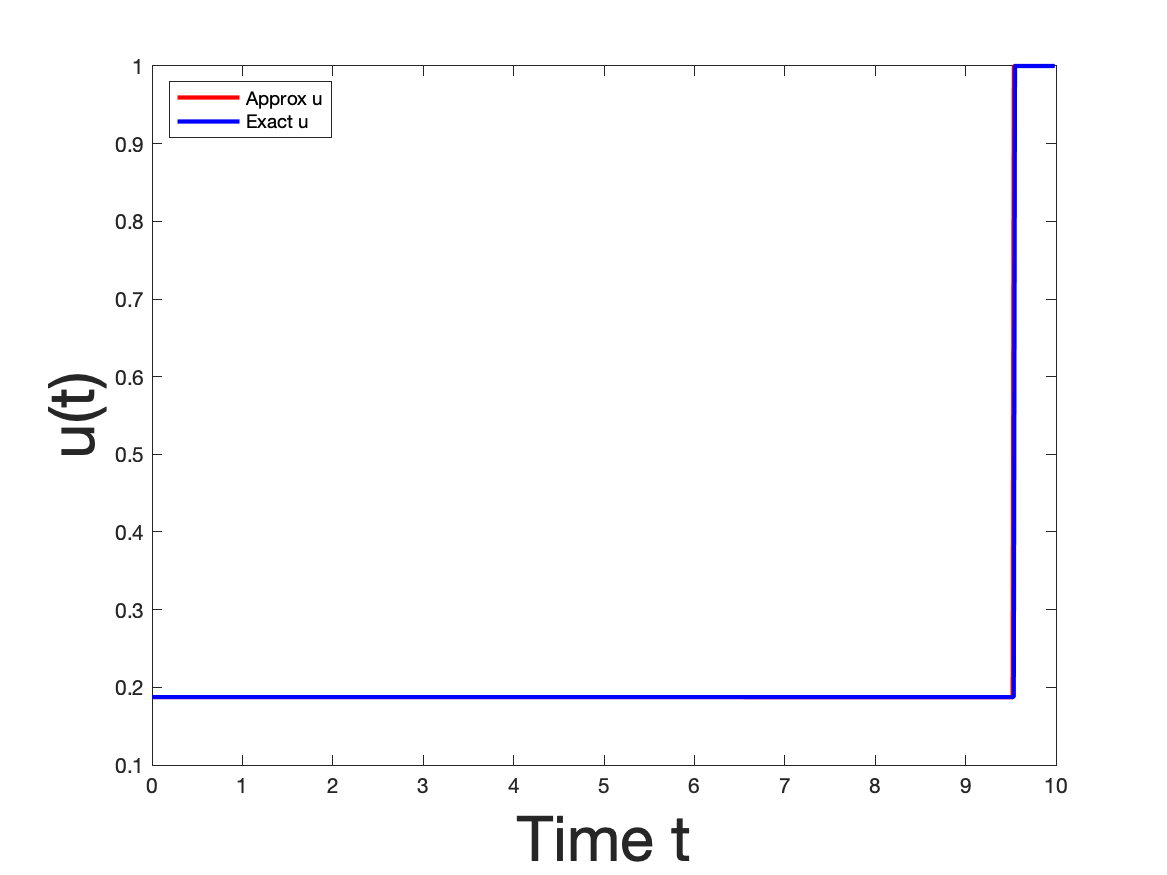}
     \caption{$u_{\rho}$ vs $u^*$ for $\rho=10^{-2}$}
     \label{fig: profitvary2}      
    \end{subfigure}\hfill
    \begin{subfigure}[t]{0.25\textwidth}
     \includegraphics[width=\linewidth]{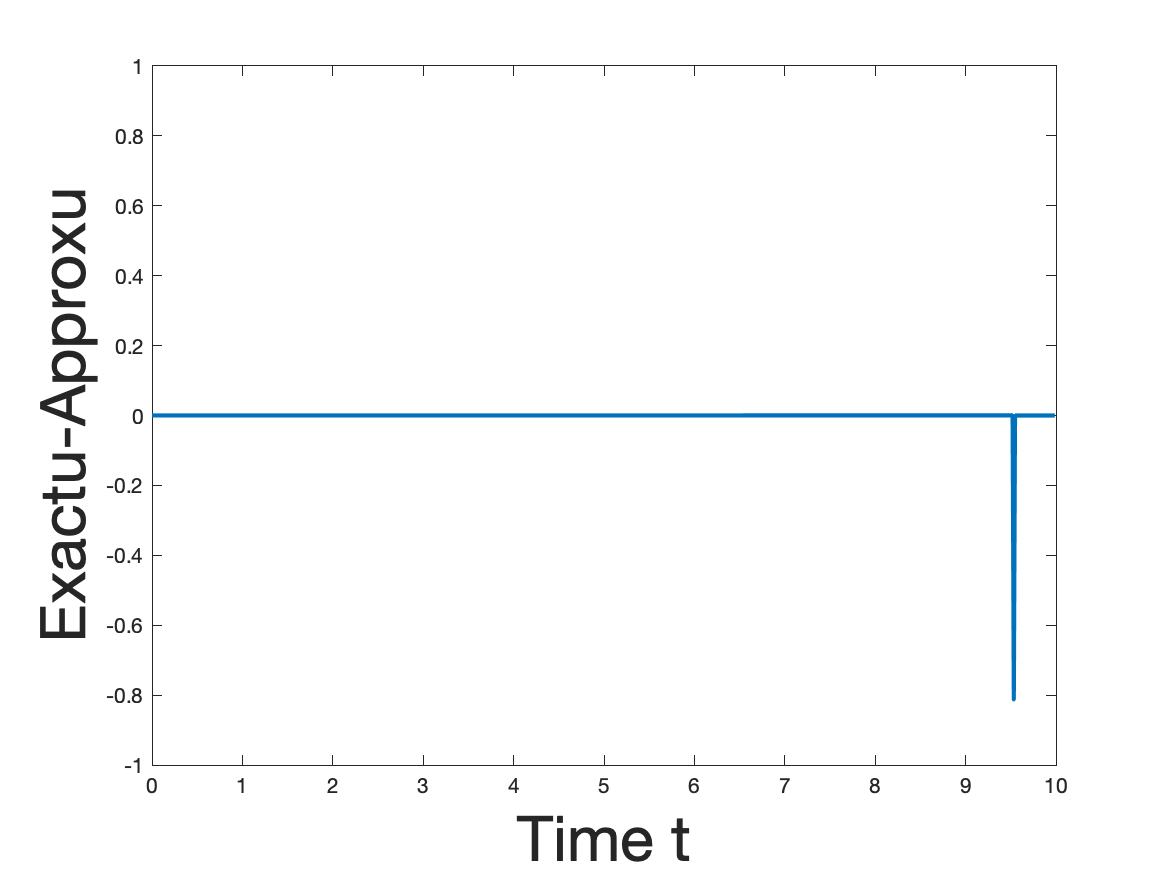}
     \caption{$u^*-u_{\rho}$ for $\rho=10^{-2}$}
     \label{fig: profitvarydiff2}      
    \end{subfigure}\hfill
    \begin{subfigure}[t]{0.25\textwidth}
     \includegraphics[width=\linewidth]{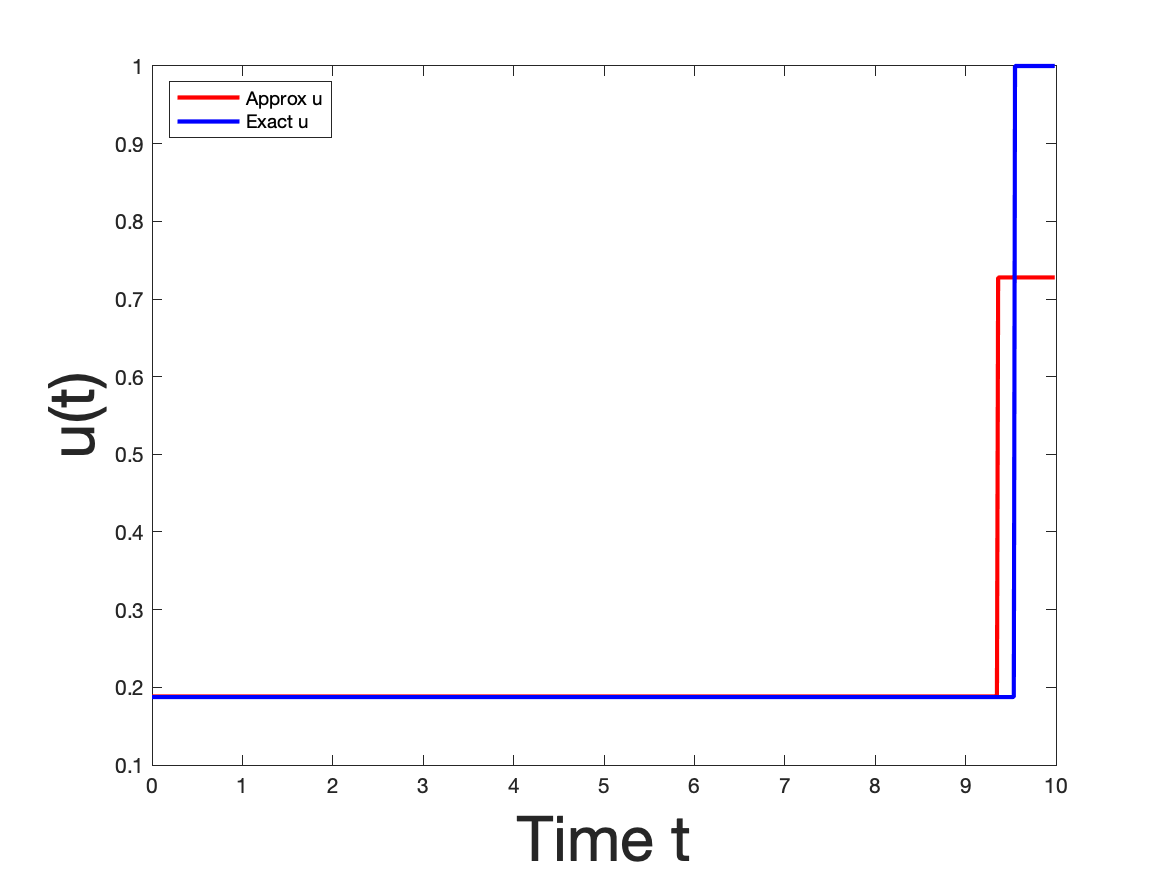}
     \caption{$u_{\rho}$ vs $u^*$ for $\rho=10^{-1}$}
     \label{fig: profitvary1}      
    \end{subfigure}\hfill
    \begin{subfigure}[t]{0.25\textwidth}
     \includegraphics[width=\linewidth]{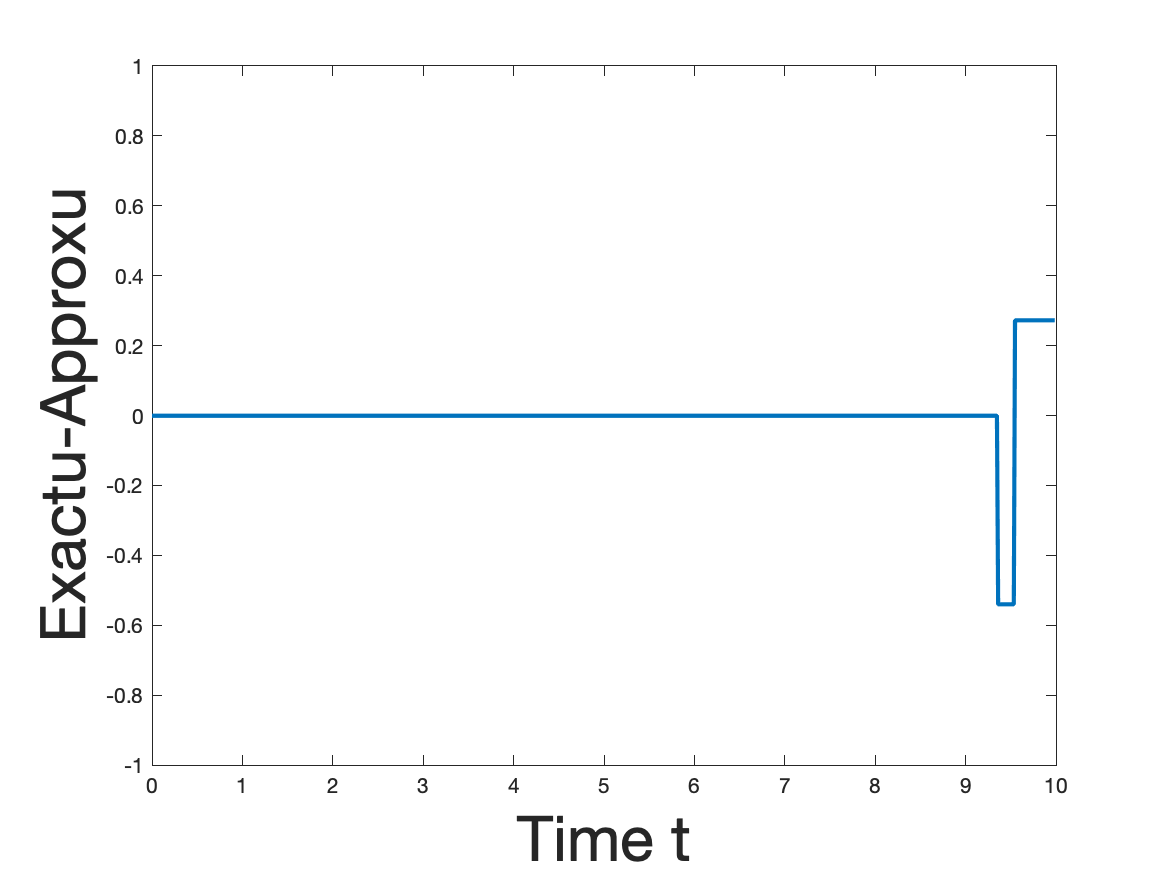}
     \caption{$u^*-u_{\rho}$ for $\rho=10^{-1}$}
     \label{fig: profitvarydiff1}      
    \end{subfigure}\hfill
    \caption{Plots of regularized control $u_{\rho}$ (red) vs optimal control $u^*$ (blue) and $u^*-u_{\rho}$ (cyan) for tuning parameter $\rho\in\{0, 10^{-8}, 10^{-5}, 10^{-4}, 10^{-2}, 10^{-1}\}$. Time interval $[0,T]$ was partitioned to have $N=750$ mesh intervals. PASA stopping tolerance was tol=$10^{-10}$.}
    \label{fig: profitvarying}
 \end{figure}
 \begin{figure}
    	\begin{subfigure}[t]{0.33\textwidth}
    		  \includegraphics[width=\linewidth]{images/profit2.png}
		\caption{$u_{\rho}$ (red) vs $u^*$ (blue) with $\rho=10^{-2}$}
   	  \end{subfigure}\hfill
	 \begin{subfigure}[t]{0.33\textwidth}
    		 \includegraphics[width=\linewidth]{images/profitdiff2.png}
    		 \caption{$u^*-u_{\rho}$ (cyan) with $\rho=10^{-2}$}
     		\label{fig: profitx2}
	\end{subfigure}\hfill
   	\begin{subfigure}[t]{0.33\textwidth}
    		 \includegraphics[width=\linewidth]{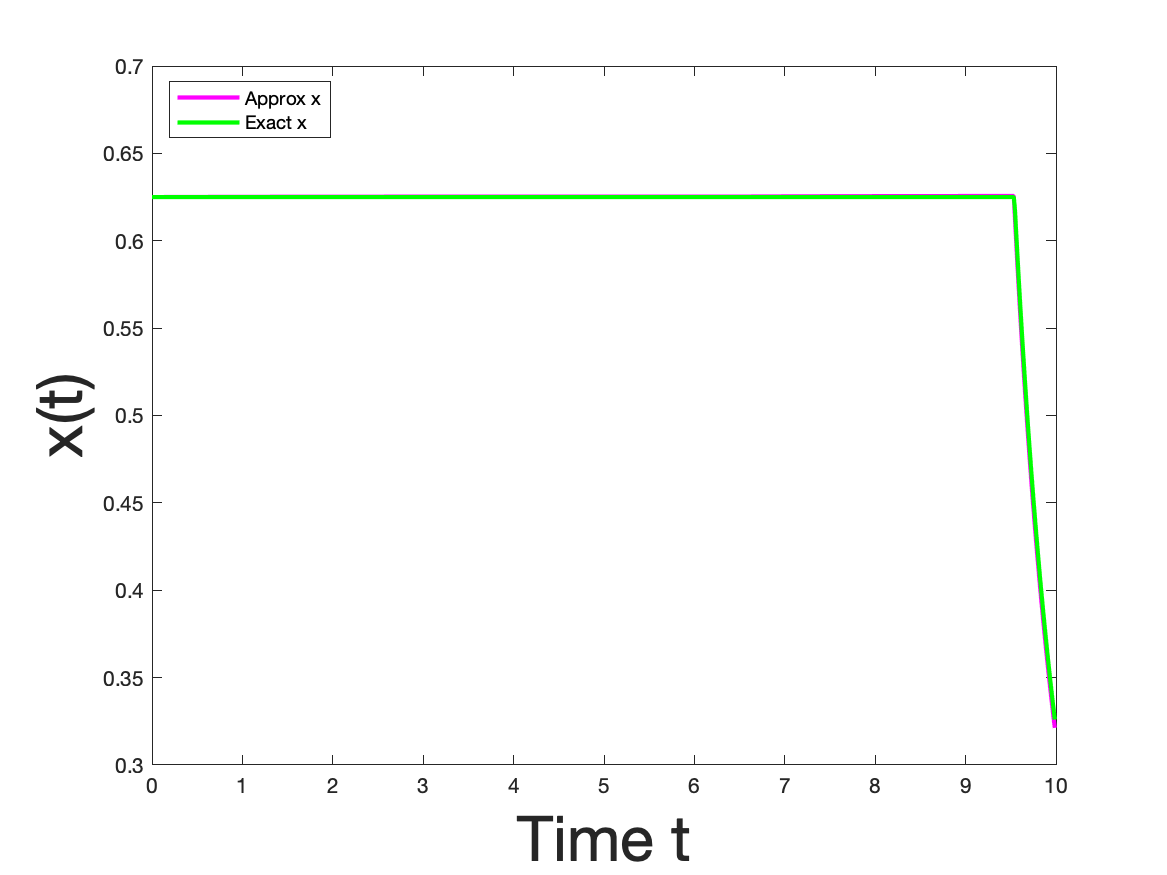}
    		 \caption{Approximate State $x_{\rho}$ (pink) vs Exact State $x^*$ (green) with $\rho=10^{-2}$}
     		\label{fig: profitx2}
	\end{subfigure}
\caption{Penalized Results with $\rho=10^{-2}$: Time interval $[0,T]$ was partitioned to have $N=750$ mesh intervals and  PASA stopping tolerance was tol=10$^{-10}$.}
 \end{figure}
 
 \begin{table}[ht]
  \centering
 \begin{tabular}{|c|c c c c c|}
 \hline
 &&&\textbf{Varying Tuning Parameter Table} &&\\
 \hline
   \textbf{Parameters}  & $\rho$ & $\norm{u^* -u_{\rho}}_{L^1}$ 
   &$\norm{u^*-u_{\rho}}_{\infty}$ & \textbf{Switch}&\textbf{Runtime (s)} \\
   \hline
       $N=750$          & 0         & 2.89548286 & 0.8125 & $0^*$ & 8.82\\
       $\text{tol}=10^{-10}$   & $10^{-9}$ & 2.89262350 & 0.8125 & $0^*$& 14.44\\ 
       $T=10$           & $10^{-8}$ & 2.89708395 & 0.8125 & $0^*$& 22.11\\
       $M=1$            & $10^{-7}$ & 2.89799306 & 0.8125 & $0^*$&28.45 \\
       $p=2$            & $10^{-6}$ & 2.88774180 & 0.8125 & $0.04^*$ & 34.88\\
       $q=2$            & $10^{-5}$ & 0.19705717 & 0.8125 & 9.0667 &6.58 \\
       $c=1$            & $10^{-4}$ & 0.05427589 & 0.8125 & 9.5200 &1.25 \\
       $x_0=0.625$      & $10^{-3}$ & 0.01628685 & 0.8125 & \textbf{9.5333}& 0.68\\
    $h=0.0133$                    & $10^{-2}$ & \textbf{0.01282480} & 0.8125 & \textbf{9.5333}&\textbf{0.38} \\
                        & $10^{-1}$ & 0.22780044 & 0.5402 &  9.9867&0.43 \\
                        \hline
 \end{tabular}
 \caption{The total variation of our exact solution is $V(u^*)=0.8125$, so if our approximated solution $u_{\rho}$ switched at a time different from $u^*$,  then $\norm{u^*-u_{\rho}}_{\infty}=V(u^*)$. 
 The switch column is indicating the first instance in time $t$ when $u_{\rho}(t)=M$. The $*$-values on the switch column indicates that the numerical solution is oscillating.}
 \label{tab: profitvary}
 \end{table}
 \begin{table}[ht]
 \centering
 \begin{tabular}{|c|l c c c |}
 \hline
 \textbf{Parameters} & $h$ & err$_{h}$ & $\frac{\text{err}_h}{\text{err}_{h/2}}$ & $\log_2\left(\frac{\text{err}_h}{\text{err}_{h/2}}\right)$\\
 \hline
 $\text{tol}=10^{-10}$ & 0.2        & 0.28091381 & & \\
 $\rho=10^{-2}$ & 0.1       		& 0.12111065 & 2.31948057     		& 1.21380176\\
                & 0.05      	       		& 0.00089762 & 134.92448804   	& 7.07600840\\
                & 0.025            		& 0.02736103 & 0.03280644     	& -4.92987718\\
                & 0.0125          		& 0.01089531 & 2.51126802    		& 1.32841601\\
                & 0.00625        		& 0.00527063 & 2.06717301     	& 1.04765914\\
                & 0.003125       		& 0.00165509 & 3.18450288     	& 1.67106818\\
                \hline
 \end{tabular}
 \caption{Convergence analysis between the penalized solution  and the exact solution to fishery problem. Here $err_{h}=\norm{u_h-u^*}_{L^1}$, where $u_h$ is the penalized solution for Fishery problem when the mesh size is $h$.}
 \label{tab: profitconverge}
 \end{table}
 \begin{figure}
       \centering
       \includegraphics[width=0.75\linewidth]{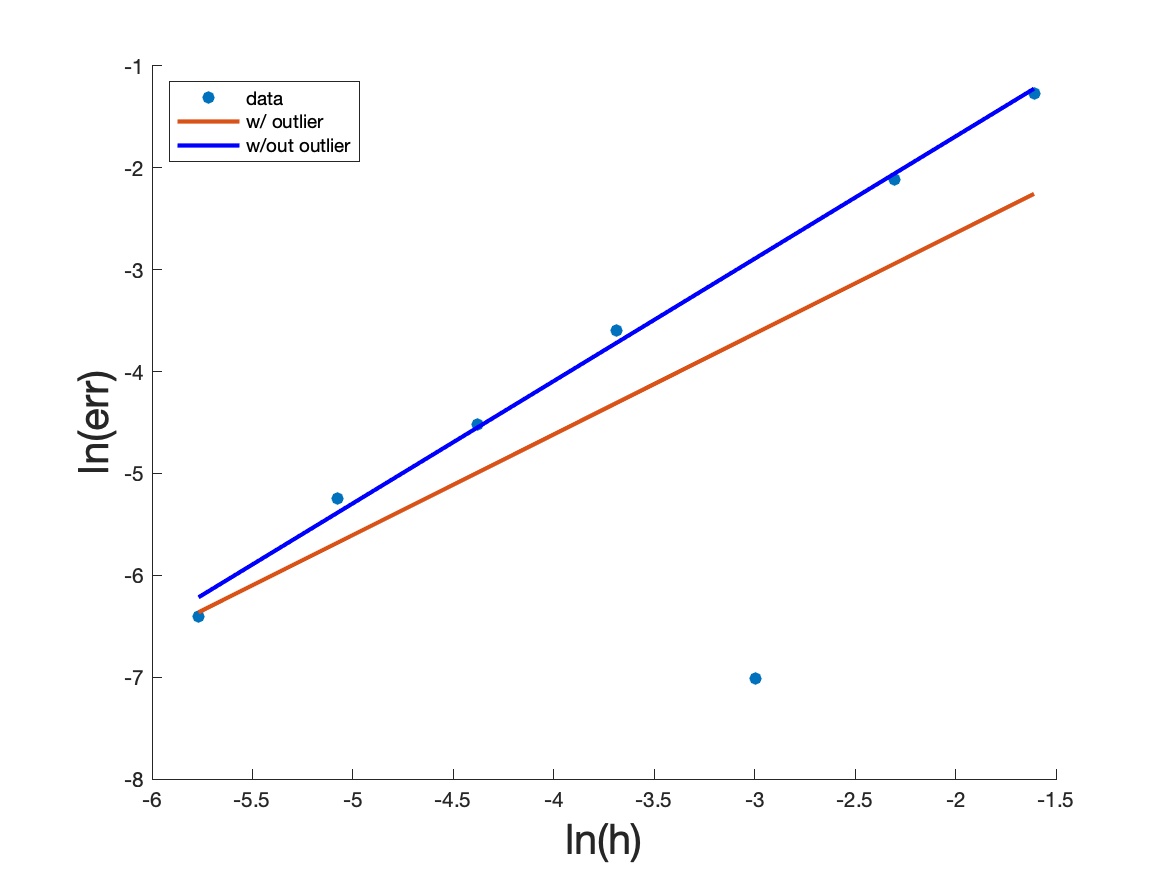}
       \caption{Linear Fit of Data Points  $(\ln{(h)},\ln{(\text{err}_h)})$ from Table \ref{tab: profitconverge}. Red line is the linear fit of the data points including outlier point  $(\ln(0.05), \ln(0.00089762))$. Blue line is the linear fit of data points excluding the outlier point. }
       \label{fig: profitlogfit}
   \end{figure}
    
    \begin{table}[ht]
    \centering
    \begin{tabular}{|c|c c c|}
    \hline
    $\ln(\text{err}_h)=m(\ln(h))+b$& \textbf{slope} $m$ &\textbf{y-intercept} $b$ &  $r^2$\\
    \hline
    \textbf{with outlier} & 0.988064 & -0.6645189 & 0.4826212\\
    \textbf{without outlier} & 1.200738 & 0.7100375 & 0.9959617\\
    \hline
    \end{tabular}
    \caption{Linear Fit Table}
    \label{tab: profitlinear}
    \end{table}

  Since we have the advantage of comparing the penalized solution $u_{\rho}$ with the exact solution $u^*$, we can determine which penalty parameter values are appropriate by comparing the plots between $u^*$ and $u_{\rho}$ as well as the plots of $u^*-u_{\rho}$.
  Additionally, we will look at the instance in time when the penalized solution $u_{\rho}$ switched from singular to non-singular to see if it is close to when the switching point should 
 occur, i.e. when {$t^*\approx9.5392$}.

  When calculating the associated switching point for $u_{\rho}$, we look at the first instance when $u_{\rho}=M$. 
  Due to the discretization of problem (\ref{eqn: profitpen}), the approximated switching points will be some mesh point value of $[0,T]$.
  We also use the $L_1$ norm difference between $u^*$ and $u_{\rho}$ for determining which penalty parameter value gives the closest approximation to $u^*$.
  We observe the $L_{\infty}$ norm difference between $u^*$ and $u_{\rho}$ as well; however, this computation may not be helpful since approximated switching points are restricted to mesh points of the partitioned time interval. 
  
   Using parameter settings given in Table \ref{tab: profitparameters}, 
   we use PASA to solve for problem (\ref{eqn: profitpen}) for varying values of penalty parameter $\rho \in\{10^{-9},10^{-8}, 10^{-7}, 10^{-6},10^{-5},10^{-4}, 10^{-3}, 10^{-2}, 10^{-1}\}$. 
   Our initial guess for problem (\ref{eqn: profitpen}) was $u(t)=0$ for all $t\in[0,T]$. 
   We partition $[0,T]$ to where there are $N=750$ mesh intervals, and we use the discretization method described in Subsection \ref{subsec: profitdiscretesummary}.
   Additionally, the stopping tolerance is set to $10^{-10}$.
   Results corresponding to the penalized solutions that PASA obtained can be found in Figure \ref{fig: profitvarying} and Table \ref{tab: profitvary}. 
   Note that from looking at the  plots of the penalized controls  of  Figure \ref{tab: profitvary} alone, we find that $\rho=10^{-2}$ to be the most appropriate penalty parameter value.  
   We believe that even if we did not know what the explicit formula for the singular case was for problem (\ref{eqn: profitpen}), the plots of the penalized controls would have led us to choose  $p=10^{-2}$  as  being the appropriate penalty parameter because no oscillations are occurring  and the end behavior of the penalized control matches the non-oscillating region from the unpenalized control.  
   This yields some evidence that we could potentially use  PASA for solving penalized control problems without apriori information of the optimal control. 
   
	We will however continue to   explain how we ruled out all cases except for $\rho=10^{-2}$. 
   In Sub-figure \ref{fig: profitvary8}, we rule out penalty parameter value $\rho=10^{-8}$, since its corresponding solution $u_{\rho}$ is oscillating in the same manner as the unpenalized solution, whose plot is given in Sub-figure \ref{fig: profitvary0}.
   We did not provide plots of penalized solution $u_{\rho}$ that corresponded to when $\rho=10^{-9},10^{-7}, 10^{-6}$ because their plots possessed many oscillations similar to Sub-figure \ref{fig: profitvary8}. 
   Additionally, one could look at the $\norm{u^*-u_{\rho}}_{L^1}$ column in Table \ref{tab: profitvary} to see that the penalized solution for $u_{\rho}$ for $\rho\in\{10^{-9}, 10^{-8},10^{-7}, 10^{-6}\}$ are negligible in comparison to the unpenalized solution. 
   In Sub-figures \ref{fig: profitvary5} and \ref{fig: profitvarydiff5}, we start to see some improvements in the singular region of $u_{\rho}$ when $\rho=10^{-5}$.
   However, due to the oscillations appearing in time interval $(8.5,9.5)$, we rule out $\rho=10^{-5}$.
   We observe in Sub-figure {\ref{fig: profitvary4}} that increasing the penalty parameter $\rho$ to $10^{-4}$ significantly reduces the number of oscillations along the singular region. 
   However, around time interval $[9.24,9.37]$, $u_{\rho}(t)\approx 0.2366$, which is a value different from the bounds of the control and the singular case solution. 
   We rule out penalty $\rho=10^{-4}$ for this reason. 
   We did not provide a figure of the solution we obtained when penalty parameter value $\rho=10^{-3}$,  but we rule out this penalty parameter value due to similar reasoning that was used when $\rho=10^{-4}$. 
   We also rule out tuning parameter $\rho=10^{-1}$ because as illustrated in Sub-figure \ref{fig: profitvary1}, $u_{\rho}$ is constant around $[9.36,10]$ with constant value  being approximately $0.7277$. 
   For this case, we say that $\rho=10^{-1}$ is over-penalizing the control along the non-singular region. 
   Based on Table \ref{tab: profitvary} and Sub-figures \ref{fig: profitvary2} and \ref{fig: profitvarydiff2}, we find that penalty parameter $\rho=10^{-2}$ as being the most appropriate penalty parameter. 
   For $\rho=10^{-2}$, we have that $u_{\rho}$ switches to the non-singular case at the node that is closest to the actual switching point. 
   In addition, $u_{\rho}$ have the lowest $L^{1}$ norm error in Table \ref{tab: profitvary} when the penalty parameter is set to $\rho=10^{-2}$.
   In Figure \ref{fig: profitx2}, we have a plot of the state solution corresponding to $u_{\rho}$ when $\rho=10^{-2}$. 
   Observe that the plot of the approximated fish population corresponding to $u_\rho$ lines up almost perfectly with the true solution. 
   
   After finding an appropriate penalty parameter value for problem (\ref{eqn: profitpen}), we would like to study the numerical convergence rate between PASA's solution to the discretized penalized problem with tuning parameter $\rho=10^{-2}$) and the exact solution to problem (\ref{eqn: minprofit}). 
   We will observe the $L^1$ norm error between the exact solution $u^*$ and the penalized solution that PASA obtained $u_h$ for varying mesh sizes $h\in\{0.2,0.1,0.05,0.025,0.0125,0.00625,0.003125\}$.
   In Table \ref{tab: profitconverge}, we have $\text{err}_h$ representing the $\norm{u^*-u_h}_{L^1}$. 
   Notice that for almost all mesh size values the $L^1$ norm error between $u_h$ and $u^*$ is relatively close to the mesh size value. 
   However, in the case when $h=0.05$, we have that the $L^1$ norm error being significantly less than the mesh size. 
   This is because the mesh size influences the discretization of time interval $[0,T]$ to where one mesh point is exceedingly close to the the true switching point value. 
   We look at the last column of Table \ref{tab: profitconverge} to determine what the rate of 
   convergence between the penalized solution and the exact solution is.
   All but two entries in the last column of Table \ref{tab: profitconverge} were values that were slightly greater than one, and the entries that were not close to one involved using the $L^1$ norm error for when $h=0.05$.
   Based on the values of the last column of Table \ref{tab: profitconverge}, we find the convergence rate as being slightly better than a linear rate. 
   Additionally, form Table \ref{tab: profitconverge}, we take the natural logarithm of values found in columns 2 and 3, and use least squares method to indicate if there were a linear relationship between $\ln{(h)}$ and $\ln{(\text{err}_h)}$.
   Since we view data point {$(\ln(0.05),\ln(0.00089762))$} as being an outlier, we also perform least squares  between $\ln(h)$ and $\ln(\text{err}_h)$ without the outlier point. 
   In Figure \ref{fig: profitlogfit} the blue line does a better job with fitting to the data than the red line. 
   In Table \ref{fig: profitlogfit}, we see that the goodness of fit associated with the blue line is very strong. 
   We look at the slope associated with the blue line to determine the rate of convergence. 
   Since the slope of the blue line is approximately $1.200738$, we have further indication that the rate of convergence between the penalized solution and the exact solution is better than linear.  




\section{Example 2: Plant Problem}
 In this section we implement a biological optimal control problem from
 David King and Jonathan Roughgarden \cite{King1982} that can be solved analytically.  
 In \cite{King1982}, King and Roughgarden use optimal control theory to study the allocation strategies that annual plants possess when distributing photosynthate to components of the plant that pertain to \emph{vegetative growth} and to components that pertain to \emph{reproductive growth}.
 The control variable $u(t)$ involved represents the fraction of photosynthate being reserved for vegetative growth. 
 The remaining photosynthate, $1-u(t)$, is assigned to aid in the reproductive growth processes.
 The annual plants that are being modeled are assumed to live in an environment where consecutive seasons are identical. 
 Since $u(t)$ represents a fraction, the control variable is assumed to be bounded in the obvious way:
 \begin{equation}\label{eqn: plantcontrolbounds}
     0\leq u(t)\leq 1  \quad \text{for all } t\in[0,T],
 \end{equation}
 where $T$ represents the maximum season length. 
 The optimal control problem has two state variables where one state $x_1$ represents the weight of the vegetative part of a plant, while the second state $x_2$ represents the weight of the reproductive part of a plant. 
 The construction of the state equations is a continuous version of a discrete time model of reproduction for annual plants that was developed by Cohen \cite{Cohen1971} and is given by 
 \begin{align}
     x_1'(t) &= u(t)x_1, \quad \quad \quad\;\;x_1(0)>0,\label{eqn: plantstate1}\\
     x_2'(t) &= (1-u(t))x_2, \quad x_2(0)\geq 0.\label{eqn: plantstate2}
 \end{align}
 
  For the construction of an objective functional,  King and Roughgarden \cite{King1982} wanted to find the strategy that maximizes total reproduction and asserted that such a strategy is ``expected in organisms if natural selection maximizes total reproduction". 
Their construction of a cost functional that measures the total reproduction produced by a plant is based on Cohen's discrete model  in \cite{COHEN} which uses the expectation of the logarithm of seed yield  to express the long-term rate of population increases observed in annual plants. 
{The objective functional was chosen to have natural logarithm of $x_2$, which helped this problem to have singular arcs.} 
  The optimal control problem is then
  \begin{equation}\label{eqn: maxplant}
     \max\limits_{u\in\mathscr{A}}\; J(u) = \int_{0}^{T} \ln{(x_2(t))} dt
  \end{equation}
  subject to the state equations (\ref{eqn: plantstate1})-(\ref{eqn: plantstate2}) and control bounds (\ref{eqn: plantcontrolbounds}).
  The set of the admissible controls is defined $\mathscr{A}=\{u\in L^1(0,T):  \text{ for all } t\in [0,T], \; u(t) \in A,\; A=[0,1]\}$.
  Existence of an optimal control follows from Filippov-Cesari Existence Theorem \cite{Filippov}. 
  An explicit solution is obtained for problem (\ref{eqn: maxplant}) that satisfies Pontryagin's Maximum Principle \cite{Pontryagin}, the Generalized Legendre Clebsch Condition or Kelley's Condition \cite{Lewis, Robbins, Kelley, Zelikin1994}, the Strengthened Generalized Legendre Clebsch Condition, and McDanell and Powers' Junction Theorem \cite{Powers}. 
  The structure of the optimal allocation strategy $u^*$ depends on the maximum season length $T$, the initial weight of the vegetative part of the plant $x_1(0)$, and the initial weight of the reproductive part of the plant $x_2(0)$. 
  
  In this section, we provide a summary of the solution to the equivalent minimization problem: 
  %
  \begin{equation}\label{eqn: plantproblem}
     \begin{array}{rl}
        \min & -J(u)=-\int_0^T\ln{(x_2(t))}dt\\
         \textrm{s.t.}& \dot x_1(t) = u(t)x_1, \\
         \textrm{   } & \dot x_2(t)=(1-u(t))x_1,\\
	\textrm{} &x_1(0)=x_{1,0}>0,\quad x_2(0)=x_{2,0}\geq 0,\\
         \textrm{  } & 0\leq u(t)\leq 1.\\
    \end{array}
  \end{equation}
  Additionally we demonstrate how to discretize the penalized version of problem (\ref{eqn: plantproblem}) 
  \begin{equation}\label{eqn: penplant}
       \begin{array}{rl}
        \min & J_{\rho}(u)=-\int_0^T\ln{(x_2(t))}dt +\rho V(u)\\
         \textrm{s.t.}& \dot x_1(t) = u(t)x_1, \\
         \textrm{   } & \dot x_2(t)=(1-u(t))x_1,\\
	\textrm{} &x_1(0)=x_{1,0}>0,\quad x_2(0)=x_{2,0}\geq 0,\\
         \textrm{  } & 0\leq u(t)\leq 1, \\
    \end{array}
  \end{equation}
  where $0\leq \rho <1$ is the bounded variation tuning parameter and $V(u)$ measures the total variation of control $u$, as defined in equation (\ref{eqn: profitV}). 
  We then use PASA to numerically solve for problems (\ref{eqn: plantproblem}) and (\ref{eqn: penplant}) when parameters $x_{1,0}$, $x_{2,0}$, and $T$ are set in such a way where $u^*$ contains a  singular subarc. 
  \subsection{ Solving for the Singular Case to the Plant Problem}\label{subsec: Plantanalytic}
    Before we present a summary of the solution to problem (\ref{eqn: plantproblem}), we would like to mention some terminology that is used in King and Roughgarden's work \cite{King1982} to describe the events when $u^*(t)=0$ and when $u^*(t)=1$. 
    We say that $u^*$ is \emph{purely reproductive} at time $t$ when $u^*(t)=0$. 
    We say that  $u^*$ is \emph{purely vegetative} at time $t$ when $u^*(t)=1$.
    Our method for obtaining a solution to the minimization problem (\ref{eqn: plantproblem}) will parallel with what was done in  \cite{King1982}. 
    Our method for solving problem (\ref{eqn: plantproblem}) involves Pontryagin's Minimum Principle \cite{Pontryagin}. 
    We compute the Hamiltonian to problem (\ref{eqn: plantproblem})
    \begin{equation}\label{eqn: planthamiltonian}
        H = -\ln{x_2} +x_1(\lambda_1-\lambda_2)u+x_1\lambda_2. 
    \end{equation}
    We take the partial derivatives of the Hamiltonian with respect to $x_1$ and $x_2$ to construct differential equations for costate variables $\lambda_1$ and $\lambda_2$: 
    \begin{align}
         \lambda_1'(t) &= -\frac{\partial H}{\partial x_1} = (\lambda_2-\lambda_1)u -\lambda_2, 
         \label{eqn: costate1plant}\\
        \lambda_2'(t) &=  -\frac{\partial H}{\partial x_2} = \frac{1}{x_2}, \label{eqn: costate2plant}
    \end{align}
    and use the transversality conditions to obtain boundary conditions for the costate equations
    \begin{equation}
        \label{eqn: plantterminal}
            \lambda_1(T)=\lambda_2(T)= 0.
    \end{equation}
    We construct the switching function $\psi(t)$ by taking the partial derivative of the Hamiltonian with respect to control $u$:
    \begin{equation}\label{eqn: plantswitch}
     \psi(t) = x_1(\lambda_1-\lambda_2)
    \end{equation}
    {Since $x_1(t)>0$ for all $t$, the sign of $\psi$ will depend on the sign of $\lambda_1-\lambda_2$. 
    By Pontryagin's Minimum Principle, the optimal allocation strategy $u^*(t)$ will be as follows }
    \begin{equation}\label{eqn: plantucases}
        u^*(t) = 
        \begin{cases}
        0 & \text{whenever }\lambda_1(t) > \lambda_2(t),\\
        \text{singular} &\text{whenever } \lambda_1(t)=\lambda_2(t),\\
        1 & \text{whenever } \lambda_1(t) < \lambda_2(t).
        \end{cases}
    \end{equation}
    
    To solve for the singular case we assume that there is an interval $I\subset[0,T]$ where $\lambda_1(t)=\lambda_1(t)$ for all $t\in I$. 
    This implies that $ \lambda_1'(t)\equiv \lambda_2'(t)$ on $I$.  
    From adjoint equations (\ref{eqn: costate1plant})-(\ref{eqn: costate2plant}), we have that the following holds on $I$
    \begin{equation*}
        (\lambda_2(t) -\lambda_1(t))u(t)-\lambda_2(t) = \frac{1}{x_2(t)} ,
    \end{equation*}
    which yields the following on interval $I$
    \begin{equation}
        \label{eqn: padj2}
        \lambda_2(t) = -\frac{1}{x_2(t)}.
    \end{equation}
    It follows from equations (\ref{eqn: costate2plant}) and (\ref{eqn: padj2}) that 
    \begin{equation*}
        \label{eqn: dlamb2}
        \dot \lambda_2 =-\lambda_2 
    \end{equation*}
    for all $t\in I$. 
    Solving for the above differential equation yields
    \begin{equation}
        \label{eqn: plantadjoint2}
        \lambda_2(t) = Ae^{-t}, 
    \end{equation}
    where $A = -\frac{1}{x_2(t_1)}e^{t_1}$ and $t_1$ is the time when the singular subarc begins. 
    We can construct a differential equation for $x_2(t)$ on interval $I$ 
    by differentiating equation (\ref{eqn: padj2})  with respect to time and rearranging terms:
    \begin{equation}
    \label{eqn: dplantstate2}
    \dot x_2 = x_2.
    \end{equation}
    Equation (\ref{eqn: plantstate2}) and (\ref{eqn: dplantstate2}) yields 
    \begin{equation*}
        (1-u(t))x_1= x_2,
    \end{equation*}
    which implies that 
    \begin{equation}
        \label{eqn: plantuform}
        \frac{x_2(t)}{x_1(t)} = 1-u(t)
    \end{equation}
    on interval $I$.
    We differentiate  (\ref{eqn: plantuform}) to get 
    \[
         -\dot u(t)  = \frac{x_1\dot x_2-\dot x_1 x_2}{x_1^2}.
    \]
    We replace $\dot x_1$ and $\dot x_2$ in the above equation with state equation (\ref{eqn: plantstate1}) and equation (\ref{eqn: dplantstate2}), respectively, and obtain the following
    \begin{equation}
         \label{eqn: plantdustep}
        -\dot u(t)=\frac{x_2}{x_1}(1-u(t)). 
    \end{equation}
    Substituting (\ref{eqn: plantuform}) into (\ref{eqn: plantdustep}) yields
     \begin{equation}
        \label{eqn: plantusep}
        \frac{\dot u(t)}{(1-u(t))^2} = -1.
     \end{equation}
     We solve for the separable equation above and obtain the following solution 
      \begin{equation}
        \label{eqn: plantusingsoln}
        u(t) = 1- \frac{1}{C-t},
    \end{equation}
     on interval $I$, where $C$ is some constant. 
     Observe that from equation (\ref{eqn: plantusingsoln}), (\ref{eqn: plantuform}) becomes 
    \begin{equation}
        \label{eqn: plantratio}
        \frac{x_2(t)}{x_1(t)} = \frac{1}{C-t},
     \end{equation}
    on interval $I$. 
    
     In order to find constant term $C$ we need to find the time, $t_2$, when the singular subarc ends.
    As mentioned in \cite{King1982}, we can be certain that $u^*$ is not singular at time $T$  because the singular case solution for adjoint variable $\lambda_2(t)$ given in equation (\ref{eqn: padj2}) would not satisfy the transversality condition (\ref{eqn: plantterminal}).
    Now from the adjoint equations (\ref{eqn: costate1plant}) and (\ref{eqn: costate2plant}) and the transversality conditions (\ref{eqn: plantterminal}), we have that $\dot \lambda_1(t) \to 0$ and $\dot \lambda_2(t) \to \frac{1}{x_2(T)}$ as $t\to T$. 
    The derivatives of $\lambda_1$ and $\lambda_2$ near $T$ and the transversality conditions imply that $\lambda_1(t) >\lambda_2(t)$ for some  interval, namely $(t_2,T]$, and so from (\ref{eqn: plantucases}), we have that $u^*(t) \equiv 0$ on $(t_2,T]$. 
    The optimal control will always be purely reproductive on $(t_2,T]$ and does not depend on $T$ nor the initial conditions $x_1(0)$ and $x_2(0)$.      
    
    In order to find $t_2$ we must first solve for the state equations and the associated adjoint equations on time interval $[t_2,T]$. 
    From the state equations (\ref{eqn: plantstate1})-(\ref{eqn: plantstate2}), $u(t)\equiv 0$ on $[t_2,T]$ implies the following
    \begin{equation}
    \label{eqn: nonsingstatesol}
        \begin{array}{rl}
         x_1(t) &= x_1(t_2) ,\\
         x_2(t) &= x_1(t_2)(t-t_2)+x_2(t_2),  
         \end{array}
    \end{equation}  
    where 
    \[
     x_1(t_2) = (C-t_2)x_2(t_2), \text{ and }
     x_2(t_2) = x_2(t_1)e^{t_2-t_1}.
    \]
    Differentiating $x_2(t)$ given in equation (\ref{eqn: nonsingstatesol}) yields 
    \begin{equation} \label{eqn: plantdx2sep}
        \dot x_2(t) = x_1(t) = x_1(t_2).
    \end{equation}
 {
 We use the differential found from the above separable equation and the transversality condition (\ref{eqn: plantterminal}) to solve for adjoint equation (\ref{eqn: costate2plant}). 
}
  
    \begin{align}
        \lambda_2(T)-\lambda_2(t) &= 
 \int\limits_{x_2(t)}^{x_2(T)} \frac{1}{x_2} \frac{dx_2}{x_1(t_2)}\nonumber\\
 \lambda_2(t) &= \frac{1}{x_1(t_2)}\ln{\left(\frac{x_2(t)}{x_2(T)}\right)}  \text{ for } t\in [t_2,T].\label{eqn: plantnonsinglam2}
    \end{align}
    We use equation (\ref{eqn: plantnonsinglam2}) and adjoint equation (\ref{eqn: costate1plant})  to obtain the following differential equation for $\lambda_1$  
    \begin{equation}
        \label{eqn: dlamb1}
        \dot \lambda_1(t)  = -\frac{1}{x_1(t_2)}\ln{\left(\frac{x_2(t)}{x_2(T)}\right)}.
     \end{equation}
   {Again we use the differential obtained from separable equation (\ref{eqn: plantdx2sep}) and the terminal condition for $\lambda_1$ (\ref{eqn: plantterminal})  to solve the above differential equation:}
  \begin{align*}
    	\lambda_1(T)-\lambda_1(t) &= -\frac{1}{(x_1(t_2))^2}\int\limits_{x_2(t)}^{x_2(T)}( \ln{x_2}-\ln{x_2(T)}) dx_2 \\
	-\lambda_1(t) &= -\frac{1}{(x_1(t_2))^2}\left[x_2(t)-x_2(T)+x_2(t)\ln\left(\frac{x_2(T)}{x_2(t)}\right)\right].
  \end{align*}
  By negating the above equality and applying logarithmic rules, we obtain a solution for $\lambda_1(t)$:
\begin{equation}
  \label{eqn: plantlambda1}
  \lambda_1(t) = -\frac{1}{(x_1(t_2))^2}\left[x_2(t)\ln{\left(\frac{x_2(t)}{x_2(T)}\right)}+x_2(T)-x_2(t)\right] \;\text{ for } t\in [t_2,T]. 
\end{equation}

    After finding the solutions for the state and adjoint equations on interval $(t_2,T]$, we are ready to solve for $t_2$. 
    Since  $\lambda_1(t)$ and $\lambda_2(t)$ are continuous on $[0,T]$ and $\lambda_1(t)\equiv \lambda_2(t)$ on the singular case, we have $\lambda_1(t_2)=\lambda_2(t_2)$ 
    We equate equations (\ref{eqn: plantlambda1}) and (\ref{eqn: plantnonsinglam2}) at $t=t_2$ and solve for $t_2$:
    \begin{equation*}
         -\frac{1}{(x_1(t_2))^2}\left[x_2(t_2)\ln{\left(\frac{x_2(t_2)}{x_2(T)}\right)}+x_2(T)-x_2(t_2)\right] = -\frac{1}{x_1(t_2)}\ln{\left(\frac{x_2(t_2)}{x_2(T)}\right)}.
    \end{equation*}
    We multiply the above equation by $x_1(t_2)$ and rearrange terms to obtain the following equation:
    \begin{equation}
        \label{eqn: x2minusx2}
         x_2(T)-x_2(t_2) = (x_1(t_2)+x_2(t_2)) \ln{\left(\frac{x_2(T)}{x_2(t_2)}\right)}.
    \end{equation}
    We use (\ref{eqn: nonsingstatesol}) to evaluate $x_2(T)$, which yields $x_2(T)= x_1(t_2)(T-t_2)+x_2(t_2)$, and apply it to equation (\ref{eqn: x2minusx2}):
\begin{equation}
  \label{eqn: planttimeeq1}
  T-t_2= \left[1+\frac{x_2(t_2)}{x_1(t_2)}\right]\left[\ln{\left(\frac{x_1(t_2)}{x_2(t_2)}(T-t_2)+1\right)}\right].
\end{equation} 
We also have $\dot \lambda_1(t_2)= \dot \lambda_2(t_2)$  so using (\ref{eqn: padj2}), (\ref{eqn: dlamb1}), and  $x_2(T)= x_1(t_2)(T-t_2)+x_2(t_2)$ we find that
\begin{equation} 
 \label{eqn: planttimeeq2}
 \frac{x_1(t_2)}{x_2(t_2)} = \ln{\left(\frac{x_1(t_2)}{x_2(t_2)}(T-t_2)+1\right)}.
\end{equation} 
 {We rewrite equations (\ref{eqn: planttimeeq1})  and (\ref{eqn: planttimeeq2}) into terms of  $y$ and $z$, where $y=T-t_2$ and $z= \frac{x_2(t_2)}{x_1(t_2)}$, to obtain a nonlinear system of equations that can be solved numerically through a nonlinear solver. Consequently, we have the following:}



\begin{align}
y&= T- t_2 \approx 2.79328213,\quad\text{ and } \label{eqn: tstarform}\\
z&=\frac{x_2(t_2)}{x_1(t_2)}\approx 0.55763674.\label{eqn: ratioformula}
\end{align}
Using $x_2(T)= x_1(t_2)(T-t_2)+x_2(t_2)$ and  $x_1(t) = x_1(t_2) $ for $t\in [t_2,T]$ yields  

\begin{equation}
  \frac{x_2(T)}{x_1(T)}= (T-t_2)+\frac{x_2(t_2)}{x_1(t_2)}.
\end{equation}
From equations (\ref{eqn: tstarform}) and (\ref{eqn: ratioformula}) we then have 
\begin{equation}
  \label{eqn: terminalratio}
  \frac{x_2(T)}{x_1(T)} \approx 3.35091887.
\end{equation}
Additionally, from equations (\ref{eqn: planttimeeq1}) and (\ref{eqn: planttimeeq2}) we have 
\begin{equation}
\label{eqn: Tmint*}
 T-t_2 = \left[1+ \frac{x_2(t_2)}{x_1(t_2)}\right]\frac{x_1(t_2)}{x_2(t_2)}= 1+\frac{x_1(t_2)}{x_2(t_2)}, 
\end{equation} 
which implies that 
\begin{equation}
\label{eqn: switchcriteria}
 \frac{x_2(t_2)}{x_1(t_2)}=\frac{1}{(T-1)-t_2}.
\end{equation}
Finally, we can evaluate equation (\ref{eqn: plantratio}) at $t=t_2$ and use the above equation to solve for constant $C$, and we have that $C=T-1$. 
\subsubsection{Summary of Singular Case Solution to Plant Problem}\label{subsec: plantsummary}
In Subsection \ref{subsec: Plantanalytic}, we found that if the optimal control $u^*$ to problem (\ref{eqn: plantproblem}) contains a singular subarc where $u^*$ becomes singular at time $t_1$, then $u^*$ will be the following on the interval $[t_1,T]$: 
\begin{equation}\label{eqn: plantsingsumu}
 u^*(t)=
 \begin{cases}
  1-\frac{1}{T-1-t} &\quad\text{ when } t\in[t_1,t_2),\\
  0,	& \quad \text{ when }t\in [t_2,T],
 \end{cases}
\end{equation}
where $t_2\approx T-2.79328213$. 
The corresponding solution to the state equations $x_1$ and $x_2$ are as follows: 
\begin{equation}\label{eqn: plantsingsumx1}
	x_1(t) = 
	\begin{cases}
		(T-1-t)x_2(t) & \quad \text{ when } t\in[t_1,t_2],\\
		x_1(t_2) & \quad \text{ when } t\in [t_2,T],
	\end{cases}
\end{equation}
\begin{equation}\label{eqn: plantsingsumx2}
	x_2(t) = 
	\begin{cases}
		x_2(t_1)e^{t-t_1}& \quad \text{ when } t\in [t_1,t_2],\\
		x_1(t_2)(t-t_2)+x_2(t_2) & \quad \text{ when } t\in [t_2, T],
	\end{cases}
\end{equation}
where $x_1(t_2) =(T-1-t_2)x_2(t_2)$ and $x_2(t_2)=x_2(t_1)e^{t_2-t_1}$. 
The corresponding solution to the adjoint variables $\lambda_1$ and $\lambda_2$ are as follows:
\begin{equation}\label{eqn: plantsingsumlamb1}
	\lambda_1(t)=
	\begin{cases}	
	\lambda_2(t) & \quad \text{ when } t\in [t_1,t_2],\\
	 -\frac{1}{(x_1(t_2))^2}\left[x_2(t)\ln{\left(\frac{x_2(t)}{x_2(T)}\right)}+x_2(T)-x_2(t)\right] & \quad \text{ when } t\in [t_2,T],
	\end{cases}
\end{equation}
\begin{equation}\label{eqn: plantsingsumlamb2}
	\lambda_2(t)=
	\begin{cases}
	  -\frac{1}{x_2(t_1)}e^{t_1-t} &  \quad \text{ when } t\in [t_1,t_2],\\
	  \frac{1}{x_1(t_2)}\ln{\left(\frac{x_2(t)}{x_2(T)}\right)} & \quad \text{ when } t\in [t_2,T].
	\end{cases}
\end{equation}

 Additionally, King and Roughgarden's \cite{King1982} use the  
 generalized Legendre-Clebsch Condition \cite{Powers, Robbins, Lewis}, also referred  as Kelley's condition \cite{Zelikin1994, Kelley}, to show that the second order necessary condition of optimality is satisfied. 
 The generalized Legendre-Clebsch Condition involves finding what is called the \emph{order} of a singular arc which is defined as the being the integer $q$ such that $\left(\frac{d^{2q}}{dt^{2q}}\frac{\partial H}{\partial u}\right)$ is the lowest order total derivative of the partial derivative of the Hamiltonian with respect to $u$, in which control $u$ appears explicitly. 
 By using equations (\ref{eqn: plantstate1}), (\ref{eqn: plantstate2}), (\ref{eqn: costate1plant}), (\ref{eqn: costate2plant}), and (\ref{eqn: plantswitch}), we found that 
 \begin{align*}
     \frac{d}{dt}\frac{\partial H}{\partial u} &=  \frac{d}{dt} \psi = -x_1\lambda_2-\frac{x_1}{x_2},\\
     \frac{d^{2}}{dt^2}\frac{\partial H}{\partial u} & =\frac{d^{2}}{dt^2}\psi = \left(-x_1\lambda_2-\frac{x_1}{x_2}-\frac{x_1^2}{x_2^2}\right)u-\frac{x_1}{x_2}+\frac{x_1^2}{x_2^2}.
 \end{align*}
 So for problem (\ref{eqn: plantproblem}), the order of singular subarc $u^*$ is 1. 
 By General Legendre Clebsch condition if $u^*$ is an optimal singular control on some interval $[t_1,t_2)$ of order $q$, then it is necessary that 
 \begin{equation*}
     (-1)^{q} \frac{\partial}{\partial u}\left[\frac{d^{2q}}{dt^{2q}}\left(\frac{\partial H}{\partial u}\right)\right]\geq 0
 \end{equation*}
 if the extremum is a minimum, and the inequality is reversed if the extremum is a maximum \cite{Powers}. 
 Additionally, if the above inequality is strict then the strengthened generalized Legendre-Clebsch condition holds. 
 On singular region $[t_1,t_2)$ we have that
 \begin{align*}
     (-1)^{1} \frac{\partial}{\partial u}\left[\frac{d^{2}}{dt^{2}}\left(\frac{\partial H}{\partial u}\right)\right] = \frac{x^2_1}{x^2_2}>0,
 \end{align*}
 so the singular arc is of order 1 and satisfies both the generalized Legendre-Clebsch Condition and the strengthened Legendre-Clebsch Condition.

What we have yet to discuss is conditions for problem (\ref{eqn: plantproblem}) that determine when the optimal control $u^*$ will be a bang-bang or concatenations of bang and singular controls.  
In \cite{King1982}, King and Roughgarden constructed a series of conditions for determining the structure of the optimal control solution $u^*$ to problem (\ref{eqn: plantproblem}). 
Most of the conditions are based upon the value of terminal time $T$ and the ratio of the reproductive to vegetative weight at time $t=0$. We provide a summary of those conditions (see \cite{King1982} for details on how these conditions were constructed):
\begin{enumerate}
    \item The optimal allocation strategy cannot be purely reproductive on the entire time interval ($u^*(t)=0$ for all $t\in[0,T]$) unless $T\leq 3.35091887$ and $\frac{x_2(0)}{x_1(0)}\leq 3.35091887-T$. 
    \item If $T>3.35091887$ and $0\leq \frac{x_2(0)}{x_1(0)}<0.55763674e^{T-2.79328213}$, the optimal allocation, $u^*$, will contain a singular subarc. 
    \begin{enumerate}[2a)]
        \item\label{itm: case1} If $\frac{x_2(0)}{x_1(0)}=\frac{1}{T-1}$, then the optimal allocation strategy $u^*$ begins with a singular subarc and switches to being purely reproductive at time $t_2=T-2.79328213$. 
        \item \label{itm: case2}If $0\leq \frac{x_2(0)}{x_1(0)}<\frac{1}{T-1}$, then the optimal control $u^*$ contributes to reproductive growth before and after the singular subarc occurs. 
        \item \label{itm: case3} If $\frac{1}{T-1}<\frac{x_2(0)}{x_1(0)}<0.55763674e^{T-2.79328213}$, then $u^*$ begins contributing to purely vegetative growth before the singular subarc occurs. At time $t_2=T-2.79328213$, $u^*$ switches from being singular to being purely reproductive. 
    \end{enumerate}
    \item If $\frac{x_2{(0)}}{x_1(0)}\geq0.55763674e^{T-2.79328213} $, then the optimal control must be bang-bang, where $u^*$ begins with purely vegetative growth and switches once to being purely reproductive. 
\end{enumerate}
The value $\frac{1}{T-1}$  corresponds to evaluating  the right hand side of equation (\ref{eqn: plantratio}) at $t=0$. 
Additionally some of the numerical values shown in the above conditions corresponded to equations (\ref{eqn: tstarform}), (\ref{eqn: ratioformula}), and (\ref{eqn: terminalratio}).

We arel only interested in finding the exact solutions to problem (\ref{eqn: plantproblem}) in the cases where the optimal control contains a singular subarc, i.e. Cases  \ref{itm: case1}, \ref{itm: case2}, and \ref{itm: case3}.
For Case \ref{itm: case1}, we use equations (\ref{eqn: plantsingsumu})-(\ref{eqn: plantsingsumlamb2}) with $t_1=0$ to construct the exact solution which is as follows: \\
\textbf{Case 2a):} If $T>3.35091887$,  $0\leq\frac{x_2(0)}{x_1(0)}<0.55763674e^{T-2.79328213}$, and $\frac{x_2(0)}{x_1(0)}=\frac{1}{T-1}$, then 
\begin{align}
 u^*(t) &=
    \begin{cases}	
     	1-\frac{1}{T-t-1} &\; 0\leq t\leq t_2,\\
	0 & \;t_2< t\leq T,
    \end{cases}\label{eqn: plantsingsumabegin}\\
x_1(t) &= 
	\begin{cases}
  		x_2(t) (T-t-1) & \;0\leq t\leq t_2,\\
		x_1(t_2)	&	\; t_2\leq t\leq T,
	\end{cases}\\
x_2(t)&=
 	\begin{cases}
	x_2(0)e^t 	& 0 \leq t\leq t_2,\\
	x_1(t_2)(t-t_2)+x_2(t_2) & t_2\leq t\leq T,
	\end{cases}\\
\lambda_1(t) &= 
	\begin{cases}
		\lambda_2(t) &  0\leq t\leq t_2,\\
		-\frac{1}{(x_1(t_2))^2}\left[x_2(t)\ln{\left(\frac{x_2(t)}{x_2(T)}\right)}+x_2(T)-x_2(t)\right] & t_2\leq t\leq T,
	\end{cases}\\
\lambda_2(t)&=
	\begin{cases}
		-\frac{1}{x_2(0)}e^{-t} & \;0\leq t\leq t_2,\\
		\frac{1}{x_1(t_2)}\ln\left(\frac{x_2(t)}{x_2(T)}\right) &\; t_2 \leq t\leq T,
	\end{cases}\label{eqn: plantsingsumaend}
\end{align}
 where $t_2=T-2.79328213$, $x_1(t_2)=(T-1-t_2)x_2(t_2)$, and $x_2(t_2)=x_2(0)e^{t_2}$. 
 
  For Case \ref{itm: case2}, we solve for state equations (\ref{eqn: plantstate1}) and (\ref{eqn: plantstate2}) on the interval $[0,t_1]$  with $u$ set to being 0, and then we use continuity of state variables to get an explicit solution for $t_1$. 
  Then we solve for adjoint equations (\ref{eqn: costate1plant}) and (\ref{eqn: costate2plant}) on $[0,t_1]$ where we use equations (\ref{eqn: plantsingsumlamb1}) - (\ref{eqn: plantsingsumlamb2}) and continuity of $\lambda_1$ and $\lambda_2$ at $t_1$ to gain the following boundary conditions: $\lambda_1(t_1)=\lambda_2(t_1)=-\frac{1}{x_2(t_1)}$.\\
 The exact solution for Case \ref{itm: case2} is as follows: \\
 \textbf{Case 2b):} If $T>3.35091887$, $0\leq\frac{x_2(0)}{x_1(0)}< 0.55763674e^{T-2.79328213}$, and $\frac{x_2(0)}{x_1(0)}< \frac{1}{T-1}$, then 
\begin{align}
   u^*(t) &= 
   \begin{cases}
	0 & 0\leq t<t_1,\\
	1-\frac{1}{T-t-1} & t_1\leq t\leq t_2,\\
   	0 & t_2< t\leq T,
   \end{cases}\label{eqn: plantsingsumbbegin}\\
  x_1(t) &= 
   \begin{cases}
	x_1(0) & 0\leq t\leq t_1,\\
	(T-t-1)x_2(t) & t_1 \leq t\leq t_2, \\
	x_1(t_2) & t_2\leq t\leq T,
   \end{cases}\\
  x_2(t)  &= 
   \begin{cases}
	x_1(0)t+x_2(0) & 0\leq t\leq t_1,\\
	x_2(t_1)e^{t-t_1} & t_1\leq t\leq t_2,\\
 	x_1(t_2)(t-t_2)+x_2(t_2) & t_2\leq t\leq T,	
   \end{cases}\\
\lambda_1(t) & = 
\begin{cases}
        \frac{x_2(t)}{(x_1(0))^2}\ln{\left(\frac{x_2(t_1)}{x_2(t)}\right)}+\frac{x_2(t)-x_2(t_1)}{(x_1(0))^2}+\frac{x_2(t)-x_2(t_1)}{x_2(t_1)x_1(0)}-\frac{1}{x_2(t_1)} & 0\leq t\leq t_1,\\
	\lambda_2(t) & t_1\leq t\leq t_2,\\
 	-\frac{x_2(t)}{(x_1(t_2))^2}\left(\ln\left (\frac{x_2(t)}{x_2(T)}\right)-1\right)-\frac{x_2(T)}{(x_1(t_2))^2} & t_2\leq t\leq T,	
\end{cases} \\
\lambda_2(t) & =
\begin{cases}
	\frac{1}{x_1(0)}\ln\left(\frac{x_2(t)}{x_2(t_1)}\right)-\frac{1}{x_2(t_1)} & 0\leq t\leq t_1,\\
	-\frac{1}{x_2(t_1)}e^{t_1-t} & t_1\leq t\leq t_2,\\
 	\frac{1}{x_1(t_2)}\ln\left(\frac{x_2(t)}{x_2(T)}\right) & t_2\leq t\leq T,	
\end{cases}\label{eqn: plantsingsumbend}
\end{align}
where $t_2\approx T-2.7932821329007607$, $x_1(t_2) = x_2(t_2)(T-t_2-1)$, and $x_2(t_1)=x_1(0)t_1+x_2(0)$. 
 Also, switch $t_1$ is the solution to the following equation
 $$\frac{x_2(t_1)}{x_1(t_1)}=\frac{1}{T-t_1-1}. $$ 
 The solution to the above equation is the following expression 
 $$t_{1,\pm} = \frac{1}{2}\left[
				T-1-\frac{x_2(0)}{x_1(0)}\pm \sqrt{T^2+2T\left(\frac{x_2(0)}{x_1(0)}-1\right)-3+\left(\frac{x_2(0)}{x_1(0)}\right)^2}
			\right],$$
	and we choose $t_1$ to be the solution that is between the values $0$ and $t_2$. 
	
  For Case \ref{itm: case3}, we solve for state equations (\ref{eqn: plantstate1}) and (\ref{eqn: plantstate2}) on the interval $[0,t_1]$  with $u$ set to being 1, and then use continuity of state variables to obtain an equation that can be used to solve for $t_1$.  
  Then we solve for adjoint equations (\ref{eqn: costate1plant}) and (\ref{eqn: costate2plant}) on $[0,t_1]$ where we used equations (\ref{eqn: plantsingsumlamb1}) - (\ref{eqn: plantsingsumlamb2}) and continuity of $\lambda_1$ and $\lambda_2$ at $t_1$ to gain the following boundary conditions: $\lambda_1(t_1)=\lambda_2(t_1)=-\frac{1}{x_2(t_1)}$.\\
 The exact solution for Case \ref{itm: case3} is as follows: \\
 \textbf{Case 2c):} If $T>3.35091887$ and 
 $\frac{1}{T-1}<\frac{x_2(0)}{x_1(0)}\leq0.55763673e^{T-2.79328213}$,
 then
\begin{align}
   u^*(t) &= 
   \begin{cases}
	1 & 0\leq t< t_1,\\
	1-\frac{1}{T-t-1} & t_1\leq t\leq t_2,\\
 	0 & t_2< t\leq T,
   \end{cases}\label{eqn: plantsingsumcbegin}\\
  x_1(t) &= 
   \begin{cases}
	x_1(0)e^{t} & 0\leq t\leq t_1,\\
	x_2(t)(T-t-1) & t_1\leq t\leq t_2,\\
 	x_1(t_2) & t_2\leq t\leq T,	
   \end{cases}\\
x_2(t)  &= 
\begin{cases}
	x_2(0) & 0\leq t\leq t_1,\\
	x_2(t_1)e^{t-t_1} & t_1\leq t\leq t_2,\\
 	x_1(t_2)(t-t_2)+x_2(t_2) & t_2\leq t\leq T,	
\end{cases}\\
\lambda_1(t) & =
\begin{cases}
	-\frac{1}{x_2(0)}e^{t_1-t} & 0\leq t\leq t_1,\\
	\lambda_2(t) & t_1\leq t\leq t_2,\\
 	-\frac{x_2(t)}{(x_1(t_2))^2}\left(\ln\left (\frac{x_2(t)}{x_2(T)}\right)-1\right)-\frac{x_2(T)}{(x_1(t_2)^2} & t_2\leq t\leq T,	
\end{cases}\\
\lambda_2(t) & = 
\begin{cases}
	\frac{1}{x_2(0)}(t-t_1)-\frac{1}{x_2(t_1)} & 0\leq t\leq t_1,\\
	-\frac{1}{x_2(0)}e^{t_1-t} & t_1\leq t\leq t_2,\\
 	\frac{1}{x_1(t_2)}\ln\left(\frac{x_2(t)}{x_2(T)}\right) & t_2\leq t\leq T,	
\end{cases} \label{eqn: plantsingsumcend}
\end{align}
where 
$t_2=T-2.79328213$, 
$x_1(t_1)=x_1(0)e^{t_1}$, $x_2(t_1)=x_2(0)$, 
$x_1(t_2)=x_2(t_2)(T-t_2-1)$, and $x_2(t_2)=x_2(t_1)e^{t_2-t_1}$.
Also, switch $t_1$ can be obtained by using a nonlinear solver for the following non-linear equation: \[\frac{x_2(0)}{x_1(0)}e^{-t_1}=\frac{1}{T-1-t_1}.\] 

\subsection{Discretization of Plant Problem}\label{subsec: plantdiscrete}
 We discretize the following penalized problem 
 \begin{equation}\label{eqn: penplant1}
       \begin{array}{rl}
        \min & J_{\rho}(u)=-\int_0^T\ln{(x_2(t))}dt +\rho V(u)\\
         \textrm{s.t.}& \dot x_1(t) = u(t)x_1, \\
         \textrm{   } & \dot x_2(t)=(1-u(t))x_1,\\
	\textrm{} &x_1(0)=x_{1,0}>0,\quad x_2(0)=x_{2,0}\geq 0,\\
         \textrm{  } & 0\leq u(t)\leq 1, \\
    \end{array}
  \end{equation}
  where $0\leq \rho <1$ is the bounded variation penalty parameter and $V(u)$ measures the total variation of control $u$,  which is defined in equation (\ref{eqn: profitV}).
  We would also like to discretize the following adjoint equations: 
  \begin{align}
         \lambda_1'(t) &= -\frac{\partial H}{\partial x_1} = (\lambda_2-\lambda_1)u -\lambda_2, 
         \label{eqn: costate1plantagain}\\
        \lambda_2'(t) &=  -\frac{\partial H}{\partial x_2} = \frac{1}{x_2}, \label{eqn: costate2plantagain}
    \end{align}
    with $\lambda_1(T)=\lambda_2(T)=0$ being the transversality conditions. 
    
    We assume the control $u$ is constant over each mesh interval.
    We partition time interval $[0,T]$, by using $N+1$ equally spaced nodes, $0=t_0<t_1<\cdots<t_N=T$. 
    For all $k=0,1,\dots,N$ we assume that $x_{1,k}=x_1(t_k)$ and $x_{2,k}=x_2(t_k)$. 
    For the control we denote $u_k=u(t)$ for all $t_k\leq t<t_{k+1}$ when $k=0,\dots N-2$ and $u_{N-1}=u(t)$ for all $t_{N-1}\leq t\leq t_N$. 
    So we have $\x_1\in\R^{N+1}$, $\x_2\in\R^{N+1}$ while $u\in\R^N$. We use a left-rectangular integral approximation for objective function $J_{\rho}$ in Problem (\ref{eqn: penplant1}).
    The discretization of problem (\ref{eqn: penplant1}) is then
    \begin{equation*}
     \begin{array}{rl}
         \min & J_{\rho}(\u) =\sum\limits_{k=0}^{N-1}(-h\ln{(x_{2,k})})+\rho\sum\limits_{k=0}^{N-2}|u_{k+1}-u_k|  \\
         & x_{1,k+1}=x_{1,k}+h(u_k x_{1,k})\;\text{ for all } k=0,\dots,N-1,\\
         & x_{2,k+1}=x_{2,k}+h((1-u_k)x_{1,k})\; \text { for all } k=0,\dots, N-1,\\
        & x_{1,0}>0,\; x_{2,0}\geq 0\\
        & 0\leq u_k\leq 1 \; \text{ for all } k=0,\dots, N-1,\\
     \end{array}
    \end{equation*}
 where $h=\frac{T}{N}$, is the mesh size and the first component of $\x_1$ and $\x_2$ is set to being the initial condition associated with the state equations. 
 Since PASA uses a gradient scheme for one of its phases, we need the cost functional to be differentiable. 
 We perform a decomposition on each absolute value term in $J_{\rho}$ to ensure that $J_{\rho}$ is differentiable. 
 We introduce two $N-1$ vectors $\bs{\zeta}$ and $\bs{\iota}$ whose entries are non-negative. 
 Each entry of $\bs{\zeta}$ and $\bs{\iota}$ will be defined as: 
{ 
\begin{equation*}
     |u_{k+1}-u_k| =\zeta_k+\iota_k \; \text{ for all } k=0,\dots, N-2.
 \end{equation*}
}
 With this decomposition in mind the discretized penalized problem will be the following:
 \begin{equation}\label{eqn: penplantdiscrete}
 \begin{array}{rl}
     \min & J_{\rho}(\u,\bs{\zeta},\bs{\iota})= \sum\limits_{k=0}^{N-1}(-h\ln{(x_{2,k})}) +\rho\sum\limits_{k=0}^{N-2}(\zeta_k+\iota_k) \\
     &x_{1,k+1} = x_{1,k}+h(u_k x_{1,k}) \; \text{ for all } k=0,\dots, N-1,\\
     &x_{2,k+1}= x_{2,k} + h(1-u_k)x_{1,k}\; \text{ for all } k=0,\dots, N-1,\\
     &x_{1,0}> 0\; x_{2,0}\geq 0,\\
     & 0\leq u_k\leq \;  \text{ for all } k=0,\dots, N-1,\\
     & u_{k+1}-u_k = \zeta_k-\iota_k \; \; \text{ for all } k=0,\dots, N-2,\\
     & \zeta_k\geq 0, \; \iota_k\geq 0 \; \; \text{ for all } k=0,\dots, N-1.\\
 \end{array}
 \end{equation}
 Notice that for problem (\ref{eqn: penplantdiscrete}), we are minimizing the penalized objective function with respect to three vectors, $\u,\bs{\zeta},$ and $\bs{\iota}$. 
 The equality constraints associated with $\bs{\zeta}$ and $\bs{\iota}$ are linear constraints that PASA can interpret. 
 The equality constraints associated with $\bs{\zeta}$ and $\bs{\iota}$ can be written like so:
 \begin{equation}\label{eqn: plantlinear}
  \left[
      \begin{array}{c | c | c}
      \bs{A} &-\bs{I}_{N-1}&\bs{I}_{N-1}
      \end{array}
   \right]
   \begin{bmatrix}
   \u\\
   \hline
   \bs{\zeta}\\
   \hline
   \bs{\iota}
   \end{bmatrix}
   =
   \bs{0},
 \end{equation}
 where $\bs{I}_{N-1}$ is the identity matrix with dimension $N-1$, 
 $\bs{0}$ is the $N-1$ dimensional all zeros vector, and {$\bs{A}$ is an $N-1\times N$ sparse matrix defined on equation (\ref{eqn: sparsematrix}).}
 For finding the gradient of $J_\rho$ in problem (\ref{eqn: penplantdiscrete}) we use Theorem \ref{thm: lagrange}. 
 In order to compute the Lagrangian of problem (\ref{eqn: penplantdiscrete}), we rewrite the discretized state equations accordingly: 
 \begin{align*}
     -x_{1,k+1}+x_{1,k}+h(u_kx_{1,k}) &=0, \text{ for all } k=0,\dots, N,\\
     -x_{2,k+1}+x_{2,k}+h(1-u_k)x_{1,k}&=0, \text{ for all } k=0,\dots, N. 
 \end{align*}
 The Lagrangian of problem (\ref{eqn: penplantdiscrete}) is then
 \begin{equation}
 \begin{split}
     \Ls{(\x_1,\x_2,\u,\bs{\zeta},\bs{\iota})} &= \sum\limits_{k=0}^{N-1} (-h\ln{(x_{2,k})}) +\rho\sum\limits_{k=0}^{N-2}(\zeta_k+\iota_k)\\ &+\sum\limits_{k=0}^{N-1}\lambda_{1,k}(-x_{1,k+1}+x_{1,k}+h(u_k x_{1,k}))+\sum\limits_{k=0}^{N-1}\lambda_{2,k}(-x_{2,k+1}+x_{2,k}+h((1-u_k)x_{1,k})),
 \end{split}
 \end{equation}
 Where $\bs{\lambda_1}, \bs{\lambda_2} \in \R^{N-1}$ are the Lagrangian multiplier vectors. 
 We take the partial derivative of $\Ls$ with respect to $u_k$ and obtain: 
 \[
 \frac{\partial\Ls}{\partial u_k} = h(x_{1,k}(\lambda_{1,k}-\lambda_{2,k})) \; \text{ for all } k=0,\dots, N-1.
 \]
 By using Theorem \ref{thm: lagrange}, we obtain the following: 
 \begin{equation}
     \grad_{\u,\bs{\zeta},\bs{\iota}}J_{\rho}=\grad_{\u,\bs{\zeta},\bs{\iota}}\Ls = 
     \begin{bmatrix}
      h( x_{1,0}(\lambda_{1,0}-\lambda_{2,0}))\\
      \vdots\\
      h (x_{1,k}(\lambda_{1,k}-\lambda_{2,k})\\
      \vdots\\
      h (x_{1,N-1})(\lambda_{1,N-1}-\lambda_{2,N-1})\\
      \hline
      \rho\\
      \vdots\\
      \rho\\
      \hline
      \rho\\
      \vdots\\
      \rho
     \end{bmatrix},
 \end{equation}
 provided that condition (\ref{eqn: thmadjoint}) in Theorem \ref{thm: lagrange} is satisfied. 
 To satisfy (\ref{eqn: thmadjoint}) we take the partial derivative of the Lagrangian $\Ls$ with respect to $x_{1,k}$ and $x_{2,k}$ for all $k=1,\dots, N$, and note that we are not taking the partial derivative of $\Ls$ with respect to $x_{1,0}$ and $x_{2,0}$ because they are known values. 
 Taking the partial derivative of $\Ls$ with respect to each state vector components yields the following expressions: 
 \begin{align}
    \frac{\partial \Ls}{\partial x_{1,k}}&= h(\lambda_{1,k}u_k+\lambda_{2,k}(1-u_k))+\lambda_{1,k}-\lambda_{1,k-1}\; \text{ for all } k=1,\dots, N-1 \label{eqn: plantpartialLx1}\\
    \frac{\partial \Ls}{\partial x_{1,N}}&=-\lambda_{1,N-1} \label{eqn: plantpartialLx1N}\\
    \frac{\partial \Ls}{\partial x_{2,k}} &= 
    -h\left(\frac{1}{x_{2,k}}\right)+ \lambda_{2,k}-\lambda_{2,k-1}
    \; \text{ for all } k=1,\dots, N-1\label{eqn: plantpartialLx2}\\\\
    \frac{\partial \Ls}{\partial x_{2,N}}&=-\lambda_{2,N-1}.\label{eqn: plantpartialLx2N}\\
 \end{align}
 For satisfying condition (\ref{eqn: thmadjoint}) from Theorem \ref{thm: lagrange}, we set the above equations equal to zeros and perform the following steps:  solve for $\lambda_{1,k-1}$ in equation (\ref{eqn: plantpartialLx1});  solve for $\lambda_{1,N-1}$ in equation (\ref{eqn: plantpartialLx1N}); solve for $\lambda_{2,k-1}$ in equation (\ref{eqn: plantpartialLx2}); and solve for $\lambda_{2,N-1}$ in equation (\ref{eqn: plantpartialLx2N}). 
 Consequently, we generate a discretization for the adjoint equations (\ref{eqn: costate1plant}) and (\ref{eqn: costate2plant}) and that produces the transversality conditions (\ref{eqn: plantterminal}): 
 \begin{align}
     \lambda_{1,k-1} &= \lambda_{1,k}+ h(\lambda_{1,k}u_k+\lambda_{2,k}(1-u_k)) \; \text{ for all } k=1,\dots , N-1,\\
     \lambda_{1,N-1}&= 0,\\
     \lambda_{2,k-1} &=\lambda_{2,k}-h\left(\frac{1}{x_{2,k}}\right) \; \text{ for all } k=1,\dots, N-1 \; \text{ and }\\
     \lambda_{2,N-1}&=0.
 \end{align}

 
\subsection{Numerical Results of Plant Problem}\label{subsec: plantnumerical}
  \begin{table}[ht]
  \centering
  \begin{tabular}{|c |l | c |c| c|}
     \hline
    \textbf{Parameter} &\textbf{Description}& \textbf{2a) Values}& \textbf{ 2b) Values}& \textbf{2c)Values}\\
    \hline
     $T$       & Terminal time                             & 5 & 5 & 5\\
     $x_{1,0}$ & Initial condition for Vegetative Weight   & 4 & 1& 1 \\
     $x_{2,0}$ & Initial condition for Reproductive Weight & 1& $10^{-4}$& 2\\
     \hline
  \end{tabular}
  \caption{Numerical Values of Parameters Used in Computations for Plant Problem}
  \label{tab: plantparameters}
  \end{table}
   \begin{table}[ht]
   \centering
   \resizebox{\textwidth}{!}{%
   \begin{tabular}{|l|l|l| l l| l l| l l|l|}
   \hline
   & $\norm{u^*-\hat{u}}_{L^1}$ & $\norm{u^*-\hat{u}}_{L^{\infty}}$ & $t_1^*$ &  $\hat{t}_1$ & $t^*_2$ & $\hat{t}_2$ & $-J(u^*)$ & $-J(\hat{u})$& \textbf{Runtime (s)}\\
   \hline
   Case 2a) & 0.017557 & 0.444444& NA& NA & 2.2067 &2.2 &11.787362 & 11.787496& 46.25  \\
   Case 2b) & 0.017519 &0.444444 & 0.2678& 0.2667 &  2.2067 & 2.2 & 3.600896  &3.600974 & 32.97 \\
   Case 2c) & 0.013588   &0.444444 & 1.5778& 1.5733 & 2.2067 &2.2 & 8.612962 & 8.613037 & 7.48 \\
   \hline
   \end{tabular}}
   \caption{Results from Unpenalized Solution to Plant Problem for Each Case. The starred notation corresponds to the exact solution while the hat notation corresponds to the approximated solution. Note $-J(u)$ is the left rectangular integral approximation of $\int_{0}^{T}\ln{(x(t))} dt$. We have that  $\norm{u^*-\hat{u}}_{L^{\infty}(0,T)}=0.444444$ for each case. This is due to the approximated switch from singular to purely reproduction, $\hat{t_2}$, being off by one node. } 
   \label{tab: PlantUnpen}
  \end{table}

  \begin{table}[ht]
   \centering
  \begin{tabular}{|c|c c c c c c|}
  \hline
   & $\rho$ & $\norm{u_{\rho}-u^*}_{L^1(0,T)}$ & $\norm{u_{\rho}-u^*}_{L^{\infty}(0,2)}$ & $u_{\rho}(0)$ & \textbf{Switch } $t_2$ & \textbf{Runtime (s)}\\
  \hline
  $N=750$    & 0          &0.01755676 & 0.25200195 &0.49799805 & 2.2&46.25\\
  tol$=10^{-10}$ &$10^{-9}$ &0.01755676 & 0.25153629 &0.49846371 & 2.2&45.81\\
  $T=5$       &$10^{-8}$ &0.01755672 & 0.24953915 &0.50046085 & 2.2&36.36\\
  $x_{1,0}=4$ &$10^{-7}$ &0.01755623 & 0.23038155 &0.51961845 & 2.2&32.01\\
  $x_{2,0}=1$ &$10^{-6}$ &\textbf{0.01755173} &0.10869044 &0.64130956 &\textbf{2.2067}&27.54\\
   $h=0.00667$  &$10^{-5}$          &0.01760049 &0.03367945 & 0.7163206 &2.2667&32.82\\
     &$10^{-4}$          &0.01785072 &\textbf{0.01809606} & \textbf{0.7319039} &2.2667&27.60\\
     &$10^{-3}$           &0.01928274 &0.02181486 & 0.7281851 &2.3667& 26.27\\
     & $10^{-2}$ & 0.02628408 &  0.03439004 &0.71560996 & 2.6067 & 23.42\\
     & $10^{-1}$ & 0.06203874 & 0.07234748 & 0.67765252 & 2.3667 & \bf{9.96}\\
  \hline
  \end{tabular}
  \caption{Varying Penalty Parameter for Plant Problem Case 2a)}
  \label{tab: plantvary}
  \end{table}
  \begin{figure}[htbp]
     \begin{subfigure}[t]{0.25\textwidth}
     \includegraphics[width=\linewidth]{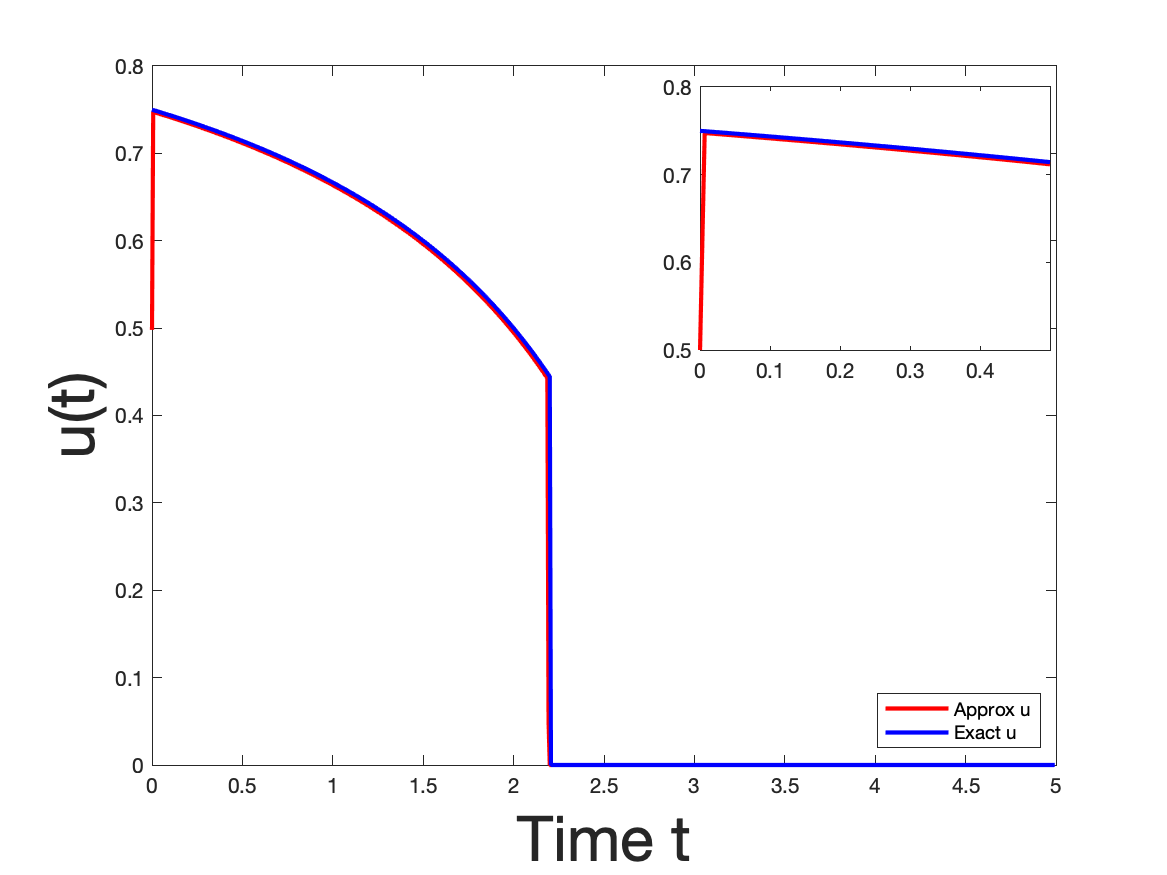}
     \caption{$u_{\rho}$ vs $u^*$ for $\rho=0$}
     \label{fig: varyplantu0}
     \end{subfigure}\hfill
     \begin{subfigure}[t]{0.25\textwidth}
      \includegraphics[width=\linewidth]{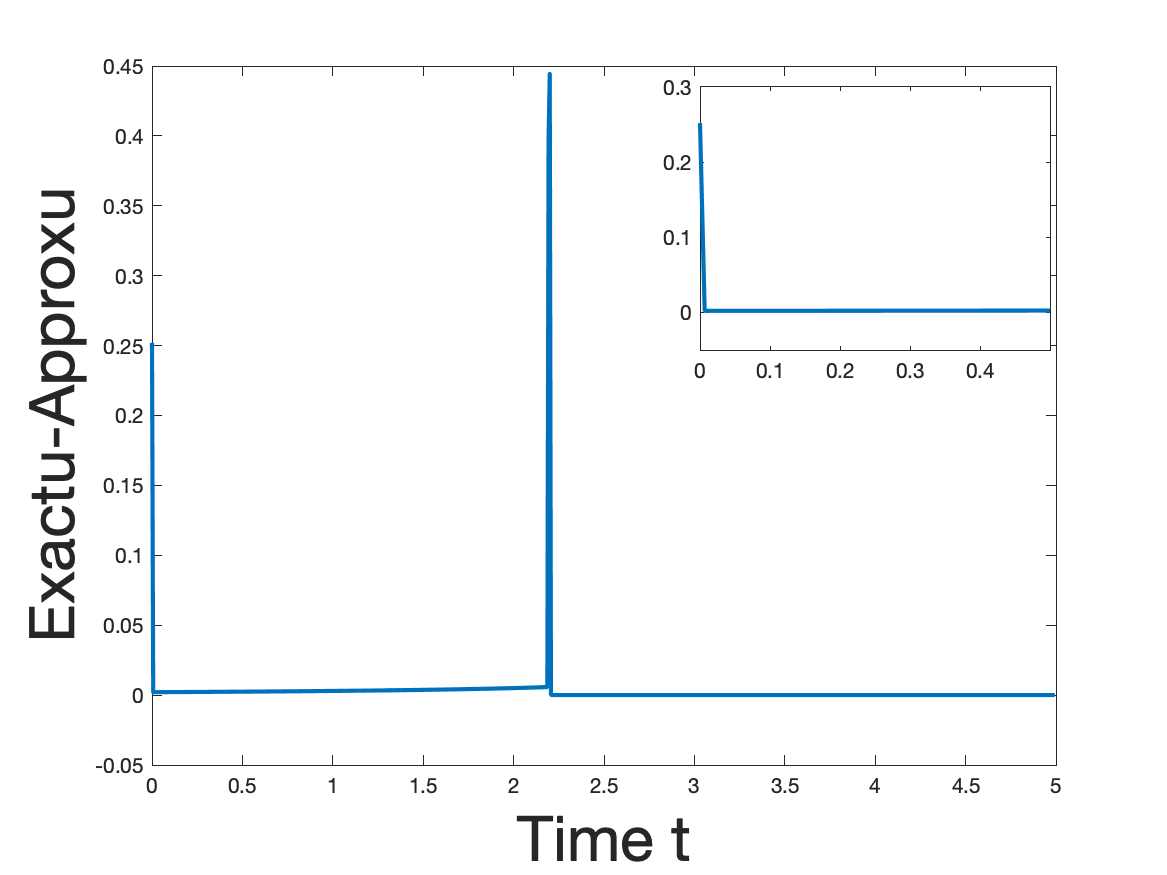}
      \caption{$u^*-u_{\rho}$ for $\rho=0$}
      \label{fig: varyplantdiff0}
      \end{subfigure}\hfill
      \begin{subfigure}[t]{0.25\textwidth}
     \includegraphics[width=\linewidth]{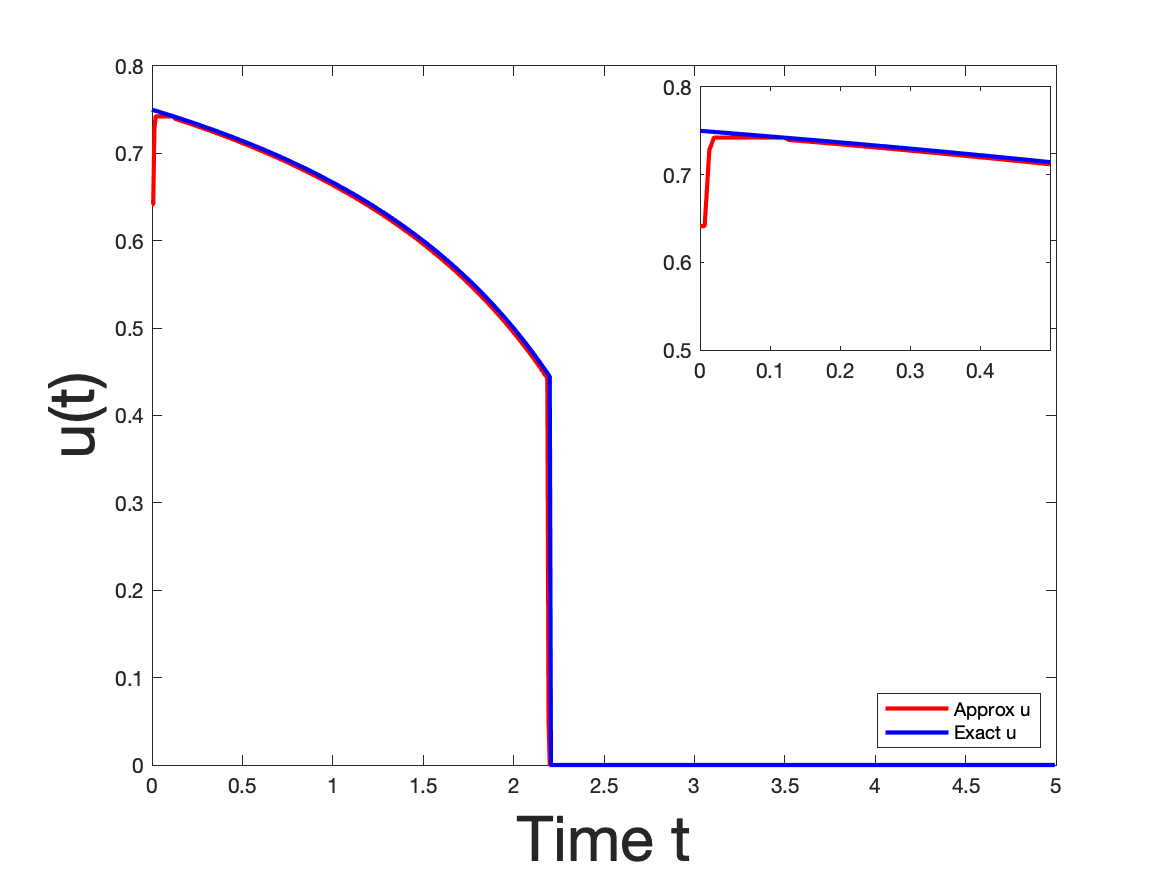}
     \caption{$u_{\rho}$ vs $u^*$ for $\rho=10^{-6}$}
     \label{fig: varyplantu6}
     \end{subfigure}\hfill
     \begin{subfigure}[t]{0.25\textwidth}
      \includegraphics[width=\linewidth]{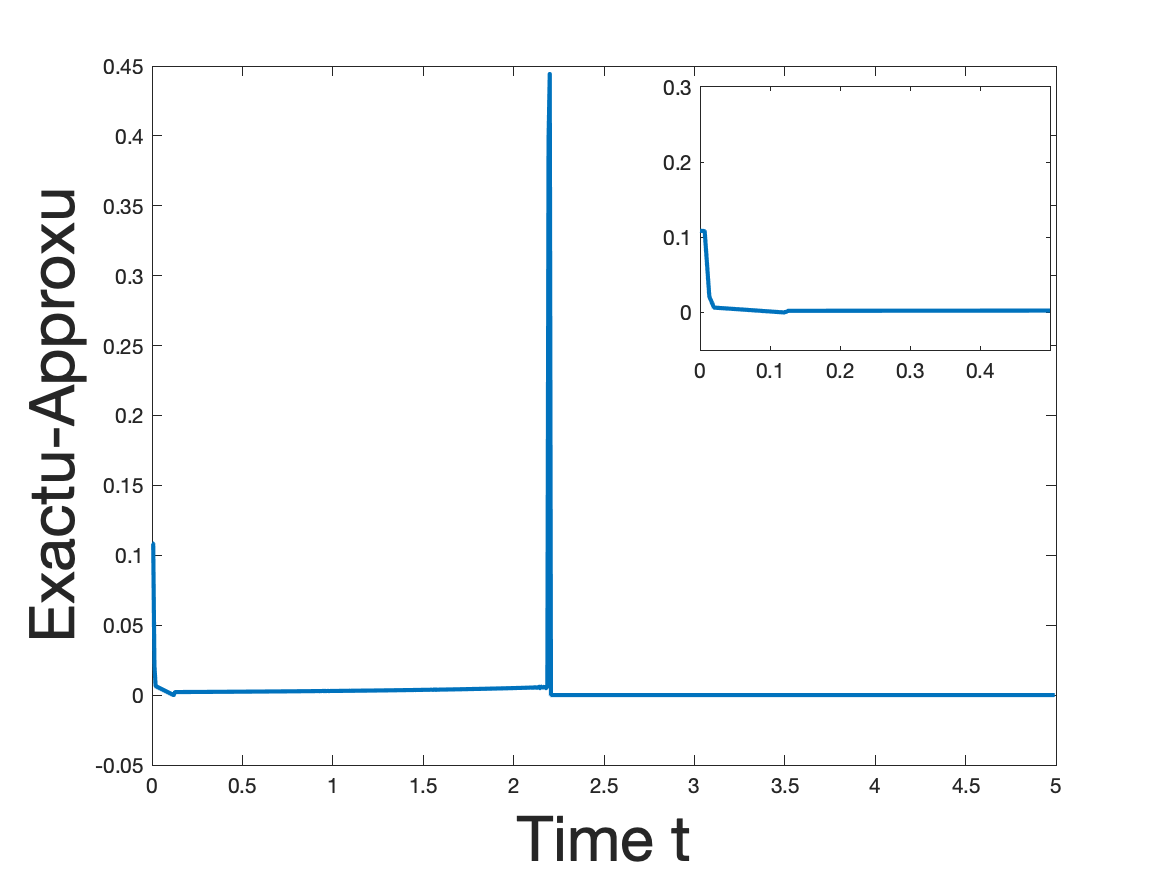}
      \caption{$u^*-u_{\rho}$ for $\rho=10^{-6}$}
      \label{fig: varyplantdiff6}
      \end{subfigure}
      \begin{subfigure}[t]{0.25\textwidth}
     \includegraphics[width=\linewidth]{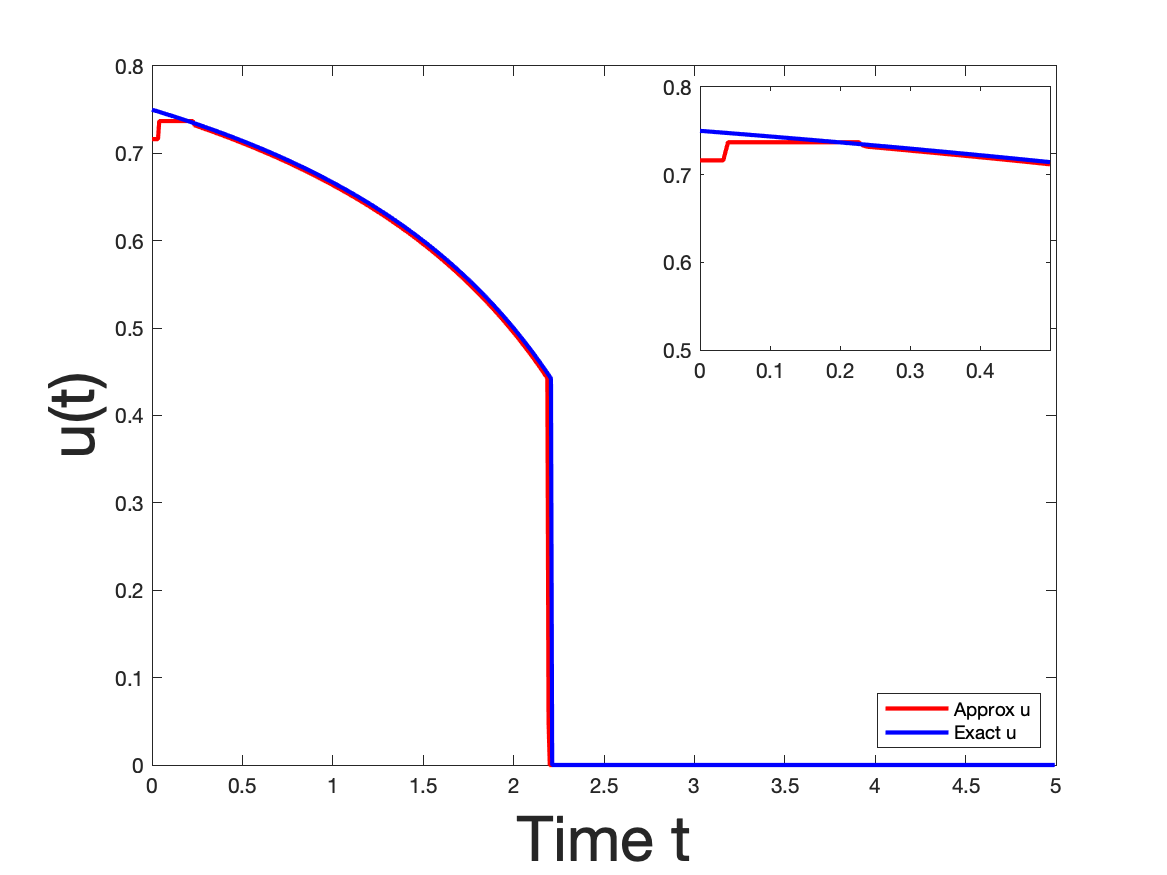}
     \caption{$u_{\rho}$ vs $u^*$ for $\rho=10^{-5}$}
     \label{fig: varyplantu5}
     \end{subfigure}\hfill
     \begin{subfigure}[t]{0.25\textwidth}
      \includegraphics[width=\linewidth]{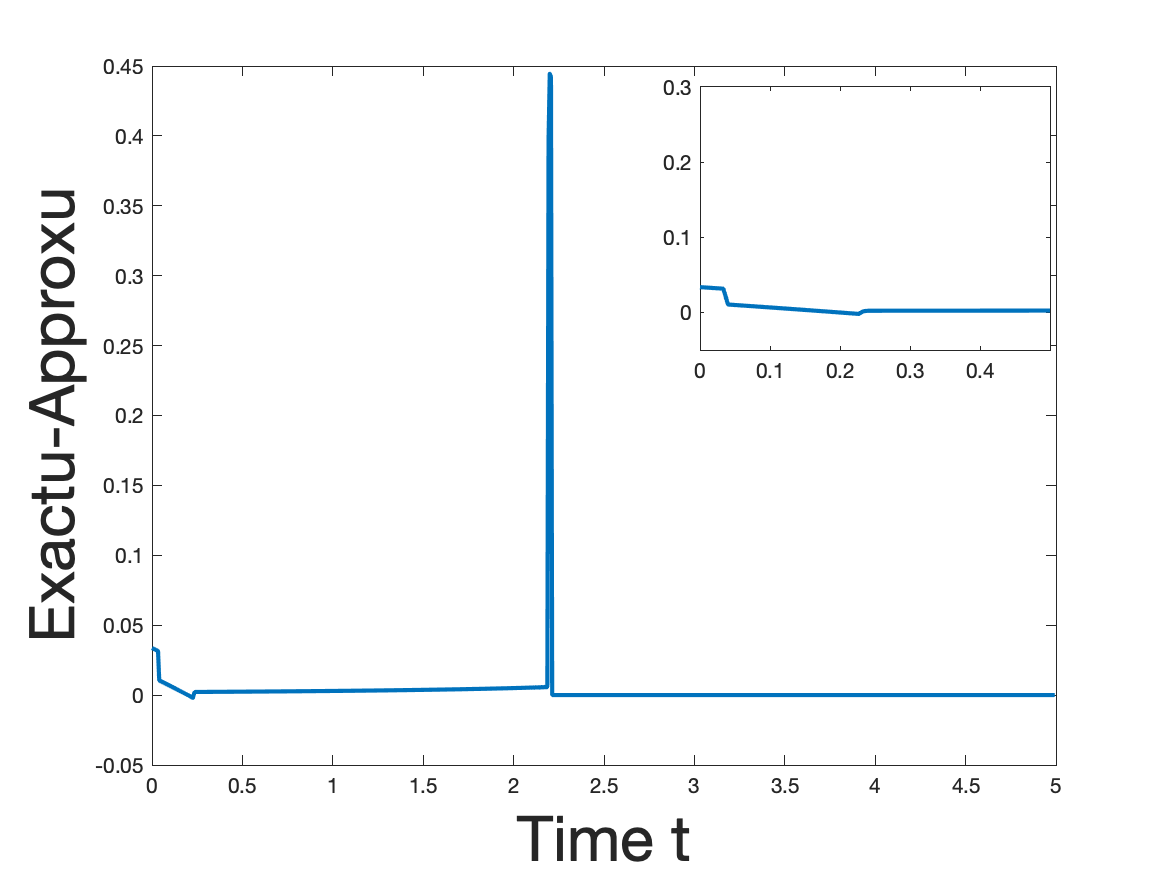}
      \caption{$u^*-u_{\rho}$ for $\rho=10^{-5}$}
      \label{fig: varyplantdiff5}
      \end{subfigure}\hfill
      %
      %
      \begin{subfigure}[t]{0.25\textwidth}
     \includegraphics[width=\linewidth]{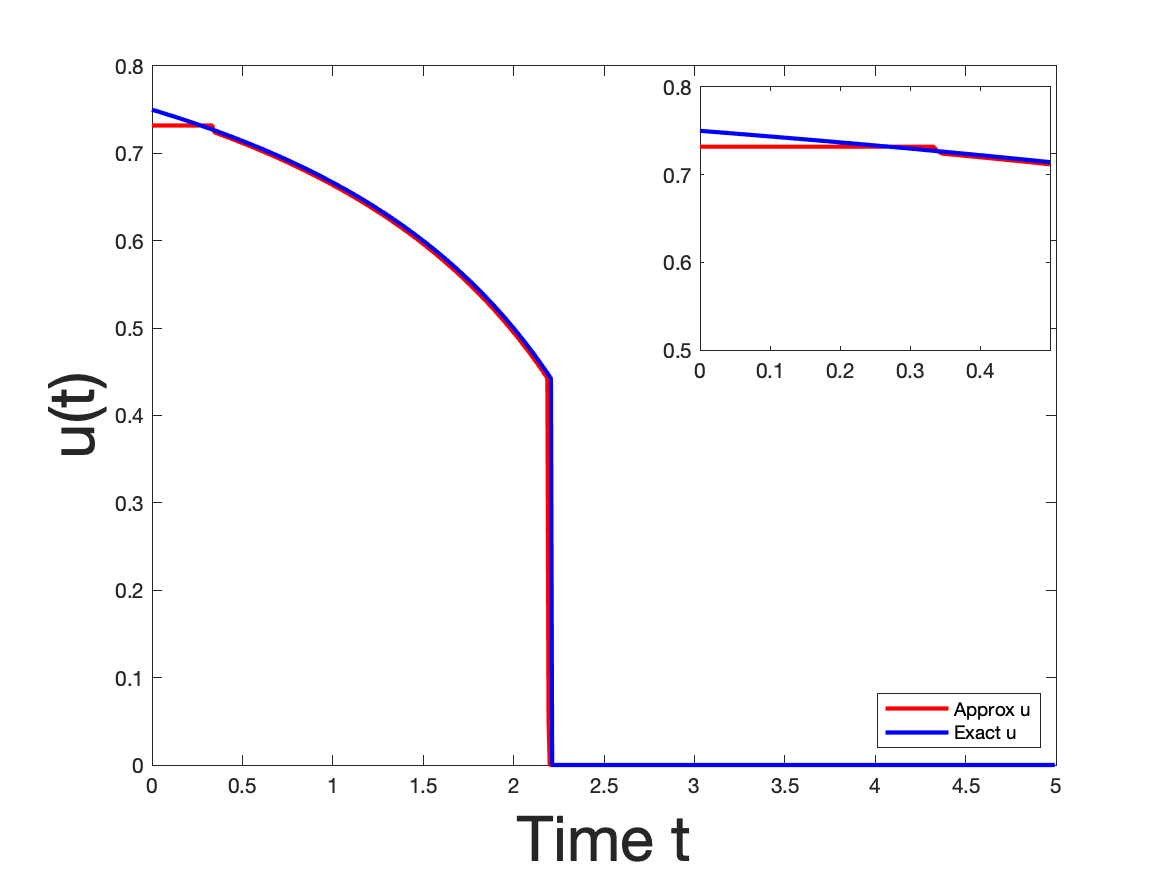}
     \caption{$u_{\rho}$ vs $u^*$ for $\rho=10^{-4}$}
     \label{fig: varyplantu4}
     \end{subfigure}\hfill
     \begin{subfigure}[t]{0.25\textwidth}
      \includegraphics[width=\linewidth]{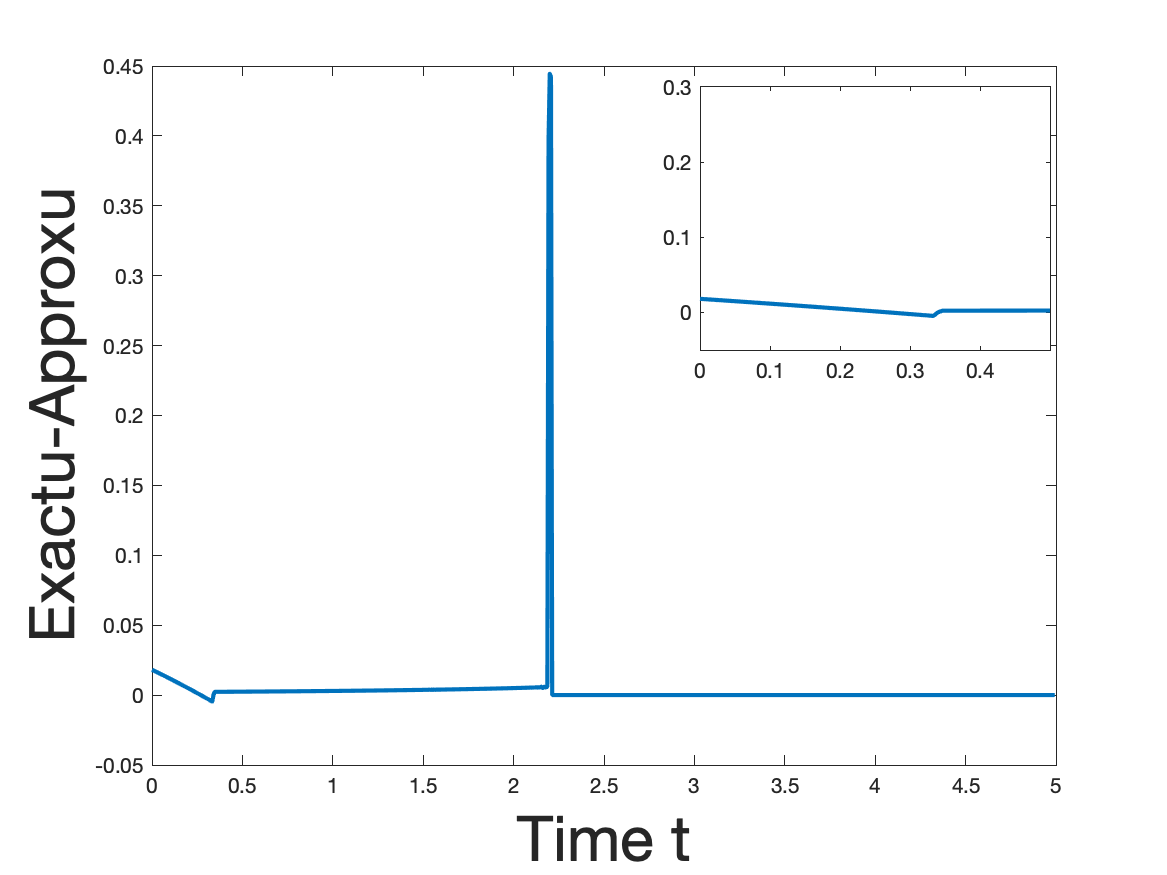}
      \caption{$u^*-u_{\rho}$ for $\rho=10^{-4}$}
      \label{fig: varyplantdiff4}
      \end{subfigure}
      \begin{subfigure}[t]{0.25\textwidth}
     \includegraphics[width=\linewidth]{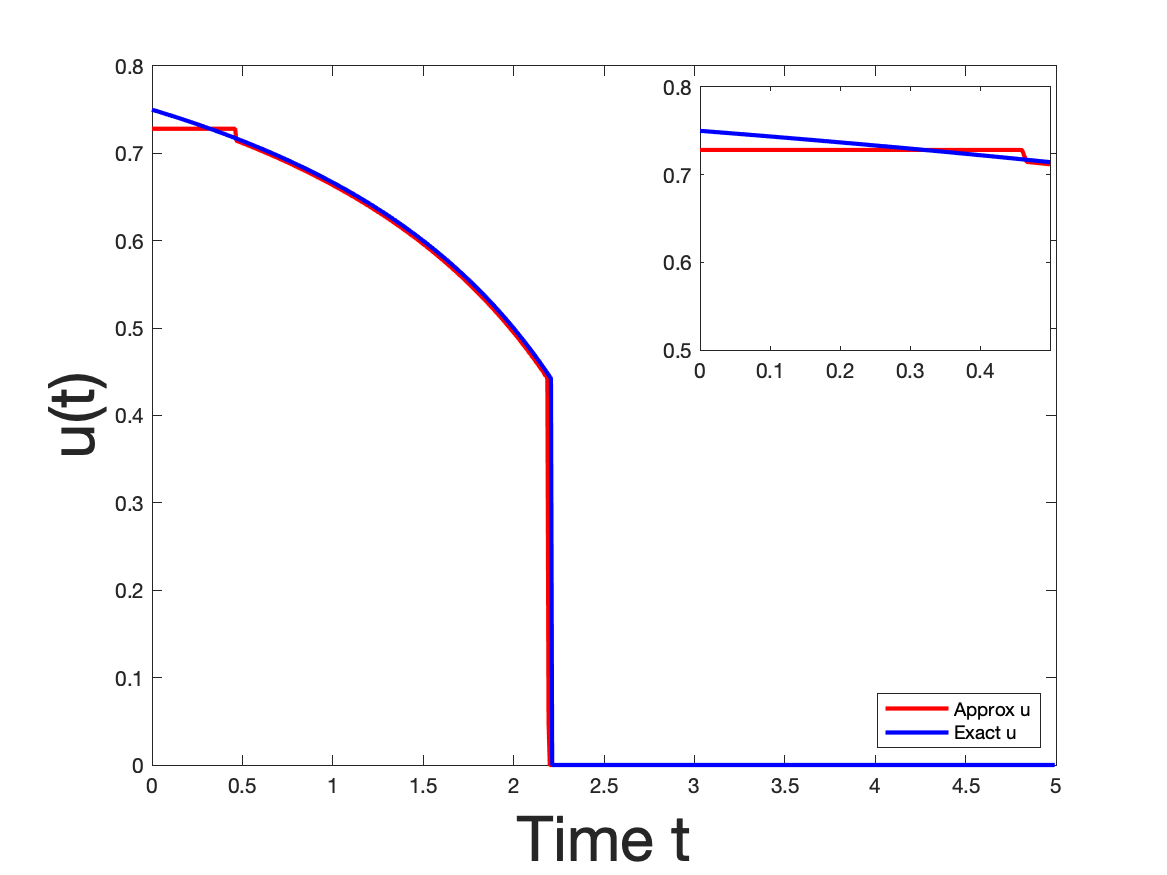}
     \caption{$u_{\rho}$ vs $u^*$ for $\rho=10^{-3}$}
     \label{fig: varyplantu3}
     \end{subfigure}\hfill
     \begin{subfigure}[t]{0.25\textwidth}
      \includegraphics[width=\linewidth]{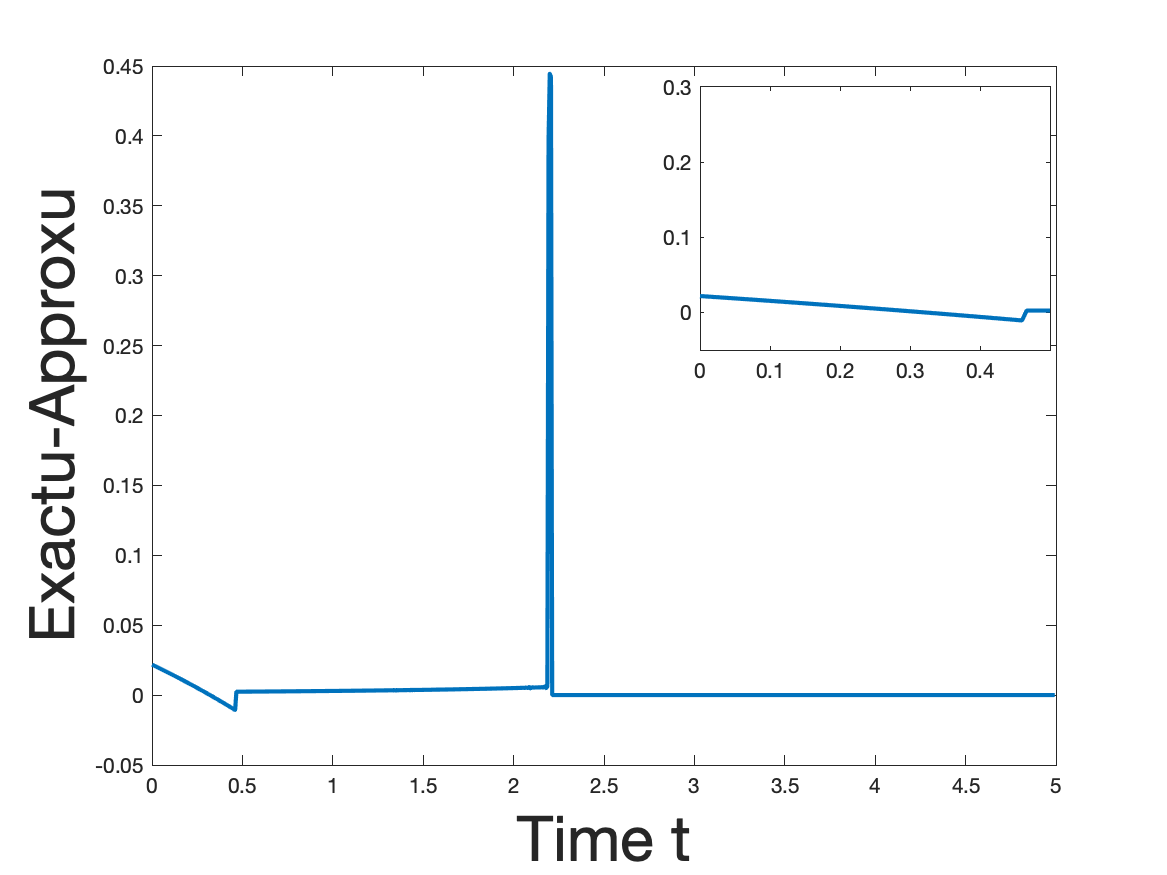}
      \caption{$u^*-u_{\rho}$ for $\rho=10^{-3}$}
      \label{fig: varyplantdiff3}
      \end{subfigure}\hfill
      \begin{subfigure}[t]{0.25\textwidth}
     \includegraphics[width=\linewidth]{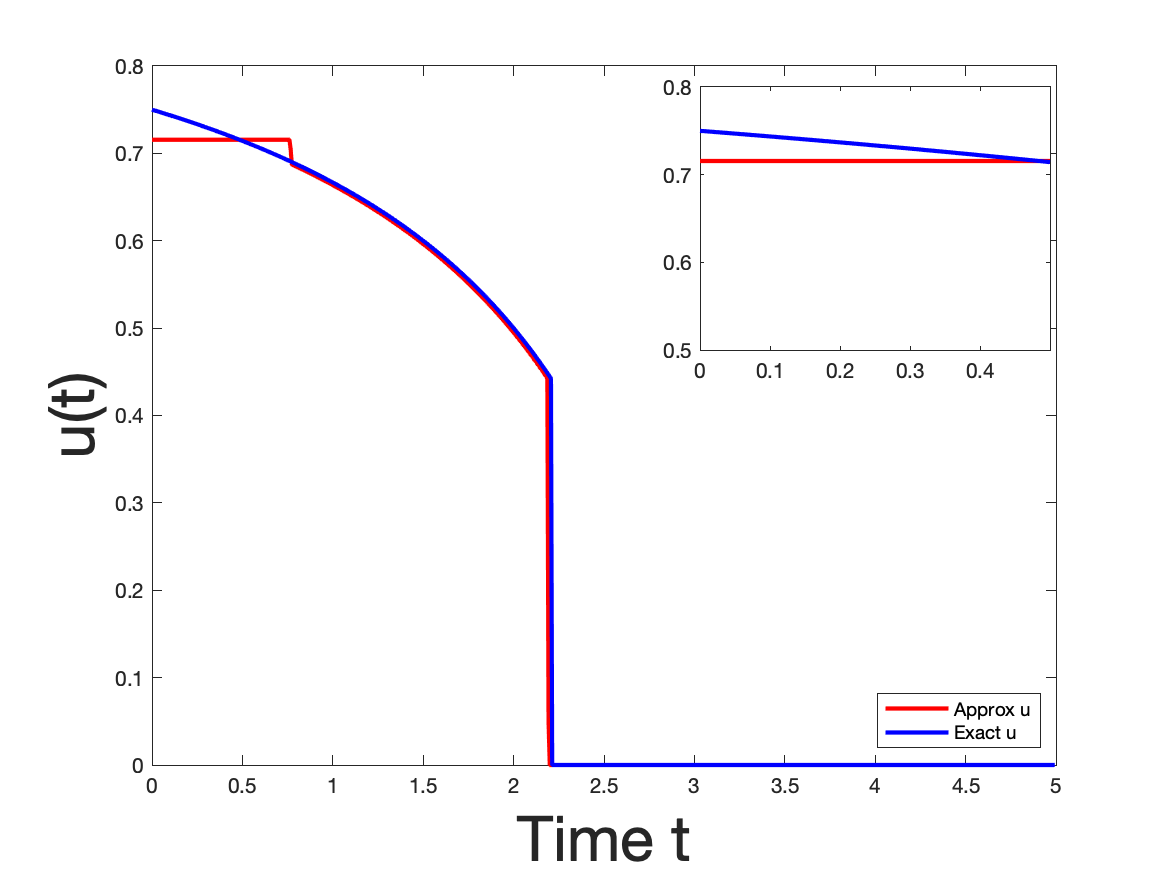}
     \caption{$u_{\rho}$ vs $u^*$ for $\rho=10^{-2}$}
     \label{fig: varyplantu2}
     \end{subfigure}\hfill
     \begin{subfigure}[t]{0.25\textwidth}
      \includegraphics[width=\linewidth]{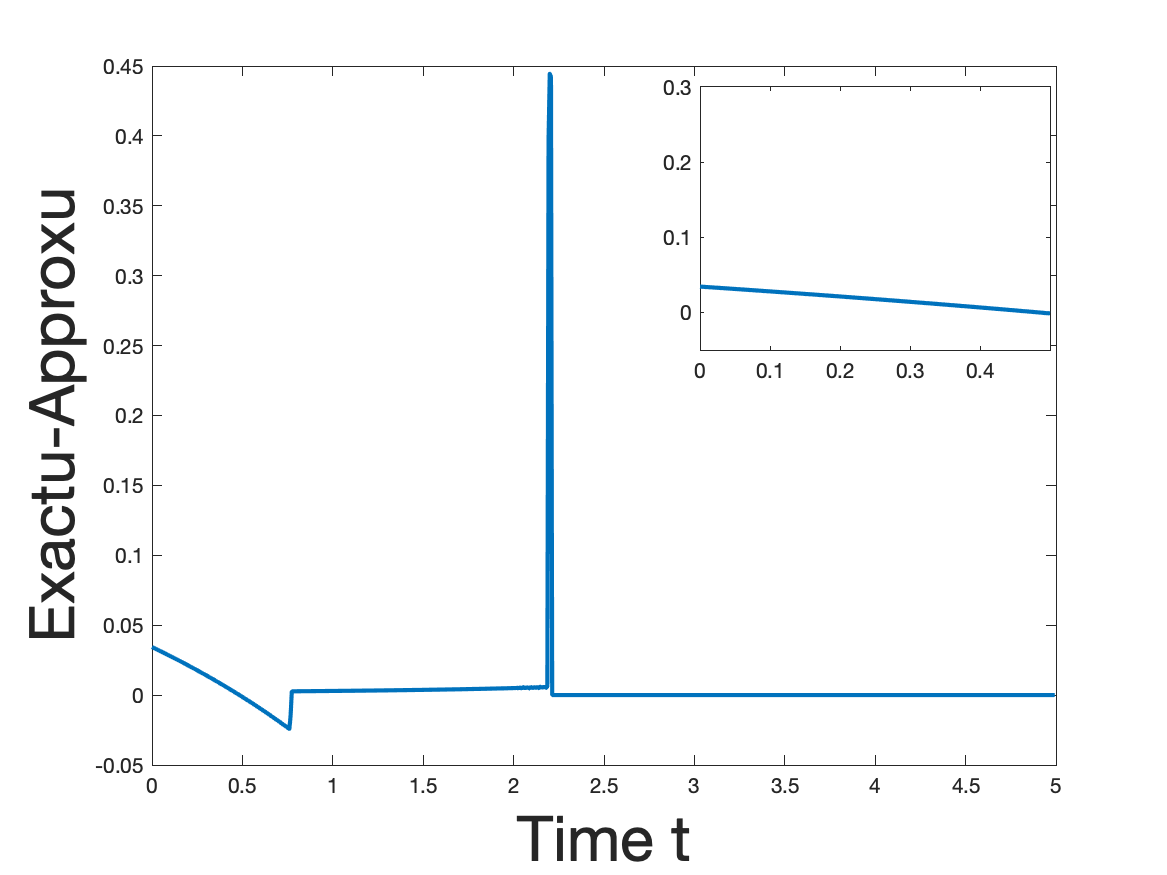}
      \caption{$u^*-u_{\rho}$ for $\rho=10^{-2}$}
      \label{fig: varyplantdiff2}
      \end{subfigure}
      \caption{Varying Penalty to Plant Problem Case 2a): Plots of regularized control $u_{\rho}$ (red) vs optimal control $u^*$ (blue) and $u^*-u$ (cyan) for tuning parameter $\rho\in\{0,10^{-6}, 10^{-5}, 10^{-4}, 10^{-3}, 10^{-2}\}$. The top right corner of each subfigure is a zoomed-in plot of the same subfigure over the time interval $[0,0.5]$. }
      \label{fig: varyingplant}
  \end{figure}
  \begin{figure}[htbp]
  \begin{subfigure}[t]{0.5\textwidth}
   \includegraphics[width=\linewidth]{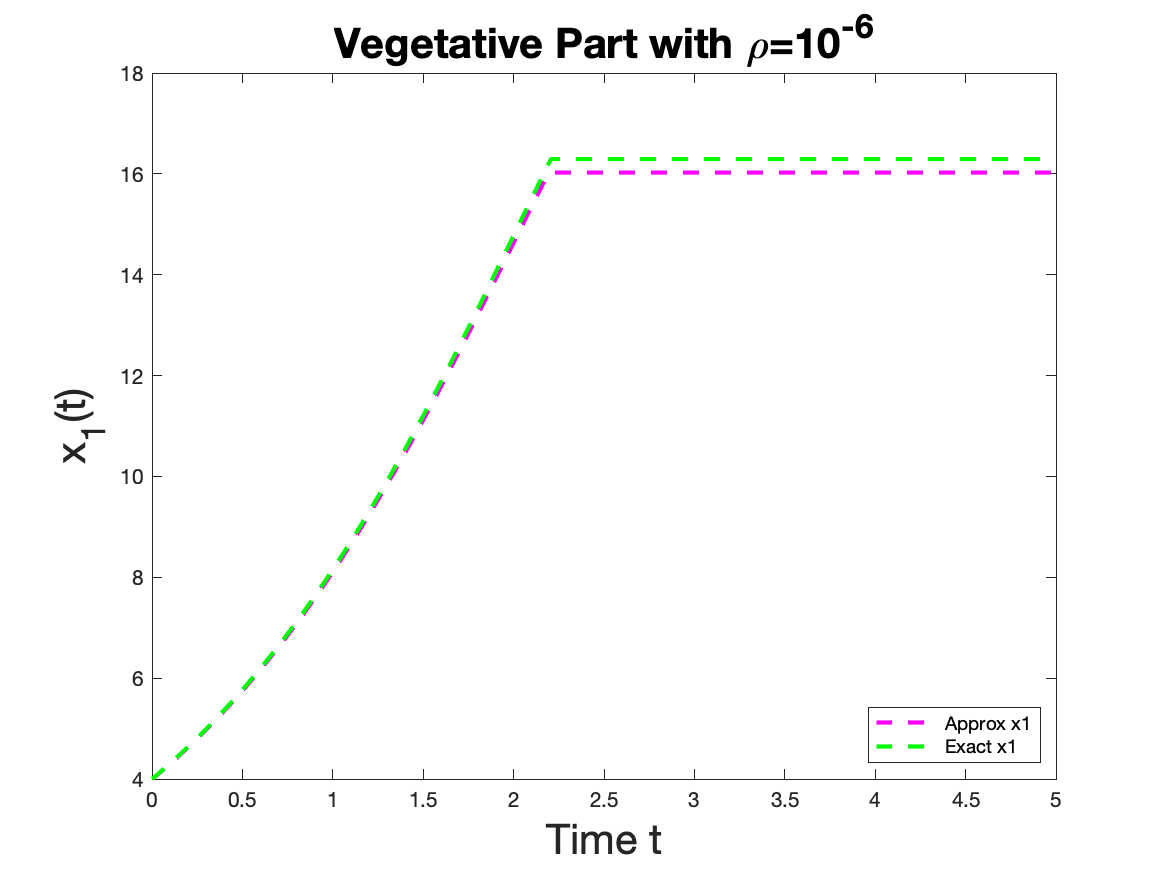}
   \caption{Exact State  $x_1$ (green) vs $x_{1,\rho}$ with $\rho=10^{-6}$ (pink)}
  \end{subfigure}\hfill
  \begin{subfigure}[t]{0.5\textwidth}
   \includegraphics[width=\linewidth]{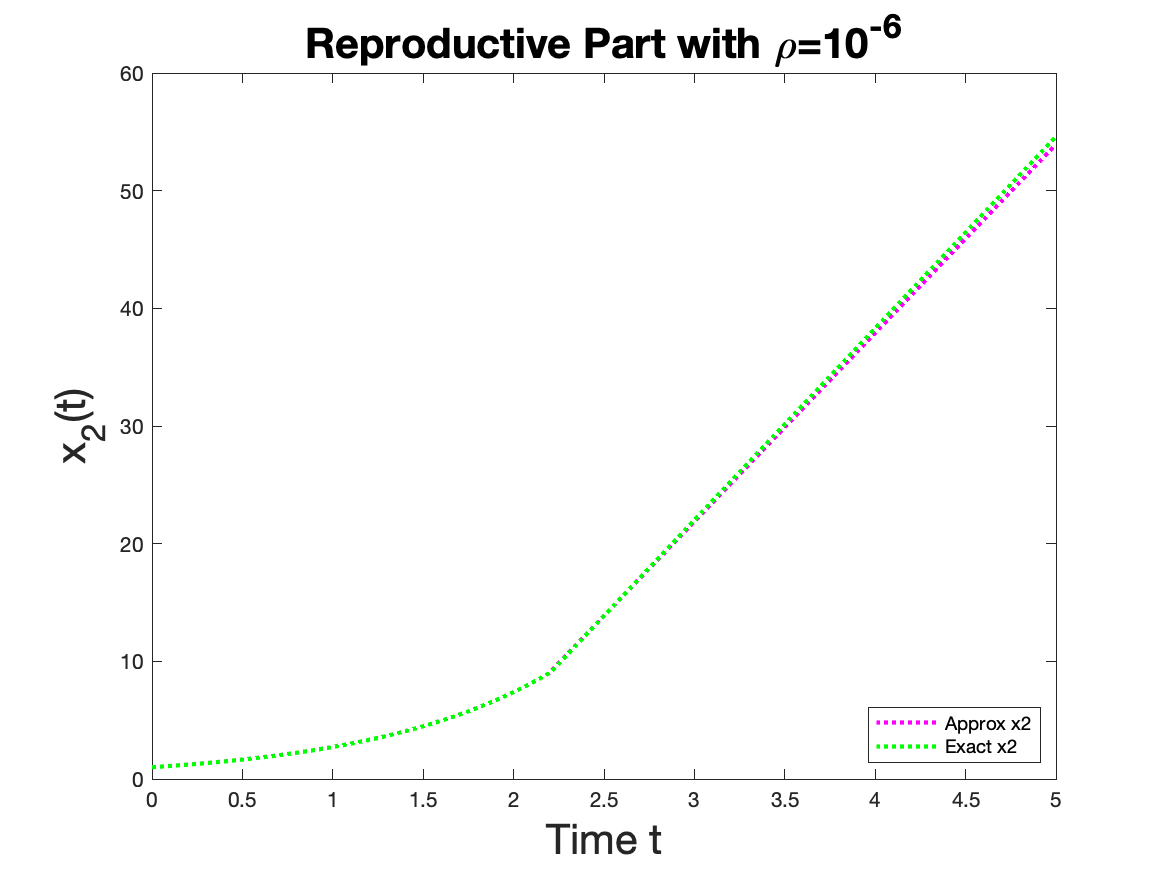}
   \caption{Exact State $x_2$ (green)  vs $x_{2,\rho}$ with $\rho=10^{-6}$ (pink)}
  \end{subfigure}\hfill
  \caption{ Corresponding Trajectories of $x_1$ and $x_2$ for Penalized Plant Problem with $\rho=10^{-6}$}
   \label{fig: plantstates6}
  \end{figure}

   
   We use PASA to numerically solve for problem (\ref{eqn: plantproblem}) with stopping tolerance set to being $10^{-10}$ and with our initial guess for the control being  $u(t)=0 $ over the entire time interval $[0,T]$. We partition the time interval to where there are $N=750$ mesh intervals with mesh size being $h=0.00667$, and the discretization process is as described in Section \ref{subsec: plantdiscrete} with $\rho=0$.
   We are interested in solving the problem when the parameters are set so that a singular case occurs. 
   King and Roughgarden \cite{King1982} mention that from a biological standpoint, the initial weight of the reproductive part of a plant, $x_2(0)$ will always be zero since germination involves vegetative growth only. 
   Consequently, the only realistic situation where the exact solution to problem (\ref{eqn: plantproblem}) contains a singular subarc is if the exact solution is of case \ref{itm: case2} with $x_{2,0}=0$. 
   However, we still would like to observe how PASA performs for all possible cases when the exact solution to problem (\ref{eqn: plantproblem}) contains a singular subarc. 
   Table  \ref{tab: plantparameters} describes the parameter settings that we used for solving problem (\ref{eqn: plantproblem}). 
   When parameters are set to the values shown in column three of Table \ref{tab: plantparameters}, the exact solution to problem  (\ref{eqn: plantproblem}) will be of Case \ref{itm: case1} where control $u^*$ begins singular and switches to the purely reproductive case at $t_2\approx2.2067$.
   When parameters are set to the values shown in column four of Table \ref{tab: plantparameters}, the exact solution to problem  (\ref{eqn: plantproblem}) will be of Case \ref{itm: case2} where control $u^*$ begins purely reproductive, switches to the singular case solution at $t_1\approx 0.2678$, and switches back to being purely reproductive at $t_2\approx 2.2067$.
   We want to emphasize that in Case 2b) parameter settings,  $x_{2,0} =10^{-4}$ because we wanted this initial value to be a close approximation of zero without having any issues in computing the integral approximation for the objective functional $J(u)=\int_0^T(-\ln{x_2(t)} )dt$.
   When parameters are set to the values shown in the last column of Table \ref{tab: plantparameters}, the exact solution to problem  (\ref{eqn: plantproblem}) will be of Case \ref{itm: case3} where control $u^*$ begins purely vegetative, switches to the singular case solution at $t_1\approx  1.5778$, and switches back to being purely reproductive at $t_2\approx 2.2067$. 
   
   For each case, we first observe the unregularized solutions that PASA obtained when solving for problem (\ref{eqn: plantproblem}) to see if it is even necessary to penalize this problem. 
   {When numerically solving for the unpenalized problem, we obtained solutions that were not oscillatory for all three cases. Figures and descriptions corresponding to the unpenalized solutions are given in Appendix Section \ref{subsec:  plantresults} and in Table \ref{tab: PlantUnpen}.}
   We find PASA's unpenalized solution to be a sufficient approximation to the true solution to problem (\ref{eqn: plantproblem}) in Cases \ref{itm: case2} and \ref{itm: case3}.
   However, in Subfigure \ref{fig: varyplantu0} the unpenalized approximation for case \ref{itm: case1}, $\hat{u}$, appears to have an unusual dip at the $t=0$. 
   Based upon explicit solution (\ref{eqn: plantsingsumabegin}) and the parameters setting given in Table \ref{tab: plantparameters}, $u^*(0)=0.75$, but for the approximated unpenalized 
   solution we have $\hat{u}(0)\approx 0.4980$.
   We use PASA to solve for Problem (\ref{eqn: plantproblem}) for Case \ref{itm: case1} with $N=1000$ and $tol=10^{-12}$ to see if such changes would provide any improvements for the initial value.
   The unpenalized solution that PASA obtained when $N=1000$ and $tol=10^{-12}$ still possessed a dip at $t=0$ where $u(0)\approx 0.4985$.
   The reason behind this initial jump is because Case \ref{itm: case1} is a degenerate case. 
   If we look further at the conditions between becoming a solution of Case \ref{itm: case1}, Case \ref{itm: case2}, and Case \ref{itm: case3}, you would notice that they all share a condition that pertains to how the ratio of initial values relates with the fraction $\frac{1}{T-1}$. 
   When having parameters set to where $\frac{x_{2,0}}{x_{1,0}}=\frac{1}{T-1}$, the discretization of problem (\ref{eqn: plantproblem}) can cause any numerical solver to converge to a solution that does not begin singular. 
   
   Although there is no chattering in Subfigure \ref{fig: varyplantu0}, we would like to add a bounded variation penalty term to the objective functional of problem (\ref{eqn: plantproblem}) to see if the approximated penalized solution, $u_{\rho}$, will not dip as significantly as the unpenalized solution did at $t=0$. 
   As before we partition time interval $[0,T]$ to where there are $N=750$ mesh intervals and use the discretization method described in Subsection \ref{subsec: plantdiscrete}. 
   We use PASA to solve for penalized problem (\ref{eqn: penplant}) for varying values of penalty parameters  $\rho\in\{0, 10^{-9}, 10^{-8}, 10^{-7}, 10^{-6}, 10^{-5}, 10^{-4}, 10^{-3}, 10^{-2}, 10^{-1}\}$, with initial guess for our control being $u(t)=0$ over the entire time interval. 
   We have parameters set to being the third column in Table \ref{tab: plantparameters}, with stopping tolerance set to being $10^{-10}$.  
   In Table \ref{tab: plantvary},  we record the $L^1$ and $L^{\infty}$ norm errors between the exact solution $u^*$ and the penalized solution that PASA obtained, $u_{\rho}$, the approximated value of $u_{\rho}(0)$, the approximated switching point $t_2$, and the runtime that was needed to solve for penalized problem (\ref{eqn: penplant}).
   In column four of Table \ref{tab: plantvary} we are finding $\norm{u^*-u_{\rho}}_{L^{\infty}(0,2)}$ because otherwise the whole column would be $0.444444$, which is what we observed in Table \ref{tab: PlantUnpen}.
   Observe from Table \ref{tab: plantvary}, the penalty parameter starts to influence the behavior of $u_{\rho}$ when $\rho\geq10^{-6}$. 
   However, it is worth noting that for penalty parameter values $\rho=10^{-9}, 10^{-8},$ $10^{-7}$, PASA obtained a solution comparable to the unpenalized solution at a faster rate. 
   Based upon Table \ref{tab: plantvary}, $u_{\rho}$ was closest to $u^*$ with respect to the $L^1$ norm when $\rho=10^{-6}$, and $u_{\rho}$ switched to being non-singular at the node that was closest to the true switching point.
   The penalized solution $u_{\rho}$ that was obtained when $\rho=10^{-4}$  had an initial value that was the closest to the true solution's initial value $u^*(0)=0.75$; however this improvement seems to be at the expense of increasing the $L^1$ norm error between $u_{\rho}$ and $u^*$. 
   In Figure \ref{fig: varyingplant}, we provide a chart of figures that pertain to $u_{\rho}$ and $u^*-u_{\rho}$ for $\rho=0,10^{-6},10^{-5}, 10^{-4},10^{-3},$ and $10^{-2}$. 
   Notice that in Figure \ref{fig: varyingplant}, the penalized solutions do not begin singular immediately; however, in comparison to the unpenalized solution, they give better approximations to $u^*(0)$.
    In Subfigure \ref{fig: varyplantu4}, the penalized solution associated with penalty parameter $\rho=10^{-4}$ begins constant with constant value approximately being $0.7319$ and switches to the singular case solution approximately at $t=0.32$.  
   Notice in Subfigures \ref{fig: varyplantu3} and \ref{fig: varyplantu2} that the these penalized solutions also begin constant at values close to 0.7, but remain constant for a longer time period. 
   Notice also in the sixth column of Table \ref{tab: plantvary} that for values $\rho\geq 10^{-5}$, $u_{\rho}$ starts to overestimate the point at which the solution should switch from singular to non-singular.
   In Subfigure \ref{fig: varyplantu6} we find that the unpenalized solution $u_\rho$ associated with penalty parameter value $\rho=10^{-6}$ did the best job in improving the approximated initial value without deviating too much from the exact solution. 
   In Figure \ref{fig: plantstates6}, we have the trajectories of $x_1$ and $x_2$ that correspond to $u_{\rho}$ when $\rho=10^{-6}$. 
   The corresponding state solutions to the penalized control are comparable to the exact state solutions. 

\section{Example 3: SIR Problem}
  Our next example is from Ledzewicz, Aghaee, and Sch\"{a}ttler's article \cite{Mahya}. 
  This example will demonstrate how PASA can be applied to problems in the absence of a true solution. 
  In \cite{Mahya}, Ledzewicz et al.  were interested in 
  using an SIR model with demography to study the 2013-2016 outbreak of Ebola in West Africa.  
  In an SIR model, the total population $N$ is divided into the following three compartments: the susceptible class ($S$), the infectious class ($I$), and the recovered class ($R$). 
  The dynamics in consideration are based on a modification of a system of equations that was first presented in Brauer and Castillo-Chavez's book \cite{Chavez}:
  \begin{align}
      \dot S &= \gamma N-\nu S -\beta \frac{IS}{N}+\rho R, & S(0) =S_0,\label{eqn: Snocontrol}\\
      \dot I &= \beta\frac{IS}{N}-(\nu+\mu)I-\alpha I, & I(0)=I_0,\label{eqn: Inocontrol}\\
      \dot R &= -\nu R-\rho R+\alpha I, & R(0)=R_0\label{eqn: Rnocontrol}.
  \end{align}
  The Greek letters in the equations above represent constant parameters. 
  In equation (\ref{eqn: Snocontrol}), the coefficient $\gamma$ is the birth rate per unit time, and it is assumed that the population of newborns are immediately classified as being susceptible to the disease.
  Standard incidence is assumed in the model with transmission rate being $\beta$, and it is assumed that infected members recover from Ebola at at rate $\alpha$. 
  Parameter $\nu$ is the natural death rate while parameter $\mu$ is the death rate due to infection. 
  We emphasize here that for this model $\rho$ is not the bounded variation penalty parameter, but instead $\rho$ is the parameter that measures the rate at which a recovered individual becomes again susceptible to Ebola.

  In \cite{Mahya}, 
  Ledzewicz, et al.  constructed an optimal control problem to gain understanding on how treatment and a theoretical vaccination of Ebola should be applied to best contribute to limiting the spread of the disease. 
  For the state equations, they incorporate controls, $u$ and $v$, into Equations (\ref{eqn: Snocontrol})-(\ref{eqn: Rnocontrol}) where $u$ represents the vaccination rate and $v$ represents the treatment rate. 
  The modified dynamics are as follows: 
  \begin{align}
      \dot S &= \gamma N-\nu S -\beta \frac{IS}{N}+\rho R-\kappa Su, & S(0) =S_0,\label{eqn: stateS} \\
      \dot I &= \beta\frac{IS}{N}-(\nu+\mu)I-\alpha I-\eta Iv, & I(0)=I_0,\label{eqn: stateI}\\
      \dot R &= -\nu R-\rho R+\kappa Su+\alpha I+\eta Iv, & R(0)=R_0,\label{eqn: stateR}
  \end{align}
  where $\kappa$ and $\eta$ are denoted as being the efficacy of vaccination and treatment respectively.  
  For constructing of an objective functional, Ledzewicz et al.  \cite{Mahya} had the following goal in mind: \\
  \emph{``Given initial population sizes of all three classes, $S_0$, $I_0$ and $R_0$, find the best strategy in terms of the combined efforts of vaccination and treatment that minimizes the number of infectious persons while at the same time also taking into account the cost of vaccination and treatment."}
  The objective functional in consideration for minimization is as follows: 
  \begin{equation}\label{eqn: SIRobjective}
      J(u,v) =\int\limits_0^T (aI(t)+bu(t)+cv(t))dt 
  \end{equation}
  The objective functional $J$ is intended to represent the weighted average of the number of infectious persons and costs of vaccination and treatment. 
  The optimal control problem is as follows 
  \begin{equation}\label{eqn: SIRproblem}
      \begin{array}{rl}
          \min\limits_{(u,v)\in \mathscr{A}}& J(u,v)=\int\limits_0^T (aI(t)+bu(t)+cv(t))dt   \\
          \textrm{subject to} &  \dot S = \gamma N-\nu S -\beta \frac{IS}{N}+\rho R-\kappa Su, \\
          & \dot I=\beta\frac{IS}{N}-(\nu+\mu)I-\alpha I-\eta Iv,\\
          & \dot R= -\nu R-\rho R+\kappa Su+\alpha I+\eta Iv, \\
          & S(0)=S_0, \; I(0)=I_0,\; R(0)=R_0,
      \end{array}
  \end{equation}
  where $\mathscr{A}$ is the set of admissible controls which is assumed as being the set of all Lebesgue measurable functions $u: [0,T]\to [0,u_{\max}]$ and $v:[0,T]\to [0,v_{\max}]$, where $u_{\max}$ is the maximum vaccination rate and $v_{\max}$ is the maximum treatment rate. 
  The assumptions on $\mathscr{A}$ ensure existence of an optimal solution to problem (\ref{eqn: SIRproblem}) which follows from Fleming and Rishel's \cite{fleming}.
  
  For problem (\ref{eqn: SIRproblem}), Ledzewicz et al. \cite{Mahya} used Pontryagin's minimum principle \cite{Pontryagin} to find the first order necessary conditions for optimality.
  By rewriting the state equations into a multi-input control-affine system of vector form and using Lie derivatives, they were able to compute the switching functions for controls $u$ and $v$ as well as multiple derivatives of the switching function, allowing  Legendre-Clebsch Condition to be checked. 
  Their methods of using Lie Derivatives is particularly useful in determining existence of a singular control when observing problems with many states and control variables involved, and we direct the reader to their article \cite{Mahya} and Sch\"{a}ttler and Ledzewicz's book \cite{geometriccontol}  for more details on this procedure. 
  
  In \cite{Mahya}, Ledzewicz et al. use the the Tomlab algorithm PROPT \cite{PROPT}, a collocation solver, to obtain an approximate solution to problem (\ref{eqn: SIRproblem}) with the numerical values of parameters set to where a singular subarc  is present in $u^*$, while $v^*$ contained no singular subarc. 
  In this section, we will discuss how to discretized the penalized Ebola problem given in Problem (\ref{eqn: SIRproblemPen}) where the penalty term involved is based on the total variation of optimal control $u$ \cite{Caponigro2018}.
  Additionally in Appendix Section \ref{subsec: matlab},  we provide the MATLAB file, demoOC.m, to demonstrate how to use PASA to solve for problem (\ref{eqn: SIRproblemPen}). 
  We will use PASA to solve for both the unpenalized and penalized problem with parameters set to being what was used in \cite{Mahya} and remark about PASA's performance on this problem.

  
  
  \subsection{ Discretization of SIR Problem} \label{sec: SIRdiscrete}
   We would like to discretize the following penalized problem 
   \begin{equation}\label{eqn: SIRproblemPen}
       \begin{array}{rl}
          \min\limits_{(u,v)\in \mathscr{A}}& J_p(u,v)=\int\limits_0^T (aI(t)+bu(t)+cv(t))dt +pV(u) \\
          \textrm{subject to} &  \dot S = \gamma N-\nu S -\beta \frac{IS}{N}+\rho R-\kappa Su, \\
          & \dot I=\beta\frac{IS}{N}-(\nu+\mu)I-\alpha I-\eta Iv,\\
          & \dot R= -\nu R-\rho R+\kappa Su+\alpha I+\eta Iv, \\
          & S(0)=S_0, \; I(0)=I_0,\; R(0)=R_0,
      \end{array}
   \end{equation}
  where $0\leq p<1$ is the bounded variation penalty parameter and function $V$ measures the total variation of the control which is defined in equation (\ref{eqn: profitV}).
  Since $v^*$ will not contain any singular subarcs, it is likely that $v$ will not need to be penalized.
  Additionally, we would like to discretize the following adjoint equations: 
  \begin{align}
      \dot\lambda_S &=\lambda_S\left(-\gamma+\nu+\beta \frac{I}{N}-\beta\frac{IS}{N^2} +\kappa u\right)+\lambda_I\left(\beta\frac{IS}{N^2}-\beta\frac{I}{N}\right)-\lambda_{R}\kappa u, \label{eqn: Sadjoint}\\
      \dot \lambda_I &= -a +\lambda_S\left(-\gamma+\beta\frac{S}{N}-\beta\frac{IS}{N^2}\right)+\lambda_I\left(-\beta\frac{S}{N}+\beta\frac{IS}{N^2}+(\nu+\mu+\alpha)+\eta v\right)-\lambda_R(\alpha+\eta v),\label{eqn: Iadjoint}\\
      \dot \lambda_{R} &= -\lambda_{S}\left(\gamma+\beta\frac{IS}{N^2}+\rho\right)+\lambda_{I}\left(\beta\frac{IS}{N^2}\right)+\lambda_{R}(\nu+\rho),\label{eqn: Radjoint}
  \end{align}
  with transversality conditions being 
  \begin{equation}\label{eqn: tranSIR}
      \lambda_S(T)=0,\quad \lambda_I(T)=0,\quad \lambda_{R}(T)=0.
  \end{equation}
  
  First we assume that controls $u$ and $v$ are constant over each mesh interval. 
  We partition time interval $[0,T]$, by using $n+1$ equally spaced nodes, $0=t_0<t_1<\cdots<t_n=T$. 
  For all $k=0,1,\dots, n$ we assume that $S_k=S(t_k)$, $I_k=I(t_k)$, and $R_k=R(t_k)$. 
  For the controls we denote $u_k=u(t)$ and $v_k=v(t)$ for all $t_k\leq t<t_{k+1}$ when $k=0,\dots, n-2$. 
  Additionally we denote $u_{n-1}=u(t)$ and $v_{n-1}=v(t)$ for all $t_{n-1}\leq t\leq t_n$. 
  We then have $\bs{S}\in\R^{n+1},\bs{I}\in\R^{n+1},$ and $\bs{R}\in\R^{n+1}$ while $\u\in\R^{n}$ and $\bs{v}\in\R^n$. 
  Furthermore, $\bs{N}\in\R^{n+1}$ where $N_k=N(t_k)=S_k+I_k+R_k$ for all $k=0,1,\dots, n$. 
  We will use a left-rectangular integral approximation for approximating objective functional $J_p(u,v)$ in problem (\ref{eqn: SIRproblemPen}), and we use forward Euler's method for discretizing the state equations.
  Since PASA uses a gradient scheme for one of its phases, we need the cost function to be differentiable. 
 We suggest a decomposition of each absolute value term arising from the total variation term in $J_{\rho}$ to ensure that $J_{\rho}$ is differentiable. 
  We introduce two $n-1$ vectors $\bs{\zeta}$ and $\bs{\iota}$ whose entries are non-negative. 
  Each entry of $\bs{\zeta}$ and $\bs{\iota}$ will be defined as: 
{  
\begin{equation*}
      |u_{k+1}-u_k| = \zeta_k+\iota_k \; \text{ for all } k=0,\dots, n-2. 
  \end{equation*}
}
  The discretization of problem (\ref{eqn: SIRproblemPen}) then becomes 
   \begin{equation}\label{eqn: pendiscreteSIR}
       \begin{array}{rl}
      \min &J_p(\u,\bs{\zeta}, \bs{\iota},\v) = \sum\limits_{k=0}^{n-1}h(aI_k+bu_k+cv_k) + p\sum\limits_{k=0}^{n-1}(\zeta_k+\iota_k)\\
      &S_{k+1}= S_k+h\left(\gamma N_k-\nu S_k-\beta \frac{I_k S_k}{N_k}+\rho R_k-\kappa S_ku_k\right)\; \text{for all }k=0,\dots, n-1,\\
      &I_{k+1}= I_k+h\left(\beta\frac{I_kS_k}{N_k}-(\nu+\mu+\alpha)I_k -\eta I_k v_k\right)\; \text{ for all } k = 0,\dots, n-1,\\
      &R_{k+1}= R_k+h(-\nu R_k-\rho R_k +\kappa S_ku_k+\alpha I_k + \eta I_kv_k)\\
      &0\leq u_k\leq u_{\max} \text{ for all } k=0,\dots, n-1,\\
      &0\leq v_k\leq v_{\max} \text{ for all } k=0, \dots, n-1,\\
      &u_{k+1}-u_k = \zeta_k-\iota_k \; \text{ for all } k=0,\dots, n-2,\\
      &\zeta_k\geq 0, \iota_k\geq 0 \; \text{ for all } k=0,\dots, n-1,
     \end{array}
  \end{equation}
  where $h=\frac{T}{n}$ is the mesh size and the first components of $S$, $I$, and $R$ are set to being the initial conditions associated with the state equations. 
  Notice that for the above problem, we are minimizing the penalized objective function with respect to vectors $\u,\bs{\zeta},\bs{\iota},$ and $\v$. 
  The equality constraints associated with $\bs{\zeta}$ and $\iota$ are linear constraints that PASA can interpret. 
  The equality constraint associated with $\bs{\zeta}$ and $\bs{\iota}$ and be written as: 
  \begin{equation}\label{ref: linearsir}
      \left[
      \begin{array}{c | c | c|c}
      \bs{A} &-\bs{I}_{n-1}&\bs{I}_{n-1}& \bs{O}_{n-1,n}
      \end{array}
   \right]
   \begin{bmatrix}
   \u\\
   \hline
   \bs{\zeta}\\
   \hline
   \bs{\iota}\\
   \hline
   \v
   \end{bmatrix}
   =
   \bs{0},
  \end{equation}
  where $\bs{I}_{n-1}$ is the identity matrix with dimension $n-1$, $\bs{O}_{n-1,n}$ is an $n-1\times n$ all zeros matrix, $\bs{0}$ is the $n-1$ all zeros vector, and {$\bs{A}$ is an $n-1\times n$ dimensional sparse matrix give in equation (\ref{eqn: sparsematrix}).}
  For finding the gradient of $J_p$ in problem (\ref{eqn: pendiscreteSIR}) we use Theorem \ref{thm: lagrange}. 
  The Lagrangian of problem (\ref{eqn: pendiscreteSIR}) is as follows
  \begin{equation}\label{eqn: SIRlagrange}
  \begin{split}
      \Ls(\bs{S},\bs{I},\bs{R},\u,\bs{\zeta},\bs{\iota},\v ) &=  \sum\limits_{k=0}^{n-1}h(a I_k+bu_k+c v_k) + p\sum_{k=0}^{n-1}(\zeta_k-\iota_k)\\
     & +\sum\limits_{k=0}^{n-1}\lambda_{S,k}\left(-S_{k+1}+S_k+h\left(\gamma N_k-\nu S_k-\beta\frac{I_k S_k}{N_k}+\rho R_k-\kappa S_k u_k\right)\right)\\
     &+ \sum\limits_{k=0}^{n-1}\lambda_{I,k}\left(-I_{k+1}+I_k+h\left(\beta\frac{I_k S_k}{N_k}-(\nu+\mu+\alpha)I_k-\eta I_k v_k\right)\right)\\
     &+ \sum\limits_{k=0}^{n-1}\lambda_{R,k}\left(-R_{k+1}+R_k+h\left(-\nu R_k-\rho R_k+\kappa S_k u_k +\alpha I_k+\eta I_k v_k\right)\right),
  \end{split}
  \end{equation}
  where $\lam_{S}\in\R^{n-1}$,$\lam_{I}\in\R^{n-1}$, $\lam_{R}\in\R^{n-1}$ are the Lagrangian multiplier vectors. 
  As suggested in Theorem \ref{thm: lagrange}, we will use the gradient of the Lagrangian  with respect to vectors $\u$, $\bs{\zeta}$, $\bs{\iota}$, and $\v$ to find $\bs{\nabla_{\u,\bs{\zeta},\bs{\iota},\v}J_p}\in\R^{2(n-1)+2(n-2)}$
  This is necessary because based upon the state equations (\ref{eqn: stateS})-(\ref{eqn: stateR}) of problem (\ref{eqn: SIRproblemPen}), $S$, $I$, and $R$ can be viewed as functions of $u$ and $v$. 
  So when computing $\bs{\nabla_{\u,\bs{\zeta},\bs{\iota},\v}J_p}\in\R^{2(n-1)+2(n-2)}$, we should consider that $\bs{S}$, $\bs{I}$, and $\bs{R}$ depend on $\u$ and $\v$.
  We find the partial derivative of $\Ls$ with respect to $u_k$, $\zeta_k$, $\iota_k$, and $v_k$, and obtain the following: 
  \begin{align}
     \frac{\partial\Ls}{\partial u_k} &= hb+h\kappa S_k(\lambda_{R,k}-\lambda_{S,k})\; \text{ for all } k=0,\dots, n-1\\
     \frac{\partial\Ls}{\partial \zeta_k} &= p \;\text{ for all } k=0,\dots, n-2\\
     \frac{\partial\Ls}{\partial \iota_k} &= p \;\text{ for all } k=0,\dots, n-2\\
     \frac{\partial\Ls}{\partial v_k} &=hc+h\eta I_k(\lambda_{R,k}-\lambda_{I,k})\; \text{ for all } k=0,\dots, n-1.
  \end{align}
  
  By Theorem \ref{thm: lagrange}, we compute $\bs{\nabla_{\u,\bs{\zeta},\bs{\iota},\v}J_p}\in\R^{2(n-1)+2(n-2)}$ as follows 
  \begin{equation}
     \bs{ \nabla_{\u,\bs{\zeta},\bs{\iota},\v}J_p }=  \bs{\nabla_{\u,\bs{\zeta},\bs{\iota},\v}\Ls} =
      \begin{bmatrix}
      hb+h\kappa S_0(\lambda_{R,0}-\lambda_{S,0})\\
      \vdots\\
      hb+h\kappa S_k(\lambda_{R,k}-\lambda_{S,k})\\
      \vdots\\
      hb+h\kappa S_{n-1}(\lambda_{R,n-1}-\lambda_{S,n-1})\\
      \hline
      p\\
      \vdots\\
      p\\
      \hline
      p\\
      \vdots\\
      p\\
      \hline
      hc+h\eta I_k(\lambda_{R,0}-\lambda_{I,0})\\
      \vdots\\
      hc+h\eta I_k(\lambda_{R,k}-\lambda_{I,k})\\
      \vdots\\
      hc+h\eta I_{n-1}(\lambda_{R,n-1}-\lambda_{I,n-1})
     \end{bmatrix},
  \end{equation}
  provided that the Lagrange multiplier vectors satisfy condition (\ref{eqn: thmadjoint}) in Theorem \ref{thm: lagrange}. 
  To satisfy condition (\ref{eqn: thmadjoint}) we take the partial derivative of the Lagrangian $\Ls$ with respect to $S_k$, $I_k$, and  $R_k$ for all $k=1,\dots, n$, and note that we are not finding the partial derivative of $\Ls$ with respect to $S_0$, $I_0$, and $R_0$ since these entries are the known values. 
  Taking the partial derivative of the Lagrangian, given in equation (\ref{eqn: SIRlagrange}), with respect to $S_k$, $I_k$, and $R_k$ yields the following expressions: 
  \begin{align}
      \frac{\partial \Ls}{\partial S_k} &= -\lambda_{S,k-1}+\lambda_{S,k}+h\lambda_{S,k}\left(\gamma-\nu-\beta\frac{I_k}{N_k}+\beta\frac{I_kS_k}{N_k^2}-\kappa u_k\right)\nonumber \\
      &\quad\;\;\;+h\lambda_{I,k}\left(\beta\frac{I_k}{N_k}-\beta\frac{I_kS_k}{N_k^2}\right)+h\lambda_{R,k}(\kappa u_k),\;\quad \text{ for all } k=1,\dots, n-1,\label{eqn: adjointSkdiscrete}\\
      \frac{\partial \Ls}{\partial S_n} &=-\lambda_{S,n-1},\label{eqn: adjointSndiscrete}\\
      \frac{\partial \Ls}{\partial I_k} &= ha-\lambda_{I,k-1}+\lambda_{I,k}+h\lambda_{I,k}\left(\beta\frac{S_k}{N_k}-\beta\frac{I_kS_k}{N_k^2}-(\nu+\mu+\alpha)-\eta v_k\right)\nonumber\\
      &\quad\;\;\; + h\lambda_{S,k}\left(\gamma-\beta\frac{S_k}{N_k}+\beta\frac{I_kS_k}{N_k^2}\right)
      +h\lambda_{R,k}(\alpha+\eta v_k),\;\quad \text{ for all } k=1,\dots, n-1,\label{eqn: adjointIkdiscrete} \\
      \frac{\partial \Ls}{\partial I_n} &= -\lambda_{I,n-1}, \label{eqn: adjointIndiscrete} \\
      \frac{\partial \Ls}{\partial R_k} &= -\lambda_{R,k-1}+\lambda_{R,k}-h\lambda_{R,k}(\nu+\rho) \nonumber\\
      &\;\quad\;\;+h\lambda_{S,k}\left(\gamma+\rho+\beta\frac{I_kS_k}{N_k^2}\right)
      -h\lambda_{I,k}\left(\beta \frac{I_kS_k}{N_k^2}\right),\; \quad\text{ for all } k=1,\dots, n-1\label{eqn: adjointRkdiscrete}\\
      \frac{\partial \Ls}{\partial R_n} &= -\lambda_{R,n-1}\label{eqn: adjointRndiscrete}. 
  \end{align}
  For finding multiplier vectors $\lam_S$, $\lam_{I}$, $\lam_{R}$ that satisfy Theorem condition (\ref{eqn: thmadjoint}),
  we set the above equations (\ref{eqn: adjointSkdiscrete})-(\ref{eqn: adjointRndiscrete}) equal to zero and do as follows:  solve for $\lambda_{S,k-1}$ in equation (\ref{eqn: adjointSkdiscrete}), solve for $\lambda_{S,n-1}$ in equation (\ref{eqn: adjointSndiscrete}),  solve for $\lambda_{I,k-1}$ in equation (\ref{eqn: adjointIkdiscrete}), solve for $\lambda_{I,n-1}$ in equation (\ref{eqn: adjointIndiscrete}),  solve for $\lambda_{R,k-1}$ in equations (\ref{eqn: adjointRkdiscrete}), and solve for $\lambda_{R,n-1}$ in equation (\ref{eqn: adjointRndiscrete}). 
  After performing the follow steps and rearranging terms we generate a discretization of equations (\ref{eqn: Sadjoint})-(\ref{eqn: Radjoint}) that satisfy the transversality conditions (\ref{eqn: tranSIR}): 
  \begin{align}
      \lambda_{S,k-1}&=\lambda_{S,k}+ h\lambda_{S,k}\left(\gamma-\nu-\beta\frac{I_k}{N_k}+\beta\frac{I_kS_k}{N_k^2}-\kappa u_k\right)\nonumber\\
      &\quad \;\; +h\lambda_{I,k}\left(\beta\frac{I_k}{N_k}-\beta\frac{I_kS_k}{N_k^2}\right)
      +h\lambda_{R,k}(\kappa u_k),\;\quad \text{ for all } k=1,\dots, n-1,\label{eqn; adjointSkdisc}\\
      \lambda_{S,n-1} &= 0,\label{eqn: adjointSndisc}\\
      \lambda_{I,k-1}&= \lambda_{I,k}+ha+h\lambda_{S,k}\left(\gamma-\beta\frac{S_k}{N_k}+\beta\frac{I_kS_k}{N_k^2}\right)\nonumber\\
       &\quad\;\; +h\lambda_{I,k}\left(\beta\frac{S_k}{N_k}-\beta\frac{I_kS_k}{N_k^2}\right)
       +h\lambda_{R,k}(\alpha+ \eta v_k),\;\quad \text{ for all } k=1,\dots, n-1,\label{eqn: adjointIkdisc}\\
       \lambda_{I,n-1} &= 0,\label{eqn: adjointIndisc}\\
       \lambda_{R,k-1} &= \lambda_{R,k}+h\lambda_{S,k}\left(\gamma+\rho+\beta\frac{I_kS_k}{N_k^2}\right)
      -h\lambda_{I,k}\left(\beta\frac{I_kS_k}{N_k^2}\right)
       -h\lambda_{R,k}(\nu+\rho),\;\quad \text{ for all } k=1,\dots, n-1,\label{eqn: adjointRkdisc}\\
       \lambda_{R,n-1}&=0. \label{eqn: adjoinRndisc}
  \end{align}



  \begin{table}[ht]
    \centering
    \begin{tabular}{|c|c|l|}
    \hline
    \textbf{Parameter}   & \textbf{Numerical Value}   & \textbf{Parameter Description}\\
    \hline
    $\gamma$    & 0.00683           & birth rate of the population\\
    $\nu$       & 0.00188           & natural death rate of the population\\
    $\beta$     & 0.2426            & rate of infectiousness of the disease\\
    $\mu$       & 0.005             & disease induced death rate\\
    $\alpha$    & 0.00002           & rate at which disease is overcome\\
    $\rho$      & 0.007             & resensitization rate\\
    $\kappa$    & 0.3               & effectiveness of vaccination\\
    $\eta$      & 0.1               & effectiveness of treatment\\
    $a$         & 5                 & weight of $I$ used in cost functional\\
    $b$         & 50                & weight for cost of vaccination\\
    $c$         & 300               & weight for cost of treatment\\
    $T$         & 50                & time horizon (weeks) \\
    $u_{\max}$  & 1                 & maximum vaccination rate\\
    $v_{\max}$  & 1                 & maximum treatment rate\\
    $S_0$       & 1000              & initial condition for $S$\\
    $I_0$       & 10                & initial condition for $I$\\
    $R_0$       & 0                 & initial condition for $R$\\
    \hline
    \end{tabular}
    \caption{Numerical Values of Parameters Used in Computations}
    \label{tab: SIRparameters}
  \end{table}
  \begin{figure}[htbp]
  \begin{subfigure}[t]{0.33\textwidth}
   \includegraphics[width=\linewidth]{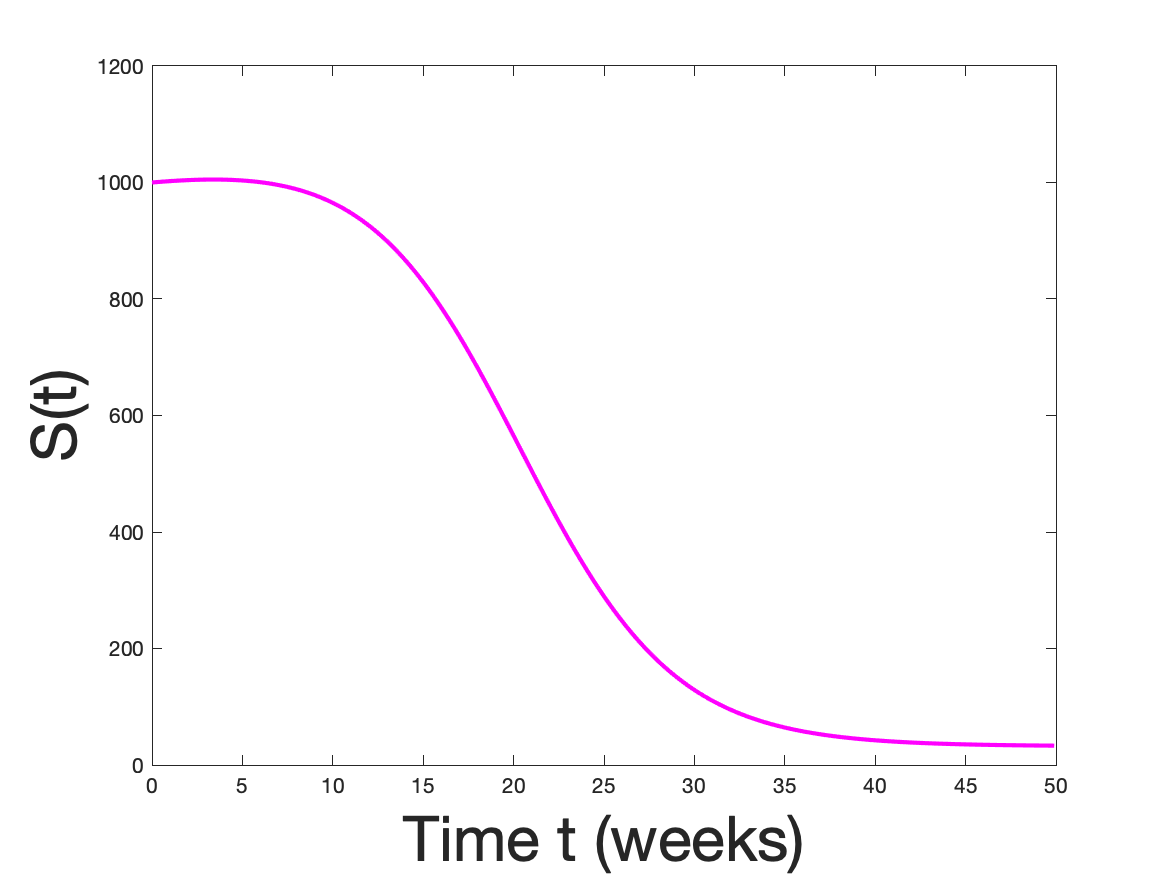}
   \caption{Susceptible Class}
  \end{subfigure}\hfill
  \begin{subfigure}[t]{0.33\textwidth}
   \includegraphics[width=\linewidth]{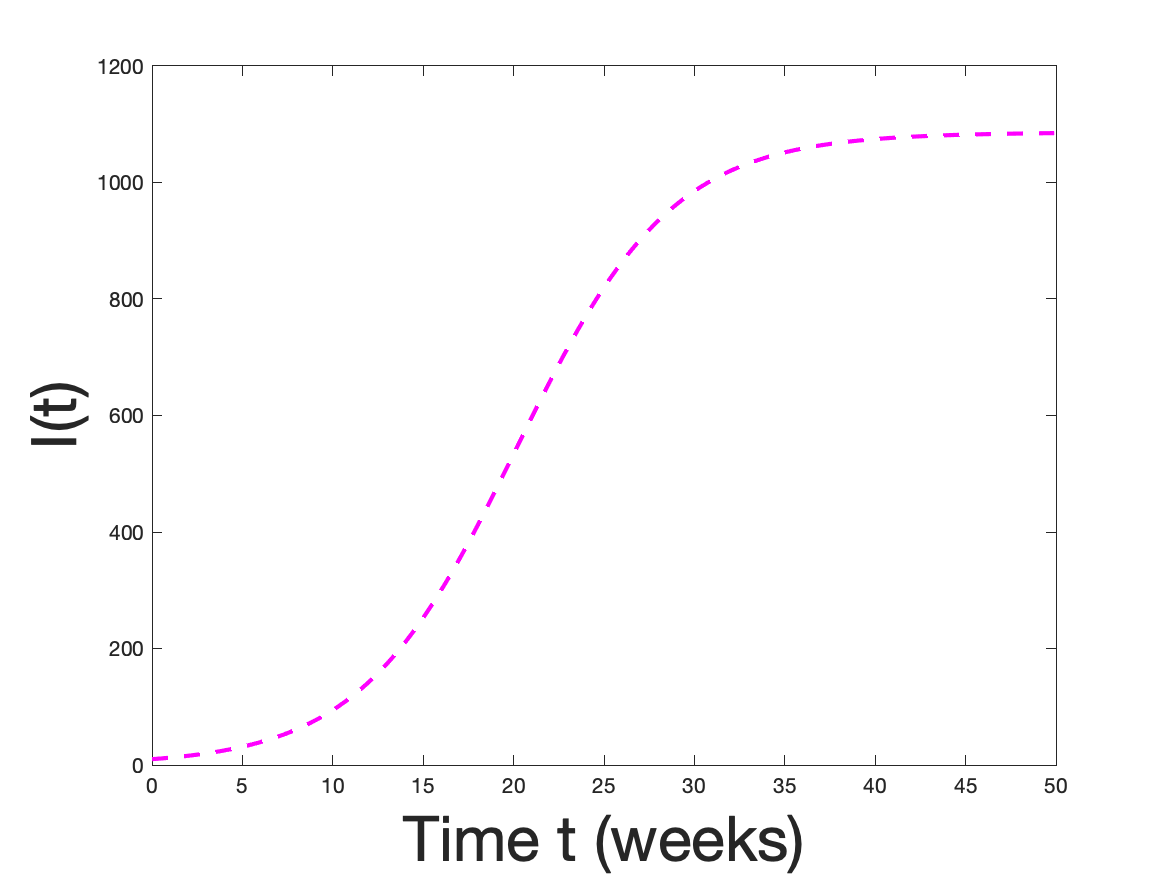}
   \caption{Infectious Class}
  \end{subfigure}\hfill
  \begin{subfigure}[t]{0.33\textwidth}
   \includegraphics[width=\linewidth]{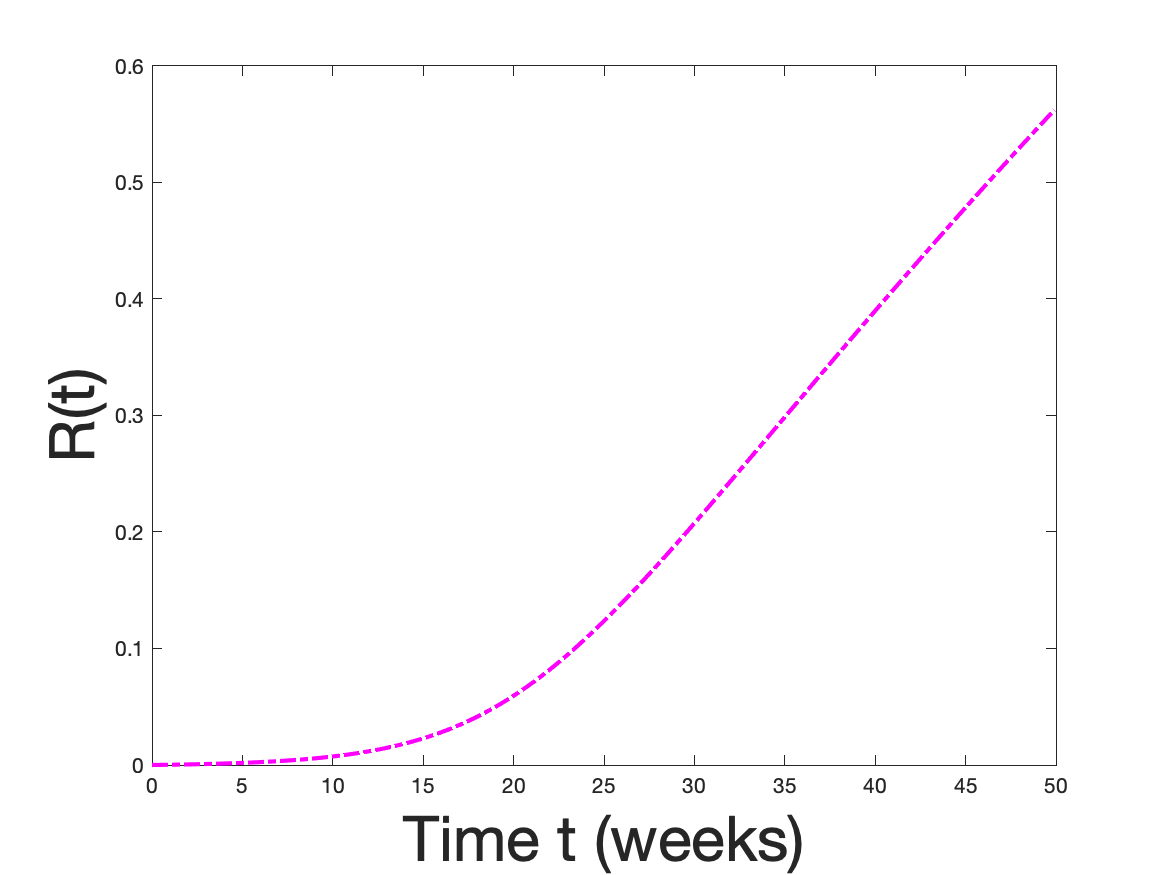}
   \caption{Recovered Class }
  \end{subfigure}\hfill
  \caption{Spread of Disease Without Vaccination and Treatment}
   \label{fig: SIRstatesnoconrols}
  \end{figure}
 
   \begin{figure}[htbp]
  \begin{subfigure}[t]{0.33\textwidth}
   \includegraphics[width=\linewidth]{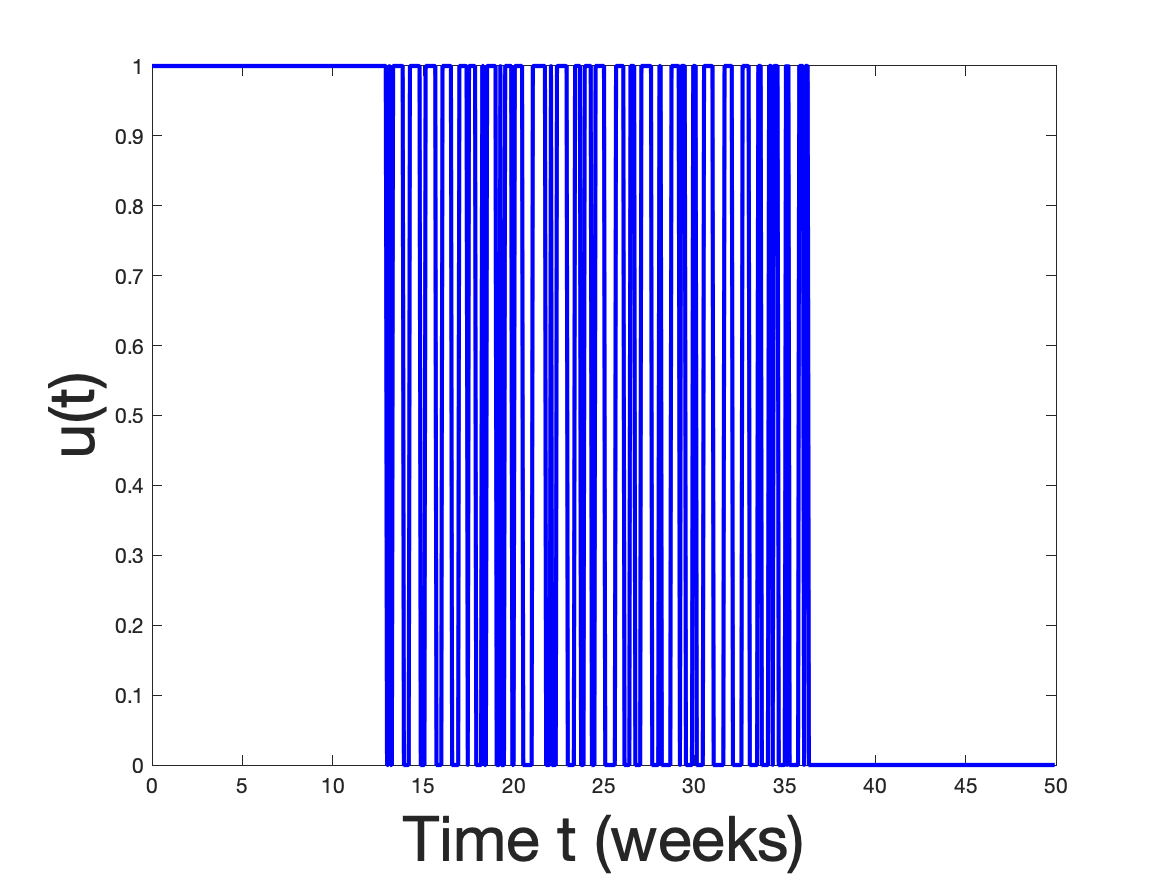}
   \caption{Unpenalized Vaccination $\hat{u}$}
   \label{subfig: nopenuSIR}
  \end{subfigure}\hfill
  \begin{subfigure}[t]{0.33\textwidth}
   \includegraphics[width=\linewidth]{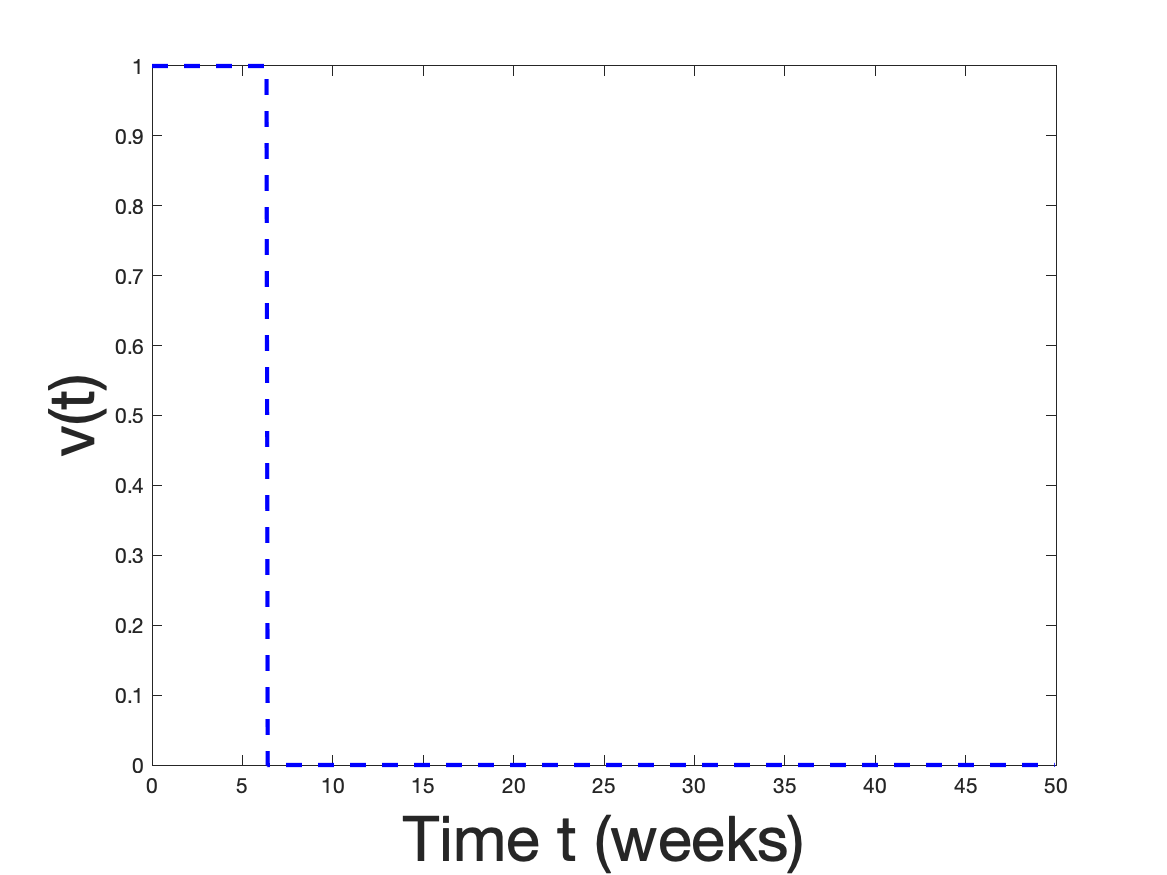}
   \caption{Unpenalized Treatment $\hat{v}$}
   \label{subfig: nopenvSIR}
  \end{subfigure}\hfill
  \begin{subfigure}[t]{0.33\textwidth}
  \includegraphics[width=\linewidth]{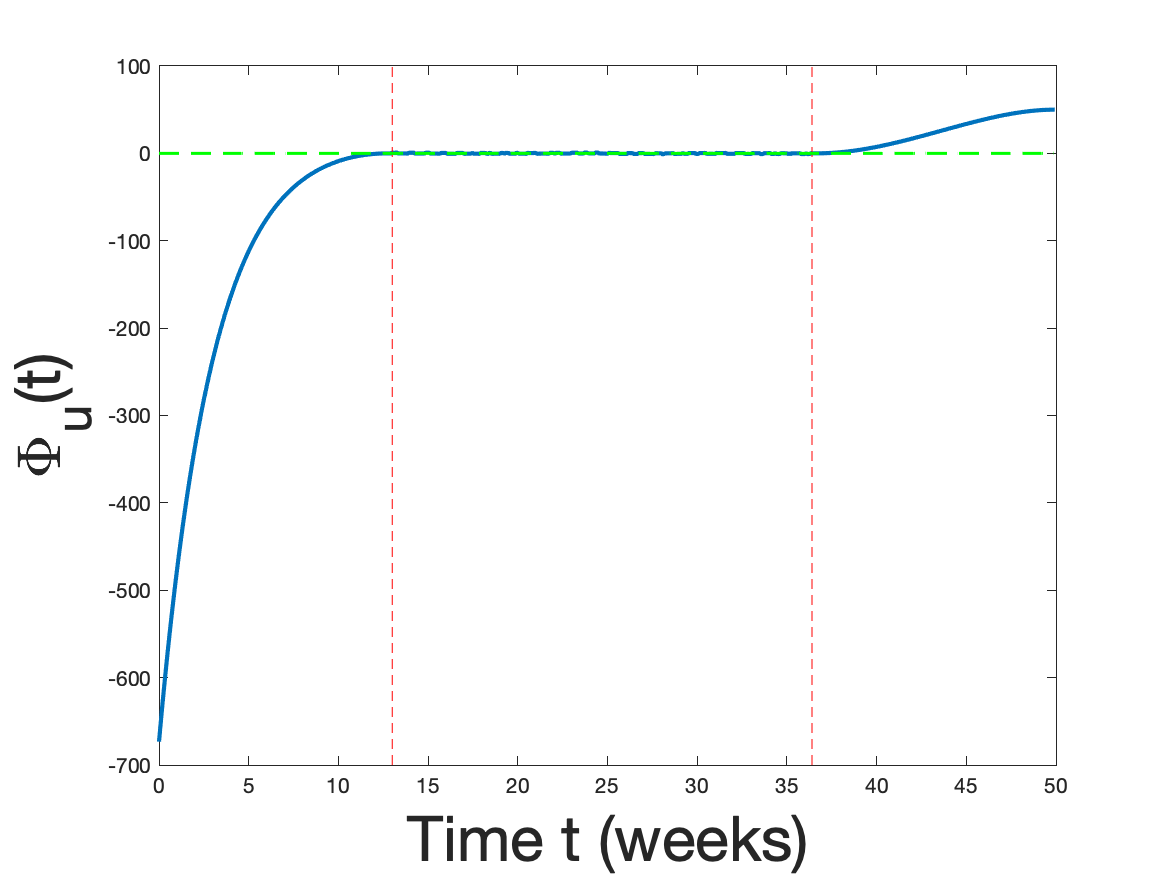}
  \caption{Switching Function $\Phi_{u}$ for $\hat{u}$}
  \label{subfig: nopenswitchuSIR}
  \end{subfigure}
   \begin{subfigure}[t]{0.33\textwidth}
   \includegraphics[width=\linewidth]{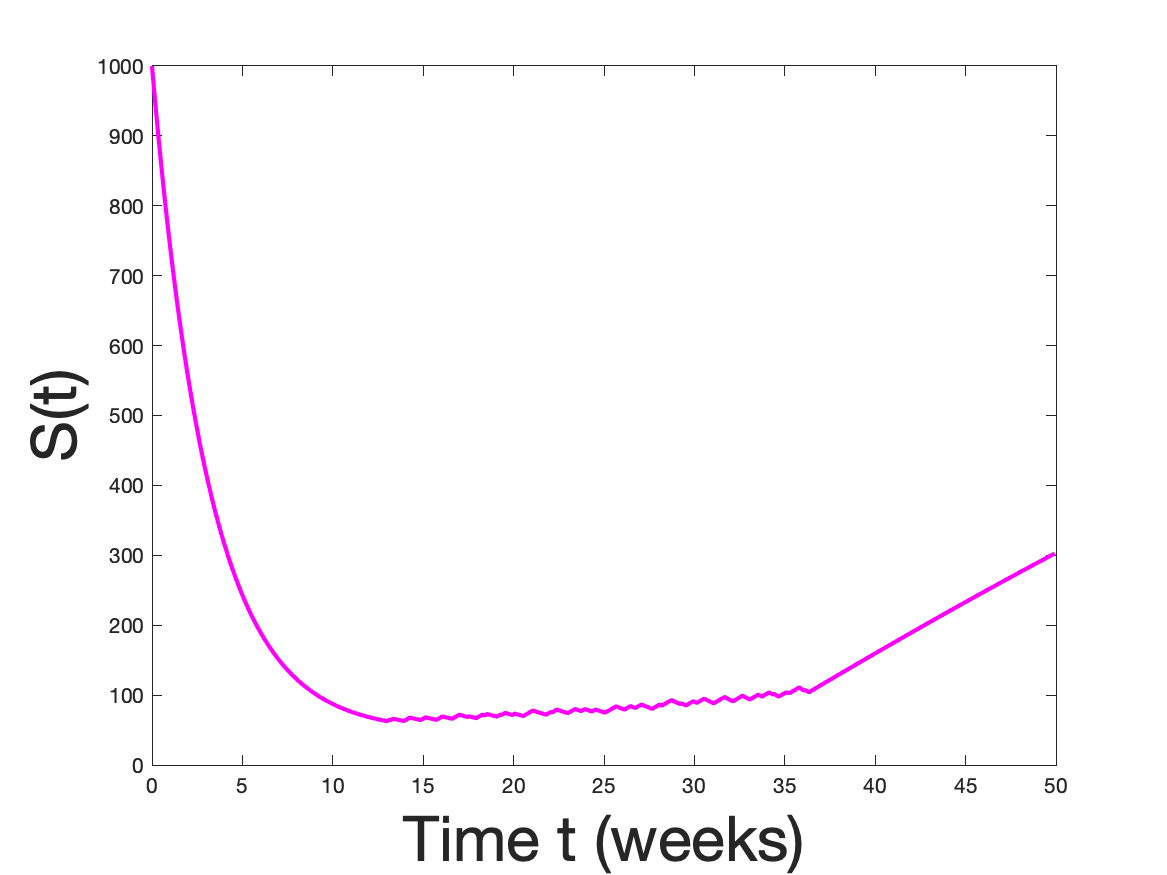}
   \caption{Susceptible Class}
   \label{subfig: nopensSIR}
  \end{subfigure}\hfill
  \begin{subfigure}[t]{0.33\textwidth}
   \includegraphics[width=\linewidth]{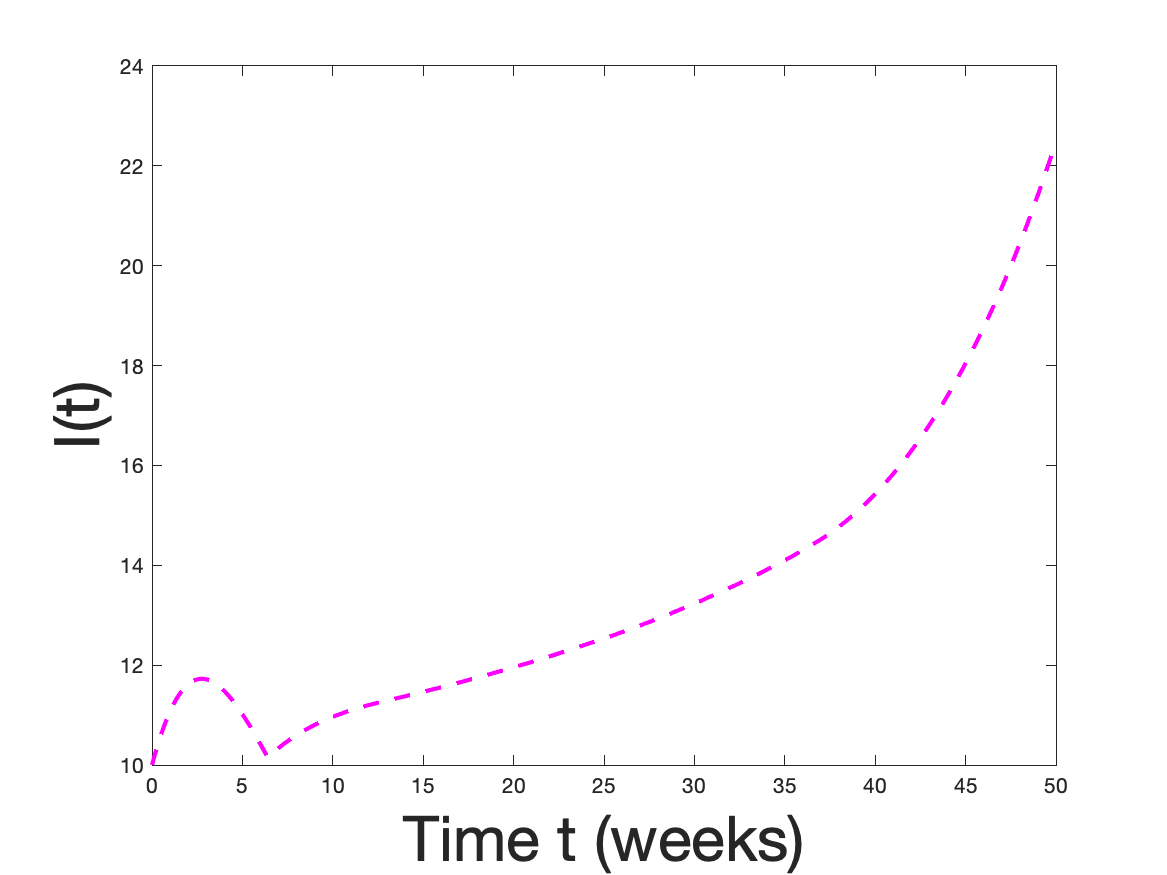}
   \caption{Infectious Class}
   \label{subfig: nopeniSIR}
  \end{subfigure}\hfill
  \begin{subfigure}[t]{0.33\textwidth}
   \includegraphics[width=\linewidth]{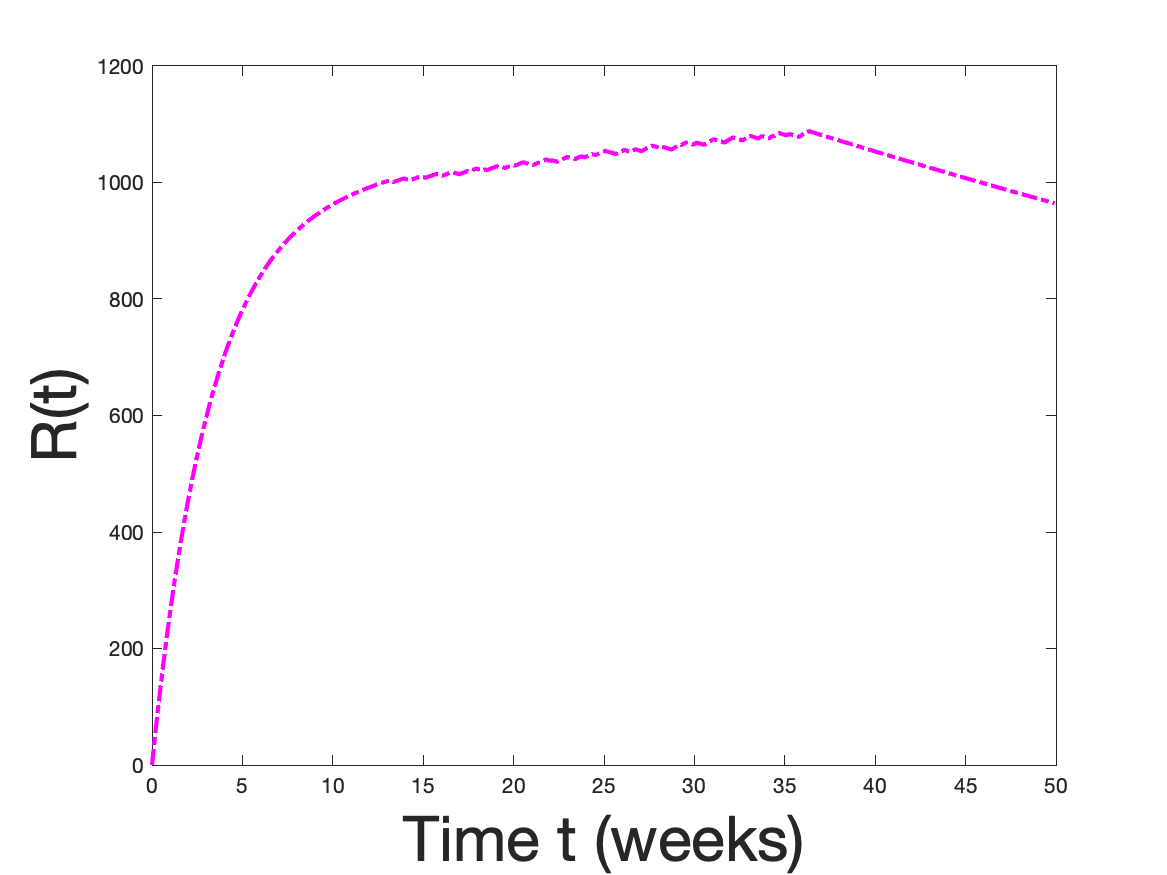}
   \caption{Recovered Class }
   \label{subfig: nopenrSIR}
  \end{subfigure}\hfill
  \caption{Spread of Disease with Unpenalized Vaccination and Treatment}
   \label{fig: SIRNoPen}
  \end{figure}
 \begin{figure}[htbp]
  \begin{subfigure}[t]{0.43\textwidth}
  \includegraphics[width=\linewidth]{images/Siru.png}
  \caption{$u_{p}$ with $p=0$}
  \label{subfig: SIRvaryu}
  \end{subfigure}\hfill
  \begin{subfigure}[t]{0.435\textwidth}
  \includegraphics[width=\linewidth]{images/SIRswitchu.png}
  \caption{$\Phi_{u_{p}}$ with $p=0$}
  \label{subfig: SIRvaryswitch}
  \end{subfigure}\hfill
  \begin{subfigure}[t]{0.435\textwidth}
  \includegraphics[width=\linewidth]{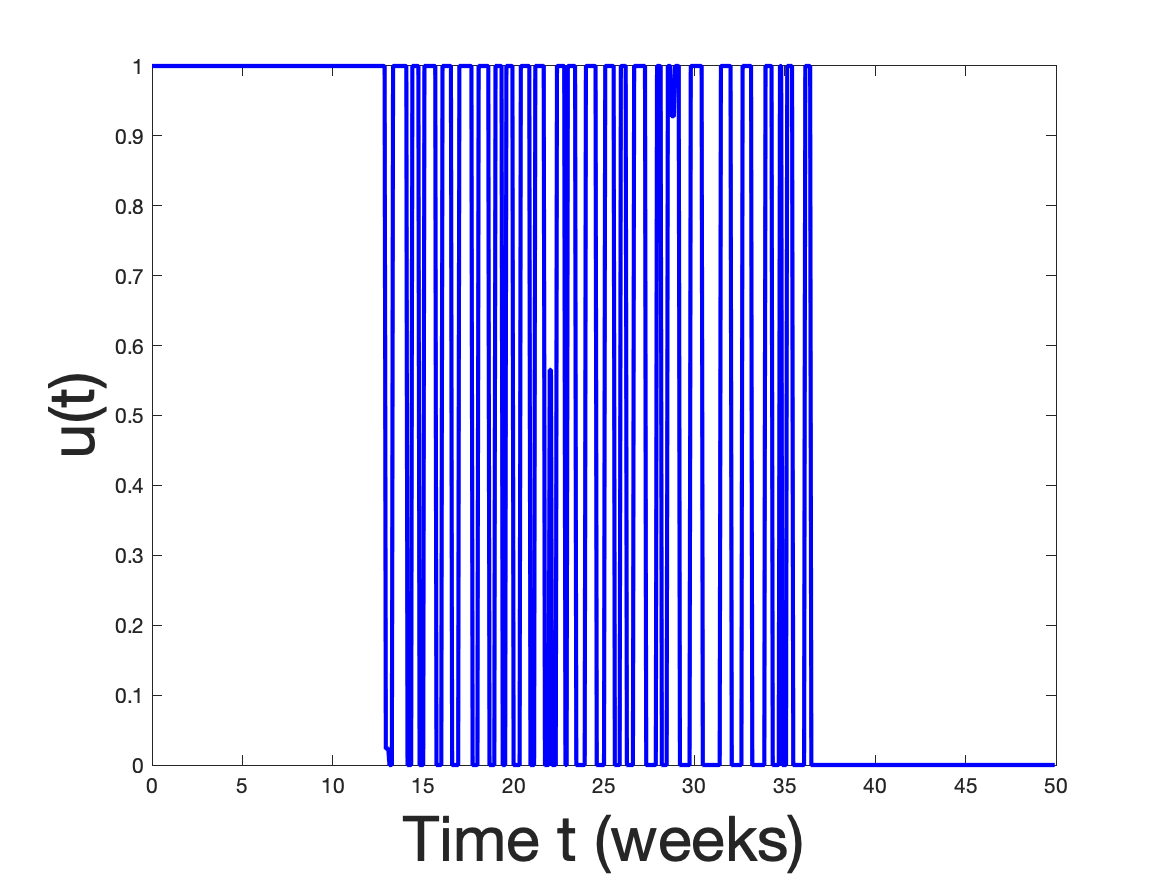}
  \caption{$u_p$ with $p=10^{-3}$}
  \label{subfig: SIRvaryu3}
  \end{subfigure}\hfill
  \begin{subfigure}[t]{0.435\textwidth}
  \includegraphics[width=\linewidth]{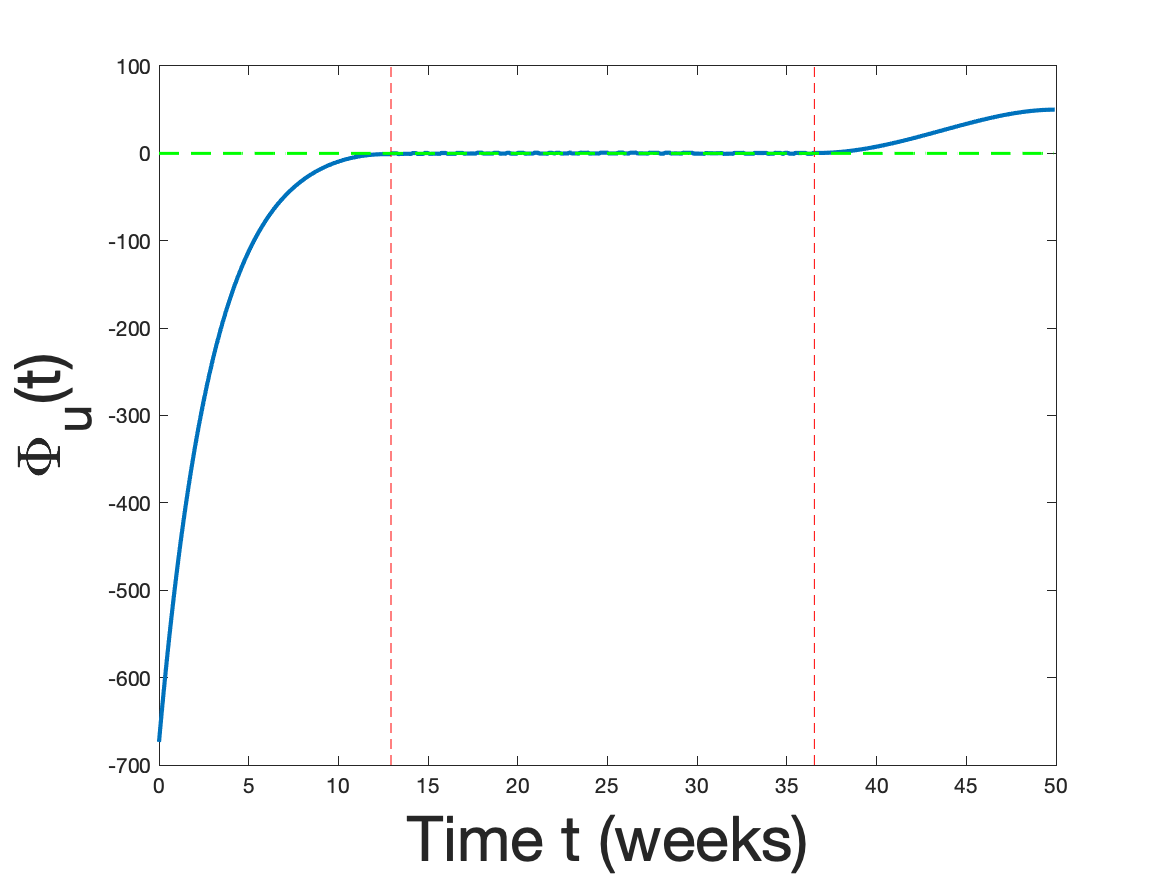}
  \caption{$\Phi_{u_p}$ with $p=10^{-3}$}
  \label{subfig: SIRvaryswitch3}
  \end{subfigure}\hfill
  \begin{subfigure}[t]{0.435\textwidth}
  \includegraphics[width=\linewidth]{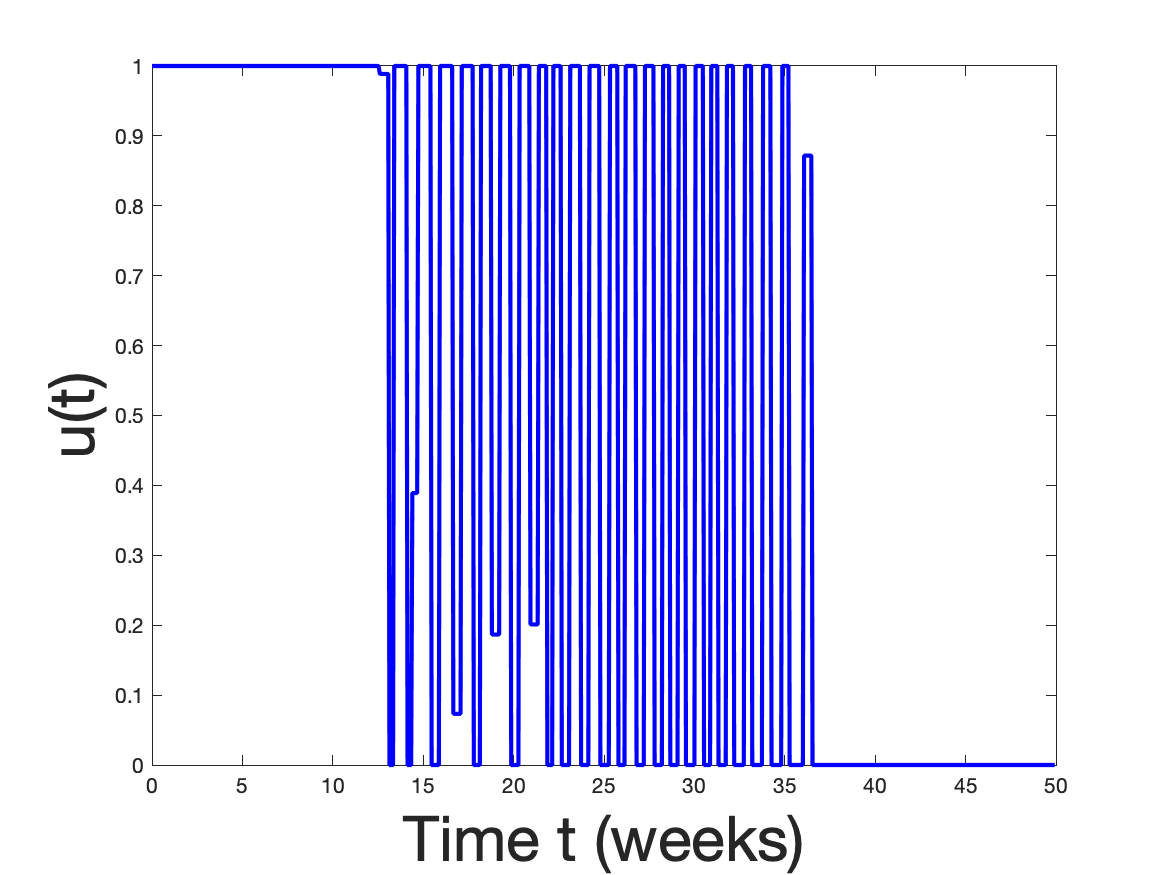}
  \caption{$u_p$ with $p=10^{-2}$}
  \label{subfig: SIRvaryu2}
  \end{subfigure}\hfill
  \begin{subfigure}[t]{0.435\textwidth}
  \includegraphics[width=\linewidth]{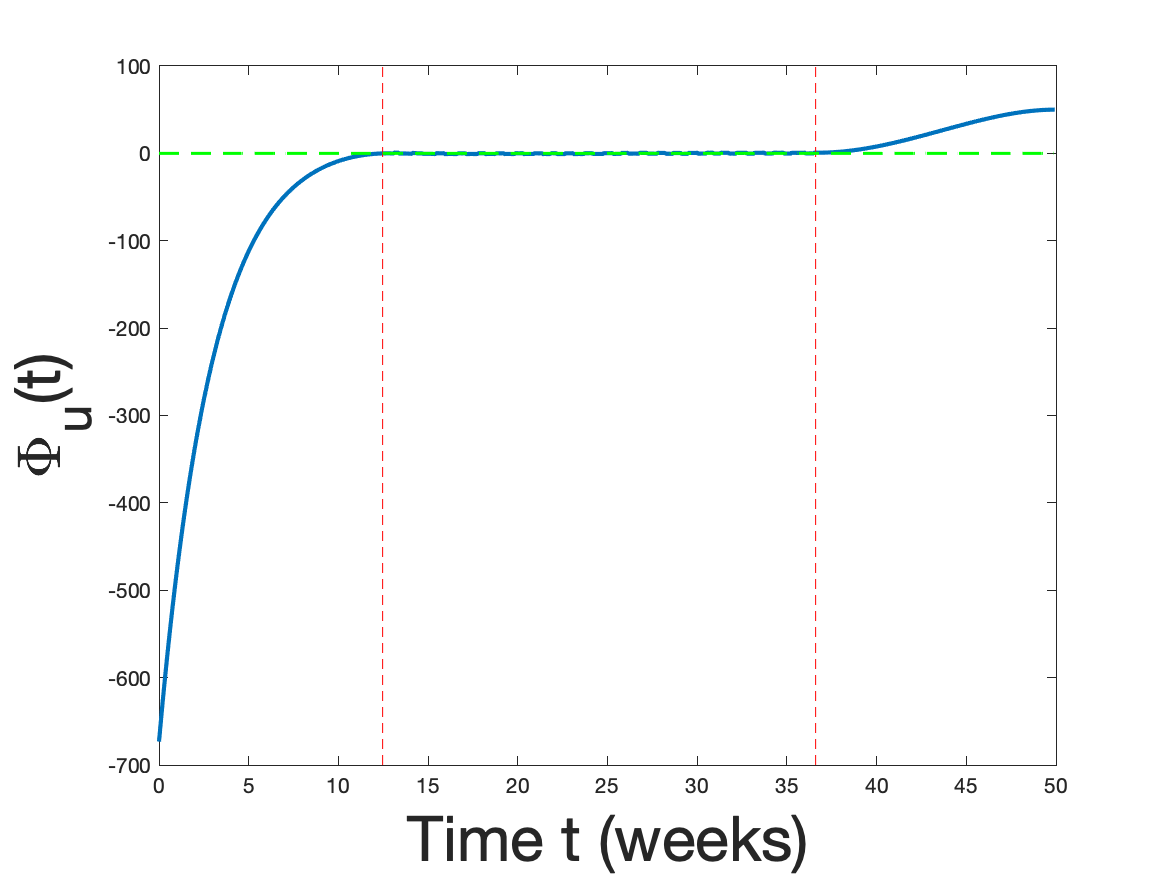}
  \caption{$\Phi_{u_p}$ with $p=10^{-2}$}
  \label{subfig: SIRvaryswitch2}
  \end{subfigure}\hfill
   \begin{subfigure}[t]{0.435\textwidth}
  \includegraphics[width=\linewidth]{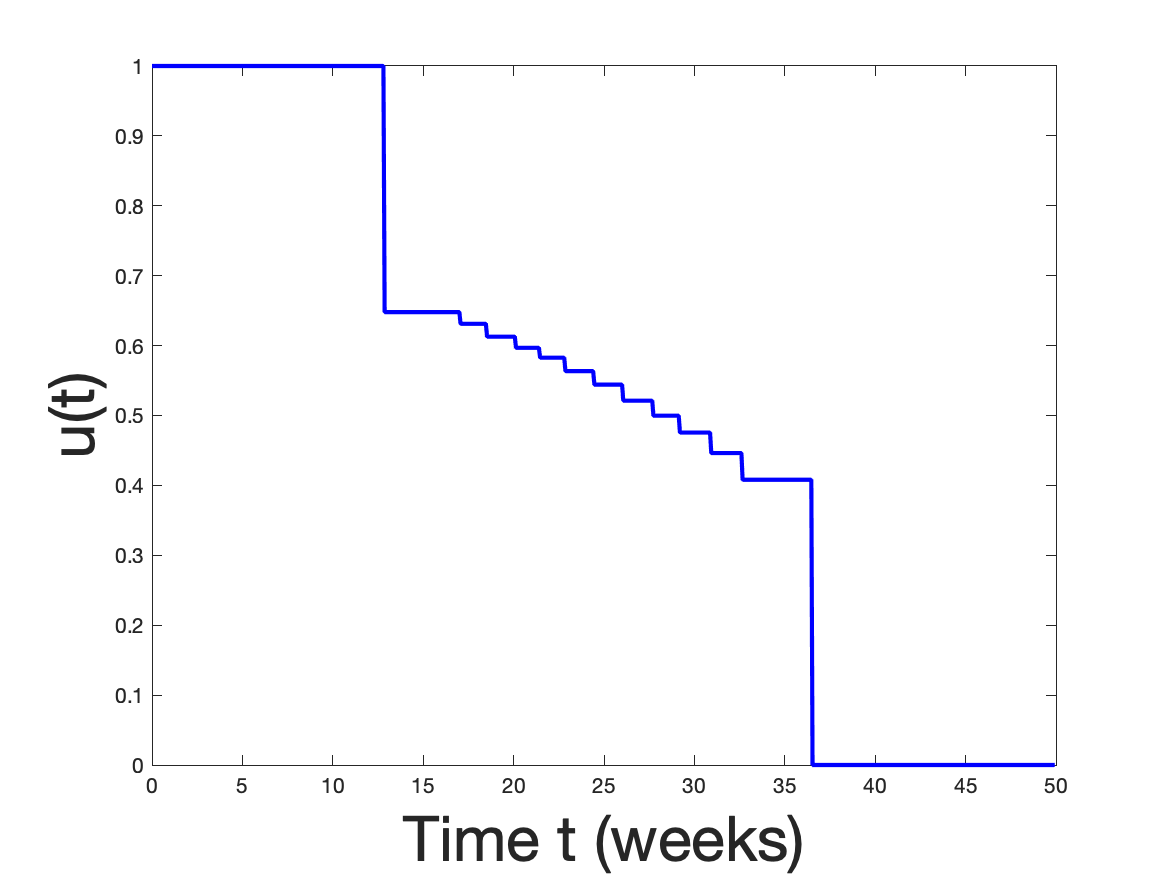}
  \caption{$u_p$ with $p=10^{-1}$}
  \label{subfig: SIRvaryu1}
  \end{subfigure}\hfill
  \begin{subfigure}[t]{0.435\textwidth}
  \includegraphics[width=\linewidth]{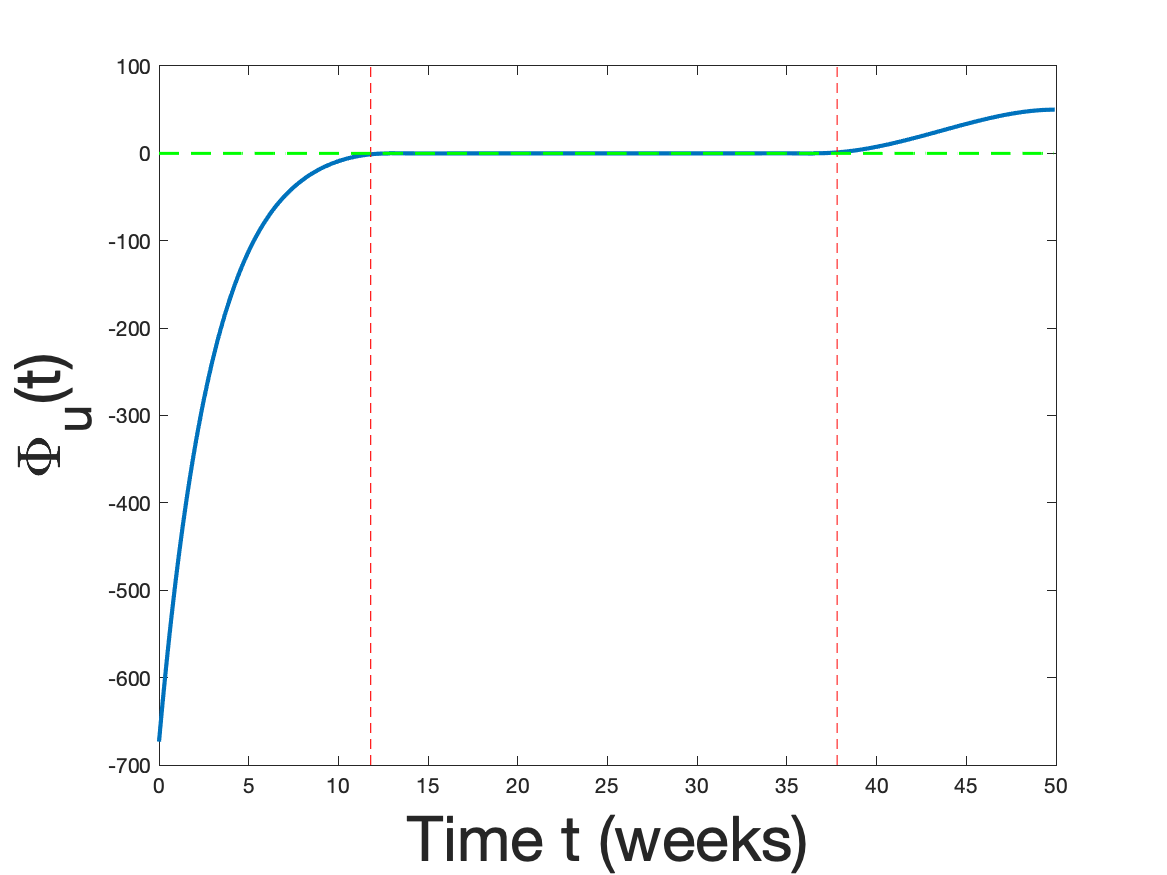}
  \caption{$\Phi_{u_p}$ with $p=10^{-1}$ }
  \label{subfig: SIRvaryswitch1}
  \end{subfigure}\hfill
  \caption{Varying Penalty Parameter for SIR Problem}
   \label{fig: SIRvary}
  \end{figure}
  
  \begin{table}[ht]
   \centering
  \begin{tabular}{|c|c c c c c |}
  \hline
  & $p$ & $J(u_p)$ & $t_{u,1}$ & $t_{u,2}$& Runtime\\
  \hline
  $N=750$       & 0          & 6572.955 & 12.933 & 36.533 & 16.77 \\
  $h=0.0667$    & $10^{-5}$  & 6572.956 & 13.067 & 36.400 & 17.07 \\
  $\text{tol}=10^{-8}$ & $10^{-4}$  & 6572.948 & 13.133 & 36.400 & 20.58 \\
                & $10^{-3}$  & 6573.018 & 12.933 & 36.533 & 19.12\\
                & $10^{-2}$  & 6573.101 & 12.467 & 36.000 & 56.58\\
                & $10^{-1}$  & 6575.429 & 12.333 & 37.067 & 10.47\\
 \hline
  \end{tabular}
  \caption{ Varying Penalty Results: $J(u_p)$ is unpenalized cost functional being approximated at $u_p$, $t_{u,1}$ and $t_{u,2}$ are the approximated switches corresponding to each solution. The Runtime column corresponds to the time (in seconds) it took to run PASA for each problem.}
  \label{tab: varySIR}
  \end{table}
  
 \begin{figure}[htbp]
  \begin{subfigure}[t]{0.33\textwidth}
   \includegraphics[width=\linewidth]{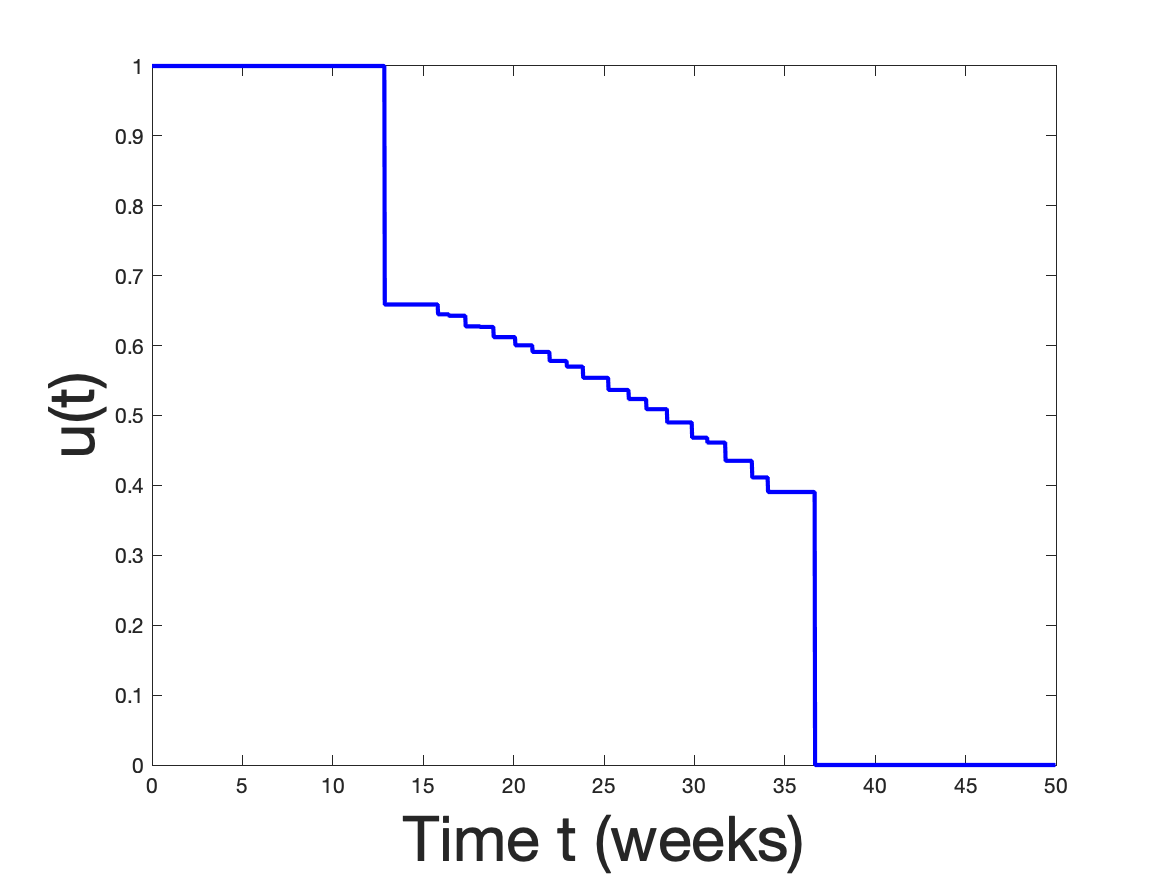}
   \caption{Penalized Vaccination $u_{p}$}
   \label{subfig: SIRu1n}
  \end{subfigure}\hfill
  \begin{subfigure}[t]{0.33\textwidth}
   \includegraphics[width=\linewidth]{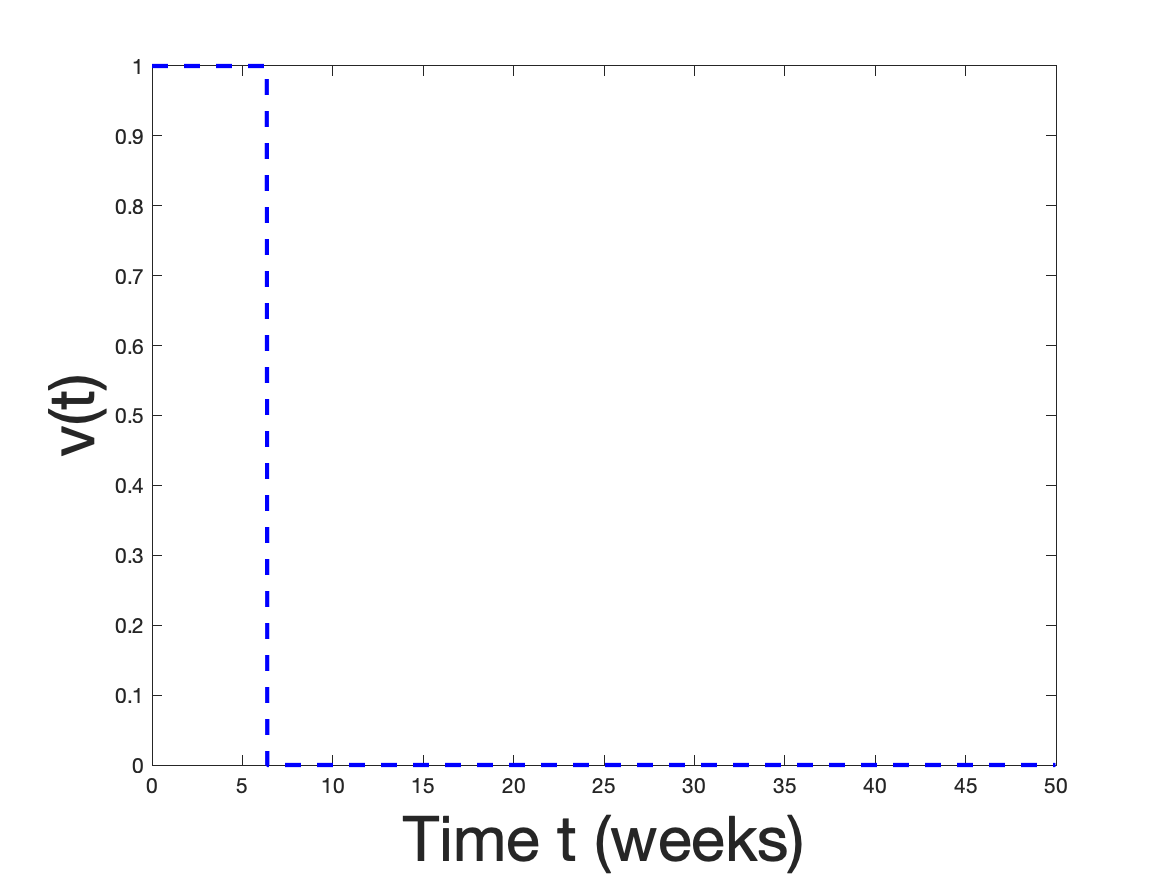}
   \caption{Treatment $v$}
   \label{subfig: SIRv1n}
  \end{subfigure}\hfill
  \begin{subfigure}[t]{0.33\textwidth}
   \includegraphics[width=\linewidth]{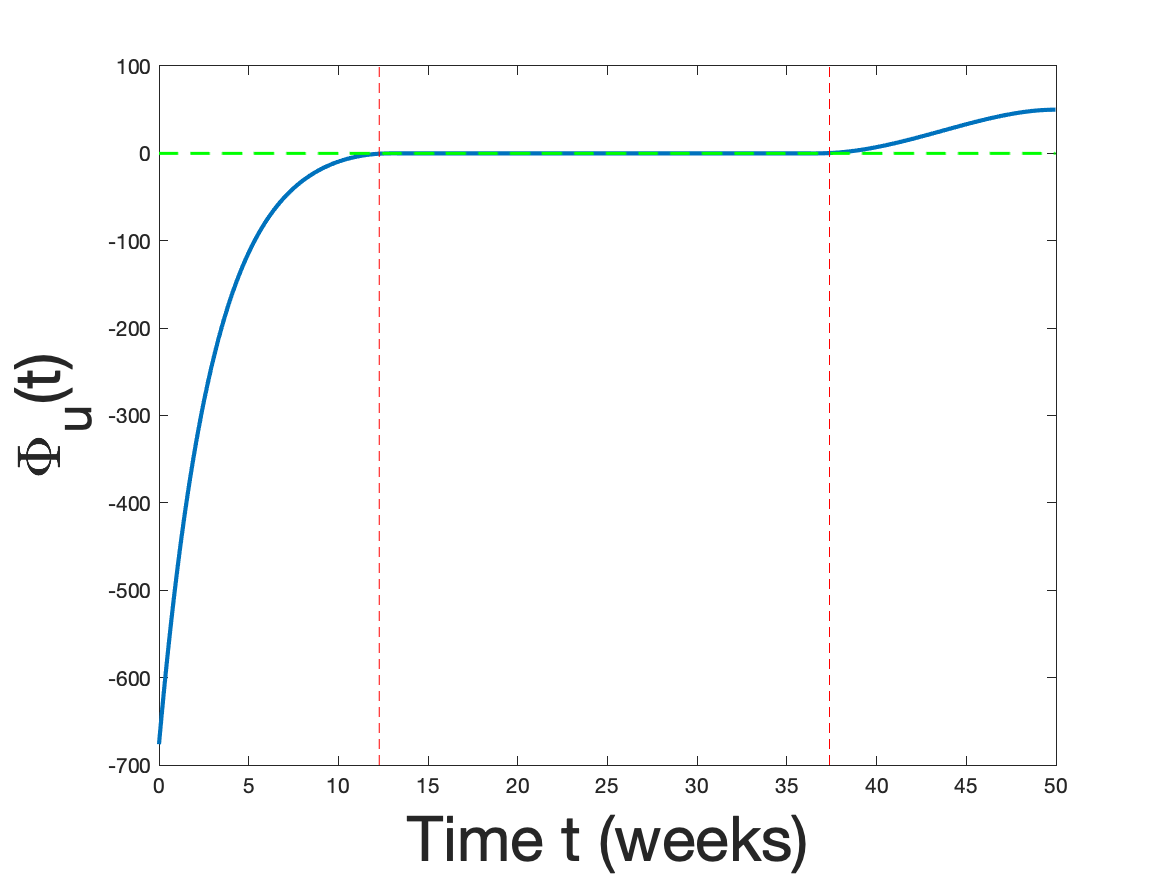}
   \caption{Switching Function, $\Phi_{u_{p}}$ for $u_p$}
   \label{subfig: SIRswitchu1n}
  \end{subfigure}\hfill
   \begin{subfigure}[t]{0.33\textwidth}
   \includegraphics[width=\linewidth]{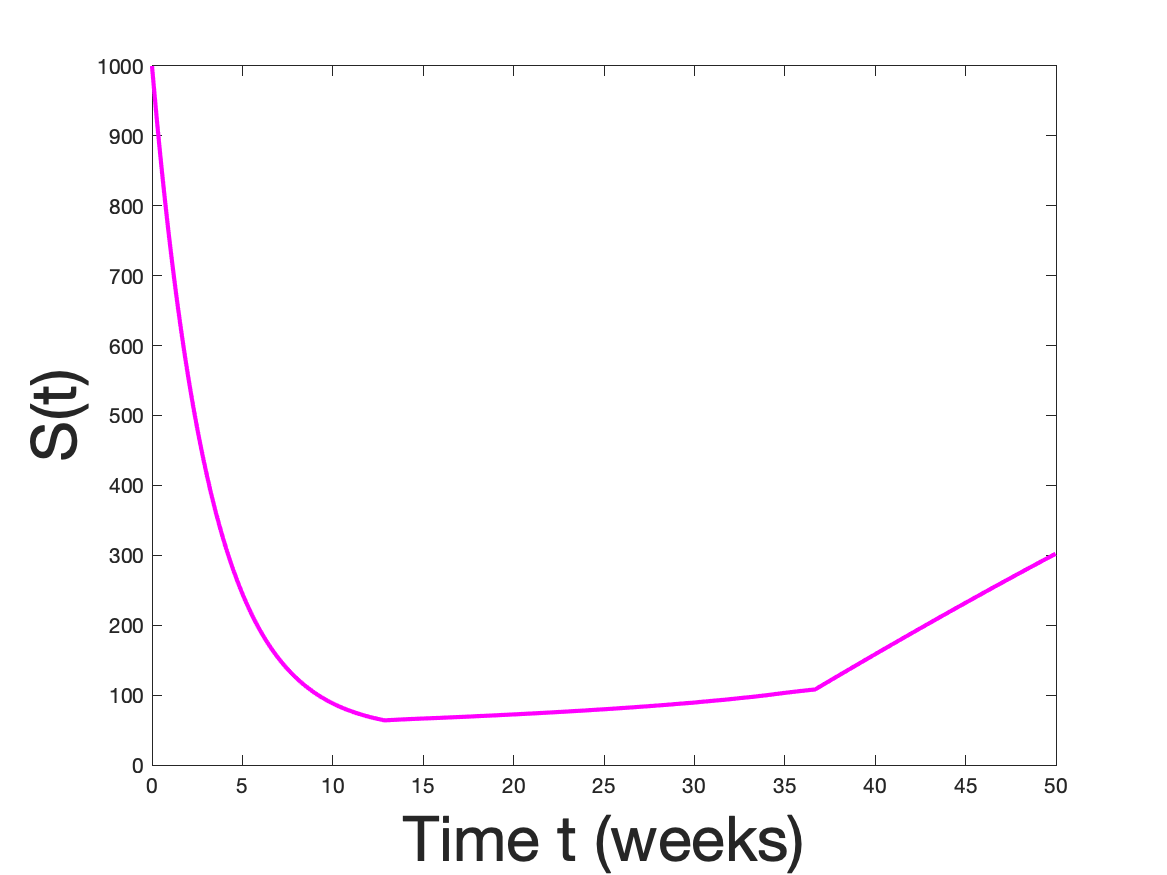}
   \caption{Susceptible Class}
   \label{subfig: SIRs1n}
  \end{subfigure}\hfill
  \begin{subfigure}[t]{0.33\textwidth}
   \includegraphics[width=\linewidth]{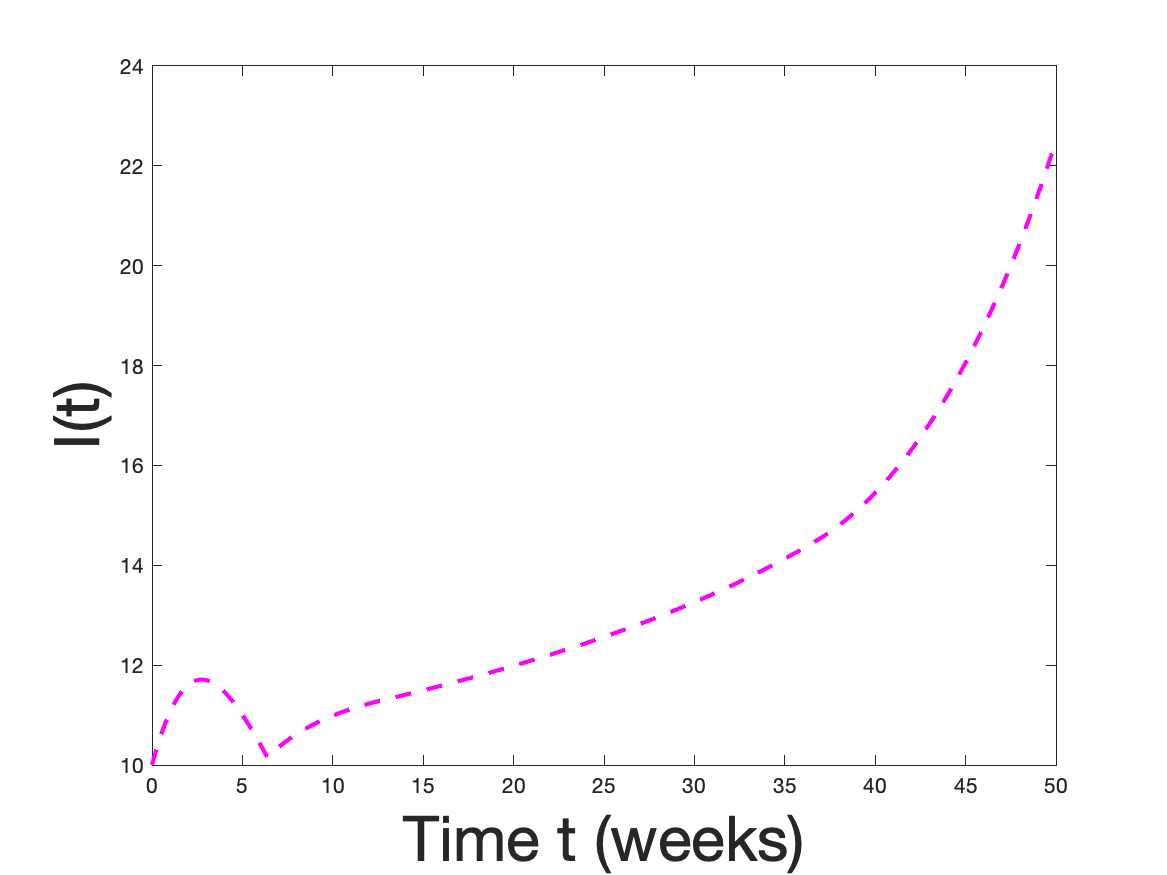}
   \caption{Infectious Class }
   \label{subfig: SIRi1n}
  \end{subfigure}\hfill
  \begin{subfigure}[t]{0.33\textwidth}
   \includegraphics[width=\linewidth]{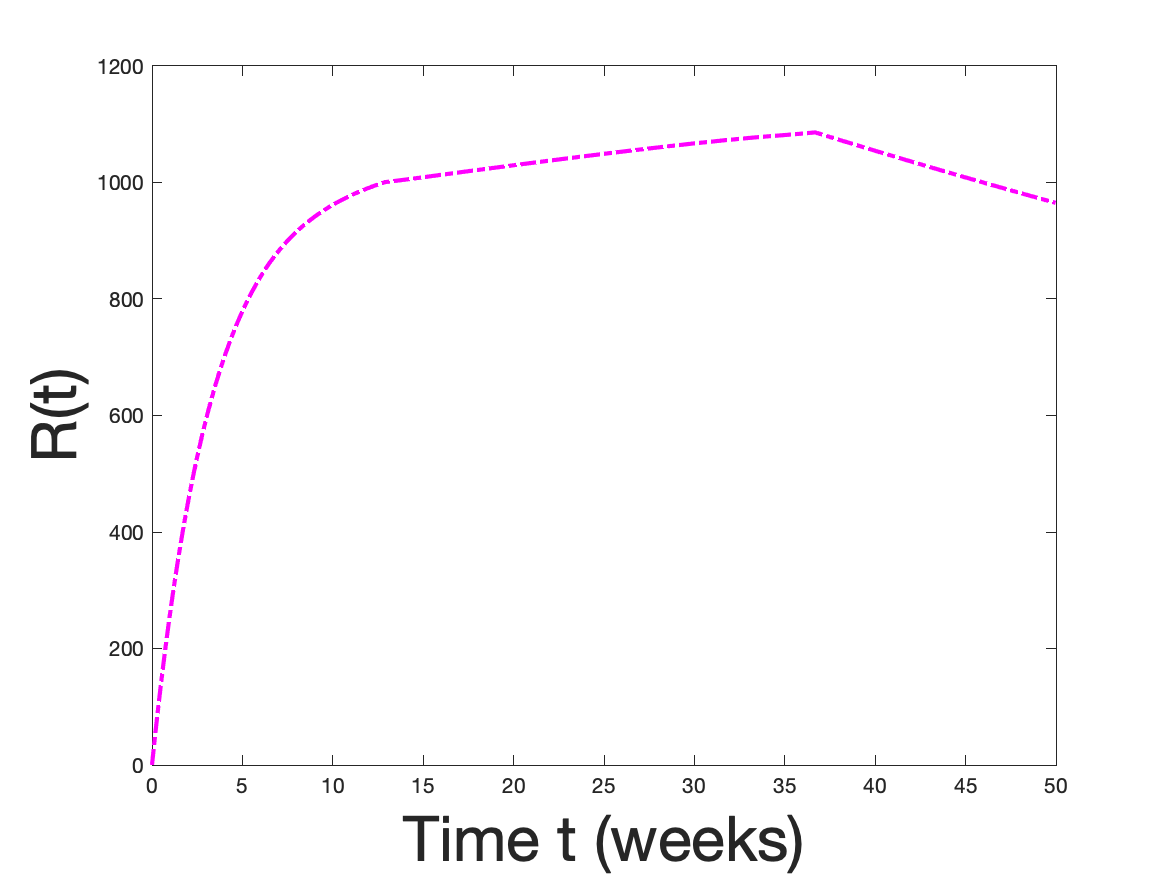}
   \caption{Recovered Class}
   \label{subfig: SIRr1n}
  \end{subfigure}\hfill
  \caption{Spread of Disease with Penalized Vaccination and Treatment. Penalty Parameter: $p=10^{-1}$, mesh intervals: $n=2000$, and Error Tolerance: $10^{-8}$ }
   \label{fig: SIRPen1n}
  \end{figure}
  
  \subsection{Numerical Results of SIR Problem}
   We first use PASA to numerically solve for the unpenalized problem given in equation (\ref{eqn: SIRproblem}) with stopping tolerance set to being $10^{-8}$ and with our initial guesses for the vaccination and treatment control being $u(t)=0$ and $v(t)=0$ over the entire time interval. 
   {The numerical parameters, which are given in Table \ref{tab: SIRparameters}, are set to the values that were used in Ledzewicz et al. \cite{Mahya} because we want to see if our numerical results are comparable to theirs.}
    {Now, Ledzewicz et al. \cite{Mahya}  are investigating optimal vaccine and treatment strategies for a theoretical vaccine and treatment for Ebola, and one can infer from their parameter settings that they were not concerned with realistic limitations in widely developing and administering these theoretical controls.}
   {For example based on these parameter values, a vaccination strategy being $u(t)=u_{\max}=1$ implies that we are vaccinating approximately 97\% of the susceptible population immediately at time $t$.} 
     {When considering realistic limitations in widely-developing and administering vaccines, one might want to consider setting $u_{\max}$ to a value that is less than one.}
   %
   %
   We partition the time interval to where there are $n=750$ mesh intervals with mesh size being $h=0.0667$, and use the discretization methods presented in Subsection \ref{sec: SIRdiscrete} with penalty parameter set to being $p=0$. 
   The results corresponding to the approximated solution that PASA obtained for problem (\ref{eqn: SIRproblem}) are given in Figure \ref{fig: SIRNoPen}.
   Observe in Subfigure \ref{subfig: nopenuSIR} that the numerical solution  associated with the vaccination control, $\hat{u}$, begins with a full dose segment (i.e. $\hat{u}=1$), followed by an interval of many oscillations, and ends with the control being turned off.  
   The many oscillations approximately start at time $t_{u,1}=12.933$ and end at time $t_{u,2}=36.533$. 
   We looked at the corresponding switching function to $\hat{u}$, which is $\Phi_u(t)=b+(\lambda_R(t)-\lambda_S(t))\kappa S(t)$ (see \cite{Mahya} for construction of switching function) to verify that the oscillating region  corresponds to a singular subarc. 
   It is necessary that the optimal control $u^*$ will be as follows 
   \begin{equation}\label{eqn: SIRucases}
     u^*(t)=
     \begin{cases}
         0          & \text{ if } \Phi_u(t)>0,\\
         u_{\max}   & \text{ if } \Phi_u(t)<0,
     \end{cases}
   \end{equation}
   and $u^*(t)$ will be singular if $\Phi_u$ is zero over an open interval of time. 
   A plot of switching function $\Phi$ is given in Subfigure \ref{subfig: nopenswitchuSIR}. 
   The green horizontal line in Subfigure \ref{subfig: nopenswitchuSIR} corresponds to the $x$-axis, while the red vertical lines correspond to the approximated times when oscillations began, $t_{u,1}=12.933$,  and ended $t_{u,2}=36.533$. 
   Observe that the switching function is on the $x$-axis on interval $(t_{u,1}, t_{u,2})$. 
   It follows that the numerical control obtained for the unpenalized problem,  $\hat{u}$, is oscillating over its singular region. 
   In Subfigure \ref{subfig: nopenvSIR}, the numerically computed treatment control $\hat{v}$ is bang-bang that begins with a full dose segment, and the control turns off approximately at time $t_{v,1}=6.4$.
    
    Although we observe in Subfigures \ref{subfig: nopensSIR}-\ref{subfig: nopenrSIR} that the dynamics of $S$, $I$, and $R$ corresponding to the approximated controls $\hat{u}$ and $\hat{v}$ give significantly favorable results in comparison to the dynamics associated with no application of vaccination or treatment, which is given in Figure \ref{fig: SIRstatesnoconrols}, we recognize that the wild oscillations found in $\hat{u}$ makes it an unrealistic strategy to implement. 
    In order to obtain a more realistic vaccination strategy, we penalize control $u$ via a total variation penalty term \cite{Caponigro2018}. 
    We will not be penalizing control $v$, since control $v$ is a piecewise constant function that can be easily interpreted.  
    We use PASA to solve for penalized problem (\ref{eqn: SIRproblemPen}) with stopping tolerance set to $10^{-8}$ and with parameter value settings given in Table \ref{tab: SIRparameters}. 
    Our initial guess for controls are $u(t)=0$ and $v(t)=0$ over the entire time interval. 
    We partitioned $[0,T]$ to $n=750$ mesh intervals of size $h=0.0667$, and the discretization procedure is given in Subsection \ref{sec: SIRdiscrete}. 
    Additionally the MATLAB file used for running PASA is given in Appendix Section \ref{subsec: matlab}. 
    We numerically solve for problem (\ref{eqn: SIRproblemPen}) for varying values of penalty parameter $p\in\{10^{-5},10^{-4},10^{-3}, 10^{-2}, 10^{-1}\}$. 
    If the approximated control, $u_p$, corresponding to penalty parameter value $p$ is no longer oscillating and if the sign of the approximated switching, $\Phi_{u_p}$ aligns with equation (\ref{eqn: SIRucases}), then we consider $p$ as being a suitable penalty parameter value to use. 
    
    In Table \ref{tab: varySIR}, we record an approximation of the unpenalized cost functional, given in equation (\ref{eqn: SIRobjective}), being evaluated at the penalized solution $u_p$. 
    The values corresponding to $J(u_p)$ are remarkably similar to the approximated objective value that corresponds to the unpenalized solution $\hat{u}$.
    We also record PASA's runtime  performance (in seconds) on solving problem (\ref{eqn: SIRproblemPen}), and we record the approximated switches of $u_{p}$, $t_{u,1}$ and $t_{u,2}$.
    When numerically solving for problem (\ref{eqn: SIRproblemPen}) we observed chattering in all cases except for when $p=10^{-1}$. 
    In Figure \ref{fig: SIRvary}, we provide the plots of the penalized solution $u_p$  and the plots of the corresponding switching function when the penalty parameter was set to being $10^{-3},$ $10^{-2}$, and $10^{-1}$.
    Subfigures \ref{subfig: SIRvaryu3}, \ref{subfig: SIRvaryu2}, and \ref{subfig: SIRvaryu1} illustrate a trend of how increasing the penalty parameter $p$ influences the penalized solution to where chattering no longer occurs. 
    We could increase the penalty parameter values larger than $10^{-1}$; however, we did not do so because the approximated switching function associated $u_p$ when $p=10^{-1}$, given  in Subfigure \ref{subfig: SIRvaryswitch1}, gave indication that the structure of $u_p$ aligns with equation (\ref{eqn: SIRucases}). 
    Meaning we have a non-oscillatory vaccination strategy that is not only more realistic to implement, but also satisfies the first order necessary conditions of optimality \cite{Pontryagin}. 
    In addition, we did not increase the penalty parameter values larger than $10^{-1}$ because of the possibility of over-penalization. 
    In terms of runtime performances, we have that PASA performed the fastest when the penalty parameter was set to being $10^{-1}$.
    
    We recognize in Subfigure \ref{subfig: SIRvaryu1} that the singular region associated with $u_p$ when the penalty parameter is set to $p=10^{-1}$ is a rather crude approximation. 
    However, if we want a smoother presentation of the singular region we recommend either increasing the number of mesh intervals used for the discretization and/or  reducing the stopping tolerance. 
    In Subfigure \ref{subfig: SIRu1n} we provide a solution that PASA obtained in solving for penalized problem (\ref{eqn: SIRproblemPen}) with $p=10^{-1}$ when time interval $[0,T]$ was partitioned to have $2000$ mesh intervals with mesh size being $h=0.025$. 
    The approximated singular subarc that is given in Subfigure \ref{subfig: SIRu1n} is comparable to the approximation that Ledezewicz et al. \cite{Mahya} obtained. 
   
    Although these results are theoretical due to no vaccination for Ebola, optimal control theory can be used to see how vaccination and treatment can influence the spread of the disease. 
    We observe in Figure \ref{fig: SIRstatesnoconrols}, that without treatment and vaccination, the susceptible class starts to decline due to many individuals becoming infectious and few individuals recover from the disease. 
    When incorporating the vaccination strategy and treatment strategy that is respectively given in Subfigures \ref{subfig: SIRu1n} and \ref{subfig: SIRv1n}, we observe in Subfigures \ref{subfig: SIRs1n}-\ref{subfig: SIRr1n} a large portion of the  vaccinated susceptible class become resistant to Ebola and hence transfer to the Recovered class.
    The positive influence that treatment has on the infectious class occurs approximately in the first 6 weeks when the treatment strategy is set to treat all infectious individuals.
    We believe that the weighted cost of treatment parameter value $c$ is influencing the optimal treatment strategy to switch to incorporating no treatment at $t_{v,1}\approx 6.4$ weeks. 
    Reducing the weighted cost of treatment parameter would alter the optimal treatment strategy to extending the region(s) when $v^*=v_{\max}$. 
    
    \section{Concluding Remarks}
     We demonstrate how Hager and Zhang's \cite{Hager2016} Polyhedral Active Set Algorithm (PASA) can be implemented for numerically solving optimal control problems that depend linearly with respect to the control. 
     If problems of this form possess a singular subarc, we recommend penalizing such problems by using a penalty term based on the total variation of the control as suggested in Caponigro et al. \cite{Caponigro2018}. 
     This penalty term allows for PASA to yield approximations that do not oscillate wildly over the singular region. 
     We provide a discretization method for a general optimal control problem that involves using Euler's method for the state variables and using the gradient of the Lagrangian of the discretized optimal control problem to approximate the gradient of the cost functional that is being used. 
     Additionally, the use of the Lagrangian consequently aids in the discretization of the adjoint equations.
     We then present in detail three optimal control problems that apply to biological models. 
     For the first two problems an explicit solution that satisfies the first order necessary conditions of optimality \cite{Pontryagin} can be obtained. 
     In the third problem, we can verify existence of a singular subarc for one of the control variables used. 
     For all three examples, we used PASA to numerically solve each problem with parameters set to where the existence of a singular subarc is possible. 
     For the example associated with the plant problem that was presented in King and Roughgarden's \cite{King1982}, three cases were investigated where each case determines the beginning behavior of the optimal allocation strategy. 
     PASA solved for the unpenalized plant problem for all three cases, and in all three cases no chattering was observed. 
     However, in the degenerative case, the case corresponding to the optimal control beginning singular, the unpenalized solution obtained had an unusual dip at $t=0$. 
     For this case, we penalized the problem via total variation and obtained a penalized solution that served to be a better approximation to the true solution. 
      When using PASA to solve for the unpenalized fishery problem \cite{Clark, Lenhart2007} and to solve for the unpenalized SIR problem that was presented in Ledzewicz et al.'s  \cite{Mahya}, we obtained numerical solutions that possessed oscillations along the singular region. 
     We then applied a total variation penalty term for varying values of the penalty parameter. 
     We found that the oscillatory solution could be controlled by increasing the penalty parameter size. 
     In these examples, we found that observing the plots of the approximated penalized solution and of the corresponding switching function to be an effective procedure for determining whether or not the penalty parameter size was over or under-penalizing the solution.  
     For the fishery problem, we used additional information such as the approximated switching points and the $L^1$  norm difference between the true solution and the penalized solution to verify that the penalty parameter value that was chosen based upon the plots of the penalized solution was the most appropriate choice in comparison to all of the other penalty parameters values being tested. 
     For the SIR problem, we found that the penalized solution that PASA obtained was comparable to the approximated solution that was presented in Ledzewicz et al. \cite{Mahya}, which was obtained via a collocation method called  PROPT \cite{PROPT}.

     Overall, the use of PASA and the use of total variation penalization make the process of numerically solving for singular control problems more tractable. 
     With the aid of PASA and total variation penalization, mathematical biologists can consider constructing an optimal control problem with a cost functional that possesses a more accurate representation of the cost in implementing the  control, which aligns more with the principle of parsimony than when using an optimal control problem with quadratic dependence. 
     Solving an optimal control problem with a more accurate representation of the cost functional ensures that the optimal control associated with the problem truly meets the criteria that are being considered for determining which control strategy was the best relative to the state dynamics being studied. 
     Additionally, using an optimal control problem where the control variables appear linearly can potentially yield optimal control solutions with regions that are both easy to interpret and implement. Incorporating a more practical, uncomplicated, and accurate optimal control into a dynamical system that is intended to model some biological phenomenon will surely provide more compelling and elucidating results that are applicable to practice and worthwhile to study. 

     \section{Acknowledgements}
	Support by the National Science Foundation (NSF) under Grants 1522629 and 1819002, and by the Office of Naval Research under Grants N00014-15-1-2048 and N00014-18-1-2100 is gratefully acknowledged. 
       S. Atkins was supported in part by the University of Florida's Graduate Student Fellowship. S. Atkins also recieved student travel support for the CMPD5 conference from NSF CMPD5 Grant 1917506, The Center for Applied Mathematics (CAM), and University of Florida's CLAS.  
      M. Martcheva was supported in part by NSF CMPD5 Grant 1917506 and from the University of Florida's Department of Mathematics.
	We also are grateful to Suzanne Lenhart for her helpful comments and suggestions of problems to explore. 
\bibliographystyle{amsplain}      
\bibliography{proceedings} 

 \section{Appendix}
  \subsection{Numercial Results Associated with Unpenalized Plant Problem}\label{subsec: plantresults}
   \begin{figure}[htbp]
   \begin{subfigure}[t]{0.5\textwidth}
    \includegraphics[width=\linewidth]{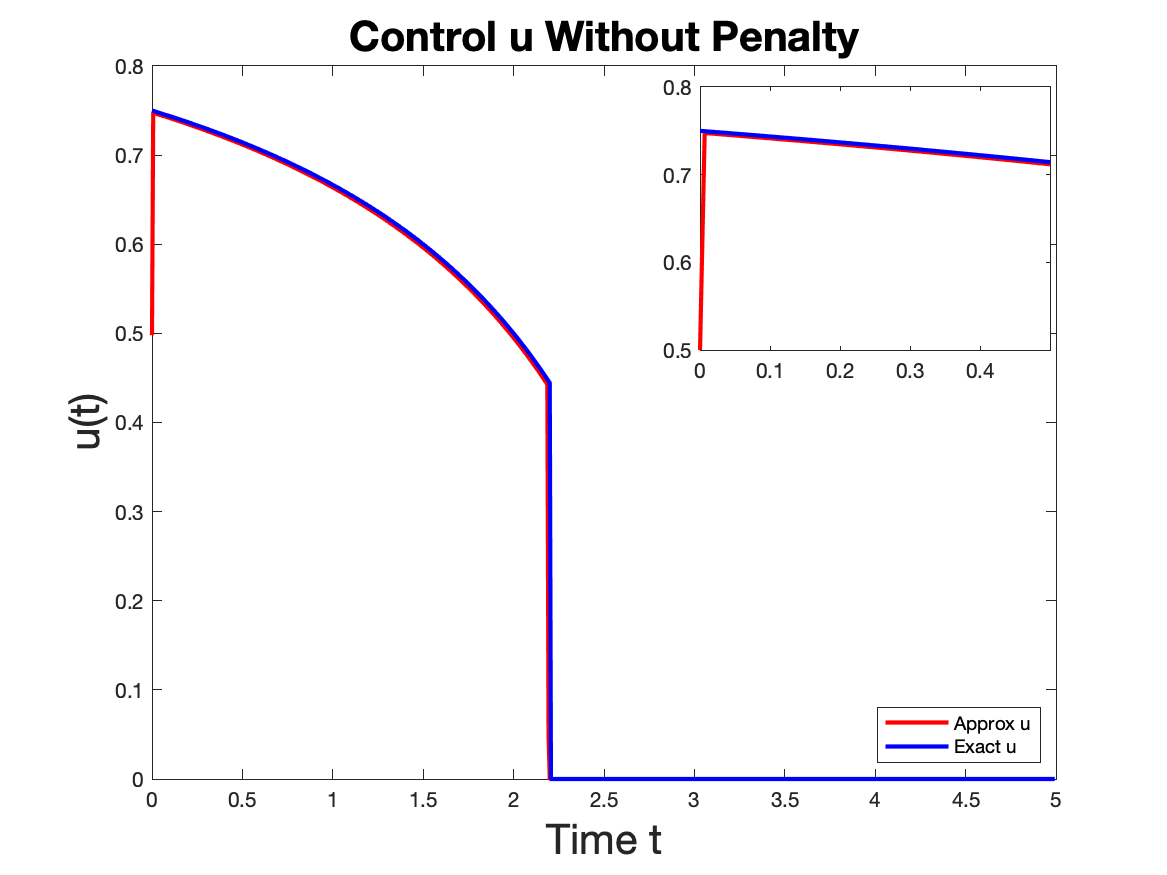}
    \caption{Approximated Allocation Strategy for Unpenalized $\;$  Problem. Inset shows zoomed in area of interest.}
    \label{fig: PlantUnopen}
  \end{subfigure}\hfill
  \begin{subfigure}[t]{0.5\textwidth}
   \includegraphics[width=\linewidth]{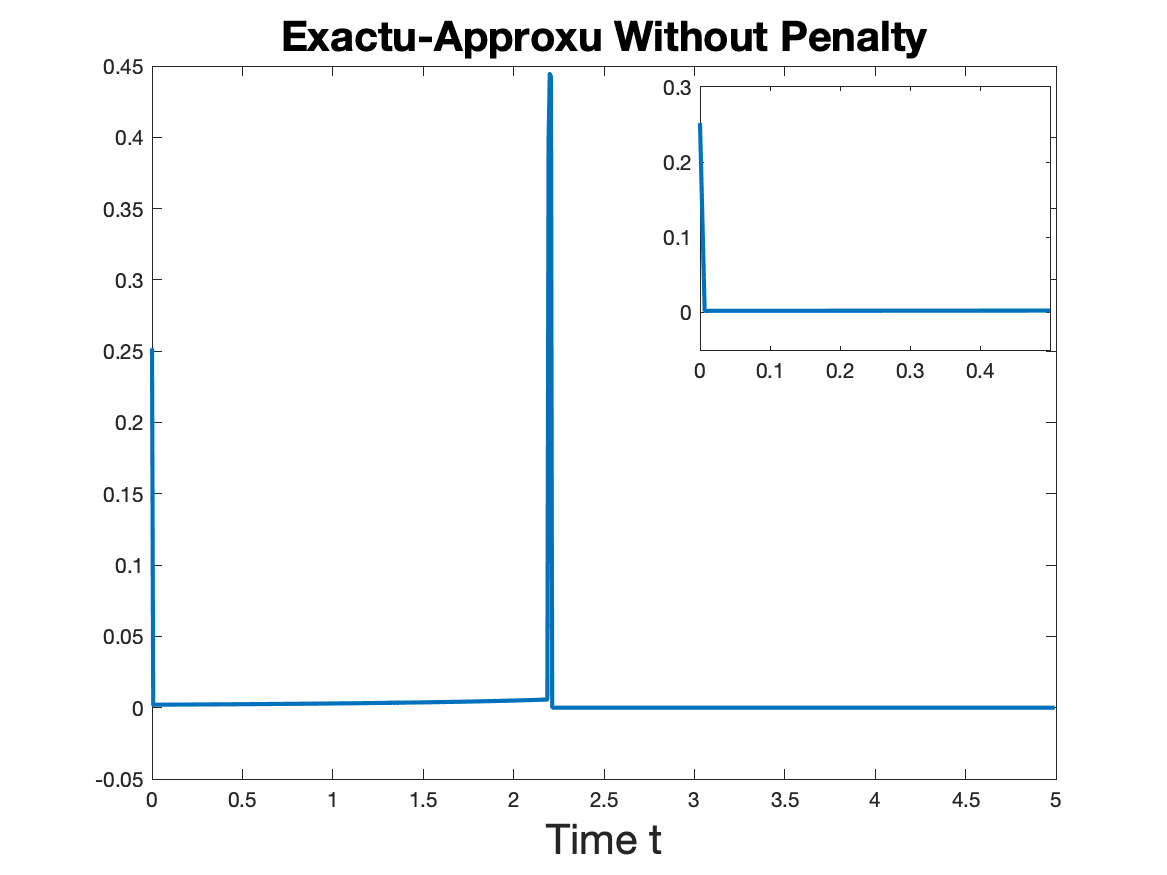}
  \caption{ $u^*-\hat{u}$ for Unpenalized Problem. Inset shows zoomed in area of interest.}
  \label{fig: plantsingdiffinset}
  \end{subfigure}
  \begin{subfigure}[t]{0.5\textwidth}
    \includegraphics[width=\linewidth]{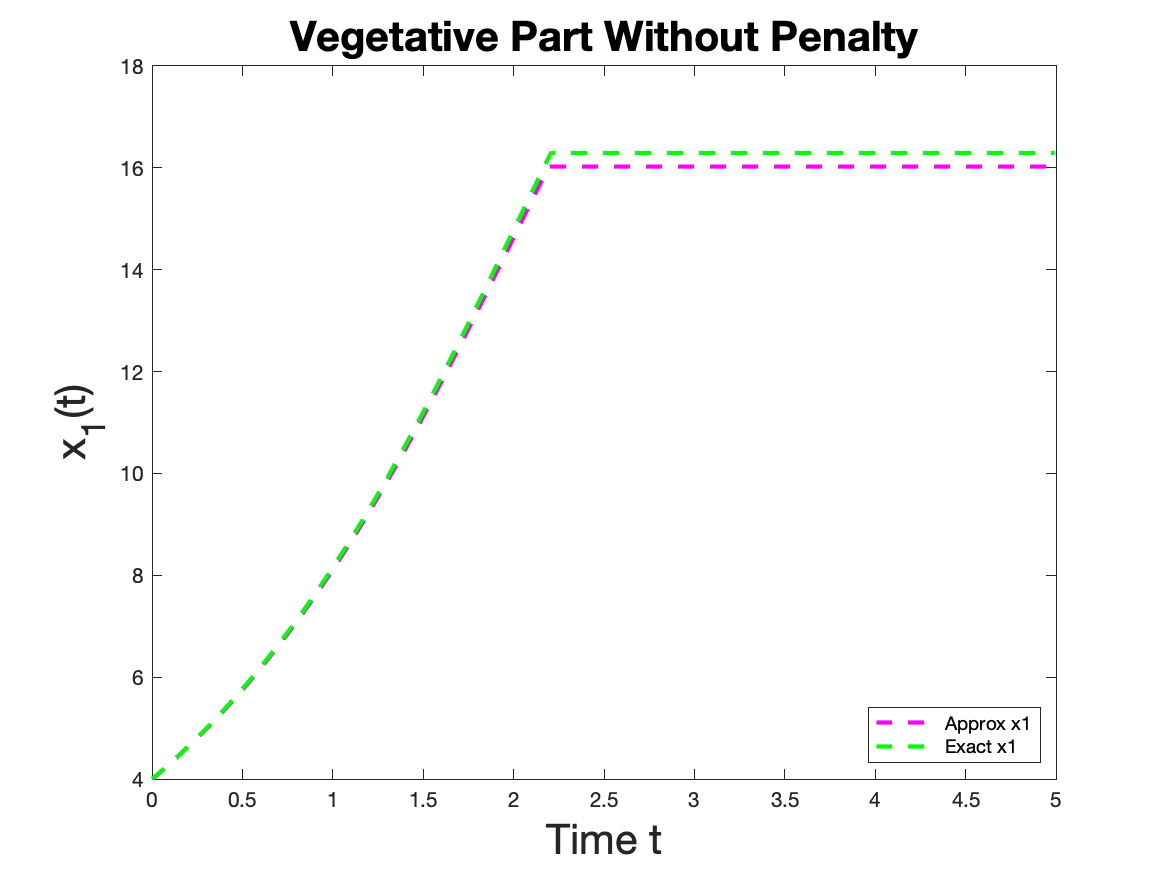}
    \caption{Vegetative Part of Plant from Unpenalized $\hat{u}$}
    \label{fig: plantsingx1}
  \end{subfigure}\hfill
  \begin{subfigure}[t]{0.5\textwidth}
    \includegraphics[width=\linewidth]{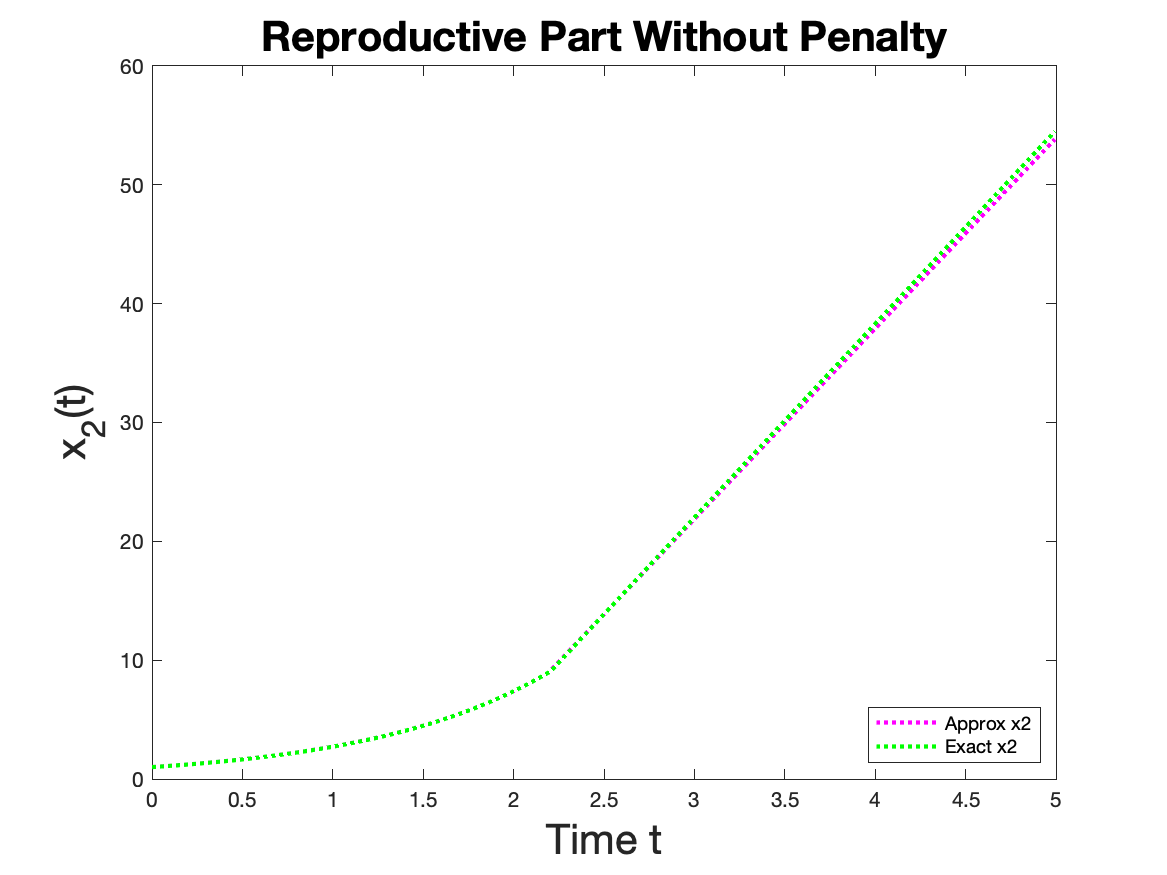}
    \caption{Reproductive Part of Plant from Unpenalized $\hat{u}$}
     \label{fig: plantsingx2}
  \end{subfigure}
  \caption{Results from Solving Unpenalized Plant Problem Case 2a)}
  \label{fig: plantsing}
  \end{figure}

   \begin{figure}[htbp]
   \begin{subfigure}[t]{0.5\textwidth}
    \includegraphics[width=\linewidth]{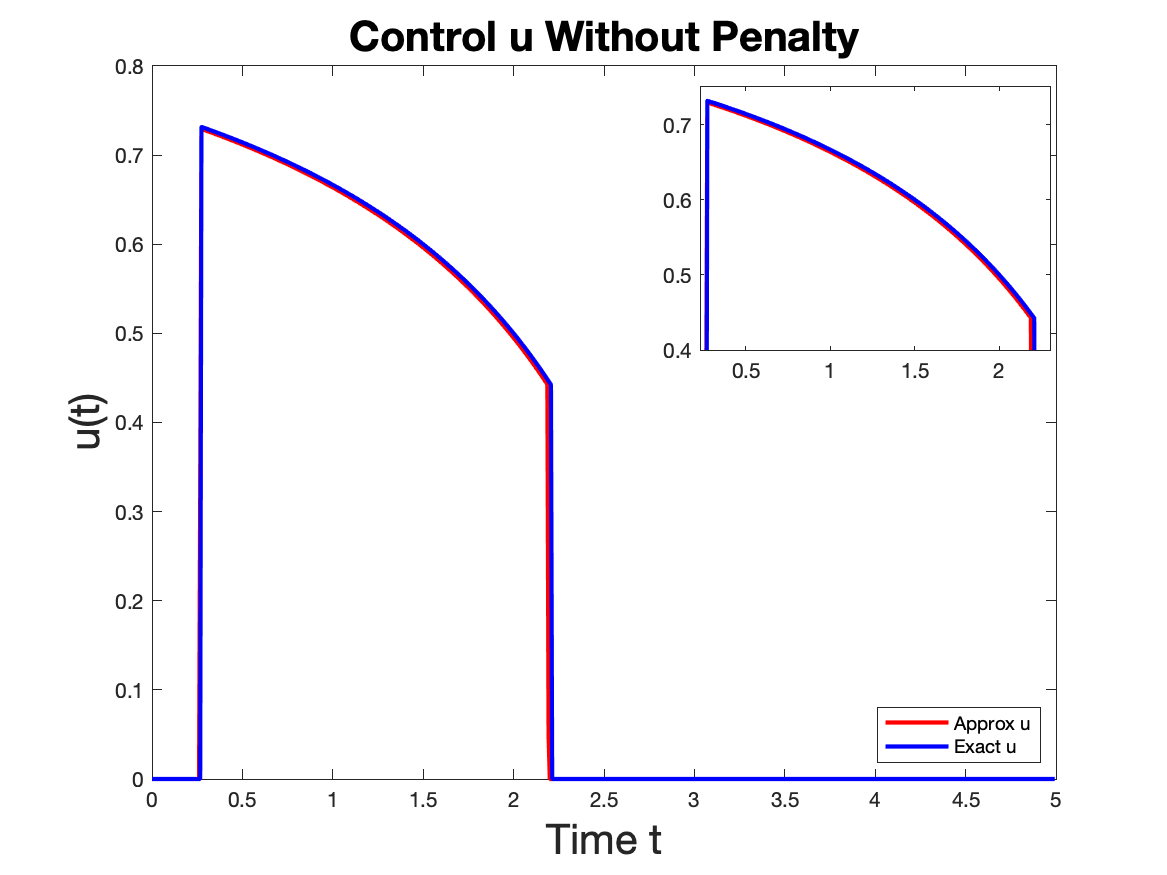}
    \caption{ Approximated Allocation Strategy for Unpenalized $\;$ Problem. Inset shows zoomed in area of interest.}
    \label{fig: PRunopen}
  \end{subfigure}\hfill
  \begin{subfigure}[t]{0.5\textwidth}
   \includegraphics[width=\linewidth]{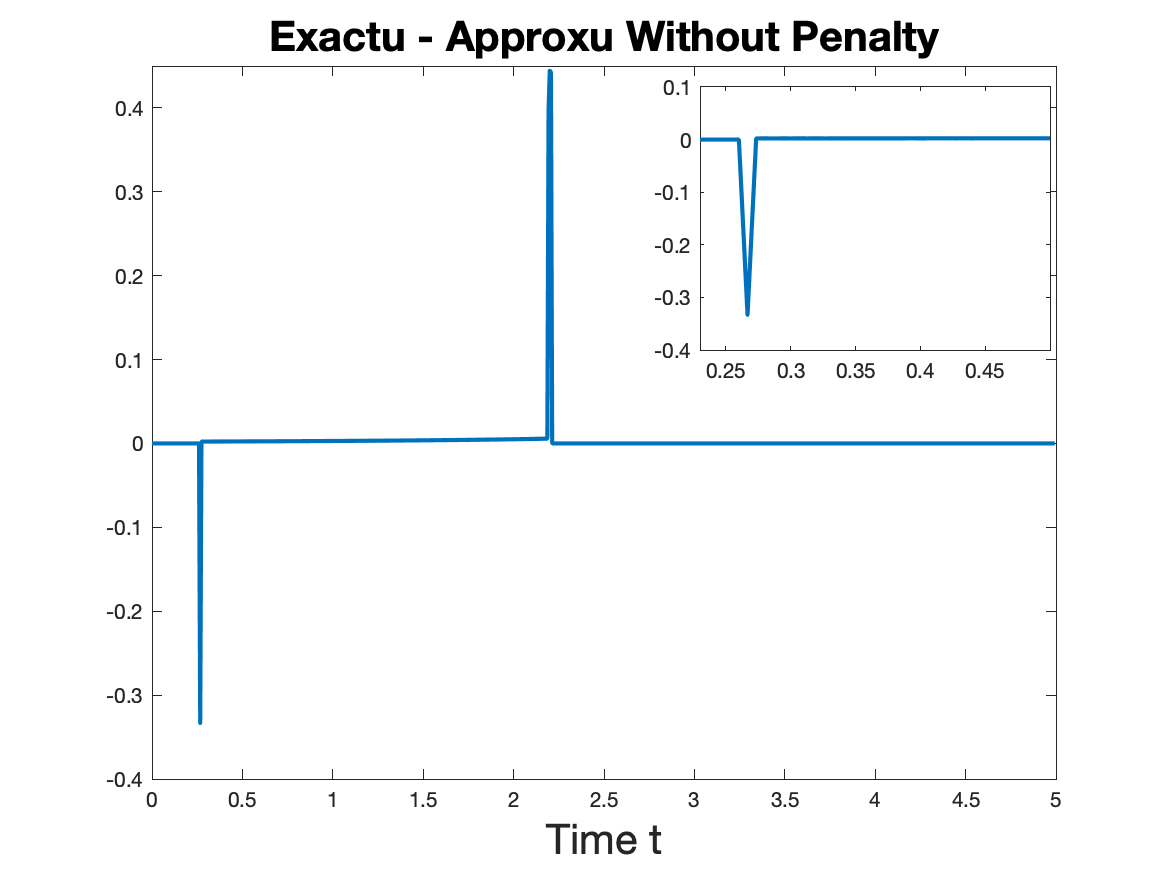}
  \caption{  $u^*-\hat{u}$ for Unpenalized Problem.  Inset shows zoomed in area of interest.}
  \label{fig: PRdiffnopen}
  \end{subfigure}
  \begin{subfigure}[t]{0.5\textwidth}
    \includegraphics[width=\linewidth]{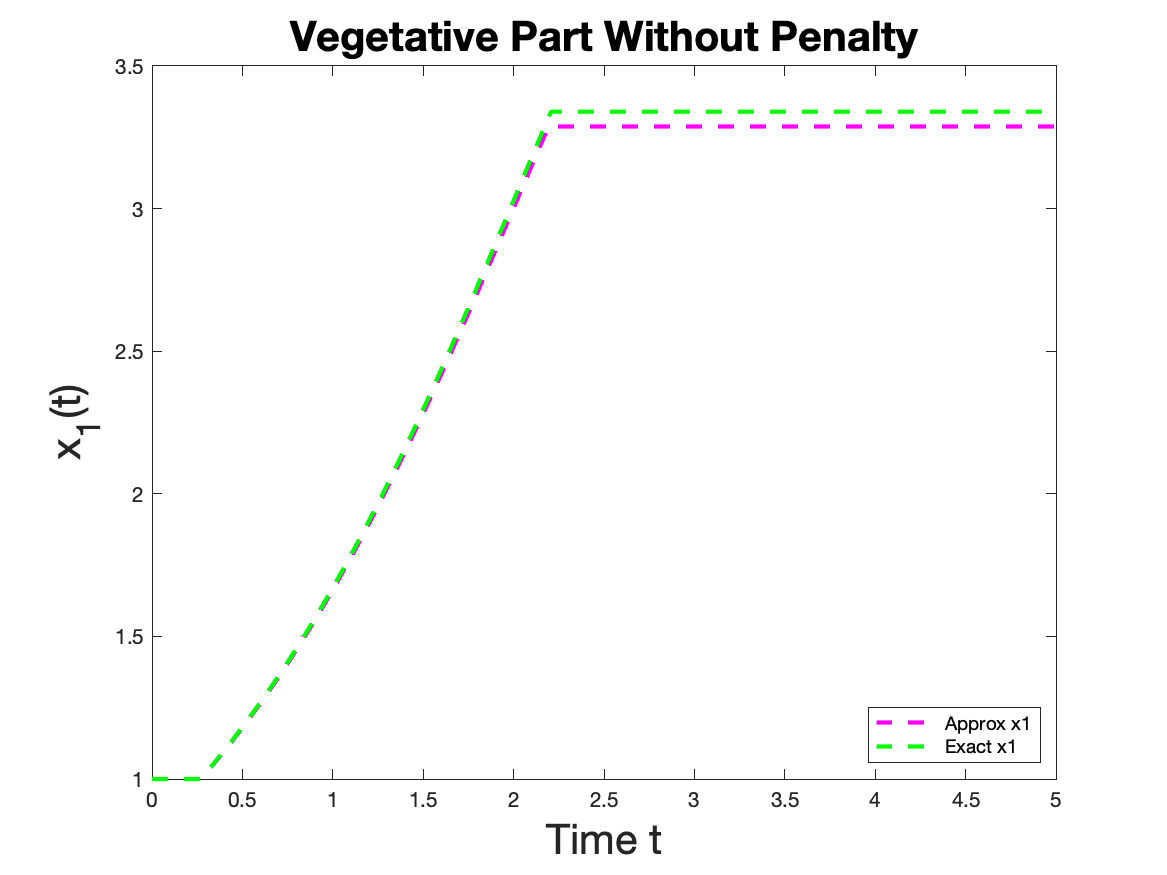}
    \caption{Vegetative Part of Plant from Unpenalized $\hat{u}$}
    \label{fig: PRsingbangx1}
  \end{subfigure}\hfill
  \begin{subfigure}[t]{0.5\textwidth}
    \includegraphics[width=\linewidth]{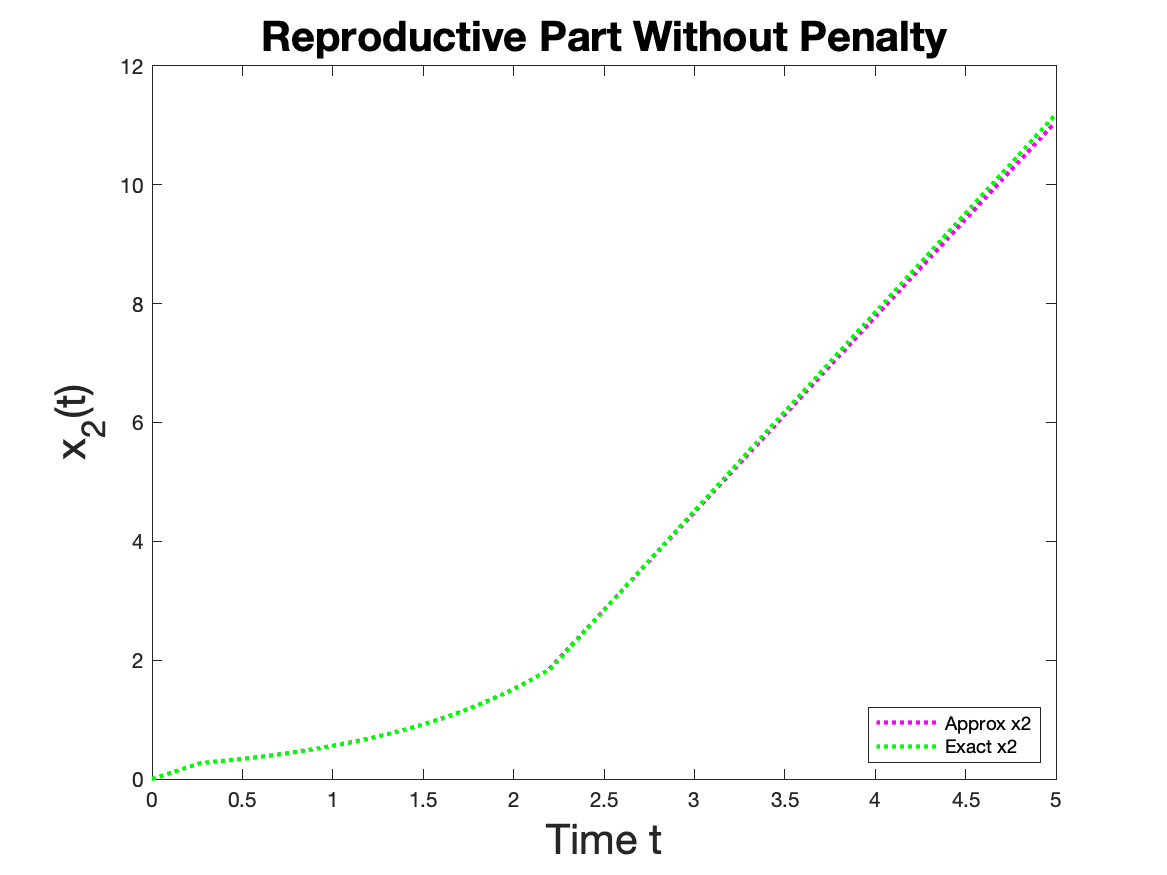}
    \caption{Reproductive Part of Plant from Unpenalized $\hat{u}$}
     \label{fig: PRsingbangx2}
  \end{subfigure}
  \caption{Results from Solving Unpenalized Plant Problem Case 2b)}
  \label{fig: PRsing}
  \end{figure}
   \begin{figure}[htbp]
   \begin{subfigure}[t]{0.5\textwidth}
    \includegraphics[width=\linewidth]{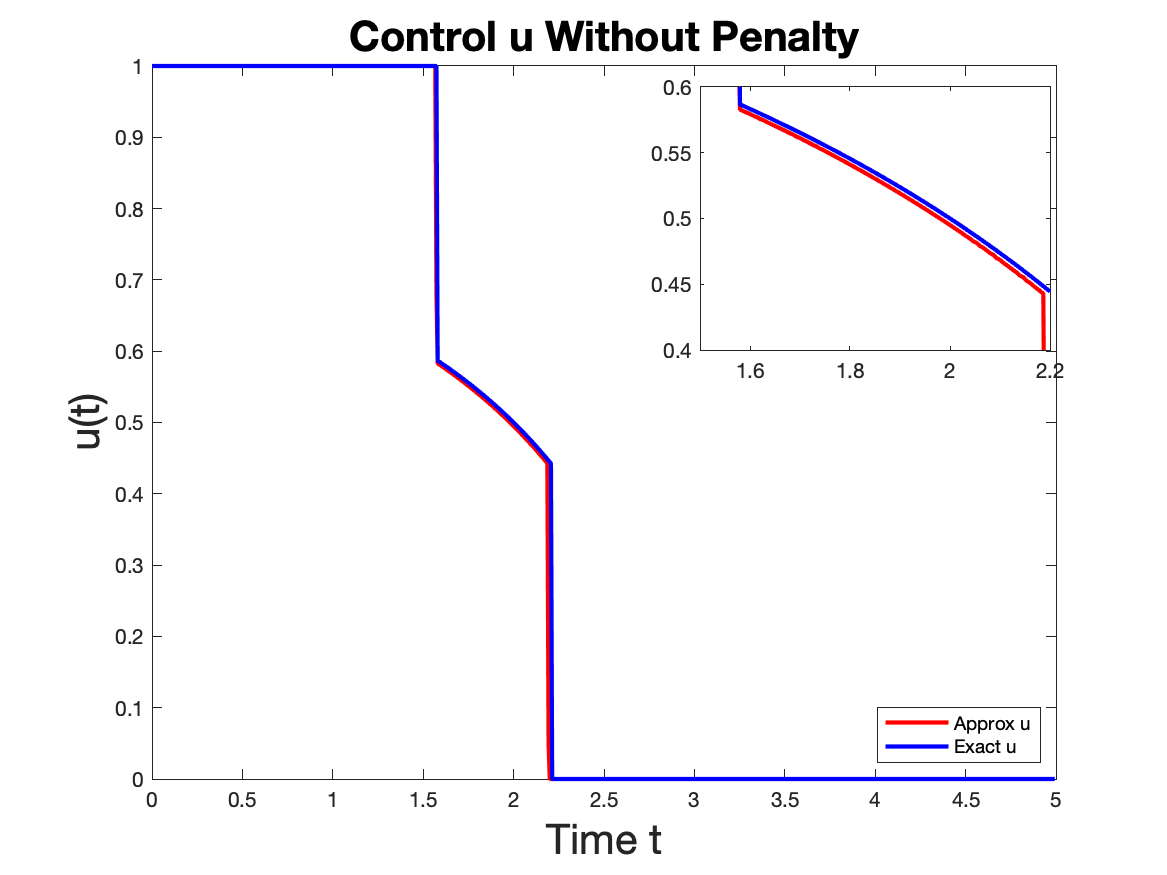}
    \caption{Approximated Allocation Strategy for Unpenalized $\;$ Problem. Inset shows zoomed in area of interest.}
    \label{fig: PVunopen}
  \end{subfigure}\hfill
  \begin{subfigure}[t]{0.5\textwidth}
   \includegraphics[width=\linewidth]{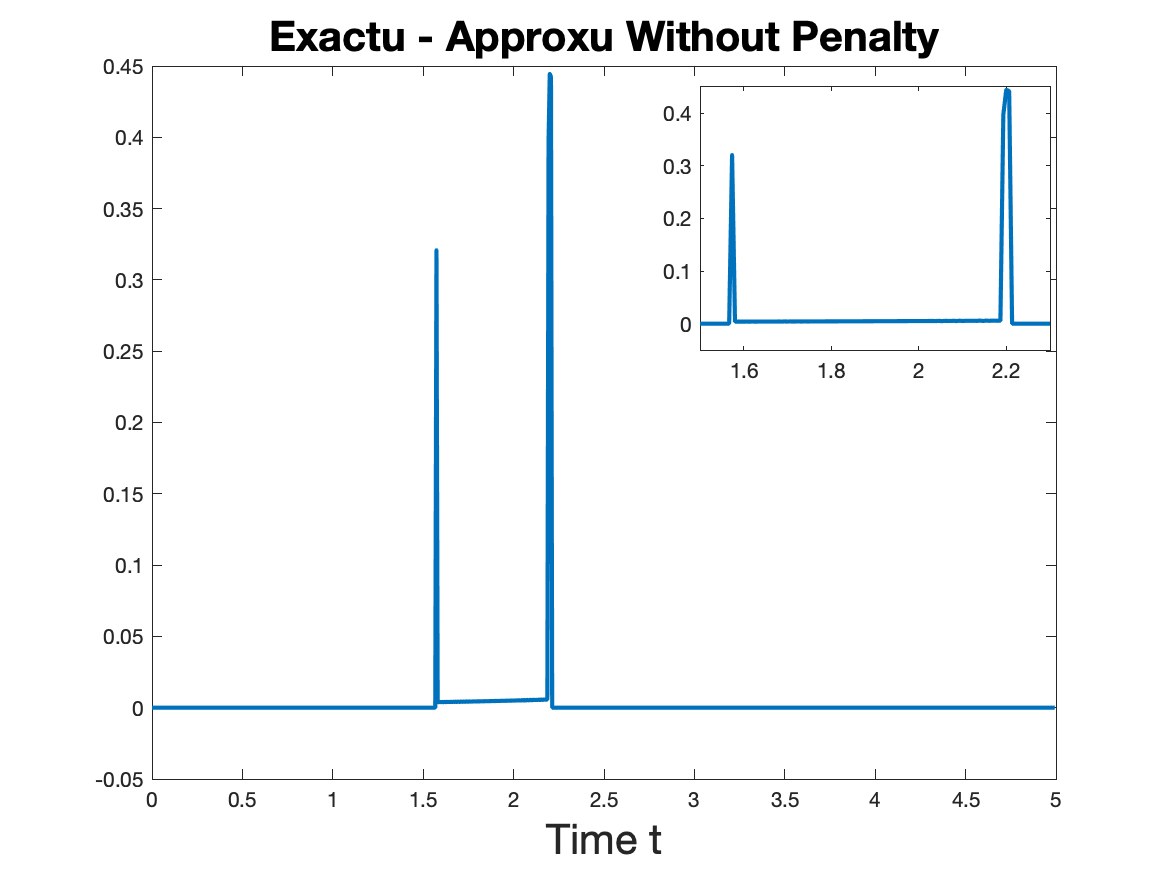}
  \caption{ $u^*-\hat{u}$ for Unpenalized Problem.  Inset shows zoomed in area of interest.}
  \label{fig: PVdiffnopen}
  \end{subfigure}
  \begin{subfigure}[t]{0.5\textwidth}
    \includegraphics[width=\linewidth]{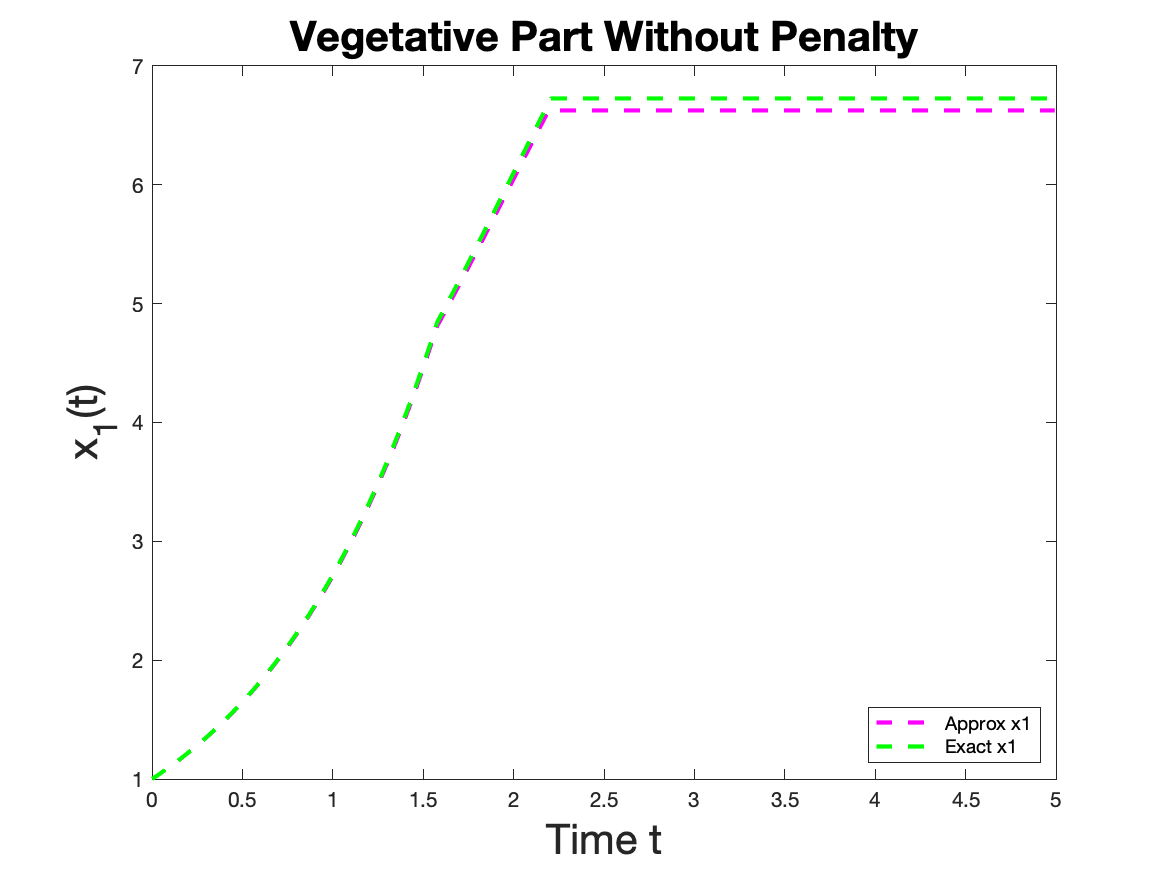}
    \caption{Vegetative Part of Plant From Unpenalized $\hat{u}$}
    \label{fig: PVsingbangx1}
  \end{subfigure}\hfill
  \begin{subfigure}[t]{0.5\textwidth}
    \includegraphics[width=\linewidth]{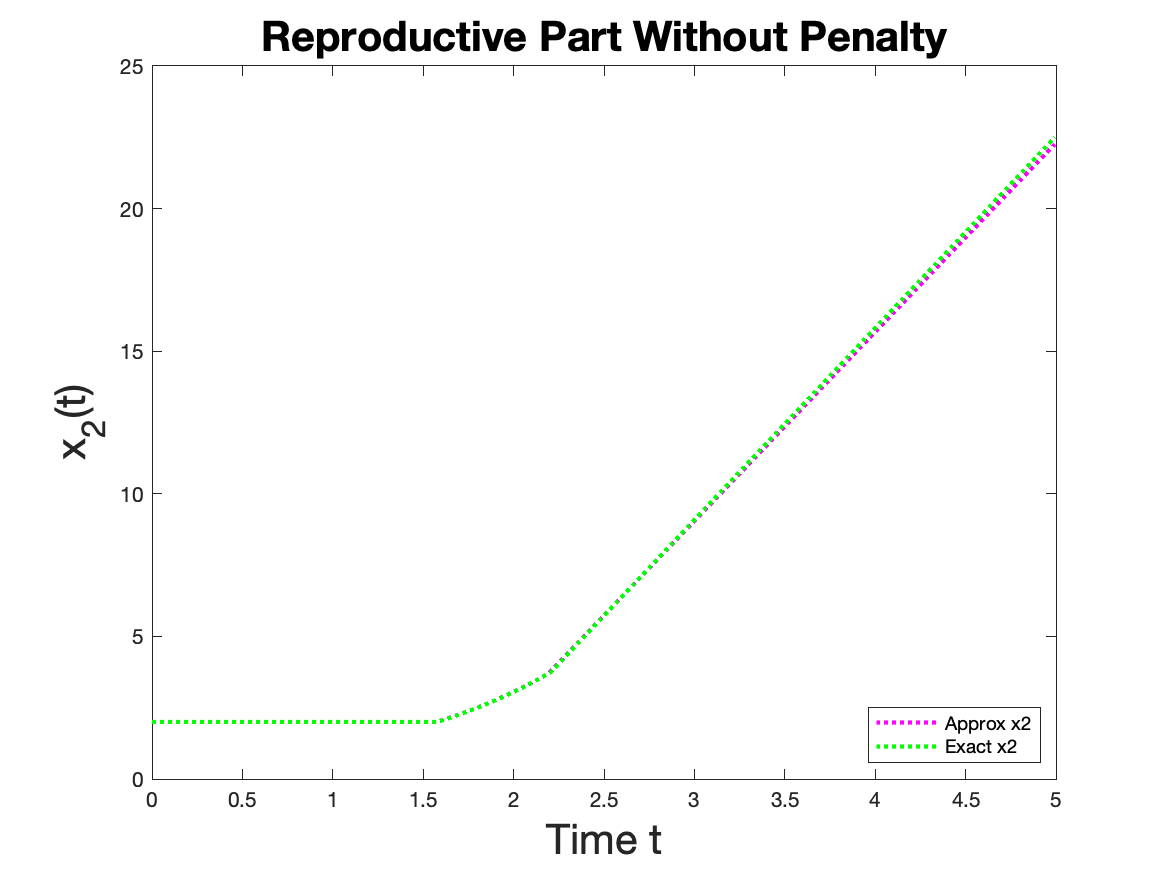}
    \caption{Reproductive Part of Plant From Unpenalized $\hat{u}$}
     \label{fig: PVsingbangx2}
  \end{subfigure}
  \caption{Results from Solving Unpenalized Plant Problem Case 2c)}
  \label{fig: PVsing}
  \end{figure}

  The unpenalized solutions that PASA obtained when parameters were set to be of case \ref{itm: case1}, case \ref{itm: case2} and case \ref{itm: case3} are shown in Figures \ref{fig: plantsing}, \ref{fig: PRsing}, and \ref{fig: PVsing} respectively; and numerical results associated with each case are provided in Table \ref{tab: PlantUnpen}. 
    As you can see in Subfigures \ref{fig: PlantUnopen}, \ref{fig: PRunopen}, and \ref{fig: PVunopen}, the unpenalized control found for each case did not present any numerical chattering. 
   We see in  Subfigures \ref{fig: PRdiffnopen} and \ref{fig: PVdiffnopen} that other than the discrepancies arising from the switching points, the unpenalized solutions for  Case 2b) and Case 2c) only differ slightly from their respective exact solution along the singular region.   
   The state solutions that correspond to the unpenalized solution for Cases 2a), 2b), and 2c), which are respectively shown in subfigures \ref{fig: plantsingx1}-\ref{fig: plantsingx2}, subfigures \ref{fig: PRsingbangx1}-\ref{fig: PRsingbangx2} and subfigures \ref{fig: PVsingbangx1}-\ref{fig: PVsingbangx2}, 
   compare relatively well with the state solutions that correspond to $u^*$. 
   In subfigures \ref{fig: plantsingx1}, \ref{fig: PRsingbangx1}, and \ref{fig: PVsingbangx1}, we see that in each case the approximate solutions associated with the vegetative parts of the plant do give slightly underestimated values along the interval $[1.5,T]$. 
   This could be due to the small discrepancies in $\hat{u}$ along the singular region. 
   Moreover, if we look at the analytic solution for $x_1$ on $[t_2,T]$ given in (\ref{eqn: plantsingsumx1}), we can see how errors between $u^*$ and $\hat{u}$ and how switching prematurely to the non-singular case could contribute errors in $x_1$ along the non-singular region $[t_2,T]$.
   Now in Subfigures \ref{fig: plantsingx2}, \ref{fig: PRsingbangx2}, and \ref{fig: PVsingbangx2} we see that solutions associated with the reproductive parts of plants give slightly underestimated values on the interval $[t_2,T]$. 
    From the explicit solution of state variable $x_2$ given in equation (\ref{eqn: plantsingsumx2}), we see that state variable $x_2$ is a line on time interval $[t_2,T]$ with slope being $x_{1}(t_2)$. 
    So for each case, the approximated solution for $x_2$ is presenting underestimated values along time interval $[t_2,T]$ because the approximation of state variable $x_1$ is yielding underestimated values near $t_2$. 
   Despite these differences, we find PASA's unpenalized solution to be a sufficient approximation to the true solution to problem (\ref{eqn: plantproblem}) in Cases \ref{itm: case2} and \ref{itm: case3}. 

  \subsection{Matlab Code for SIR Problem}\label{subsec: matlab}
  The following is the MATLAB file called demoOC.m. It is used for solving problem (\ref{eqn: SIRproblemPen}). 
  This file is included when downloading PASA package, and it is located in the following in directory \\SuiteOPT/PASA/MATLAB . 
  To access PASA software that can be used on MATLAB for Linux and Unix operating systems, download the SuiteOPT Version 1.0.0 software given on \url{http://users.clas.ufl.edu/hager/papers/Software/ }. For future reference, any updates to the software will be uploaded to this link. 
  We recommend the reader to read MATLAB file readme.m which is located in the same directory as demoOC.m. 
 The readme.m file will go into detail about the inputs that are used for running PASA. 
  Note that we have state variables initialized as $n$-vectors rather than $n+1$-vectors because from the discretization of problem (\ref{eqn: SIRproblemPen}) given in subsection \ref{sec: SIRdiscrete},  the $(n+1)^{th}$ mesh point of $S$, $I$, and $R$ is not present in the following: the approximation of the penalized objective functional, the approximation of gradient of the objective functional, and the discretized adjoint equations. 
  We only used the $(n+1)^{th}$ mesh points to generalize the transversalty conditions for the adjoint equations. 
  We still are discretizing $[0,T]$ with $n$ mesh intervals, which is why the mesh size is set to being $h=T/n$.

\begin{lstlisting}
% This test problem is an optimal control problem that arises in the
% treatment of the ebola virus (see *). An important feature of the problem is
% that it is singular. If the dynamics were discretized and the cost
% was optimized, then the solution oscillates wildly in the singular
% region, so the optimal control in that region would not be determined.
% The oscillations are removed using a penalty term based on the total
% variation of the optimal control. The dynamics are discretized using
% Euler's method with a constant control on each mesh interval. If u_i
% is the control on the i-th interval, then we add the constraint
%
% u_{i+1} - u_i = iota_i - zeta_i
%
% where zeta_i and iota_i >= 0. The penalty term that we add to the objective
% function is p * sum_{i = 1}^n iota_i + zeta_i, which is p times the total
% variation in the control. Besides the control u, which can be singular,
% there is another control v which is bang/bang. The optimization in this
% problem is over x = [u, zeta, iota, v] where u and v are both constrained
% to the interval [0, 1], while the only constraint on zeta and iota is the
% nonnegativity constraint. If there are n intervals in the mesh, then the
% dimension of x is 4*n - 2 since u and v both have n components while
% zeta and iota have n-1 components. For pasa, the constraints should be
% written in the form bl <= A*x <= bu and lo <= x <= hi. Thus a row of the A
% matrix corresponds to the constraint:
%
%          0 <= u_i - u_{i+1} + iota_i - zeta_i <= 0,  1 <= i <= n-1
% 
% The state variable in the control problem is the triple (S, I, R)
% corresponding to the susceptible, infected, and recovered individuals.
% The control u is vaccination rate while v is referred to as treatment rate.
% The controls are linear in both the dynamics and the cost, which leads
% to the singularity of the solution.
%
%   * Optimal Control for a SIR Epidemiological Model with Time-varying
%   Populations by Urszula Ledzewicz, Mahya Aghaee, and Heinz Schaettler,
%   2016 IEEE Conference on Control Applications (CCA), pp. 1268-1273,
%   DOI: 10.1109/CCA.2016.7587981

function demoOC
    % ------ Initialize constant parameters ------ %
    global gamma beta nu rho kappa mu alpha eta T n a b c p
    clf ;
    gamma = 0.00683 ;  %birth rate of the population
    nu = 0.00188 ;     %natural death rate of the population
    beta = 0.2426 ;    %rate of infectiousness of the disease
    rho = 0.007 ;      %resensitization rate
    kappa = 0.3 ;      %effectiveness of vaccination
    mu = 0.005 ;       %disease induced death rate
    alpha = 0.00002 ;  %rate at which disease is overcome
    eta = 0.1 ;        %effectiveness of treatment
    T = 50 ;           %time horizon (weeks)
    n = 750 ;          % Dimension of u and v, the number of mesh intervals
    a = 5 ;            % Constant in cost function
    b = 50 ;           % Constant in cost function
    c = 300 ;          % Constant in cost function
    p = 1e-1 ;         % penalty parameter in cost function
    umax = 1 ;         % maximum vaccination rate
    vmax = 1;          % maximum treatment rate
    h = T/n ;          % Step size
    S = zeros (n, 1) ;
    I = zeros (n, 1) ;
    R = zeros (n, 1) ;

%% --- Store problem description in a structure which we call pasadata --- %
    % ------------ Setup sparse matrix A ----------- %
    A1 = spdiags([ones(n-1,1) -ones(n-1,1)], [0,1], n-1, n) ;
    A2 = speye(n-1) ;
    A = sparse (n-1, 4*n-2) ;
    % A1 corresponds to the u_i - u_{i+1} term in the contraint while
    % A2 is used for the zeta and iota terms in the constraint
    A(:, 1:3*n-2) = [A1 -A2 A2] ;

    % put A in the structure
    pasadata.A = A ;

    % If the constraint bl <= A*x <= bu is present, then pasa uses the
    % dimensions of A to determine the number of linear constraints
    % and the number of components in x. Thus if the constraint matrix lies
    % in the upper left nrow by ncol submatrix of a larger matrix Afull,
    % set pasadata.A = Afull(1:nrow, 1:ncol). In the example above, the
    % statement A = sparse (n-1, 4*n-2) specified the dimensions of A

    % ------- store the bounds for A*x ------- %
    pasadata.bl = zeros (n-1, 1) ;
    pasadata.bu = zeros (n-1, 1) ;

    % -------- store the bounds for x -------- %
    pasadata.lo = zeros(4*n-2,1) ;
    pasadata.hi = ...
            [umax*ones(n,1); inf*ones(n-1,1); inf*ones(n-1,1); vmax*ones(n,1)] ;

    % The codes to evaluate the cost function and its gradient appear below.
    % Store the name of the codes in the pasadata structure.
    pasadata.grad  = @grad ;  % objective gradient
    pasadata.value = @cost ;  % objective value

%% ---------------- User defined parameter values for pasa ------------- %%
    % Type "pasa readme" for discussion of the parameters.
    % By default, there is no printing of statistics.
    pasadata.pasa.PrintStat = 1 ;    % print statistics for used routines
    pasadata.pasa.grad_tol = 1.e-8 ; % PASA stopping tolerance (1.e-6 default)

    % --------------- Call pasa to determine optimal x -------------- %
    [x, stats] = pasa (pasadata) ;

    % Since pasadata.pasa.PrintStat = 1, the statistics are displayed at
    % the end of the run. Since the stats structure was included as an
    % output, the corresponding numerical entries can be found in the
    % structures stats.pasa and stats.cg

    % ---------------------- Plot the states ---------------------- %
    u = x (1:n) ;          % extract the control u from the returned solution x
    v = x (3*n-1:4*n-2) ;  % extract the control v from the returned solution x
    [S, I, R] = state (u, v) ; % the state associated with the optimal controls

    t = linspace (0, T-h, n) ;
    subplot(3,2,1)
    plot (t, S,'linewidth',2) ;
    xlabel('Time')
    ylabel('S')
    subplot(3,2,2)
    plot (t, I, 'linewidth',2) ;
    xlabel('Time')
    ylabel('I')
    subplot(3,2,3)
    plot (t, R, 'linewidth',2) ;
    xlabel('Time')
    ylabel('R')

    % ----------------- Plot the controls u and v ----------------- %
    subplot(3,2,5);
    plot (t, u,'linewidth',2) ;
    xlabel('Time')
    ylabel({'Control u','(Vaccination)'})
    subplot(3,2,6);
    plot (t, v, 'linewidth',2) ;
    xlabel('Time')
    ylabel({'Control v','(Treatment)'})

%% ------------------ User defined functions for pasa ------------------ %%
    % ---- Objective function ---- %
    function J = cost(x)
        h = T/n ;
        u = x(1:n) ;
        zeta = x(n+1:2*n-1) ;
        iota = x(2*n:3*n-2) ;
        v = x(3*n-1: 4*n-2) ;
        [S, I, R] = state (u, v) ;
        J = h*(a*sum(I) + b*sum(u) + c*sum(v)) ;
        J = J + p*sum(zeta + iota) ;
    end

    % ---- Gradient of objective function ---- %
    function g = grad(x)
        h = T/n ;
        u = x(1:n) ;
        v = x(3*n-1: 4*n-2) ;

        % Compute state and costate
        [S, I, R] = state (u, v);
        [lS, lI, lR] = costate (S, I, R, u, v) ;

        % Update gradient values
        for i=1:n
            Fu(i) = h*(b - lS(i)*kappa*S(i) + lR(i)* kappa*S(i)) ;
            Fv(i) = h*(c - lI(i)*eta*I(i) + lR(i)*eta*I(i)) ;
        end

        % Store gradient values in array g to return
        g = [Fu, p*ones(1,n-1), p*ones(1,n-1), Fv] ;
    end

    % ---- State ---- %
    function [S, I, R] = state (u, v)
        h = T/n;
        S(1) = 1000 ;
        I(1) = 10 ;
        R(1) = 0 ;

        for i = 1:n-1
            N = S(i) + I(i) + R(i) ;
            S(i + 1) = S(i) + h*(gamma*N - nu*S(i) - (beta*S(i)*I(i))/N ...
                +  rho*R(i) - kappa*S(i)*u(i)) ;
            I(i + 1) = I(i) + h*(beta*S(i)*I(i)/N - (nu + mu + alpha)*I(i)...
                - eta*I(i)*v(i)) ;
            R(i + 1) = R(i) + h*( - nu*R(i) - rho*R(i) ...
                +  kappa*S(i)*u(i) + alpha*I(i) + eta*I(i)*v(i)) ;
        end
    end

    % ---- Costate ---- %
    function [lS, lI, lR] = costate (S, I, R, u, v)
        h = T/n ;
        lS(n) = 0 ;
        lI(n) = 0 ;
        lR(n) = 0 ;

        for i=n:-1:2
            N = S(i) + I(i) + R(i) ;
            lS(i-1) = lS(i)+ h*lS(i)*(gamma - nu - beta*(I(i)/N) ...
                        + beta*(S(i)*I(i)/(N^2)) - kappa*u(i))...
                        + h*lI(i)*(beta*(I(i)/N) - beta*(S(i)*I(i)/(N^2)))...
                        + h*lR(i)*(kappa*u(i)) ;
            lI(i-1) = lI(i)+ h*a +h*lS(i)*(gamma - (beta*S(i))/N ...
                        + (beta*S(i)*I(i))/N^2)...
                        + h*lI(i)*(( beta*S(i))/N - (beta*S(i)*I(i))/(N^2) ...
                        - (nu + mu + alpha) - eta*v(i))...
                        + h*lR(i)*(alpha + eta*v(i)) ;
            lR(i-1) = lR(i)+ h*lS(i)*(gamma + (beta*S(i)*I(i))/(N^2) + rho)...
                        + h* lI(i)*(- (beta*S(i)*I(i))/(N^2) ) ...
                        + h* lR(i)*(- nu - rho) ;
        end
    end
end
\end{lstlisting}

\end{document}